\documentclass[12pt]{amsart}
\usepackage{graphicx}
\usepackage{amssymb,amstext,amsfonts,amscd}
\newtheorem{theorem}{Theorem}

\newtheorem{lemma}[theorem]{Lemma}

\newtheorem{proposition}[theorem]{Proposition}
\newtheorem{corollary}[theorem]{Corollary}

\numberwithin{equation}{section}

\newcommand{\bn}{\mathbb{N}}
\newcommand{\bz}{\mathbb{Z}}
\newcommand{\br}{\mathbb{R}}
\newcommand{\bc}{\mathbb{C}}
\newcommand{\bp}{\mathbb{P}}
\newcommand{\bq}{\mathbb{Q}}
\newcommand{\bff}{\mathbb{F}}
\newcommand{\M}{\mathcal{M}}
\newcommand{\Ha}{\mathcal{H}}
\newcommand{\La}{\mathcal{L}}

\newcommand{\Ta}{\mathcal {T}}
\newcommand{\da}{\mathcal{D}}
\newcommand{\Ka}{\mathcal{K}}
\newcommand{\bee}{\mathbb{E}}
\newcommand{\bhh}{\mathbb{H}}

\newcommand{\alg}{\mathrm{alg}}

\newcommand{\rk}{\mathop\mathrm{rk}}

\newcommand{\Pic}{\mathop\mathrm{Pic}}

\newcommand{\Proof}{\noindent{\it{Proof.}}\ \ }
\newcommand{\QED}{\ \ $\Box$}

\begin{document}
\pagestyle{plain}

\begin{abstract}
Classification of real K3 surfaces $X$ with a non-symplectic involution 
$\tau$ is considered. 
For some exactly defined and one of the weakest possible type of  
degeneration (giving the very reach discriminant), we show that the 
connected component of their moduli is defined by the isomorphism class of 
the action of $\tau$ and the anti-holomorphic involution 
$\varphi$ in the homology lattice. (There are very few similar cases known.)   
For their classification we apply invariants of integral lattice involutions 
with conditions which were developed by the first author in 1983.  
As a particular case, we describe connected components of 
moduli of real non-singular curves $A\in |-2K_V|$ for the classical 
real surfaces: $V=\bp^2$, hyperboloid, ellipsoid, $\bff_1$, $\bff_4$. 

As an application, we describe all real polarized K3 surfaces which are  
deformations of general real K3 double rational scrolls (surfaces 
$V$ above). There are very few exceptions. For example, any non-singular 
real quartic in $\bp^3$ can be constructed in this way.         
\end{abstract}

\title{Real K3 surfaces with non-symplectic involutions and  
applications} 

\author{Viacheslav V. Nikulin and Sachiko Saito}

\address{Department of Pure Mathematics, The University of Liverpool, 
Liverpool, L69 3BX, United Kingdom\\ 
Steklov Mathematical Institute, ul. Gubkina 8, Moscow, GSP-1, Russia} 
\email{vnikulin@liv.ac.uk, vvnikulin@list.ru} 

\address{
Department of Mathematics Education, Hakodate Campus,
Hokkaido University of Education,
1-2 Hachiman-cho, Hakodate 040-8567, Japan}
\email{sachi63@cc.hokkyodai.ac.jp}

\maketitle

\tableofcontents

%%%%%%%%%%%%%%%%%%%%%%%%%%%%%%%%%
\section{Introduction}

In this paper, a K3 surface is a simply-connected 
K\"alerian compact complex surface with the trivial canonical class. 
Thus, there exists a 2-dimensional holomorphic differential form 
$\omega_X$ with zero divisor. The form $\omega_X$ is unique up to 
$\lambda \omega_X$ where $\lambda \in \bc$ and $\lambda\not=0$. 

A K3 surface may have two very different kinds 
of holomorphic involutions $\tau$: 
symplectic when $\tau^\ast(\omega_X)=\omega_X$ and non-symplectic 
when $\tau^\ast(\omega_X)=-\omega_X$. 

There exists only one type of 
symplectic involutions \cite{Nikulin79.1} and many types of non-symplectic 
involutions \cite{Nikulin81}, \cite{Nikulin86}. The types of non-symplectic 
involutions are classified by invariants which are all presented in 
Figure \ref{Sgraph} below. They describe irreducible components of moduli 
of complex K3 surfaces $X$ with non-symplectic involutions $\tau$. 
See more delicate classification in \cite{AlexeevNikulin88} and 
\cite{AlexeevNikulin89}. 

In this paper we consider K3 surfaces with non-symplectic involutions 
$(X,\tau)$. The quotient surface $Y=X/\{1,\tau\}$ is either Enriques 
surface when $\tau$ has no fixed points, or a rational surface when 
$\tau $ has fixed points. The set of fixed points is a non-singular curve 
$A=X^\tau$, and its image $A\subset Y$ is an element of the linear system 
$|-2K_Y|$ where $K_Y$ is the canonical system of $Y$. Then $X$ is 
the double covering of $Y$ ramified in $A$ and $\tau$ the involution of 
the double covering. Pairs $(A,Y)$ when $Y$ is either Enriques or a rational 
surface (they are characterized by the condition $q(Y)=0$ where $q$ is 
the irregularity) and $A\in |-2K_Y|$ is a non-singular curve are 
called right DPN-pairs. Classification of K3 surfaces with non-symplectic 
involutions $(X,\tau )$ is equivalent to classification of right DPN-pairs 
$A\subset Y$. 

In this paper we consider description of connected components of moduli
of real K3 surfaces with non-symplectic involutions and real right 
DPN-pairs. Equivalently, we consider classification of 
triplets $(X,\tau,\varphi)$ where $X$ is a K3 surface, $\tau$ its holomorphic 
non-symplectic involution and $\varphi$ its anti-holomorphic involution. 
The involutions $\tau$ and $\varphi$ must commute: $\varphi\tau=\tau\varphi$. 
Equivalently, we consider classification of real right DPN-pairs 
$A\subset Y$ with an anti-holomorphic involution $\theta$ of $Y$ where  
$\theta(A)=A$. 

This classification depends on the definition of degeneration of 
$(X,\tau,\varphi)$. We choose one of {\it the weakest types of degeneration}  
(i. e. {\it giving one of the reachest discriminants}). 
Geometrically (over $\bc$) it is equivalent to appearance of 
an exceptional curve $F$ with square $-2$ on 
the quotient $Y$. Over $\br$ this curve should be either real: 
$\theta(F)=F$, or $\theta(F)\cdot F \le 1$. See Section \ref{modinvolutions}. 

We mention that our definition of degeneration is weaker (it gives 
a much reacher discriminant) then the one used in \cite{DIK2000} for 
description of connected components of moduli of real Enriques surfaces. 
Thus, even for real Enriques surfaces our results are 
stronger (more exact). Difference of our results from 
\cite{DIK2000} is in the appearance of the group $G$.  

On the one hand, our type of degeneration gives {\it a reach discriminant.} 
For example, it is enough to describe 
connected components of moduli of real bidegree $(4,4)$ curves on 
the hyperboloid which was one of the main purposes of this paper. 
A degeneration of this case gives non-singular curves  
$A\in |-2K_{\bff_2}|$ which we must exclude. 
Here we denote as $\bff_n$, $n\ge 0$ a relatively 
minimal rational surface with 
the section $s$ having $s^2=-n$. See Proposition 
\ref{factorU(2)}.  On the other hand, using Global Torelli Theorem  
\cite{PS71} and epimorphicity of Torelli map \cite{Kulikov77} for K3 
surfaces we prove that connected components of 
moduli are described by sufficiently simple invariants which are 
{\it isomorphism classes of integral involutions with conditions.} 
They are triplets $(L,\varphi,(S,\theta))$. Here  $L=H_2(X, \bz)$ is the 
even unimodular homology lattice with the 
action of $\varphi$, $S=L^\tau$ and $\theta=\varphi|S$. The pair 
$(S,\theta)$ determines the type of real K3 surfaces $(X,\tau,\varphi)$ 
with non-symplectic involution whose connected components of moduli 
we want to describe. 
See Theorems \ref{genmod}, \ref{gnmoddpn+}, \ref{gnmoddpn}. 
All genus invariants of integral involutions with condition were 
classified in \cite{Nikulin83}. See Sect. \ref{genclassification} 
about these invariants.  They give very universal invariants which 
work in all cases we consider. This paper can be considered as 
one of applications of this classification and the invariants which it 
uses.

We mention that similar arithmetic description of connected components of 
moduli had been obtained in \cite{Nikulin79} for real polarized K3 surfaces. 
Similar result was obtained in \cite{DIK2000} for real Enriques surfaces 
with the much stronger (or with the much smaller discriminant) 
type of degeneration used in \cite{DIK2000}. Thus, this paper adds some  
more cases to these beautiful series. It seems there are very few 
geometrically significant cases when similar arithmetic results are valid.  

Our results are valid for any type $(S,\theta)$ of real K3 
surfaces with non-symplectic involutions (or the real right DPN-pairs).   
In Sections \ref{P2} --- \ref{F4} we specify them in the most simple  
and general cases (by the dimension $20-\rk S$ of moduli).   
We consider all cases when $\rk S\le 2$. They give the following 
real DPN-pairs and the corresponding real K3 surfaces with non-symplectic 
involutions as double coverings:

\medskip

\noindent
$\bp^2$: Real non-singular curves $A$ of degree 6 on the real projective 
plane $\bp^2$; see Sect. \ref{P2}.

\medskip

\noindent
$\bhh$ (hyperboloid): Real non-singular curves $A$ of bidegree $(4,4)$ on 
hyperboloid and on $\bp^1\times \bp^1$ with all possible spinor 
(i. e. with empty set of real points) real structures; see Sect. 
\ref{hyperboloid}.

\medskip

\noindent
$\bee$ (ellipsoid): Real non-singular curves $A$ of bidegree 
$(4,4)$ on ellipsoid, see Sect. \ref{ellipsoid}.  

\medskip

\noindent
$\bff_1$: Real non-singular curves $A\in |-2K_{\bff_1}|$ on 
$\bff_1$ with real structure; see Sect. \ref{F1}.   

\medskip

\noindent
$\bff_4$: Real non-singular curves $A\in|-2K_{\bff_4}|$ on $\bff_4$ with a  
real structure; see Sect. \ref{F4}.

\medskip

In all these cases we get classification of connected components of 
moduli of the non-singular real curves $A$ on the real surfaces up to the   
action of the automorphism group of the real surfaces. 

The case $\bp^2$ had been considered in \cite{Gudkov69} (isotopy 
classification), \cite{Rokhlin78} (classification of divideness and 
complex orientations), \cite{Nikulin79} (classification of 
connected components of moduli). 

For hyperboloid $\bhh$, the isotopy classification had been obtained  
by Gudkov \cite{Gudkov79}. Zvonilov \cite{Zvonilov92} added description of 
divideness by $A(\br)$ of $A(\bc)$ and 
possible complex schemes (i. e. complex orientations of $A(\br)$ which are 
induced by $A(\bc)$ in the dividing case). We show that these 
visual invariants: {\it isotopy type, divideness and complex orientations (in some  
very few dividing cases), define the connected component of moduli} up to the 
action of the automorphism group which has 8 connected components for 
hyperboloid.   

For ellipsoid $\bee$, the isotopy classification had been obtained  
by Gudkov and Shustin \cite{GudkovShustin80}. 
Zvonilov \cite{Zvonilov92} added description of divideness by 
$A(\br)$ of $A(\bc)$ and possible complex schemes which 
are prescribed by the isotopy type 
in this case. We show that these visual invariants: {\it isotopy type and 
divideness, define the connected component of moduli,} up to the 
action of the automorphism group of ellipsoid which has 2 connected 
components. (Similar result is valid for $\bp^2$.)

For $\bff_1$ we show that also {\it the connected component of moduli is 
defined by the isotopy type and by divideness.} The related result 
had been obtained by Itenberg \cite{Itenberg92}, \cite{Itenberg94} who 
gave classification of connected components of moduli of real curves of 
degree 6 with one non-degenerate quadratic singular point. It seems,  
our classification should follow from Itenberg's results, but it would be 
hard for us to deduce it from \cite{Itenberg92}, \cite{Itenberg94}. 

For $\bff_4$ we also show that {\it the connected component of 
moduli is defined by the isotopy type and by divideness.} 
This case corresponds to very important class of K3 surfaces which 
are elliptic K3 surfaces with section. The involution $\tau$ is the  
inverse map of the corresponding elliptic pencil.

In Sect. \ref{polK3} we consider an important application of the above 
results to description of real polarized K3 surfaces which are 
{\it deformations of general real K3 double rational scrolls.}  
We mean the following. Let $(X,P)$ be a polarized K3 surface with a  
primitive ample $H$ where $P^2=n$. For a {\it general K3 double 
rational scroll} 
the complete linear system $|P|$ or $|2P|$ gives the double covering 
$\pi:X\to Y=X/\{1,\tau \}$ where $Y$ is one of surfaces $Y=\bp^2$, 
$\bp^1\times \bp^1$, $\bff_1$ or $\bff_4$ above. Moreover, $|kP|$ is 
the preimage $|kP|=\pi^\ast |\widetilde{P}|$ of some standard linear system 
$\widetilde{P}$ on $Y$. General K3 double rational scrolls give 
codimension 1 moduli subspace (if $n\ge 4$) of the connected moduli 
space $\M_n$ of the polarized K3 surfaces.  
Thus, any polarized K3 surface $(X,P)$ over $\bc$ is a deformation of some 
general K3 double rational scroll. It is interesting to ask, if the 
same is valid over $\br$.  

Using classification results of Sects. \ref{P2} --- \ref{F4} 
and the classification results about real polarizes K3 
surfaces \cite{Nikulin79}, in Sect. \ref{polK3} we show that 
a real polarized K3 surface $(X,P)$ of degree $n=P^2\ge 6$ 
is not a deformation of a general real K3 double rational scroll, 
if the real part $X(\br)$ either consists of 10 spheres or it consists of 
$k$ spheres and $X(\br )\sim P\mod 2$ in $H_2(X,\bz)$. ``All other''   
(in some sense) real polarized K3 surfaces $(X,P)$ are  
deformations of general real K3 double rational scrolls. 
In \cite{Nikulin79} all possible genus invariants 
(for the action of the anti-holomorphic involution $\varphi$ in homology) 
$$
(n; r,a,\delta_P, \delta_{\varphi},\delta_{\varphi P})
$$ 
of real polarized K3 surfaces were defined and classified.  
Excluding few (possibly extremely few) exceptions the genus 
invariants determine the connected component of moduli of 
real polarized K3 surfaces. ``All other'' above means all possible genus 
invariants different from genus invariants of the shown above exceptions.   
Additional considerations show that for $P^2=2$ (double coverings of 
$\bp^2$) and $P^2=4$ (quartics in $\bp^3$) any real K3 surface 
is a deformation of a general real K3 double rational scroll. 
For $P^2=2$ it is obvious. For quartics, $P^2=4$, we can additionally use 
double coverings of ellipsoid $\bee$ which one cannot use for $P^2\not=4$. 
This gives a very effective method of construction of all (i. e.  
representatives of all connected components of moduli) 
real quartic surfaces. Using results of Shah \cite{Shah81} (over $\bc$),  
one can describe all these deformations very explicitly. It gives 
an effective method of construction of real quartic surfaces.  

The first author wants to thank University of Liverpool, 
Steklov Mathematical Institute, Max-Planck Institut f\"ur Mathematik and 
Hakodate Campus of Hokkaido University of Education for hospitality.

%%%%%%%%%%%%%%%%%%%%%%%%%%%%%%%%%

\bigskip

\section{Connected components of
moduli of real K3 surfaces with non-symplectic involutions
and real DPN-pairs}
\label{modrealK3}

\subsection{Connected components of moduli of non-degenerate real K3
surfaces with non-symplectic involutions}
\label{modinvolutions}

We remind that a K3 surface $X$ is a non-singular projective algebraic
surface over $\bc$ (more generally a compact K\"ahlerian surface)
such that it has a  holomorphic 2-dimensional differential form
$0\not=\omega_X$ with zero divisor (or the canonical class $K_X=0$) and
$X$ is simply-connected. The form $\omega_X$ is
defined uniquely up to $\lambda \omega_X$, $\lambda \in \bc^\ast$.
E.g. see \cite{Shafarevich65}, \cite{PS71} and \cite{Saint-Donat74} 
about K3 surfaces.

\medskip

An involution $\tau$ of a K3 surface $X$ can be {\it symplectic,} if
$\tau^\ast \omega_X=\omega_X$, and {\it non-symplectic,} if
$\tau^\ast \omega_X=-\omega_X$. Symplectic involutions on K3
surfaces were classified in \cite{Nikulin79.1}, for this case
it is natural to consider K\"ahlerian K3 surfaces, and there exists only
one type of such involutions. Non-symplectic involutions on K3 surfaces
were classified in \cite{Nikulin81}, \cite{Nikulin86}; all K3 surfaces
with  non-symplectic involutions are algebraic.

The main invariant of a K3 surface with a non-symplectic involution
$\tau$ is the isomorphism class of the lattice $S$ which is
the fixed part
$H_2(X,\bz)^\tau=S$. Here and in what follows we denote
$$
M^\xi=\{x\in M\ |\xi(x)=x\},
$$
and by
$$
M_\xi=\{x\in M\ |\xi(x)=-x\}
$$
for an action of an involution $\xi$ on a module $M$.
By the Hodge decomposition, the lattice $S$ is hyperbolic (i. e. it 
has exactly one positive square), 
it is contained in the Picard lattice $N(X)$ of $X$
(which is generated by classes of all algebraic curves on $X$).
All possible $S$ are known (up to isomorphisms) and are described in
\cite{Nikulin81} (see also \cite{Nikulin86}).

We can additionally fix a half-cone (the {\it light-cone})
$V^+(S)$ of the cone
$$
V(S)=\{x\in S\otimes \br\ |\ x^2>0\}.
$$
We can also fix a fundamental chamber $\M\subset V^+(S)$ for the group
$W^{(-2)}(S)$ generated by reflections in all elements with square $(-2)$
in $S$. This is equivalent to fixing {\it a fundamental subdivision}
$\Delta(S)=\Delta(S)_+\cup -\Delta(S)_+$ of all elements with
square $-2$ in $S$. The $\M$ and $\Delta(S)_+$ define each other by
the condition $(\M, \Delta_+)\ge 0$. These additional structures
$\M\subset V^+(S)$ of the hyperbolic lattice $S$
are defined uniquely up to the action of the group
$\{\pm 1\}W^{(-2)}(S)$.

We can restrict considering K3 surfaces $X$ with a non-symplectic
involution $\tau$, i.e. pairs $(X,\tau)$, such that
$V^+(S)$ contains a hyperplane section of $X$ and the set
$\Delta(S)_+$ contains only classes of effective curves with square
$-2$ of $X$. I. e. $V^+(S)$ and the fundamental subdivision
$\Delta(S)_+$ are prescribed by the geometry of the K3 surface $X$.

\medskip

If a pair $(X, \tau)$ is general, then $S=N(X)$ is the Picard
lattice of $X$ (remind that it is generated by classes of 
algebraic curves) and $\M$ gives the $nef$ cone (or K\"ahlerian cone)
of $X$. The weakest condition of degeneration
(i. e. giving the most reach discriminant) is the
following condition:

\medskip

($\da$): {\it We say that $(X,\tau)$ of the type $H_2(X,\bz)^\tau=S$
is degenerate} if there exists $h\in \M$ such that $h$ is not
$nef$ for $X$. I.e., there are elements of $S$ which are
$nef$ in general (when $S=N(X)$),
but they are not $nef$ for $X$.
It is not difficult to see \cite{Nikulin86}, that this is equivalent
to existence of an exceptional curve with square $-2$ on the quotient
$Y=X/\{1,\tau\}$.  This is also equivalent to have
an element $\delta \in N(X)$ with $\delta^2=-2$ such that
$\delta=(\delta_1+\delta_2)/2$ where $\delta_1\in S$, $\delta_2\in
S^\perp_{N(X)}$ and $\delta_1^2=\delta_2^2=-4$. Remark that
$(\delta_1,S)\equiv 0\mod 2$ and $(\delta_2, S^\perp)\equiv 0\mod 2$.
I.e. $\delta_1$ and $\delta_2$ are roots with square $-4$ for lattices
$S$ and $S^\perp$ respectively.

\medskip

Let us consider examples of lattices $S$ of a small rank, when
$\rk S\le 2$. All these cases will be considered in details in 
Sect. \ref{typesStheta} below.  
We use notations: $\langle A\rangle$ denote a lattice
defined by a symmetric integral matrix $A$. We denote
$U=
\left(\begin{array}{cc}
0 & 1\\
1 & 0
\end{array}\right).
$
For a lattice $K$, we denote by $K(m)$ a lattice obtained from
the lattice $K$ by multiplication of the form of $K$ by $m\in \bq$.
By $\oplus$ we denote the orthogonal sum of lattices.

If $\rk S=1$, then $S\cong \langle 2 \rangle$,
$Y=X/\{1,\tau\}=\bp^2$ and all $X$ are non-degenerate,
because $S$ has no elements with square $-4$.

If $\rk S=2$,
then the lattice $S$ is isomorphic to one of three lattices:
$S\cong U(2)$, $\langle 2 \rangle \oplus \langle -2 \rangle$ or $U$. 

If $S\cong U(2)$, then
non-degenerate K3 surfaces $(X,\tau)$ give
$X/\{1,\tau\}\cong \bp^1\times \bp^1=\bff_0$, but degenerate K3 surfaces
$(X,\tau)$ give $X/\{1,\tau\}\cong \bff_2$ where
we denote by $\bff_n$ a relatively minimal rational surface
with a section having the square $-n$, $n\ge 0$.

If $S\cong \langle 2 \rangle \oplus \langle -2 \rangle$, then
$Y=X/\{1,\tau\}=\bff_1$. 
All $X$ are non-degenerate because $S$ has no elements
with square $-4$.

If $S\cong U$, then $Y=X/\{1,\tau\}=\mathbb{F}_4$. 
All $X$ are non-degenerate because $S$ has no elements $\delta_1$
with $\delta_1^2=-4$ and $(\delta_1, S)\equiv 0\mod 2$.

We shall consider all these cases in Sect \ref{typesStheta} below. 
\medskip

Further we fix one of possible lattices $S$ together with
$\M\subset V^+(S)$ and
consider K3 surfaces with non-symplectic involutions $(X,\tau)$
and a fixed identification $S=H_2(X,\bz)^\tau$ such that
$\M\subset V^+(S)$ are defined by the geometry of $X$.
In complex case we can consider two pairs $(X,\tau)$ and
$(X^\prime,\tau^\prime)$ to be isomorphic if there exists
an isomorphism $f:X\to X^\prime$ of K3 surfaces which is
identical on $S$. Then $f\tau =\tau^\prime f$ (it is sufficient to
check this identity for the action on second homology)
and $f$ defines the isomorphism of K3 with involutions.
The moduli space of pairs $(X,\tau)$
is connected, and its subspace of
degenerate pairs $(X,\tau)$ has codimension one over $\bc$
(if it is not empty). Thus, the space of non-degenerate pairs is
also connected. This follows from the fact that there exists
a unique (up to automorphisms) primitive embedding
$S\to L$ where $L\cong H_2(X,\bz)$ is an unimodular even lattice
of signature $(3,19)$.

\medskip

For real K3 surfaces with non-symplectic involutions $(X,\tau)$,
we should additionally fix the type $\theta$ of the action of the
anti-holomorphic involution $\varphi$ of $X$ on $S$. The $\theta$ should
satisfy the following properties:
$\theta(V^+(S))=-V^+(S)$ and $\theta(\Delta(S)_+)=-\Delta(S)_+$.
It follows that the lattice $S_+=S^\theta$ is negative definite
and it has no elements with square $-2$. Moreover, the linear subspace
$S_-\otimes \br$ where $S_-=S_\theta$ must intersect the interior of the
$nef$ cone $\M$.
For the fixed type $(S,\theta)$ we consider K3 surfaces $X$
with holomorphic involutions $\tau$ of the type $S$ and
an anti-holomorphic involution $\varphi$, such that $\varphi|S=\theta$,
in particular, $\varphi(S)=S$.
It follows that the anti-holomorphic involution $\varphi$ commutes with
$\tau$. The $\varphi$ defines then the real structure of the pair $(X,\tau)$.
The triplets $(X,\tau,\varphi)$ of the type $(S,\theta)$ give 
{\it real K3 surfaces with non-symplectic involutions of the 
type $(S,\theta)$.}

We consider the following real analogy for real K3 surfaces
with non-symplectic involutions of the type $(S,\theta)$
of the degeneration ($\da$) above. An element $h\in S$ is called
{\it real} if $\theta(h)=-h$, i.e. $h\in S_-$. For a general real $X$ we have
$N(X)=S$, and all real $nef$ elements are elements of $S_-\cap \M$.

\medskip

($\da\br$):
{\it A real K3 surface $(X,\tau,\varphi)$ with a non-symplectic 
involution of the type $(S,\theta)$ is called degenerate} if there exists a real
element $h\in S_-\cap \M$ which is not nef for $X$. I. e. $h$ is $nef$
for general real K3 surfaces with non-symplectic involutions of the type
$(S,\theta)$, but it is not
$nef$ for the triplet $(X,\tau,\varphi)$ itself. This is equivalent
to have an element $\delta \in N(X)$ with $\delta^2=-2$ such that
$\delta=(\delta_1+\delta_2)/2$ where $\delta_1\in S$, $\delta_2\in
S^\perp_{N(X)}$ and $\delta_1^2=\delta_2^2=-4$ (i. e. $(X,\theta)$
is degenerate in the sense of ($\da$) as a complex surface).
Additionally, $\delta_1$  must be orthogonal to an element
$h\in S_-\cap \text{int}(\M)$ with $h^2>0$. Here $\text{int}(\M)$ denote
the interior part of $\M$, i. e. the polyhedron $\M$ without its faces.

\medskip

The condition ($\da\br$) for $(X,\tau,\varphi)$ implies the
condition $\da$ for $(X,\tau)$. Thus, condition ($\da\br$) is
stronger for $(X,\tau)$ than ($\da$), and it is not difficult to
give examples of possible $(S,\theta)$ when ($\da$) does not
imply ($\da\br$).

On the other hand, for some types $(S,\theta)$
the conditions ($\da$) implies ($\da\br$). E. g. it is obviously
true if $S=S_-$.
It is easy to see that ($\da$) implies ($\da\br$)
for all lattices $S$ of $\rk S\le 2$ above and all possible
$\theta$ for these lattices. Only for $S=U(2)$, it is possible to
have $S_-\not=S$ (this is the case when $Y=X/\{1,\theta\}$ is an ellipsoid
as a real surface). Then $S_-=\langle 4\rangle $,
$S_+=\langle -4\rangle$, and all elements of
$S$ with square $-4$ belong to $S_+$. Then ($\da$) implies
($\da\br$). But, to formulate a result about connected
components of moduli of non-degenerate real K3 surfaces with
non-symplectic involutions (Theorem \ref{genmod}) in general (for
arbitrary $S$), we have to consider the type ($\da\br$) of degeneration.

\medskip

Let $L\cong H_2(X,\bz)$ where $H_2(X,\bz)$ is the homology lattice of K3
with intersection pairing. Then $L$ is an even unimodular lattice of
signature $(3,19)$.
Let us fix a primitive embedding $S\subset L$, it is unique up to
automorphisms of $L$. Let $\Delta(S,L)^{(-4)}$ be the set of all elements
$\delta_1$ in $S$ such that $\delta_1^2=-4$ and there exists
$\delta_2\in S^\perp_L$ such that $(\delta_2)^2=-4$ and
$\delta=(\delta_1+\delta_2)/2\in L$. Then $\delta^2=-2$.
All elements $\delta_1\in \Delta(S,L)^{(-4)}$ are roots
of $S$ since $-4=(\delta_1)^2|2(\delta_1, S)$ because
$(\delta_1,  S)\subset 2\bz$. Let $W^{(-4)}(S,L)\subset O(S)$ be the group
generated by reflections in all roots from $\Delta(S,L)^{(-4)}$, 
and $W^{(-4)}(S,L)_\M$ the stabilizer subgroup of $\M$ in
$W^{(-4)}(S,L)$. The set $\Delta(S,L)^{(-4)}$ is
invariant with respect to $W^{(-2)}(S)$. It follows that the
$W^{(-4)}(S,L)_\M$ is generated by reflections $s_{\delta_1}$ in
$\delta_1\in \Delta(S,L)^{(-4)}$ such that the the hyperplane
$(\delta_1)^\perp_S\otimes \br$ in $S\otimes \br$
intersects the interior of $\M$.

The real analogy of the group $W^{(-4)}(S,L)_\M$ is {\it the subgroup}
$G\subset W^{(-4)}(S,L)_\M$ which is generated by reflections
$s_{\delta_1}$ in all elements
$\delta_1\in \Delta(S,L)^{(-4)}$ which are contained either in
$S_+$ or in $S_-$ (i. e., $s_{\delta_1}$ should commute with $\theta$)
and such that $s_{\delta_1}(\M)=\M$. Obviously,
$G\subset W^{(-4)}(S,L)_\M$.

\medskip

We consider {\it two real K3 surfaces $(X,\tau,\varphi)$ and
$(X^\prime,\tau^\prime,\varphi^\prime)$ with non-symplectic
involutions of the type $(S,\theta)$ to be
isomorphic with respect to the group $G$,}
if there exists an isomorphism $f:X\to X^\prime$ such that
$f\tau=\tau^\prime f$, $f\varphi=\varphi^\prime f$ and
$f|S\in G$.

\medskip

We correspond to a real K3 surface$(X,\tau,\varphi)$ 
with a non-symplectic involution of the type $(S,\theta)$
{\it an integral involution $(L,\varphi,S)$
with the condition $(S,\theta)$ on the sublattice} $S$
(see \cite{Nikulin83}).
Here $L=H_2(X;\bz)$, $\varphi$ is the action of $\varphi$ on
$L$ such that $\varphi|S=\theta$. The integral involution
$(L,\varphi,S)$ has the properties (e.g. see \cite{Kharlamov76} or
Sect. 3.10 in  \cite{Nikulin83}):

\medskip

(RSK3) {\it The lattice $L$ is even unimodular of signature $(3,19)$.
The lattice $L^\varphi$ is hyperbolic (it has the signature $(1,t_{(-)})$).}

\medskip

Clearly, we can consider abstract integral involutions $(L,\varphi,S)$
of the type $(S,\theta)$ satisfying (RSK3).
We consider {\it two integral involutions
$(L,\varphi,S)$ and $(L^\prime,\varphi^\prime,S)$ with the
condition $(S,\theta)$ to be isomorphic for the group $G$}
if there exists an isomorphism $\xi:L\to L^\prime$ of lattices such that
$\xi\varphi=\varphi^\prime\xi$ and $\xi|S$ belongs to the group $G$.
We {\it denote by $\text{In}(S,\theta,G)$ the set} of isomorphism classes of
integral involutions $(L,\varphi,S)$ satisfying (RSK3) with the condition
$(S,\theta)$ and the group $G$.

\medskip

By monodromy consideration, two triplets $(X,\varphi,\tau)$ and
$(X^\prime,\varphi^\prime,\tau^\prime)$ of the type $(S,\theta)$ which
belong to one connected component of moduli give isomorphic integral
involutions with the condition $(S,\theta)$. Thus, we have {\it the natural
map from the set of connected components of moduli of triplets
$(X,\varphi,\tau)$ to the set $\text{In}(S,\theta,G)$.}

Using Global Torelli Theorem for K3 surfaces \cite{PS71} and
epimorphicity of Torelli map for K3 surfaces \cite{Kulikov77},
we can prove the main result which is similar to
Theorem 3.10.1 from \cite{Nikulin79} about moduli of real polarized
K3 surfaces. This type of statements is very interesting because
it reduces the problem of description of connected components of
moduli of real algebraic varieties to a purely arithmetic problem.

\begin{theorem}
\label{genmod}
The natural map above gives the one to one correspondence between
connected components of moduli of ($\da\br$)-non-degenerate
real K3 surfaces $(X,\tau,\varphi)$ with non-symplectic
involutions of the type $(S,\theta)$ and the set $\text{In}(S,\theta,G)$ 
of isomorphism classes for $G$ of integral 
involutions $(L,\varphi,S)$ with the condition $(S,\theta)$ satisfying (RSK3). 
\end{theorem}

\Proof
Let $(L,S,\varphi)$ be an isomorphism class from
$\text{In}(S,\theta,G)$. We have the orthogonal decomposition
up to finite index
\begin{equation}
S_+\oplus S_-\oplus L_{\tau}^\varphi\oplus L_{\tau,\varphi}\subset L
\label{orthdecom}
\end{equation}

Here is the key statement (compare with
conditions ($\da$) and ($\da\br$) above).

\begin{lemma}
\label{keylemma}
Let $\delta\in L$ has $\delta^2=-2$. Assume that
$\delta=(\delta_1+\delta_2)/2$ where $\delta_1\in S$,
$\delta_2\in L_\tau$ and $(\delta_1)^2=(\delta_2)^2=-4$.

If additionally $\delta_2\in L_\tau^\varphi \cup L_{\tau,\varphi}$
and $(\delta_1,h)=0$ for $h\in S_-$ with $h^2>0$, then
$\delta_1\in S_+\cup S_-$.
\end{lemma}

\Proof
We have  $(\delta_1)^\ast=\delta_1/2=
(\delta_1)^\ast_+ + (\delta_1)^\ast_-$ where
$(\delta_1)^\ast_\pm \in (S_\pm )^\ast$. Since $(\delta_1,h)=
((\delta_1)^\ast_-,h)=0$ and the lattice
$S_-$ is hyperbolic, then $((\delta_1)^\ast_-)^2<0$
if $(\delta_1)^\ast_-\not=0$.

Assume that $\delta_2\in L_{\tau,\varphi}$. Then
$(\delta_1)^\ast_-+\delta_2/2\in (L_\varphi)^\ast$. The lattice
$L_\varphi$ is 2-elementary (since $L$ is unimodular) which means that
$2L_\varphi^\ast\subset L_\varphi$.
It follows, $\widetilde{\delta}=2(\delta_1)^\ast_- + \delta_2\in L_\varphi$
and ${\widetilde{\delta}}^2=(2(\delta_1)^\ast_-)^2-4<-4$ if
$(\delta_1)^\ast_-\not=0$. Since
$L$ is even, it follows that $\widetilde{\delta}^2\le -6$.
Thus, we have
$\delta=(\delta_1)^\ast_+ +\widetilde{\delta}/2$.
We then get $(\delta_1)^\ast_+=(\delta_1)_+/2$ where
$(\delta_1)_+\in S_+$. Since $S_+$ is negative definite, even and does
not have elements with square $-2$, then $((\delta_1)_+)^2\le -4$ if
$(\delta_1)_+^\ast\not=0$. If both elements $(\delta_1)^\ast_\pm\not=0$,
we then get $\delta^2\le (-4-6)/2<-2$ which is a contradiction.
It follows $\delta_1\in S_+\cup S_-$.

If $\delta_2\in L_\tau^\varphi$, one should replace $\varphi$ by
$\tilde{\varphi}=\tau\varphi$ and argue similarly. We have
$L^{\tilde{\varphi}}_\tau =L_{\tau,\varphi}$ and
$L_{\tau,\tilde{\varphi}}=L^\varphi_\tau$.
\QED

\medskip

We shall use the lemma later.
The rest considerations are similar to the proof of Theorem 3.10.1 in
\cite{Nikulin79} and are now standard.

Both lattices $L^\varphi_\tau$ and
$L_{\tau,\varphi}$ are hyperbolic and define hyperbolic spaces
$\La^\varphi_\tau$ and $\La_{\tau,\varphi}$. E.g. to define
$\La^\varphi_\tau$, one should choose one of half-cones
$V^+(L^\varphi_\tau)$
of the cone $V(L^\varphi_\tau)=\{x\in L^\varphi_\tau
\otimes \br\ |\ x^2>0\}$ and
consider the set $\La^\varphi_\tau=V^+(L^\varphi_\tau)/\br^+$ of rays in
$V^+(L^\varphi_\tau )$. Similarly,
$\La_{\tau,\varphi}=V^+(L_{\tau,\varphi})/\br^+$.

By global Torelli Theorem
\cite{PS71}, the triplets
\begin{equation}
\text{($\da\br$)-non-degenerate triplets}\
(X,\tau, \varphi)\ \text{of the type}\ (S,\theta)
\label{triplets}
\end{equation}
are described by their
periods. Since $S\subset N(X)$,
periods should be orthogonal to $S$.
It follows that periods of data \eqref{triplets} are pairs
\begin{equation}
(\br^+\omega_+, \br^+\omega_-)\in \La^\varphi_\tau \times
\La_{\tau,\varphi}
\label{period}
%\tag{7}
\end{equation}
(see Sect. 3.10 in \cite{Nikulin79} for details).
(We identify cohomology with homology using Poincar\'e duality.)

Let $O(L,S,G,\varphi)$ be the group of all automorphisms
$\alpha \in O(L)$ such that $\alpha(S)=S$, $\alpha|S\in G$ and
$\alpha\varphi=\varphi\alpha$.
The group $O(L,S,G,\varphi)$ acts in
$\La^\varphi_\tau \times \La_{\tau,\varphi}$ as follows.
For its element $\beta$
we set $\beta(\br^+\omega_+,\br^+\omega_-)=
(\br^+(\pm \beta(\omega_+)), \br^+(\pm \beta(\omega_-)))$ where
$\pm \beta(\omega_+)$ and $\pm \beta(\omega_-)$ are chosen to be in
$V^+(L^\varphi_\tau)$ and $V^+(L_{\tau,\varphi})$ respectively.
By global Torelli Theorem
\cite{PS71}, periods \eqref{period} define isomorphic data
\eqref{triplets}, if
and only if they belong to one orbit of $O(L,S,G,\varphi)$. These
periods we shall call equivalent.

Let $W^{(-2)}(L^\varphi_\tau)$ and $W^{(-2)}(L_{\tau,\varphi})$ be the groups
generated by all reflections in elements with square $-2$ of
lattices $L^\varphi_\tau$ and $L_{\tau,\varphi}$ respectively. They are
subgroups of $O(L,S,G,\varphi)$. Both these groups
are discrete in hyperbolic spaces $\La^\varphi_\tau$ and
$\La_{\tau,\varphi}$ and have fundamental chambers
$\Omega^\varphi_\tau \subset \La^\varphi_\tau$ and
$\Omega_{\tau,\varphi}\subset \La_{\tau,\varphi}$ respectively. They are
connected. We get all periods of data \eqref{triplets} up to equivalence
if consider periods

\begin{equation}
(\br^+\omega_+, \br^+\omega_-)\in
\Omega^\varphi_\tau \times \Omega_{\tau,\varphi}.
\label{period_omega}
%\tag{8}
\end{equation}
Further we consider only periods \eqref{period_omega}.

Consider
$(\br^+\omega_+, \br^+\omega_-)\in
\Omega^\varphi_\theta\times \Omega_{\tau,\varphi}$.
Consider the set
\begin{equation}
\Delta^{(-2)}(\br^+\omega_+, \br^+\omega_-)=
\{\delta\in L_\tau\ |\ \delta^2=-2,\
(\delta, \omega_+)=(\delta, \omega_-)=0\}.
\label{Delta_omega}
%\tag{9}
\end{equation}
A point $(\br^+\omega_+, \br^+\omega_-)$ in
$\Omega^\varphi_\theta\times \Omega_{\tau,\varphi}$  %%%\thetag{8}
does not correspond to periods of data \eqref{triplets}
if \\
$\Delta(\br^+\omega_+, \br^+\omega_-)\not=
\emptyset$. Really, assume
that $\delta\in \Delta(\br^+\omega_+, \br^+\omega_-)$. Conditions
$(\delta, \omega_+)=(\delta, \omega_-)=0$ mean that $\delta\in N(X)$.
Since $\delta^2=-2$, by Riemann-Roch Theorem for K3
surfaces $\delta$ or $-\delta$ is effective.
Since $\delta\in L_\tau$, then $\tau(\delta)=-\delta$
which is impossible for an effective class $\pm \delta$ and
a holomorphic involution $\tau$.

Thus, we should delete from the set
$\Omega^\varphi_\tau\times \Omega_{\tau,\varphi}$
%%%\thetag{8}
of possible periods of data \eqref{triplets} all points
$(\br^+\omega_+, \br^+\omega_-)$ where $\omega_+$ or $\omega_-$
belongs to a face of the polyhedron $\Omega^\varphi_\tau$ or
$\Omega_{\tau,\varphi}$ since they are orthogonal to elements of
the corresponding lattices with square $-2$ and define elements
of $\Delta(\br^+\omega_+, \br^+\omega_-)$.

Thus, all possible periods of data \eqref{triplets}
are points
\begin{equation}
(\br^+\omega_+, \br^+\omega_-)\in
(\Omega^\varphi_\tau)^\prime\times (\Omega^-_{\tau,\varphi})^\prime
\label{period_omega_prime}
%\tag{10}
\end{equation}
where $^\prime$ means that we deleted all faces of the fundamental chambers
(we can call them as {\it open fundamental chambers}). The set is connected since
it is the product of open convex sets.
For points $(\br^+\omega_+, \br^+\omega_-)$
from $(\Omega_\varphi^+)^\prime\times
(\Omega^-_{\tau,\varphi})^\prime$,    %%%\thetag{10}
elements $\delta\in \Delta^{(-2)}(\br^+\omega_+, \br^+\omega_-)$
cannot belong to $L^\varphi$ or $L_{\tau,\varphi}$.

Any $\delta\in \Delta^{(-2)}(\br^+\omega_+, \br^+\omega_-)$
has the property that $\delta=\delta_+\oplus \delta_-$ where
$\delta_+\in (L^\varphi_\tau)^\ast$ and $\delta_-\in L_{\tau,\varphi}^\ast$.
Moreover $(\delta_+)^2<0$ if $\delta_+\not=0$,
and $(\delta_-)^2<0$ if $\delta_-^2<0$ since they are elements of
hyperbolic lattices orthogonal to their elements $\omega_+$
and $\omega_-$ with positive square. Since for points
\eqref{period_omega_prime}
elements of $\Delta^{(-2)}(\br^+\omega_+, \br^+\omega_-)$
cannot belong to $L_{\tau,\varphi}$ or $L^\varphi$, we can also assume that
$\delta_+\not=0$ and $\delta_-\not=0$.

Let
\begin{equation}
\begin{split}
\Delta^{(-2)}=&\{\delta=\delta_+\oplus \delta_-\in L_\tau\ |\
\delta_+\in (L^\varphi_\tau)^\ast,\
\delta_-\in (L_{\tau,\varphi})^\ast,\\
&(\delta_+)^2<0,\ (\delta_-)^2<0,\
\delta^2=(\delta_+)^2+(\delta_-)^2=-2\}.
\end{split}
\label{Delta}
%\tag{11}
\end{equation}
Let $\delta=\delta_+\oplus \delta_-\in \Delta^{(-2)}$. Let
$\Ha_{\delta_+}=\{\br^+\omega_+\in \La^\varphi_\tau\ |\
(\omega_+, \delta_+)=0\}$ be the hyperplane in $\La^\varphi_\tau$ which is
orthogonal to $\delta_+$. Similarly, let
$\Ha_{\delta_-}=\{\br^+\omega_-\in \La_{\tau,\varphi}\ |\
(\omega_-, \delta_-)=0\}$ be the hyperplane in $\La_{\tau,\varphi}$
which is orthogonal to $\delta_-$. Points
$(\br^+\omega_+, \br^+\omega_-)\in \Ha_{\delta_+}\times \Ha_{\delta_-}$
are exactly the points such that
$\Delta(\br^+\omega_+, \br^+\omega_-)$ contains $\delta$, and they should
be deleted from the set of possible periods of data \eqref{triplets}.
Thus, the set of possible periods of data \eqref{triplets} is
\begin{equation}
(\br^+\omega_+, \br^+\omega_-)\in
\Omega_1=(\Omega^\varphi_\tau)^\prime\times
(\Omega_{\tau,\varphi})^\prime\,-\,
\left(\bigcup_{\delta \in \Delta^{(-2)}} \Ha_{\delta_+}\times \Ha_{\delta_-}
\right).
\label{Omega_1}
%\tag{12}
\end{equation}
The space $\Omega_1$ is connected because
$\Ha_{\delta_+}\times \Ha_{\delta_-}$ has codimension 2. Moreover,
the set of hyperplanes $\Ha_{\delta_+}$, $\delta\in \Delta^{(-2)}$,
is locally finite in $\La^\varphi_\tau$ because all
$\delta_+$ belong to a lattice
$(L^\varphi_\tau)^\ast$ and $0>(\delta_+)^2>-2$.
Hyperplanes $\Ha_{\delta_-}$,
$\delta \in \Delta^{(-2)}$,
have the same property in $\La_{\tau,\varphi}$.

{\it The set \eqref{Omega_1} gives periods of all triplets
$(X,\tau, \varphi)$ of the type $(S,\theta)$ without additional
condition that they are ($\da\br$)-non-degenerate (it follows that
Theorem \ref{genmod} is valid if we delete the words
($\da\br$)-non-degenerate).} Let us add the last condition.

Let us denote by
\begin{equation}
\widetilde{\Delta}^{(-4)}\subset L_\tau
\label{degenerate-4}
\end{equation}
the set of all elements $\delta_2\in L_\tau$ such that
${\delta_2}^2=-4$ and there exists
$\delta_1 \in S$ with conditions: ${\delta_1}^2=-4$,
$\delta=(\delta_1+\delta_2)/2\in L$, there exists an element
$h\in S_-\cap \text{int}(\M)$ such that $(\delta_1,h)=0$.

For $(\br^+\omega_+, \br^+\omega_-)\in \La^\varphi_\tau \times
\La_{\tau,\varphi}$, we denote by
$\Omega^{(-4)}(\br^+\omega_+, \br^+\omega_-)$ the set of all
$\delta_2\in \widetilde{\Delta}^{(-4)}$ such that
$(\omega_+,\delta_2)=(\omega_-,\delta_2)=0$. A triplet \eqref{triplets}
is ($\da\br$)-degenerate, if and only if its periods
$(\br^+\omega_+, \br^+\omega_-)$ have the non-empty set
$\Omega^{(-4)}(\br^+\omega_+, \br^+\omega_-)$. Thus, we should delete
from $\Omega_1$ in \eqref{Omega_1} all these points.

Like for $\Delta^{(-2)}$ above, the set of these points has codimension one
only if there exists
$\delta_2\in \Omega^{(-4)}(\br^+\omega_+,\br^+\omega_-)$
such that $\delta_2\in L^\varphi_\tau\cup L_{\tau,\varphi}$.
By Lemma \ref{keylemma}, then the corresponding $\delta_1$
belongs to $S_-\cup S_+$ (see definition of $\widetilde{\Delta}^{(-4)}$ in
\eqref{degenerate-4}).
Let $\delta=(\delta_1+\delta_2)/2\in L$ and
$\delta^\prime =(-\delta_1+\delta_2)/2\in L$.
Since $\delta^2=(\delta^\prime)^2=-2$, there exist reflections
$s_\delta$, $s_{\delta^\prime}$ in these elements which belong to $O(L)$.
Their composition
$s=s_\delta s_{\delta^\prime}$ acts in $S=L^\tau$ as
the reflection $s_{\delta_1}$ in $\delta_1$, and it acts in
$L_\tau$ as the reflection $s_{\delta_2}$ in $\delta_2$.
Moreover, the reflection $s_{\delta_1}\in G$ since
$\delta_1\in S_-\cup S_+$ and $(h,\delta_1)=0$ for
$h\in \text{int}(\M)$. It follows that $s\in O(L,S,G,\varphi)$.

Thus, we get all equivalent
periods of data \eqref{triplets}, if we replace $\Omega^\varphi_\tau$,
$(\Omega^\varphi_\tau)^\prime$ and $\Omega_{\tau,\varphi}$,
$(\Omega_{\tau,\varphi})^\prime$ by the corresponding fundamental
chambers (and open fundamental chambers)
$\widetilde{\Omega}^\varphi_\tau$,
$(\widetilde{\Omega}^\varphi_\tau)^\prime$ and
$\widetilde{\Omega}_{\tau,\varphi}$,
$(\widetilde{\Omega}_{\tau,\varphi})^\prime$
for the groups $W(L_\tau^\varphi)$ and
$W(L_{\tau,\varphi})$ generated by reflections in all elements with
square $-2$ and in all elements from
$\widetilde{\Delta}^{(-4)}$ which are contained in lattices
$L^\tau_\varphi$ and  $L_{\tau,\varphi}$ respectively.

\medskip

Let $\Delta^{(-4)}=\widetilde{\Delta}^{(-4)}-
(L^\varphi_\tau\cup L_{\tau,\varphi})$. If $\delta_2\in \Delta^{(-4)}$,
then we write $\delta_2=(\delta_2)_+ + (\delta_2)_-$ where
$(\delta_2)_+\in (L^\varphi_\tau)^\ast$ and
$(\delta_2)_-\in (L_{\tau,\varphi})^\ast$.

The set of all periods of data \eqref{triplets} is then
\begin{equation}
(\br^+\omega_+, \br^+\omega_-)\in
\Omega=(\widetilde{\Omega}^\tau_\varphi)^\prime\times
(\widetilde{\Omega}_{\tau,\varphi})^\prime\,-\,
\left(\bigcup_{\delta \in \Delta\cup \Delta^{(-4)}}
\Ha_{\delta_+}\times \Ha_{\delta_-}
\right)
\label{Omega}
%\tag{12}
\end{equation}
where
$(\widetilde{\Omega}^\tau_\varphi)^\prime$ and
$(\widetilde{\Omega}_{\tau,\varphi})^\prime$ are open
locally-finite polyhedra in hyperbolic spaces; the products
$\Ha_{\delta_+}\times \Ha_{\delta_-}$ have codimension two and their
set is locally finite in
$(\widetilde{\Omega}^\tau_\varphi)^\prime\times
(\widetilde{\Omega}_{\tau,\varphi})^\prime$.
It follows that $\Omega$ is connected.

By epimorphicity of Torelli map for K3
surface \cite{Kulikov77} applied to real K3 surfaces, there exist
data \eqref{triplets} with any periods from \eqref{Omega}.
Thus, by Global Torelli Theorem \cite{PS71} and epimorphicity of
Torelli map for K3 surfaces \cite{Kulikov77},
the moduli space of data \eqref{triplets} is covered
by the connected space $\Omega$ defined in \eqref{Omega}.
More exactly, the moduli space of data \eqref{triplets} is the quotient of
$\Omega$ by the stabilizer subgroup of $\Omega$ in $O(L,S,G,\varphi)$,
and it is connected.

This finishes the proof of the Theorem.

\medskip

Let us add the standard details to that. Assume
$(\br^+\omega_+, \br^+\omega_-)\in \Omega$. We can choose
$\omega_+$ and $\omega_-$ to satisfy conditions $\omega_+^2=\omega_-^2>0$.
Let $\omega=\omega_+ + i\omega_-$. Then $\omega^2=0$ and
$(\omega,\overline{\omega})>0$. We have $\tau(\omega)=-\omega$ and
$\varphi(\omega)=\overline{\omega}$. Let $H^{2,0}=\bc \omega\subset
L\otimes \bc$. We then have $\tau(H^{2,0})=H^{2,0}$ and
$\varphi(H^{2,0})=H^{0,2}=\overline{H^{2,0}}$.

We choose $h\in S_-$ such that $h$ belongs to the interior of the
fundamental chamber $\M$. This is possible because of our conditions on
$\theta$.

Let $N=\{x\in L\ |\ (\omega,x)=0\}=\{x\in L\ |\
(\omega_+,x)=(\omega_-,x)=0\}$. The lattice $N$ is hyperbolic.
Let $\delta\in N$ and $\delta^2=-2$. We claim that $(h,\delta)\not= 0$.

We can write $\delta=(\delta_1+\delta_2)/2$ where $\delta_1\in S$
and $\delta_2\in L_\tau$. We have ${\delta_2}^2<0$ if $\delta_2\not=0$
because its orthogonal complement in $L_\tau \otimes \br$ contains,
$\omega_+\perp \omega_-$ with positive squares
and $L_\tau$ has exactly two positive squares.

We have $\delta_1\not=0$. Otherwise,
$\delta\in L_\tau$ and $\Delta^{(-2)}(\br^+\omega_+, \br^+\omega_-)\not=
\emptyset$ which is not the case for points of $\Omega$ by our
considerations above.

Assume $\delta_1\not=0$ and ${\delta_1}^2\ge 0$. Then
$(h,\delta)=(h,\delta_1)\not=0$ because $S$ is hyperbolic and
$h^2>0$.

Assume $\delta_1\not=0$ and ${\delta_1}^2<0$. Since $-2=\delta^2=
({\delta_1}^2+{\delta_2}^2)/4$ and $L$ is even, then we have three
possibilities: (0) $\delta_2=0$ and
$\delta=\delta_1/2 \in S$; (1) ${\delta_1}^2=-2$ and  ${\delta_2}^2=-6$.
(2) ${\delta_1}^2={\delta_2}^2=-4$. Let us consider all of them.

(0) and (1). Since $h$ is in interior of $\M$ and $\M$ is the fundamental
chamber for the group generated by reflections in all elements of $S$
with square $-2$, we then have $(h, \delta)\not=0$ for (0), and
$(h,\delta)=(h, \delta_1)/2\not=0$ for (1).

(2) If $(h,\delta_1)=0$, we get a contradiction with our construction
of $\Omega$ in \eqref{Omega}. Thus, $(h,\delta)=(h,\delta_1)/2\not=0$.

It proves the claim: $(h,\delta)\not=0$ for any $\delta\in N$
with $\delta^2=-2$.

By epimorphicity of the period map for K3 surfaces \cite{Kulikov77},
there exists a complex K3 surface $X$ and an isomorphism
$\xi:H_2(X,\bz)\to L$ such that $\xi^{-1}(H^{2,0})=H^{2,0}(X)$ and
$\xi^{-1}(h)$ is a polarization of $X$. By Global Torelli Theorem
for K3 surfaces \cite{PS71}, the choice of $(X,\xi)$ is unique up to
isomorphisms. The pair $(X, \tau\xi)$ satisfies the same
conditions. Thus, $\xi^{-1}\tau\xi$ is induced by a unique automorphism
of $X$ which we denote by the same letter $\tau$ (for simplicity).
The $\tau$ acts as $-1$ on $H^{2,0}(X)$ and as an involution in
$H_2(X,\bz)$. It follows that $\tau$ is a non-symplectic involution of
$X$ of the type $S$.

The pair $(\overline{X},\varphi\xi)$ is
another pair with the same
properties as $(X,\xi)$  where $\overline{X}$ is the surface
$X$ with the conjugate complex structure.
It follows that $\xi^{-1}\varphi\xi$ is
induced by an anti-holomorphic involution of $X$ which we denote by
the same letter $\varphi$.

Finally we get a unique, up to isomorphisms,  a
triplet $(X,\tau,\varphi)$ from \eqref{triplets}
with periods which are equivalent to
$(\br^+\omega_+, \br^+\omega_-)\in \Omega$.

\QED

\medskip

\subsection{Connected components of moduli of non-degenerate real
right DPN-pairs}
\label{modDPN}

Below we consider a description of K3 surfaces with non-symplectic
involutions using the corresponding quotients and double coverings.
See \cite{Nikulin81}, \cite{Nikulin86}, \cite{AlexeevNikulin88},
\cite{AlexeevNikulin89} and also \cite{DIK2000} for details.

Let $(X,\tau)$ be a K3 surface with a non-symplectic involution $\tau$
and $A=X^\tau$ the set of fixed points of $\tau$. Since
$\tau^\ast(\omega_X)=-\omega_X$ and $\omega_X$ has zero divisor,
it follows that $A$ is a non-singular curve of $X$.
It follows that the surface $Y=X/\{1,\tau \}$ is non-singular, the
quotient map $\pi: X\to Y$ is a double covering ramified along a
non-singular curve $A\subset Y$. Simple considerations show that
$Y$ is Enriques surface if $A=\emptyset$, and $Y$ is a rational surface
if $A\not=\emptyset$. The non-singular curve $A\in |-2K_Y|$ and
the irregularity $q_Y=\dim H^0(Y, \Omega^1)=0$. The opposite
statement is also valid. If $Y$ is a non-singular surface,
$A\in |-2K_Y|$ a non-singular curve, and $q_Y=0$,
then the double covering of $Y$ ramified along $A$ gives a K3 surface
$X$ with a non-symplectic involution
of the double covering. A pair $(Y,A)$ with these properties is called
a {\it right DPN pair.} The surface $Y$ is called
{\it a right DPN-surface.} We follow the terminology of \cite{Nikulin86},
\cite{AlexeevNikulin88} and \cite{AlexeevNikulin89}.

We mention that a general (not necessarily right) DPN-pair is a pair
$(Y,A)$ where $Y$ is a non-singular surface with $q_Y=0$,
$A\in |-2K_Y|$ and $A$ has only ADE-singularities;
the surface $Y$ is called a DPN-surface.
Any DPN-pair gives a right one after a sequence of blow-ups in singular
points of $A$. This sequence is canonical, and studying of general
DPN-pairs can be reduced to studying right ones. Any general DPN-pair
can be obtained from a right one by a sequence of contractions of
exceptional curves of the 1st kind. Further we consider only
right DPN-pairs and right DPN-surfaces.

\medskip

Like K3 surfaces with non-symplectic involutions, {\it we call a right
DPN-pair $(Y,A)$ to be ($\da$)-degenerate if the corresponding K3 surface
with non-symplectic involution $(X,\tau)$ is ($\da$)-degenerate.}
Similarly, we call {\it a real (i.e. with an anti holomorphic involution
$\theta$) DPN-pair $(Y,A)$ to be ($\da\br$)-degenerate
if the corresponding real K3 surface with non-symplectic involution
$(X,\tau,\varphi)$ is ($\da\br$)-degenerate.} Below we discuss that.
In particular, we shall see that the property to be degenerate depends only
on the quotient $Y=X/\{1,\tau \}$. It does not depend on $A$.

\medskip

We follow notations from
\cite{AlexeevNikulin88} which are a little different from \cite{Nikulin86}.
We use the following facts about right DPN-pairs $(Y,A)$ and
the corresponding K3 surfaces with non-symplectic involutions $(X,\tau)$
of the type $S$.

Any exceptional curve $F$ of $Y$ (i. e. an irreducible algebraic curve $F$
with negative square $F^2<0$) is a non-singular rational curve
with square $-1$, $-2$ or $-4$. If $F^2=-4$ {\it (the type I),}
then $F$ is a component of $A$ and its preimage
$E=\pi^{-1}(F)$ is a component of $X^{\tau}$ with $E^2=-2$.
If $F^2=-1$ {\it (the type II)}, then $(F,A)=2$. If $A$ intersects
$F$ in two different points
{\it (the type IIa),} then $\pi^{-1}(F)=E$ where $E$ is a
non-singular rational curve with square $-2$.
If $F^2=-1$ and $A$ is tangent to $F$ with multiplicity 2
{\it (the type IIb),} then $\pi^{-1}(F)=E+\tau(E)$ where $E$ is a non-singular
rational curve with square $-2$ and $(E,\tau(E))=1$
 (remark that $(E+\tau(E))^2=-2$ like for $E^2$ of the type  IIa).
If $F^2=-2$ {\it (the type III),}
then $\pi^{-1}(F)=E+\tau(E)$ is disjoint union
of two non-singular rational curves of $X$ with square $-2$.

We have the following description of exceptional curves on $Y$
in terms of $X$.
Let $\Delta(S)^{(-2)}=\{\delta\in S\ |\ \delta^2=-2\}$. Let
$\Delta(S,X)^{(-4)}$ be the set of all $\delta_1\in S$
such that ${\delta_1}^2=-4$ and there exists $\delta_2\in S^\perp_{N(X)}$
such that ${\delta_2}^2=-4$ and $\delta=(\delta_1+\delta_2)/2\in N(X)$
(then $\delta^2=-2$). All elements of $\Delta(S,X)^{(-4)}$ are
roots of $S$ with square $-4$. Let $W(S,X)\subset O(S)$ be the group
generated by reflections in all elements of
$\Delta(S,X)=\Delta(S)^{(-2)}\cup \Delta(S,X)^{(-4)}$.
The set $\Delta(S,X)$ is invariant with respect to
$W(S,X)$. The group $W(S,X)$ is a discrete reflection group in
$V^+(S)$ and in the corresponding hyperbolic space $V^+(S)/\br^+$.
All reflections from $W(S,X)$ are reflections with respect to elements
from $\Delta(S,X)$. All
faces of a fundamental chamber of $W(S,X)$ are orthogonal to
$\Delta(S,X)$. Let $\Ka(X)\subset N(X)\otimes \br$ be the K\"ahlerian
(or the $nef$ cone) of $X$. The crucial fact is that
{\it $\M(S,X)=(S\otimes \br)\cap \Ka(X)$ is a fundamental chamber for the
group $W(S,X)$. Its faces are orthogonal exactly to preimages
$\pi^\ast(F)$ of exceptional curves $F$ of $Y$.}

Let us take a face of $\M(S,X)$. Then it is orthogonal to
$\delta\in \Delta(S,X)$. We have:

\noindent
$\delta=\pi^\ast(F)$ where $F^2=-1$, if $\delta^2=-2$ and
$\delta$ is not the class of a component of $X^\tau$;

\noindent
$2\delta=\pi^\ast(F)$ where $F^2=-4$, if $\delta^2=-2$ and
$\delta$ is the class of a component of $X^\tau$.

\noindent
$\delta=\pi^\ast(F)$ where $F^2=-2$, if $\delta^2=-4$,
i.e. $\delta\in \Delta(S,X)^{(-4)}$.

By the projection formula, we also have
$\M(S,X)=\pi^\ast(\M(Y))$ where $\M(Y)$ is the
K\"ahlerian cone of $Y$.

Using these results, we get

\begin{proposition}
\label{degendpn}
A right DPN-pair $(Y,A)$ (or the corresponding
non-symplectic involution $(X,\tau)$) is ($\da$)-degenerate,
if and only if $Y=X/\{1, \tau \}$ has an exceptional curve
with square $-2$.
\end{proposition}

\Proof
We have $\M(S,X)\subset \M$ where $\M$ is the fundamental chamber for
the group $W^{(-2)}(S)$ generated by reflections in $\Delta^{(-2)}(S)$.

If $X$ is ($\da$)-degenerate, then $\M(S,X)\not=\M$. It follows that
at least one face of $\M(S,X)$ is not orthogonal to an element $\delta\in S$
with $\delta^2=-2$. Then it is orthogonal to
$\delta \in \Delta(S,X)^{(-4)}$ where $\delta^2=-4$.
It follows that $\delta=\pi^\ast(F)$ where $F$ is an exceptional curve of
$Y$ with $F^2=-2$.

Let $F$ be an exceptional curve on $Y$ with $F^2=-2$.
Then $\delta=\pi^\ast(F))\in S$ where $\delta^2=2F^2=-4$.
It follows that at least one face of $\M(S,X)$ is orthogonal to
an element of $S$ with square $-4$. Then it cannot be orthogonal to an
element of $S$ with square $-2$. Otherwise these two elements are
$\bq$-proportional which is impossible for an integral lattice.
It follows that $\M(S,X)\not=\M$ and $Y$ is ($\da$)-degenerate.
\QED

\medskip

As an example, let us apply cited above theory of DPN-surfaces to $S=U(2)$. 
The lattice $S$ has no elements with square $-2$. Then 
$Y$ has no exceptional curves with square $-1$. Moreover, 
$\rk H^2(Y,\br)=\rk S=2$. Since $Y$ can be only a rational or 
Enriques surface, it follows that $Y$ is a relatively minimal 
rational surface $Y\cong \bff_n$ where $n=0$, $2$ 
or $4$ since a right DPN-surface may have only exceptional curves with 
squares $-1$, $-2$ or $-4$. For $S=U(2)$, the set of fixed points 
$X^\tau = C_{9}\in |-K_Y|$ is a non-singular irreducible curve of 
genus $9$ (see the formula \eqref{fixedtau1} for $X^\tau$ below). 
This is only possible for $\bff_0$ or $\bff_2$ since any curve 
$C\in |-2K_{\bff_4}|$ contains the exceptional section (rational) of 
$\bff_4$ (see Sect. \ref{F4} below). The case $Y=\bff_2$ gives ($\da$)-degenerate case because the exceptional 
section of $\bff_2$ gives an exceptional curve with square $(-2)$. 
In the  ($\da$)-non-degenerate case, $Y\cong \bp^1\times \bp^1$. 
Thus, we get 

\begin{proposition}
\label{factorU(2)}
For a non-symplectic involution $\tau$ of a K3 surface $X$ with the 
lattice $S=U(2)$ the surface $Y=X/\{1,\tau\}$ is $Y\cong \bp^1\times \bp^1$ 
in the $\da$-non-degenerate case, and $Y\cong \bff_2$ in the 
$\da$-degenerate case. 
\end{proposition}

%1020 

We have a similar to Proposition \ref{degendpn} statement over $\br$.

 \begin{proposition}
\label{degendpnoverR}
A right real DPN-pair $(Y,A,\theta)$ (or the corresponding
real non-symplectic involution $(X,\tau,\varphi)$ where
$\theta=\varphi\mod\{1,\tau\}$) is
($\da\br$)-degenerate, if and only if $Y=X/\{1, \tau\}$
has an exceptional curve $F$ with square $-2$ such that
$(F, \theta(F))\le 1$ (equivalently, either $\theta(F)=F$ or
the intersection matrix of $F$ and $\theta(F)$ is negative definite).
\end{proposition}

\Proof
If $\overline{Y}$ is the complex conjugate surface to $Y$, then
$\overline{Y}$ has the same exceptional curves as $Y$ with the complex
structures on them changed by the conjugate. It follows that for any
exceptional curve $F$ of $Y$ we have $\theta(F)=F_1$ where
$F_1$ is another exceptional curve of $Y$ (called $\theta$-conjugate)
and $\theta^\ast([F])=-[F_1]$ where $[C]$ denotes the homology class
(or cohomology class since we always identify homology and cohomology)
of a holomorphic curve $C$.

The pair $(Y,A,\theta)$ is ($\da\br$)-degenerate, if and only if
the inclusion
$$
\M(S,X)\cap (S_-\otimes\br) \subset \M\cap (S_-\otimes \br)
$$
is strict.
By considerations for the proof of Proposition \ref{degendpn},
it is equivalent to existence of an exceptional curve $F$ of $Y$ with
$F^2=-4$ such that $\pi^\ast[F]=\delta\in \Delta(S,X)^{(-4)}$ and
the inequality $(x,\delta)\ge 0$ considered in $\M\cap (S_-\otimes \br)$
gives its proper subset. Let $\delta=(\delta_- + \delta_+)/2$ where
$\delta_-\in {S_-}^\ast$ and $\delta_+\in {S_+}^\ast$.

The inequality $(x,\delta)\ge 0$
is equivalent to $(x, \delta_-)\ge 0$ for $x\in S_-\otimes\br$.
If ${\delta_-}^2\ge 0$, the inequality is trivial on $\M\cap (S_-\otimes \br)$
since
$S_-$ is hyperbolic and $\M$ has only elements with non-negative square.

Thus, we can assume that ${\delta_-}^2<0$ if the pair $(Y,A,\theta)$ is
($\da\br$)-degenerate. Let $\widetilde{F}$ be
one of components of $\pi^{-1}(F)$. We have
$$
\delta=[\widetilde{F}]+[\tau(\widetilde{F})],
$$
$$
\delta_-=[\widetilde{F}]+[\varphi(\widetilde{F})]+
[\tau(\widetilde{F})]+[\tau(\varphi(\widetilde{F})],
$$
$$
\delta_+=[\widetilde{F}]-[\varphi(\widetilde{F})]+
[\tau(\widetilde{F})]-[\tau(\varphi(\widetilde{F})].
$$
It follows, $0> {\delta_-}^2=(\pi^\ast(F)+\pi^\ast(\theta(F)))^2=
2(F+\theta(F))^2=2(F^2+\theta(F)^2+2(F, \theta(F))=
2(-4+2(F, \theta(F))$. Then $(F, \theta(F))\le 1$.

Let us prove the opposite statement. Let $F$ be an exceptional curve
of $Y$ with $F^2=-2$ and $(F, \theta(F))\le 1$. If $(F, \theta(F))<0$,
then $\theta(F)=F$ because both curves $F$ and $\theta(F)$ are
irreducible. If $(F,\theta(F))=0$ or $1$, then $\theta(F)\not=F$,
and the intersection matrix of the curves $F$, $\theta(F)$ is negative
definite. It follows that $F$ and $\theta(F)$ (where $\theta(F)=F$ for
the first case) define a face of the $nef$ cone $\M(Y)$ of $Y$:
there exists a $nef$ element $H$ of $Y$
with $H^2>0$ such that $(H,F)=(H,\theta(F))=0$ and $(H,E)>0$
for any other exceptional curve $E$ of $Y$ which is different from
$F$ and $\theta(F)$. The elements $(-\theta^\ast(H))$
and $h=H-\theta^\ast(H)$ are also $nef$. We have $\theta^\ast (h)=-h$, thus
$h$ is a real $nef$ element. Since
$\theta^\ast [F]=-[\theta(F)]$ , it follows that
$(h,F)=(H,F)-(\theta^\ast(H),F)=(H,F)-(H,\theta^\ast[F])=
(H,F)+(H,\theta(F))=0$. Similarly, $(h,\theta(F))=0$. Moreover,
$(h,E)=(H,E)+(-\theta^\ast(H),E)\ge (H,E)>0$ for any
exceptional curve $E$ of $Y$ which is
different from the curves $F$ and $\theta(F)$.
It follows that the real $nef$ element $h$ belongs to a face of
$\M(Y)$ which is orthogonal to $F$ and $\theta(F)$, and $h$ does not
belong to a face of $\M(Y)$ which is orthogonal to any other
exceptional curve $E$ of $Y$ different from $F$ and $\theta(F)$.

We have $\M(X,S)=\pi^\ast(\M(Y))$ and $\pi^\ast(h)\in \M(X,S)$ is
a real $nef$ element which belongs to one or two faces of
$\M(X,S)$ which are orthogonal to elements $\pi^\ast[F]$ and
$\pi^\ast[\theta(F)]$ with square $-4$ of $S$,
and $\pi^\ast(h)$ does not belong to any other face of $\M(X,S)$.
It follows that $\pi^\ast(h)$ is inside of $\M$ because faces of $\M$
are orthogonal to elements with square $-2$. It finishes the proof
of the statement.

\QED

\medskip

Let $(Y,\theta,A)$ be {\it a real right DPN-pair,} i. e.
$(Y,\theta)$ is a non-singular projective algebraic surface over
$\bc$ with an anti-holomorphic involution $\theta$ (i. e. a real
projective algebraic surface $Y/\br$), and $A\in |-2K_Y|$ is a
non-singular curve such that $\theta (A)=(A)$ (i. e. a real curve
$A/\br$). Later we shall
also use notation $(Y,A)/\br$ (or just $(Y,A)$) for a pair defined
over $\br$. It is natural to consider two such pairs $(Y,A,\theta)$ and
$(Y^\prime, A^\prime, \theta^\prime)$ to be isomorphic if there
exists an isomorphism 
$f:Y\to Y^\prime$ of algebraic surfaces over $\bc$
such that $f(A)=A^\prime$ and $f\theta=\theta^\prime f$. The corresponding
moduli space should parameterize isomorphism classes of such pairs.
But we consider a more delicate isomorphism relation between real
right DPN-pairs. Studying connected components of moduli for
this more delicate isomorphism relation gives an additional information
about monodromy.

There exist two real double coverings of $Y$ ramified along $A$
which give two real K3 surfaces with non-symplectic involutions
$(X,\tau,\varphi)$ and $(X,\tau,\widetilde{\varphi})$ where
$\widetilde{\varphi}=\tau \varphi$. They define the same
real right DPN-pair $(Y=X/\{1,\tau\},\,A=X^\tau,\,\theta=\varphi\mod \{1,\tau\})$.
{\it The real K3 surfaces with non-symplectic involutions
$(X,\tau,\varphi)$ and $(X,\tau,\widetilde{\varphi}=\tau\varphi)$
are called related.} A real right DPN-pair $(Y,A,\theta)$
together with a choice of one, between two, real double
coverings of $Y$ ramified along $A$ (equivalently together with a
choice of one, between two, real K3 surfaces with non-symplectic
involutions $(X,\tau,\varphi)$ such that its quotient by $\tau$ gives
$(Y,\theta,A)$) is called {\it a positive real right DPN-pair}. We shall
denote a positive DPN-pair as $(Y,\theta,A)^+$ (or $(Y,A)^+$ for
$(Y,A)/\br$). Then a {\it related DPN-pair} will be denoted by
$(Y,\theta, A)^-$. If $(Y,\theta,A)^+$ is given by $(X,\tau,\varphi)$,
then
$(Y,\theta,A)^-$ is given by
$(X,\tau, \widetilde{\varphi}=\tau\varphi)$. Obviously, considering of
real right DPN-pairs $(Y,\theta,A)$ is equivalent to considering of
sets of related positive real right DPN-pairs $\{(Y,\theta,A)^+,(Y,\theta,A)^-\}$.
 {\it An isomorphism of two positive real right DPN-pairs}
is equivalent to an isomorphism of the corresponding real K3 surfaces with
non-symplectic involutions. Thus, we can define {\it positive real right
DPN-pairs of the type $(S,\theta)$} and {\it isomorphism of
real right DPN-pairs with respect to the group $G$ defined by $(S,\theta)$}.
Obviously, an isomorphism of real non-symplectic involutions
$(X, \tau, \varphi)$ and $(X^\prime,\tau^\prime, \varphi^\prime)$
defines the corresponding isomorphism of the related non-symplectic
involutions $(X,\tau, \widetilde{\varphi})$ and $(X^\prime,\tau^\prime,
\widetilde{\varphi^\prime})$. Moreover, the type $(S,\theta)$ and the
group $G$ don't change for related real non-symplectic involutions.
Moreover, by Proposition \ref{degendpnoverR} related pairs are
$(\da\br)$-degenerate simultaneously.

Since positive real right DPN-pairs are equivalent to
real non-symplectic involutions of K3, we get from Theorem \ref{genmod} an
equivalent

\begin{theorem}
\label{gnmoddpn+}
The natural map gives the one to one correspondence between
connected components of moduli of ($\da\br$)-non-degenerate
positive real right DPN-pairs $(Y,\theta,A)^+$
of the type $(S,\theta)$ and the set $\text{In}(S,\theta,G)$
of isomorphism classes for $G$ of integral
involutions $(L,\varphi,S)$ with the condition $(S,\theta)$
satisfying (RSK3).
\end{theorem}

Integral involutions $(L,\varphi, S)$ and $(L,\widetilde{\varphi},S)$
corresponding to related pairs $(Y,\theta,A)^+$ and $(Y,\theta,A)^-$
are related as $\widetilde{\varphi}=\tau\varphi$. We remind that
$\tau$ acts as $+1$ on $S$ and as $-1$ on $S^\perp_L$. Naturally,
{\it integral involutions $(L,\varphi,S)$ and $(L,\widetilde{\varphi},S)$
are called related}. We denote the involution
$(L,\varphi,S)\to (L,\widetilde{\varphi},S)$ in $\text{In}(S,\theta,G)$
as $\tau$ too. From Theorem \ref{gnmoddpn+} we get

\begin{theorem}
\label{gnmoddpn}
The natural map gives the one to one correspondence between
connected components of moduli of ($\da\br$)-non-degenerate
real right DPN-pairs $(Y,\theta, A)$
of the type $(S,\theta)$ and the set $\text{In}(S,\theta,G)/\{1,\tau\}$
of isomorphism classes for $G$ of pairs
$
\{(L,\varphi,S), (L,\widetilde{\varphi}=\tau\varphi, S)\}
$
of related integral involutions with the condition $(S,\theta)$
satisfying (RSK3).
\end{theorem}

\subsection{General classification results about real K3 surfaces
with non-symplectic involutions and real right DPN-pairs}
\label{genclassification}

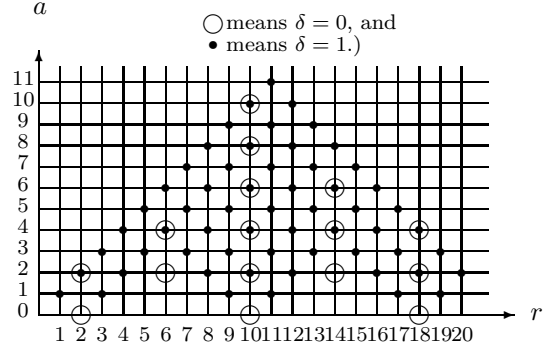
\begin{figure}
\begin{picture}(200,140)
%% All real K3
\put(66,110){\circle{7}}
\put(71,108){{\tiny means $\delta=0$, and}}
\put(66,102){\circle*{3}}
\put(71,100){{\tiny means $\delta=1$.)}}

\multiput(8,0)(8,0){20}{\line(0,1){94}}
\multiput(0,8)(0,8){11}{\line(1,0){170}}
\put(0,0){\vector(0,1){100}}
\put(0,0){\vector(1,0){180}}
\put(  6,-10){{\tiny $1$}}
\put( 14,-10){{\tiny $2$}}
\put( 22,-10){{\tiny $3$}}
\put( 30,-10){{\tiny $4$}}
\put( 38,-10){{\tiny $5$}}
\put( 46,-10){{\tiny $6$}}
\put( 54,-10){{\tiny $7$}}
\put( 62,-10){{\tiny $8$}}
\put( 70,-10){{\tiny $9$}}
\put( 76,-10){{\tiny $10$}}
\put( 84,-10){{\tiny $11$}}
\put( 92,-10){{\tiny $12$}}
\put(100,-10){{\tiny $13$}}
\put(108,-10){{\tiny $14$}}
\put(116,-10){{\tiny $15$}}
\put(124,-10){{\tiny $16$}}
\put(132,-10){{\tiny $17$}}
\put(140,-10){{\tiny $18$}}
\put(148,-10){{\tiny $19$}}
\put(156,-10){{\tiny $20$}}

\put(-8, -1){{\tiny $0$}}
\put(-8,  7){{\tiny $1$}}
\put(-8, 15){{\tiny $2$}}
\put(-8, 23){{\tiny $3$}}
\put(-8, 31){{\tiny $4$}}
\put(-8, 39){{\tiny $5$}}
\put(-8, 47){{\tiny $6$}}
\put(-8, 55){{\tiny $7$}}
\put(-8, 63){{\tiny $8$}}
\put(-8, 71){{\tiny $9$}}
\put(-10, 79){{\tiny $10$}}
\put(-10, 87){{\tiny $11$}}
\put( -2,114){{\footnotesize $a$}} %vertical
\put(186, -2){{\footnotesize $r$}} %horizontal

%% \delta_\varphi = 0
%%  ( 8r,8a)
\put( 16, 0){\circle{7}}
\put( 80, 0){\circle{7}}
\put(144, 0){\circle{7}}
\put( 16,16){\circle{7}}
\put( 48,16){\circle{7}}
\put( 80,16){\circle{7}}
\put(112,16){\circle{7}}
\put(144,16){\circle{7}}
\put( 48,32){\circle{7}}
\put( 80,32){\circle{7}}
\put(112,32){\circle{7}}
\put( 80,48){\circle{7}}
\put( 112,48){\circle{7}}
\put( 80,64){\circle{7}}
\put( 80,80){\circle{7}}    %%14data

%% \delta_\varphi = 1
%%  ( 8r,8a)
\put(  8, 8){\circle*{3}}
\put( 24, 8){\circle*{3}}
\put( 72, 8){\circle*{3}}
\put( 88, 8){\circle*{3}}
\put(136, 8){\circle*{3}}
\put(152, 8){\circle*{3}}
\put( 16,16){\circle*{3}}
\put( 32,16){\circle*{3}}
\put( 64,16){\circle*{3}}
\put( 80,16){\circle*{3}}
\put( 96,16){\circle*{3}}
\put(128,16){\circle*{3}}
\put(144,16){\circle*{3}}
\put( 24,24){\circle*{3}}
\put( 40,24){\circle*{3}}
\put( 56,24){\circle*{3}}
\put( 72,24){\circle*{3}}
\put( 88,24){\circle*{3}}
\put(104,24){\circle*{3}}
\put(120,24){\circle*{3}}
\put(136,24){\circle*{3}}
\put( 32,32){\circle*{3}}
\put( 48,32){\circle*{3}}
\put( 64,32){\circle*{3}}
\put( 80,32){\circle*{3}}
\put( 96,32){\circle*{3}}
\put(112,32){\circle*{3}}
\put(128,32){\circle*{3}}
\put( 40,40){\circle*{3}}
\put( 56,40){\circle*{3}}
\put( 72,40){\circle*{3}}
\put( 88,40){\circle*{3}}
\put(104,40){\circle*{3}}
\put(120,40){\circle*{3}}
\put( 48,48){\circle*{3}}
\put( 64,48){\circle*{3}}
\put( 80,48){\circle*{3}}
\put( 96,48){\circle*{3}}
\put(112,48){\circle*{3}}
\put( 56,56){\circle*{3}}
\put( 72,56){\circle*{3}}
\put( 88,56){\circle*{3}}
\put(104,56){\circle*{3}}
\put( 64,64){\circle*{3}}
\put( 80,64){\circle*{3}}
\put( 96,64){\circle*{3}}
\put( 72,72){\circle*{3}}
\put( 88,72){\circle*{3}}
\put( 80,80){\circle*{3}}   %%49data
\put( 88,88){\circle*{3}}
\put( 96,80){\circle*{3}}
\put( 104,72){\circle*{3}}
\put( 112,64){\circle*{3}}
\put( 120,56){\circle*{3}}
\put( 128,48){\circle*{3}}
\put( 136,40){\circle*{3}}
\put( 144,32){\circle*{3}}
\put( 144,32){\circle{7}}
\put( 152,24){\circle*{3}}
\put( 160,16){\circle*{3}}
\end{picture}

\caption{All possible $(r, a, \delta)=(r(S),a(S),\delta(S))=
(r(\tau),a(\tau),\delta(\tau))$}
\label{Sgraph}
\end{figure}

The first, 
the most important invariant of a K3 surface with a
non-symplectic involution $(X,\tau)$ 
is the isomorphism class of the lattice $S=H_2(X,\bz)^\tau$. 
{\it The lattice $S$  defines
the type of $(X,\tau)$ and of the corresponding right DPN-pair
$(Y=X/\{1,\tau\},A=X^\tau)$.} 
The lattice $S$ is any hyperbolic 
(i. e. of signature $(1,t_{(-)}(S))$ 
and of the rank $r(S)=1+t_{(-)}(S)$) even lattice 
having a primitive embedding $S\subset L$ to the even
unimodular lattice $L\cong H_2(X,\bz)$ of signature $(3,19)$
and such that there exists an involution $\tau$ of $L$ with
$L^\tau=S$. The last property is equivalent for $S$ to
be $2$-elementary, i. e. $S^\ast/S\cong (\bz/2\bz)^{a(S)}$ where
$a(S)\ge 0$ is
an  invariant of $S$. 
Another invariant $\delta(S)$ of $S$ is equal to $0$ or $1$. 
The $\delta(S)=0$
if and only if $(x^\ast)^2\in \bz$ for any $x^\ast\in S^\ast$; equivalently, 
$(z,\tau(z))\equiv 0\mod 2$ for any $z\in L$. 
(Namely, $\delta(S)$ means $\delta(\tau)$ like $\delta(\varphi)$ below.) 
The triplet $(r(S),a(S),\delta(S))$ is
the complete list of invariants of the isomorphism class of $S$.
All possible triplets $(r,a,\delta)=(r(S),a(S),\delta(S))$ were classified in 
\cite{Nikulin79}, see also \cite{Nikulin81} and \cite{Nikulin86}. 
They are presented in Figure  \ref{Sgraph} (it is the same as Figure 2 in 
\cite{Nikulin86}). 

%1380

We have 
\begin{equation}
\begin{split}
X^\tau=&\emptyset,\ \text{if}\ \ (r(S),a(S),\delta(S))=(10,10,0);\\
X^\tau=&C_1+{C_1}^\prime\  \text{if}\ (r(S),a(S),\delta(S))=(10,8,0);\\
X^\tau =& C_{g(S)}+E_1+\cdots +E_{k(S)},\ \ g(S)=(22-r(S)-a(S))/2,\ 
k(S)=(r(S)-a(S))/2.
\end{split}
\label{fixedtau}
\end{equation} 
Here $C_g$ denote a curve of genus $g$ and $E_i\cong \bp^1$. 
We have 
\begin{equation}
X^\tau \sim 0\mod2\ \text{in}\ H_2(X,\bz)
\label{fixedtau1}
\end{equation}
if and only if $\delta(S)=0$. 
The dimension of moduli of pairs $(X,\tau)$ and the corresponding
DPN-pairs $(Y,A)$ is equal to $20-r(S)$. See \cite{Nikulin81} and
\cite{Nikulin86}.

%real from here

{\it The type $(S,\theta)$ of a real K3 surface with a non-symplectic
involution $(X,\tau,\varphi)$ is defined by the isomorphism class $\theta$ of
the action of $\varphi$ on $S$.} Any $\theta$ with the conditions
{\it $S^\theta$ is negative definite and does not have elements
$x\in S^\theta$ with $x^2=-2$} can be taken as a type. Equivalently,
$(-\theta)\in O^+(S)$ and $(-\theta)$ gives an involution of a
fundamental chamber $\M$ of $W^{(-2)}(S)$.
Finding of all possible types $(S,\theta)$ is a purely arithmetic and
not difficult problem related with the automorphism group of
the hyperbolic lattice $S$ and
its subgroup $W^{(-2)}(S)$.
In this paper, we shall consider lattices $S$ with $r(S)\le 2$
when the problem of finding of possible types $(S,\theta)$ is very simple.
It is why we don't consider this problem in this paper any more. See 
\cite{DIK2000} about some results in this direction. 

\medskip

Assume that the type $(S,\theta)$ is fixed. We denote $S_+=S^\theta$
and $S_-=S_\theta$. We shall use invariants
\begin{equation}
s=\rk S,\ p=\rk S_+
\label{invariants(S,theta)}
\end{equation}
Then $S$ has the signature $(s_{(+)},s_{(-)})=(1,s-1)$ and
$S_+$ has the signature $(p_{(+)},p_{(-)})=(0,p)$.

First we calculate the group $G$ which we had uses in
Theorems \ref{genmod}, \ref{gnmoddpn+} and \ref{gnmoddpn}.
Let $\alpha\in S$ and $\alpha^2=-4$. If $\alpha\in \Delta(S,L)^{(-4)}$,
then $\alpha$ is a root of $S$, i. e. $(\alpha, S)\equiv 0\mod 2$.
Let $\langle -4 \rangle$ be a lattice generated by $\beta$ with $\beta^2=-4$.
Let $S\oplus \langle -4\rangle$ be the orthogonal sum of $S$ and
$\langle -4 \rangle$. Let $S\oplus \langle -4 \rangle \subset S(\alpha)$
be its overlattice of index two which is generated by
$(\alpha+\beta)/2$.

\begin{lemma}
\label{DeltaSL}
The element $\alpha\in \Delta(S,L)^{(-4)}$
if and only if there exists a primitive embedding $S(\alpha)\subset L$.

In particular, this is true if $r(S)+a(S)<20$.
It is also valid if $r(S)<10$. Thus under one of these conditions,
$\Delta(S,L)^{(-4)}$ consists of all roots $\alpha\in S$ with
$\alpha^2=-4$.

\end{lemma}

\Proof
The first statement is trivial.

We have $\rk S(\alpha)=r(S)+1$, the minimal number
$l(A_{S(\alpha)})$ of generators of
$A_{S(\alpha)}=S(\alpha)^\ast/S(\alpha)$ is not more
than $a(S)+1$. Thus, $\rk S(\alpha)+l(A_{S(\alpha)})\le
r(S)+1+a(S)+1<22=\rk L$.
By Corollary 1.12.3 in \cite{Nikulin79}, there exists a primitive
embedding $S(\alpha) \subset L$.

The last statement follows from the obviouse inequality $a(S)\le r(S)$.

\QED

Using results of \cite{Nikulin79} more heavily, one can give
a necessary and sufficient condition of $\alpha\in \Delta (S,L)^{(-4)}$.
We have restricted here the simple sufficient condition which is enough
for this paper. 

Let us introduce
$$
\Delta(S_\pm, S)^{(-4)}=\{\alpha\in S_\pm\ |\ \alpha^2=-4,\
(\alpha,S)\equiv 0\mod 2\} .
$$
Let $W^{(-4)}(S_\pm,S)$ be the group generated
by reflections in all elements of
$\Delta(S_\pm, S)^{(-4)}$. From Lemma \ref{DeltaSL}, we get

\begin{lemma}
\label{GroupG}
If $a(S)+r(S)<20$ (in particurlar, if $r(S)<10$),
then the group $G=W^{(-4)}(S_+,S)\times W^{(-4)}(S_-,S)$.
\end{lemma}

In general, the group $G$ is a subgroup of
$W^{(-4)}(S_+,S)\times W^{(-4)}(S_-,S)$. In Theorems \ref{genmod},
\ref{gnmoddpn+} and \ref{gnmoddpn} one can replace $G$ by
the bigger  group $W^{(-4)}(S_+,S)\times W^{(-4)}(S_-,S)$
to simplify considerations. Of course, the results can be weaker then.

\medskip

We say that {\it two integral involutions $(L,\varphi,S)$ and
$(L^\prime,\varphi^\prime, S)$ of the type $(S,\theta)$ have the
same genus with respect to the group $G$ of $(S,\theta)$}
if there exists an automorphism $\xi:S\to S$ from $G$
which can be continued to an isomorphism
$(L,\varphi,S)\otimes \br\to
(L^\prime,\varphi^\prime, S)\otimes \br$ over
$\br$, and an isomorphism
$(L,\varphi,S)\otimes \bz_p\to (L^\prime,\varphi^\prime, S)\otimes \bz_p$
over the ring $\bz_p$ of $p$-adic integers for any prime $p$.
Other speaking, the integral involutions
$(L,\varphi,S)$ and $(L^\prime,\varphi^\prime, S)$ cannot be distinguished
considering over $\br$ and modulo any integer $N$ using isomorphisms
on $S$ which belong to the group $G$.
All the genus invariants (for an arbitrary
even lattice $S$ with an involution $\theta$) were found in
\cite{Nikulin83} together with necessary and sufficient conditions
of existence. In many cases the genus has only one isomorphism class.
Then the genus invariants give complete invariants of the isomorphism class.
Below we remind these results adapting them to our case.

We assume that the integral involution $(L,\varphi,S)$ satisfies the
condition (RSK3). I. e. the lattice $L$ is even unimodular of signature
$(l_{(+)},l_{(-)})=(3,19)$. The lattice $L^\varphi$ is hyperbolic of
the signature $(t_{(+)},t_{(-)})=(1,t_{(-)})$.
Then the only real invariant of $(L,\varphi,S)$ is
\begin{equation}
r=\rk L^\varphi=1+t_{(-)}.
\label{inv-r}
\end{equation}

Below we describe {\it the genus invariants} of the integral evolution
$(L,\varphi,S)$ of the type $(S,\theta)$. To simplify notations, we
temporarily denote $L_+=L^\varphi$ and $L_-=L_\varphi$.

Since $L$ is unimodular, we have
\begin{equation}
A_{L_{\pm}} = {L_{\pm}}^*/L_{\pm}\cong L/(L_+ \oplus L_-)\cong
(\bz/2\bz)^a
\label{inv-a}
\end{equation}
where $a\ge 0$ is an integer. It is one of the most important
genus invariants.

We also have {\it genus invariants}
\begin{equation}
\delta_\varphi = \left\{
\begin{array}{cl}
0 & \mbox{if}\ (x,\varphi (x)) \equiv 0 \mod {2} \ \forall x \in L \\
1 & \mbox{otherwise}
\end{array}
\right. ,
\label{inv-deltaphi}
\end{equation}
and
\begin{equation}
\delta_{\varphi S} = \left\{
\begin{array}{cl}
0 & \mbox{if}\ (x, \varphi (x)) \equiv (x, s_{\varphi}) \mod 2 \
                                                        \forall x \in L\\
  & \mbox{for some element}\ s_{\varphi}\ \mbox{in}\ S\\
1 & \mbox{otherwise}
\end{array}
\right.  .
\label{inv-deltaSphi}
\end{equation}
If $\delta_{\varphi S}=0$, then the element $s_\varphi$ occurring in the
definition of $\delta_{\varphi S}$ is uniquely defined modulo $2S$.
It is called {\it the characteristic element of the involution $\varphi$.}

Depending on these invariants,
we divide $(L,\varphi,S)$ into the following 3 {\it types}.
\begin{quote}
{\it Type\,0:}\ \ \ \ $\delta_{\varphi S} = 0$ and $\delta_\varphi = 0$;\\
{\it Type\,Ia:} \ \   $\delta_{\varphi S} = 0$ and $\delta_\varphi = 1$;\\
{\it Type\,Ib:} \ \   $\delta_{\varphi S} = 1.$
\end{quote}

For $x_\pm \in S_\pm $ we put ,
\begin{equation}
\delta_{x_\pm } = \left\{
\begin{array}{cl}
0 & \mbox{if}\ \ (x_\pm, L_\pm ) \equiv 0 \mod 2 ,\\
1 & \mbox{otherwise}.
\end{array}
\right.
\label{deltax1}
\end{equation}
Equivalently, $\delta_{x_\pm}=0$ iff $\frac{1}{2}x_\pm\in 2L_{\pm}^\ast$.

Since $L$ is unimodular,
\begin{equation}
\delta_{x_\pm}=0 \Leftrightarrow \exists x_{\pm}^\prime\in L_{\mp}\,:\
\frac{1}{2}(x_{\pm} + x_{\pm}^\prime)\in L.
\label{deltax}
\end{equation}
The elements $x_{\pm}^\prime$ are defined by the elements $x_{\pm}$
uniquely modulo $2L_{\mp}$; this enables us, for elements
$x_+\in S_+$ and $x_-\in S_-$ for which $\delta_{x_+}=\delta_{x_-}=0$,
to define the invariant
$$
\rho_{x_+x_-}=
{1\over 2}(x_+,x_-^\prime)\mod 2=-{1\over 2}(x_+^\prime, x_-)\mod 2
\in \bz/2\bz.
$$
Thus, we get a list of genus invariants
\begin{equation}
(r,a;\delta_{x_+},\delta_{x_-}, \rho_{x_+x_-};
\delta_\varphi, \delta_{\varphi S},s_\varphi)
\label{genusinvariants}
\end{equation}
where they are defined. E. g. the element $s_\varphi\in S\mod 2S$ is
defined only if $\delta_{\varphi S}=0$;
the invariant $\rho_{x_+x_-}$ is defined only if
$\delta_{x_+}=\delta_{x_-}=0$.
From \cite{Nikulin83}, it follows that:  

{\it Two integral
involutions with the condition
$(S,\theta)$ have the same genus with respect to
$G$, if and only if the invariants \eqref{genusinvariants}
are conjugate by the group $G$ of the condition $(S,\theta)$.}

\medskip

To formulate conditions of existence of an integral involution with the 
genus invariants \eqref{genusinvariants}, we need to reformulate
the invariants \eqref{genusinvariants}.

We have the function $\delta_{\pm}: S_\pm\to \bz/2\bz$ where
$x_\pm \mapsto \delta_{x_\pm}$. We set
$$
H_\pm = \delta_{\pm}^{-1}(0)/2S_\pm\subset S_\pm/2S_\pm.
$$
The subgroups are equivalent to the invariants $\delta_{x_\pm}$.
We have a more exact range for the subgroups $H_\pm$:
$$
\Gamma_\pm=(2S)_\pm/2S_\pm \subset H_\pm\subset
2({S_\pm}^\ast\cap ({1\over 2}S_\pm))/2S_\pm \cong
({S_\pm}^\ast\cap ({1\over 2}S_\pm))/S_\pm =A_{S_\pm}^{(2)}
\subset A_{S_\pm}.
$$
Here $(2S)_\pm$ are the orthogonal projections of
$2S\subset S_+\oplus S_-$ to $S_\pm$ respectively. This projections
also give the graph $\Gamma$ of the isomorphism $\gamma$ of the groups
$\Gamma_+$ and $\Gamma_-$. The $A_{S_\pm}^{(2)}$ denote the subgroup of
$A_{S_\pm}$ generated by all elements of order two. Let
$$
H=H_+\oplus_\gamma H_-=(H_+\oplus H_-)/\Gamma.
$$
For simplicity we identify $H_\pm=H_\pm \mod \Gamma\subset H$.
We can define {\it a finite quadratic form
$q_\rho:H\to {1\over 2}\bz/2\bz$} as follows
$$
q_\rho|H_+=q_{S_+}|H_+\text{;\ }q_\rho|H_-=-q_{S_-}|H_-\text{;\ }
q_\rho(x_+,x_-)={\rho_{x_+x_-}\over 2} \mod \bz
\text{\ for\ } x_\pm\in H_\pm.
$$
The finite quadratic form $q_\rho$ is defined correctly and is equivalent to
the invariants $\rho_{x_+x_-}$. In $H_+\times \Gamma_-$ and
$\Gamma_+\times H_-$ the invariant $\rho_{x_+x_-}$
is defined by the discriminant quadratic
forms $q_{S_+}$ and $q_{S_-}$ of lattices $S_+$ and $S_-$ respectively.
We shall also denote by
$$
\rho:H_+\times H_-\to \bz/2\bz
$$
the pairing $\rho(x_+,x_-)=\rho_{x_+x_-}$. Thus
$q_\rho(x_+,x_-)=\rho(x_+,x_-)/2\mod 1$.

Moreover (see \cite{Nikulin83}), {\it the characteristic element}
$$
v=s_\varphi \in H=H_+\oplus_\gamma H_-\subset (S_+\oplus S_-)/2S,
\text{\ if\ } \delta_{\varphi S}=0.
$$
In particular, $v=0$ if $\delta_\varphi=0$. The element $v$ should be  
{\it characteristic for the quadratic form $q_\rho$ on $H$} which means 
that $q_\rho(x,v)\equiv q_\rho(x,x)\pmod 1$ for any $x\in H$. The element 
$v$ is zero, if $\delta_\varphi=0$, and $v$ is not zero, if 
$\delta_\varphi=1$.  

\medskip

Thus, {\it the genus invariants \eqref{genusinvariants}
are equivalent to the described above genus invariants}
\begin{equation}
(r,a;H_+, H_-,q_\rho;\delta_\varphi, \delta_{\varphi S},v) .
\label{geninv2}
\end{equation}

\medskip

Now we introduce some numerical invariants of the data \eqref{geninv2}
which are important for their existence.

Let $a_{M}$ be the rank over $(\bz/2\bz)$ of a
$2$-elementary group $M$. Thus, we have {\it the numerical invariants
$a_{H_\pm}$, $a_{\Gamma_\pm}=a_{\Gamma_+}=a_{\Gamma_-}$ and
$a_H=a_{H_+}+a_{H_-}-a_{\Gamma_\pm}$.}

Any finite 2-elementary quadratic form
is an orthogonal sum $\oplus$ of elementary forms.
The elementary forms are:

\noindent
$z$ on $\bz/2\bz\xi$ with $z(\xi)=0\mod 2$;

\noindent
$w$ on $\bz/2\bz\xi$ with $w(\xi)=1\mod 2$;

\noindent
$u_+(2)$ on $\bz/2\xi_1+\bz/2\bz\xi_2$ with $u_+(\xi_1)=u_+(\xi_2)=0\mod 2$,
$u_+(\xi_1,\xi_2)={1\over 2} \mod 1$, \newline $\sigma(u_+(2))\equiv 0\mod 8$;

\noindent
$v_+(2)$ on $\bz/2\xi_1+\bz/2\bz\xi_2$ with $v_+(\xi_1)=v_+(\xi_2)=1\mod 2$,
$v_+(\xi_1,\xi_2)={1\over 2} \mod 1$; \newline $\sigma(v_+(2))\equiv 4\mod 8$;

\noindent
$q_\alpha(2)$, $\alpha =\pm 1 \mod 4$, on $\bz\xi$ with
$q_\alpha(2)(\xi)=\alpha/2\mod 2$. By definition,
$\sigma(q_\alpha(2))\equiv \pm 1\mod 8$.

Any finite 2-elementary quadratic form $f$ with a non-degenerate bilinear
form is an orthogonal sum of $u_+(2)$,
$v_+(2)$ and $q_\alpha(2)$ and has the invariant $\sigma\mod 8$ which is
equal to the sum of $\sigma$ for all its summands. It has an
additional invariant $\delta$. The $\delta=0$ if $f$ is orthogonal sum of
only $u_+(2)$ and $v_+(2)$ (equivalently, $f$ takes values in $\bz$).
Otherwise, $\delta=1$. Up to isomorphisms, the $f$ is defined by its rank $k$,
$\sigma \mod 8$ and $\delta\in \{0,1\}$. Then we denote $f$ as
$f=q(k,\delta,\sigma)$.

The quadratic form $q_\rho$ on $H$ is one of forms of the rank
$a_H=a_{H_+}+a_{H_-}-a_{\Gamma_\pm}$ with the
invariants $\delta_H=0$ or $1$; $\mu_\rho=0$ or $1$:
$$
\delta_H=0:\  q_\rho=z^{k_\rho-\mu_r}
\oplus w^{\mu_\rho}\oplus v_+(2)^{\sigma_\rho/4}
\oplus u_+(2)^{(a_H-k_\rho)/2-\sigma_\rho/4},
$$
where $\mu_\rho=0$ or $1$, $\sigma_\rho\equiv 0$ or $4 \mod 8$ and
$\sigma_\rho\ge 0$,
and $\mu_\rho+\sigma_\rho/4\le 1$;
$$
\delta_H=1,\ \mu_\rho=0:\  q_\rho=z^{k_\rho}\oplus q(a_H-k_\rho,1,
\sigma_\rho),
$$
$$
\delta_H=1,\ \mu_\rho=1:\  q_\rho=z^{k_\rho-1}\oplus w\oplus
q_1(2)^{a_H-k_\rho}
$$
where $a_H>k_\rho$.

The invariant $\delta_H=0$ iff $q_\rho(x)\in \bz$ for any $x\in H$.
We similarly introduce {\it the invariants}
$\delta_{H_\pm}$ for the form $q_\rho|H_\pm$. We have
$$
\delta_H=\max\{{\delta_{H_+},\delta_{H_-}}\}.
$$

If $\delta_{\varphi S}=0$ and $v=s_\varphi \in H$
is the characteristic element of $\varphi$,
we get {\it the invariant}
$$
c_v\mod 4\ {where}\ c_v/2=q_\rho(v)\mod 2.
$$
Thus, we get invariants
\begin{equation}
(a_{H_\pm},a_{\Gamma_\pm}, \delta_{H_\pm};k_\rho,\mu_\rho,\sigma_\rho,c_v)
\label{genusinvariants3}
\end{equation}
of the data \eqref{geninv2}.

We need some more invariants of data \eqref{geninv2}.
An element $\tilde{v}\in A_{S_\pm}^{(2)}$ is called {\it characteristic} if
$q_{S_\pm}(x)\equiv q_{S_\pm}(\tilde{v},x)\mod 1$
for any $x\in A_{S_\pm}^{(2)}$. We remind that $A_{S_\pm}^{(2)}$ is
generated by all elements of order 2 in $A_{S_\pm}$.

We say that $\delta_{\varphi S_{\pm}}=0$ if $\delta_{\varphi S}=0$,
the element $v=s_\varphi\in H_\pm$ and $v$ is equal to a characteristic element
of $A_{S_\pm}^{(2)}$. Otherwise $\delta_{\varphi S_{\pm}}=1$.

Assume that  $\delta_{\varphi S_{\pm}}=0$. We consider a non-degenerate
finite quadratic form 
$$
\gamma_\pm=
\left\{ \begin{array}{ll}
0, & \text{if $\delta_\varphi=0$, equivalently $v=0$},\\
q_1(2)\oplus q_{-1}(2), &\text{if $\delta_\varphi=1$ and 
$c_v\equiv 0\mod 4$,}\\
q_1(2)^2, &\text{if $c_v\equiv 2\mod 4$,}\\
q_{\pm 1}(2), &\text{if $c_v\equiv \pm 1\mod 4$.}
\end{array}
\right.
$$
Let $v_{\gamma_\pm}$ be the characteristic element of $\gamma_\pm$ (it is
unique). We have  $\gamma_\pm(v_{\gamma_\pm})=c_v/2\mod 2$.
We denote $(q_{S_\pm})_v=(v\oplus v_{\gamma_\pm})^\perp_
{q_{S_\pm}\oplus (\mp\gamma_{\pm})}/[v \oplus v_{\gamma_\pm}]$.
There exists a unique even $2$-adic lattice $K((q_{S_\pm})_v)$ having
the discriminant quadratic form $(q_{S_\pm})_v$ and the same rank as
the form $(q_{S_\pm})_v$ (i. e. the minimal number of
generators of the $2$-group where the form is defined). We denote by
$A_q$ the group where a finite form $q$ is defined and by
$|A_q|$ its order.
{\it The invariant
$$
\varepsilon_{v_\pm}\in (\bz/2\bz)
$$
is defined by}
$$
5^{\varepsilon_{v_\pm}}\equiv {\det K((q_{S_\pm})_v)\over
|A_{S_\pm}|\cdot |A_{\gamma_\pm}|} \pmod {\pm {\bq_2^\ast}^2}
$$
where $\bq_p$ is the field of $p$-adic integers.

Finally we get invariants
\begin{equation}
(\delta_{\varphi S_+}, \delta_{\varphi S_-}, \varepsilon_{v_+}, \varepsilon_{v_-})
\label{genusinvariants4}
\end{equation}

By our construction, the finite quadratic form $q_\rho$ on $H$ has
an embedding to the discriminant quadratic form $q_{L_+}$ which
is $2$-elementary, and the characteristic element $v$ 
(if $\delta_{\varphi S}=0$) goes to the characteristic element 
of $q_{L_+}$. The following conditions from \cite{Nikulin83}
are sufficient and necessary for its existence (we numerate these condition
as in \cite{Nikulin83}). They use invariants
\eqref{genusinvariants3} and invariants $r$, $a$, $\delta_{\varphi S}$,
and $\delta_\varphi$.

\medskip

\noindent
CONDITIONS 1.8.1 (from \cite{Nikulin83}).

{\it
\noindent
{\bf Type 0} ($\delta_\varphi=\delta_{\varphi S}=0$):

\medskip

\noindent
Conditions on $H_{\pm}$.\  $\delta_{H_+}=\delta_{H_-}=0.$

\noindent
Inequality. $a\ge {a_{H_+}}+{a_{H_-}}-{a_{\Gamma_\pm}}+k_\rho$.

\noindent
Congruences. 1) $a+r\equiv 0\pmod 2$; 2) $2-r\equiv 0\pmod 4$.

\noindent
Boundary condition. If ($a = {a_{H_+}}+{a_{H_-}}-{a_{\Gamma_\pm}}+k_\rho,\
\mu_\rho=0$), then $2-r\equiv \sigma_\rho\pmod 8$.

\noindent
{\bf Type Ia} ($\delta_\varphi=1$, $\delta_{\varphi S}=0$):

\noindent
Conditions on $H$, $\rho$ and $c_v$.
1) $c_v\equiv a_{H_+}+a_{H_-}-a_{\Gamma_\pm}-k_\rho\pmod 2$.
2) If $\mu_\rho=0$, then $c_v\equiv \sigma_\rho \pmod 4$.
3) $a_{H_+}+a_{H_-}-a_{\Gamma_\pm}>0$.
4) If $\delta_{H_+}=\delta_{H_-}=0$, then $k_\rho\ge 1$.
5) If ($\delta_{H_+}=\delta_{H_-}=0$, $k_\rho=\mu_\rho=1$), then
$c_v\equiv 2\pmod 4$.

\noindent
Inequality.  $a\ge {a_{H_+}}+{a_{H_-}}-{a_{\Gamma_\pm}}+k_\rho$.

\noindent
Congruences. 1) $a+r\equiv 0\pmod 2$; 2) $2-r\equiv c_v\pmod 4$.

\noindent
Boundary condition. If ($a = {a_{H_+}}+{a_{H_-}}-{a_{\Gamma_\pm}}+k_\rho,\
\mu_\rho=0$), then $2-r\equiv \sigma_\rho\pmod 8$.

{\bf Type Ib} ($\delta_\varphi=1$, $\delta_{\varphi S}=1$):

\noindent
Inequality. $a\ge {a_{H_+}}+{a_{H_-}}-{a_{\Gamma_\pm}}+k_\rho+1$.

\noindent
Congruences. $a+r\equiv 0\pmod 2$.

\noindent
Boundary condition. 1) If ($a = {a_{H_+}}+{a_{H_-}}-{a_{\Gamma_\pm}}+
k_\rho+1,\
\mu_\rho=0$), then $2-r\equiv \sigma_\rho\pm 1\pmod 8$.
2) If ($a = {a_{H_+}}+{a_{H_-}}-{a_{\Gamma_\pm}}+k_\rho+2,\
\mu_\rho=0$), then $2-r\not\equiv \sigma_\rho+4 \pmod 8$.
}

\medskip

The following conditions from \cite{Nikulin83}
are necessary and sufficient for existence
of the lattices $K_\pm=(S_\pm)^\perp_{L_\pm}$ with the local invariants
prescribed by the data \eqref{genusinvariants} or \eqref{geninv2}.
They also use the invariants \eqref{invariants(S,theta)} and
\eqref{genusinvariants4}. The $l(A)$ denotes the minimal number
of generators of a finite Abelian group $A$. We numerate the conditions as in
\cite{Nikulin83} and drop unnecessary (for our case) conditions.

\medskip

\noindent
CONDITION 1.8.2 (from \cite{Nikulin83}).
{\it

\noindent
Conditions on the invariants $\delta_{\varphi S_\pm}$. If
$\delta_{\varphi S}=1$, then $\delta_{\varphi S_+}=\delta_{\varphi S_-}=1$.

\noindent
Inequalities.
1) $p+1\le r\le 21-s+p$;
2) $r-a\ge -2a_{H_+}+p+l(A_{S_+})$;
3) $r+a\le 2a_{H_-}+p-s-l(A_{S_-})+22$.
\medskip

\noindent
Boundary conditions.
1) If $\delta_{\varphi S_+}=0$ and
$r-a=-2a_{H_+}+p+l(A_{S_+})$, then
$2-r\equiv 4\varepsilon_{v_+}+c_v\pmod 8$.
 
2) If $\delta_{\varphi S_-}=0$ and
$r+a=2a_{H_-}+p-s-l(A_{S_-})+22$, then
$2-r\equiv 4\varepsilon_{v_-}+c_v\pmod 8$.  }

\medskip

By \cite{Nikulin83} both these conditions are necessary and sufficient
for existence of the integral involution. Thus, finally we get

\begin{theorem} (\cite{Nikulin83}). The invariants \eqref{genusinvariants} or
\eqref{geninv2} give the complete genus invariant of integral involutions
$(L,\varphi, S)$ of the type $(S,\theta)$ from the set
$\text{In}(S,\theta,G)$. The conditions 1.8.1 and 1.8.2 above
are necessary and sufficient for existence
of the involution from $\text{In}(S,\theta,G)$ with the genus
invariants.
\label{geninvtheorem}
\end{theorem}

\medskip

Obviously, the genus of the integral involution
$(L,\varphi,S)$ defines the genus of the
related involution $(L,\widetilde{\varphi}=\tau\varphi,S)$. Thus,
invariants \eqref{genusinvariants} or
\eqref{geninv2} of $(L,\varphi,S)$ define similar invariants of the
related involution. We denote by
\begin{equation}
\left(r(\varphi),a(\varphi);H(\varphi)_+, H(\varphi)_-,q_{\rho(\varphi)};
\delta_\varphi, \delta_{\varphi S},v(\varphi)\right)  
\label{geninv2phi}
\end{equation}
and
\begin{equation}
\left(r(\tau\varphi),a(\tau\varphi);
H(\tau\varphi)_+, H(\tau\varphi)_-,q_{\rho(\tau\varphi)}; 
\delta_{\tau\varphi},\delta_{\tau\varphi S},v(\tau\varphi)\right) 
\label{geninv2tauphi}
\end{equation}
the genus invariants \eqref{geninv2} for $\varphi$ and
$\widetilde{\varphi}=\tau\varphi$. Similarly we mark by $\varphi$ or
$\tau\varphi$ the equivalent genus invariants \eqref{genusinvariants}.

\begin{theorem} We have the following relations between
genus invariants of the related involutions
$(L,\varphi,S)$ and $(L,\widetilde{\varphi}=\tau\varphi,S)$:

\begin{equation}
r(\varphi)+r(\tau\varphi)=22-s+2p,
\label{related-r}
\end{equation}

\begin{equation}
a(\tau\varphi)-a(\varphi)=a(S)-2a_{H(\varphi)}+2\rk \rho(\varphi),
\label{related-a}
\end{equation}

\begin{equation}
s_\varphi+s_{\tau\varphi}\equiv s_\theta\in S\mod 2L,
\label{related-chara}
\end{equation}
where 
$\rk \rho$ means the rank of the matrix which gives 
$\rho$ in some bases of $H_+$ and $H_-$ over the field $\bz/2\bz$,  
and 
$s_\theta\in S$ is the characteristic element of $(S,\theta)$
defined by the property 
$2(x,\theta(x))\equiv (x,s_\theta)\mod 2$ 
for any $x\in S^\ast$ (the $s_\theta$ is defined $\mod 2S$).
In particular, 
\begin{equation}
\delta_{\varphi S}=\delta_{(\tau\varphi) S}.
\label{related-delta}
\end{equation}

The $H(\varphi)_\pm$ and $H(\tau\varphi)_\pm$ are orthogonal with respect to 
the discriminant bilinear form $b_{S_\pm}$ on $A_{S_\pm}$ respectively. 

Moreover,
\begin{equation}
l(A_{S_+})+a(\varphi)-2a_{H(\varphi)_+}=
l(A_{S_-})+a(\tau\varphi)-2a_{H(\tau\varphi)_-}.
\label{related-H+-}
\end{equation}

%  $$
%  l(A_{S_-})+a(\varphi)-2a_{H(\varphi)_-}=
%  l(A_{S_+})+a(\tau\varphi)-2a_{H(\tau\varphi)_+}.
%  $$
%  (This relation is equivalent to the previous one.)

\label{relinvtheorem}
\end{theorem}

\Proof
We have the orthogonal decomposition up to a finite index
\begin{equation}
S_+\oplus S_-\oplus L^\varphi_\tau\oplus L_{\tau,\varphi}\subset L.
\label{orthogonal1}
\end{equation}
We have $r(\tau\varphi)=\rk S_+ + \rk L_{\tau,\varphi}=
\rk S_+ + (\rk L- \rk L^\varphi-\rk S_-)=p+22-r(\varphi)-(s-p)$.
It follows the first relation \eqref{related-r}.

The proof of second relation \eqref{related-a} follows the proof of Lemma 3.3.2 in
\cite{NikulinSujatha}. There was shown that
$$
2^{a(\tau\varphi)}=\# A_{L^\tau_\varphi}\cdot \# A_{L^{\tau,\varphi}}
\cdot \# A_{L_\varphi}/
(2^{2\rk H(\varphi)_+}\cdot 2^{2\rk {(H(\varphi)_+}^\perp \cap H(\varphi)_-)})
$$
where ${H(\varphi)_+}^\perp\subset H(\varphi)_-$
is the orthogonal complement to $H(\varphi)_+$ 
with respect to the pairing $\rho(\varphi)$. We have
$\rk {(H(\varphi)_+}^\perp \cap H(\varphi)_-)=\rk H(\varphi)_- -
\rk \rho(\varphi )$. Moreover, $\# A_{L_\varphi}=2^{a(\varphi)}$,
$\# A_{L^\tau_\varphi}\cdot \# A_{L^{\tau,\varphi}}=
\# A_{L^\tau}\cdot 2^{\rk \Gamma_\pm}=2^{a(S)+2\rk \Gamma_\pm}$. Thus, we
get
$$
a(\tau\varphi)=a(\varphi)+a(S)+2\rk \Gamma_\pm -2\rk H(\varphi)_+
-2\rk H(\varphi)_-
+2\rk \rho(\varphi).
$$
It follows \eqref{related-a}.

We can write up $x\in L$ as 
$x={s_+}^\ast\oplus {s_-}^\ast\oplus {x_\tau^\varphi}^\ast
\oplus {x_{\tau,\varphi}}^\ast$ where 
${s_+}^\ast\in {S_+}^\ast$, 
${s_-}^\ast\in {S_-}^\ast$, ${x_\tau^\varphi}^\ast \in 
{L_\tau^\varphi}^\ast$, ${x_{\tau,\varphi}}^\ast\in 
{L_{\tau,\varphi}}^\ast$. We have 
$$
(x,\varphi(x))=({s_+}^\ast)^2- ({s_-}^\ast)^2+
({x_\tau^\varphi}^\ast)^2-({x_{\tau,\varphi}}^\ast)^2 
\equiv (x,s_\varphi)\mod 2;
$$
$$
(x,\tau\varphi(x))=({s_+}^\ast)^2- ({s_-}^\ast)^2
-({x_\tau^\varphi}^\ast)^2+({x_{\tau,\varphi}}^\ast)^2
\equiv (x,s_{\tau\varphi})\mod 2. 
$$
Taking sum, we get 
$$
(x,s_\varphi)+(x,s_{\tau\varphi})\equiv 
2\left(({s_+}^\ast)^2- ({s_-}^\ast)^2\right)\equiv (p_S(x),s_\theta)\mod 2.
$$
where $p:L\to S^\ast$ the orthogonal projection. Here $p(L)=S^\ast$ 
since $L$ is unimodular. Since $L$ is unimodular, it follows that 
$s_\varphi+s_{\tau\varphi}\equiv s_\theta \mod 2L$. It proves the third 
relation \eqref{related-chara}.

Let $u\in S_-$ has the invariant $\delta_u=0$ for $\varphi$ 
and $v\in S_-$ has the invariant $\delta_v=0$ for $\tau\varphi$. 
Then (see \cite{Nikulin83}, (7.6)) 
there exist $u_1\in L^\varphi$ and
$v_1\in L^{\tau\varphi}=L_{\tau,\varphi}$ such that
$u/2+u_1/2\in L$ and $v/2+v_1/2 \in L$. We then have
$(u/2,v/2) = (u/2+u_1/2,v/2+v_1/2) \in \bz$. Thus
$u/2 \in A_{S_-}$ and $v/2 \in A_{S_-}$ are orthogonal with respect to
the discriminant bilinear form $b_{S_-}$ on $A_{S_-}$. It proves that
the groups $H(\varphi)_-$ and $H(\tau\varphi)_-$
are orthogonal to each other. The proof for $H_+$ is similar.

By the relation (8.9) from \cite{Nikulin83}, we have 
$$
l(A_{L^{\varphi}_\tau})=l(A_{S_+})+a(\varphi)-2a_{H(\varphi)_+};
$$
$$
l(A_{L_{\tau,\varphi}})=l(A_{S_-})+a(\varphi)-2a_{H(\varphi)_-}.
$$
Similar relation for $\tau\varphi$ gives 
$$
l(A_{L^{\varphi}_\tau})=l(A_{S_-})+a(\tau\varphi)-2a_{H(\tau\varphi)_-};
$$
$$
l(A_{L_{\tau,\varphi}})=l(A_{S_+})+a(\tau\varphi)-2a_{H(\tau\varphi)_+}.
$$
It follows the last relation \eqref{related-H+-} of the theorem.
\QED 

\medskip 

{\it Theorem \ref{relinvtheorem} permits to choose one from two 
related involutions by some conditions on its invariants.} It helps to 
classify the pairs of related involutions which is important for 
a classification of the corresponding DPN-pairs. 

\medskip

One should follow \cite{Nikulin83} {\it to give the exact construction of 
the genus invariants \eqref{geninv2tauphi} of the related involution 
$\tau\varphi$.} The idea is as follows. Using construction of 
\cite{Nikulin83}, we can calculate the genus of the integral involution 
from the invariants \eqref{geninv2phi} of $\varphi$. 
From the genus, we can find the invariants \eqref{geninv2tauphi} of 
$\tau\varphi$. We don't know a simpler procedure in general. 
Below we give this construction in details. 

Let $(A_q, q)=(A_{q(\varphi)},q(\varphi))$ 
be the discriminant quadratic form of the lattice $L^\varphi$. 
It is defined by the invariants 
$(r(\varphi),a(\varphi),\delta(\varphi))$. 
Let $\xi^{(+)}:H=H_+\oplus_\gamma H_-\to A_q$ be an 
embedding of $(H,q_\rho)$ to $(A_q,q)$ such that $\xi^{(+)}(v(\varphi))$ 
is the characteristic element of $q$, if $\delta_{\varphi S}=0$, 
and $\xi^{(+)}(H)$ does not contain a characteristic element of $q$ if 
$\delta_{\varphi S}\not=0$. Such an embedding is defined 
uniquely up to automorphisms of $q$.  

We denote by $(A_{-q},-q)$ the form $(A_q,q)$ with the the opposite sign and 
by $\Delta_{q}\subset A_q\oplus A_{-q}$ the corresponding isotropic diagonal 
for $q\oplus (-q)$. 
The $\xi^{(+)}$ gives $\xi^{(-)}:(H,-q_\rho) \to (A_{-q},-q)$ where 
$\xi^{(-)}=\xi^{(+)}$ on the groups. Let 
$\Gamma(H_+)\subset A_{S_+}\oplus A_{-q}$ be the graph of 
$\xi^{(-)}|H_+$. We then get the non-degenerate quadratic form  
\begin{equation}
(A_{k_+}=\Gamma(H_+)^\perp_{q_{S_+}\oplus (-q)}/\Gamma(H_+),
k_+=-q_{S_+}\oplus q|A_{k_+})
\label{kplus}
\end{equation}
with the 
corresponding projection 
$\pi_{k_+}:\Gamma(H_+)^\perp_{q_{S_+}\oplus (-q)}\to A_{k_+}$ 
and the graph $\Gamma_{\pi_{k_+}}\subset 
A_{S_+}\oplus A_{-q}\oplus A_{k_+}$.    
The form $k_+$ can be identified with the discriminant form of the 
lattice $L^\varphi_\tau$. As a result, we get the 
quadratic form 
$$
(A_{S_+}\oplus A_{-q}\oplus A_{k_+}, 
q_{S_+}\oplus A_{-q}\oplus {k_+})
$$
with the isotropic subgroup 
$\Gamma_{(S_+,-q,k_+)}=\Gamma(H_+)\oplus \{0\}+
\Gamma_{\pi_{k_+}}$. Further we identify 
$A_{S_+}$, $A_{-q}$, $A_{k_+}$ with the corresponding subgroups 
in $A_{S_+}\oplus A_{-q}\oplus A_{k_+}$.  
By the construction, 
$\Gamma_{(S_+,-q,k_+)}^\perp= \Gamma_{(S_+,-q,k_+)}$ is the maximal 
isotropic subgroup. Moreover, 
$\Gamma_{(S_+,-q,k_+)}\cap A_{S_+}=
\Gamma_{(S_+,-q,k_+)} \cap A_{-q}=
\Gamma_{(S_+,-q,k_+)}\cap A_{k_+}=\{0\}$,  
$\Gamma_{(S_+,-q)}=\Gamma_{(S_+,-q,k_+)}\cap (A_{S_+}\oplus A_{-q})=
\Gamma(H_+)$.  
Let $\Gamma_{(S_+,k_+)}=\Gamma_{(S_+,-q,k_+)}
\cap (A_{S_+}\oplus A_{k_+})$. By our 
construction, we then have the natural identification 
\begin{equation}
\left(A_q=(\Gamma_{(S_+,k_+)})^\perp_{q_{S_+}\oplus k_+}/
\Gamma_{(S_+,k_+)},\ q=q_{S_+}\oplus k_+|A_q\right)  
\label{Aqphi}
\end{equation}
where $\Gamma_{(S_+,-q,k_+)}/\Gamma_{(S_+,k_+)}$ gives the 
graph of $\Delta_q\subset A_q\oplus A_{-q}$ above. Let \newline  
$\pi_q:(\Gamma_{(S_+,k_+)})^\perp_{q_{S_+}\oplus k_+}\to A_q$ 
be the corresponding projection.  
Let $\Gamma(H_-)\subset A_q\oplus A_{S_-}$ be the graph of 
$\xi^{(+)}|H_-$ and 
$\Gamma_{(S_+,S_-,k_+)}=\widetilde{\Gamma(H_-)}=
(\pi_q\oplus \text{id}(A_{S_-}))^{-1}(\Gamma(H_-))\subset 
A_{S_+}\oplus A_{S_-}\oplus A_{k_+}$. Let 
\begin{equation}
\left(A_{k_-}= (\Gamma_{(S_+,S_-,k_+)})^\perp_
{q_{S_+}\oplus q_{S_-}\oplus q_{k_+}}/\Gamma_{(S_+,S_-,k_+)},\ 
k_-=(-q_{S_+})\oplus (-q_{S_-})\oplus (-q_{k_+})|A_{k_-}\right).
\label{kmin}
\end{equation} 
The form $k_-$ can be identified with the discriminant form of 
the lattice $L_{\tau,\varphi}$. We remark that one gets the same form as 
\begin{equation}
(A_{k_-}=\Gamma(H_-)^\perp_{q_{S_-}\oplus q}/\Gamma(H_-),
k_-=-q_{S_-}\oplus (-q)|A_{k_-}). 
\label{kmin1}
\end{equation} 
Thus, finally we get the form 
\begin{equation}
\left(A_{S_+}\oplus A_{S_-}\oplus A_{k_+}\oplus A_{k_-}, 
q_{S_+}\oplus q_{S_-}\oplus k_+\oplus k_-\right)
\label{orthogonal2}
\end{equation}
which corresponds to the orthogonal decomposition \eqref{orthogonal1}.
Further we use this form. 
Let 
\begin{equation}
\Gamma=\Gamma_{(S_+,S_-,k_+,k_-)}\subset A_{S_+}\oplus A_{S_-}
\oplus A_{k_+}\oplus A_{k_-}
\label{Gamma}
\end{equation}
be the graph of the projection 
$\pi_{k_-}: \Gamma_{(S_+,S_-,k_+)})^\perp_
{q_{S_+}\oplus q_{S_-}\oplus q_{k_+}}\to A_{k_-}$. 
By our construction, the group 
$\Gamma\subset  A_{S_+}\oplus A_{S_-}
\oplus A_{k_+}\oplus A_{k_-}$ is the maximal isotropic subgroup: 
$\Gamma^\perp=\Gamma$. We have 
$$
\Gamma=L/(S_+\oplus S_-\oplus L^\varphi_\tau\oplus L_{\tau,\varphi})
$$
for the orthogonal decomposition \eqref{orthogonal1}. The equality  
$\Gamma^\perp=\Gamma$ is equivalent for the lattice $L$ to be unimodular.

The orthogonal decomposition \eqref{orthogonal2} and the isotropic 
subgroup $\Gamma$ (together with the corresponding real invariants) are 
equivalent to the genus of the integral involution. Thus, we have calculated 
the genus from the invariants \eqref{geninv2phi} of $\varphi$. Using 
this genus, now we can calculate the invariants \eqref{geninv2tauphi} of 
$\tau\varphi$. For $\widetilde{\varphi}=\tau\varphi$, the lattices 
$L^\varphi_\tau$ and $L_{\tau,\varphi}$ change places:  
$L^{\widetilde{\varphi}}_\tau=L_{\tau,\varphi}$ and 
$L_{\tau,\widetilde{\varphi}}=L_\tau^\varphi$.
 
Let $\alpha_i\in \{S_+,S_-,k_+,k_-\}$.  We denote  
$\Gamma_{(\alpha_1,...,\alpha_n)}=
\Gamma\cap (A_{\alpha_1}\oplus \cdots \oplus A_{\alpha_n})$ 
and by $p_{(\alpha_1,...,\alpha_n)}$ 
the corresponding orthogonal projection into the group 
$A_{\alpha_1}\oplus \cdots \oplus A_{\alpha_n}$. 

By definition of the invariants \eqref{geninv2phi} and 
\eqref{geninv2tauphi}, we have: 
\begin{equation}
H(\varphi)_+=p_{S_+}(\Gamma_{(S_+,S_-,k_-)})
\label{Hphipl}
\end{equation}
with the natural embedding $H(\varphi)_+\subset A_q$ (use \eqref{Aqphi}) since 
$H(\varphi)_+\perp \Gamma_{(S_+,k_+)}$. Similarly, we have
\begin{equation}
H(\varphi)_-=p_{S_-}(\Gamma_{(S_+,S_-,k_+)})
\label{Hphimin}
\end{equation}
with the natural embedding $H(\varphi)_-\subset A_q$ 
defined by the graph  
$\Gamma_{(S_+,S_-,k_+)}/\Gamma_{(S_+,k_+)}$. The form $q$ on 
$A_q$ defines then the pairing $\rho(\varphi)$ between  
$H(\varphi)_+$ and $H(\varphi)_-$. 

For $\tau\varphi$ we similarly get:  
\begin{equation}
H(\tau\varphi)_+=p_{S_+}(\Gamma_{(S_+,S_-,k_+)})
\label{Htauvarphipl}
\end{equation}
with the natural embedding into  
\begin{equation}
\left(A_{q(\tau\varphi)}=(\Gamma_{(S_+,k_-)})^\perp_
{q_{S_+}\oplus k_-}/
\Gamma_{(S_+,k_-)},\ q(\tau\varphi)=q_{S_+}
\oplus k_-|A_{q(\tau\varphi)}\right);    
\label{qtauphi}
\end{equation}
\begin{equation}
H(\tau\varphi)_-=p_{S_-}(\Gamma_{(S_+,S_-,k_-)})
\label{Htauphimin}
\end{equation}
with the natural embedding into $A_{q(\tau\varphi)}$ 
defined by the graph 
$\Gamma_{(S_+,S_-,k_-)}/\Gamma_{(S_+,k_-)}$. 
The form $q(\tau\varphi)$ defines then the pairing 
$\rho(\tau\varphi)$ in between 
$H(\tau\varphi)_+$, $H(\tau\varphi)_-$. 

This calculates the genus invariants \eqref{geninv2tauphi}.  
For these calculations, the trivial (from the construction above) 
and crucial observation is:   

\begin{theorem} Assume that $q=q(\varphi)=q_1\oplus q^\prime$ and 
$\xi^{(+)}$ gives an 
embedding of $q_\rho$ into $q_1\oplus \{0\}$.  Then  
the calculations above of the genus invariants 
$H(\tau\varphi)_+$, $H(\tau\varphi)_-$ and 
$\rho(\tau\varphi)$ give the same result if one replaces 
$q=q(\varphi)$ by $q_1$.
\label{calculHtauphi}
\end{theorem} 

The Theorem \ref{calculHtauphi} reduces calculations of the 
invariants $H(\tau\varphi)_\pm$ and $\rho(\tau\varphi)$ 
to some finite and simple calculations. 
 
\medskip

There is a case when these calculations are very simple 

\begin{corollary} Assume that $S_+=0$, $S_-=S$. 
Then $H(\varphi)_+=H(\tau\varphi)_+=0$, $\rho(\varphi)=\rho(\tau\varphi)=0$ 
and $H(\varphi)_+\perp H(\tau\varphi)_+$ are orthogonal complements to 
each other with respect to $q_S=q_{S_-}$.  
\label{calculHtauphi0}
\end{corollary}

\medskip

In many cases the genus invariants \eqref{geninv2}  
have a unique isomorphism class of integral involutions.   
By Theorems \ref{genmod}, \ref{gnmoddpn+} and \ref{gnmoddpn}, then    
the genus invariants define the connected component 
of moduli of real K3 surfaces with a non-symplectic involution 
or the connected component of moduli of the corresponding 
real right DPN-pairs. We shall use the invariant 
$\kappa(A)=m$ of a finite Abelian group 
$A=(\bz/2)^m\oplus (\bz/4\bz)^n$. Using results of  
\cite{Nikulin83} and \cite{Nikulin79.1}, \cite{Nikulin79}, we get the 
following statement.   

\begin{theorem}
\label{genus-isomor} 
The genus invariants \eqref{geninv2} have a unique isomorphism class from 
$\text{In}(S,\theta,G)$ of 
integral involutions $(L,\varphi,S)$ of the 
type $(S,\theta)$ satisfying (RSK3) 
(in particular, by Theorems \ref{genmod}, \ref{gnmoddpn+} and \ref{gnmoddpn}, 
the moduli of real K3 surfaces with a non-symplectic involution, 
the moduli of the corresponding positive real right DPN-pairs, 
the moduli of the corresponding real right DPN-pairs having these invariants 
are connected), 
if both conditions $(a_\pm)$ below are satisfied: 

\medskip 

$(a_+)$ an even lattices $K_+$ with the discriminant quadratic form 
$k_+$ (see \eqref{kplus}) and the signature $(1,r-p-1)$ is unique up to 
isomorphisms, and the canonical homomorphism $O(K_+)\to O(k_+)$ is 
epimorphic;

$(a_-)$ an even lattice $K_-$ with the discriminant quadratic form $k_-$ 
(see \eqref{kmin1}) and the signature $(1,21-r-s+p)$ is unique up to 
isomorphisms, and the canonical homomorphism $O(K_-)\to O(k_-)$ is 
epimorphic.

\medskip 

The condition $(a_\pm)$ is valid if the corresponding condition $(b_\pm)$ 
below satisfies: 

\medskip

$(b_+)$ either 
$r-a>-2a_{H_+}+p+l(A_{S_+})$, or 
($r-a=-2a_{H_+}+p+l(A_{S_+})$, $a\ge 2a_{H_+}-\kappa(A_{S_+})+3$), or   
($r-a=-2a_{H_+}+p+l(A_{S_+})$, $a = 2a_{H_+}-\kappa(A_{S_+})+2$, 
$\delta_{\varphi S_+}=0$); 

$(b_-)$ either 
$r+a<2a_{H_-}+p-s-l(A_{S_-})+22$, or 
($r+a=2a_{H_-}+p-s-l(A_{S_-})+22$, $a\ge 2a_{H_-}-\kappa(A_{S_-})+3$), or 
($r+a=2a_{H_-}+p-s-l(A_{S_-})+22$, $a = 2a_{H_-}-\kappa(A_{S_-})+2$, 
$\delta_{\varphi S_-}=0$).  
\label{isomclass}
\end{theorem}

\Proof 

The first statements follows from Theorem 1.3.1 in \cite{Nikulin83} 
(see also Remark 1.6.2 there). 

The lattice $K_+$ is indefinite and even, 
its discriminant form $k_+$ is a $2$-form. Moreover, the discriminant 
group $A_{k_+}\cong (\bz/2\bz)^m\oplus (\bz/4\bz)^n$ since the lattice 
$S$ is $2$-elementary and the discriminant groups $A_{S_\pm}$ are annulated 
by $4$. 
Applying Theorem 1.14.2 from \cite{Nikulin79} to this case, we see 
that $(a_+)$ is valid  
if either $\rk K_+> l(A_{k_+})$ or $\rk K_+=l(A_{k_+})$, 
but the discriminant form $k_+\cong u_+(2)\oplus k_1$, $v_+(2)\oplus k_1$.  

The relation (8.9) in \cite{Nikulin83} gives that 
$$
l(A_{k_+})=l(A_{S_+})+a-2a_{H_+}.
$$
Thus, the condition $\rk K_+>l(A_{k_+})$ means  
$r-p>l(A_{S_+})+a-2a_{H_+}$. Equivalently, $r-a>-2a_{H_+}+p+l(A_{S_+})$. 

We have $A_{k_+}\cong (\bz/2\bz)^m\oplus (\bz/4\bz)^n$. Let us calculate  
$m$. Assume \newline $A_{S_+}\cong (\bz/2\bz)^\alpha\oplus (\bz/4\bz)^\beta$ 
where $\alpha=\kappa (A_{S_+})$.    
By \eqref{kplus},  
$\# A_{k_+}=2^{m+2n}=\#A_{S_+}\cdot \# A_{q}/(\#H_+)^2=
2^{\alpha+2\beta+a-2a_{H_+}}$. It follows, 
$m+2n=\alpha+2\beta+a-2a_{H_+}$. 
Moreover, the relation 
$l(A_{k_+})=l(A_{S_+})+a-2a_{H_+}$ used above gives 
$m+n=l(A_{k_+})=l(A_{S_+})+a-2a_{H_+}=\alpha+\beta+a-2a_{H_+}$. 
It follows, $m=\alpha+a-2a_{H_+}=\kappa(A_{S_+})+a-2a_{H_+}$. 

From classification of finite quadratic forms, we get that 
$k_+\cong u_+(2)\oplus k_1$ or $v_+(2)\oplus k_2$, if 
$m\ge 3$. Moreover, it is shown in (\cite{Nikulin83}, page 116), that 
$k_+\cong u_+(2)\oplus q_1$,  $v_+(2)\oplus q_2$, if $m=2$ and 
$\delta_{\varphi S_+}=0$. Thus, cited above conditions of 
Theorem 1.14.2 from \cite{Nikulin79} are equivalent to $(b_+)$. 
Thus, $(b_+)$ implies $(a_+)$.  

The proof of $(b_-)$ implies $(a_-)$ is the same. 

\bigskip

We mention the useful formulae which follow from the proof: 
\begin{equation}
A_{k_\pm}\cong (\bz /2\bz)^{m_{\pm}+a-2a_{H_\pm}}
\oplus (\bz/4\bz)^{n_{\pm}}, 
\text{\ if\ } A_{S_\pm}\cong (\bz/2\bz)^{m_\pm}\oplus (\bz/4\bz)^{n_{\pm}}.
\label{invariantsAk}
\end{equation} 

We also mention the following useful statement which follows 
from (Theorem 3.6.3 in \cite{Nikulin79}): 

\begin{proposition} The condition  $(a_\pm)$ 
of Theorem \ref{isomclass} is valid, if $A_{k_\pm}$ 
is $2$-elementary. 
In particular, this is true (by \eqref{invariantsAk}), if the lattice 
$S_{\pm}$ is $2$-elementary.  
\label{genus-isomor2}
\end{proposition}

\medskip

At last, we mention the main geometric interpretation of the 
invariants $(r(\varphi),a(\varphi),\delta(\varphi))=(r,a,\delta)$. 
We denote by $S_g$ an oriented surface of the genus $g$. 

We have (see Theorem 3.10.6 in \cite{Nikulin79}) for the real part 
of $X_\varphi(\br)=X^\varphi$ of $X$ with the real structure defined by 
$\varphi$ the same result as for the holomorphic non-symplectic 
involution $\tau$: 

\begin{equation}
X_\varphi(\br)=
\begin{cases}
\emptyset, &\text{if $(r,a,\delta)=(10,10,0)$}\\
T_1\amalg T_1, &\text{if $(r,a,\delta)=(10,8,0)$}\\
T_g\amalg (T_0)^k, &\text{otherwise, where}\\
              &g=(22-r-a)/2,\ k=(r-a)/2 
\end{cases};
\label{realcomponents}
\end{equation}
and 
\begin{equation}
X_\varphi(\br)\sim s_\varphi \mod 2 \text{\ in\ } H_2(X,\bz).
\label{realmod2}
\end{equation}

We have the same formula for $X_{\tau\varphi}(\br)$ using invariants 
$(r(\tau\varphi), a(\tau\varphi),\delta(\tau\varphi))$ of $\tau\varphi$.

\subsection{Possible types of $(S,\theta)$ with $\rk S\le 2$}
\label{typesStheta}

By classification in \cite{Nikulin79} and also in  
\cite{Nikulin81}, \cite{Nikulin86} (see the beginning of 
Sect. \ref{genclassification} and Figure \ref{Sgraph}), 
there are the following and 
only the following possibilities for $S$ with $r(S)=\rk S\le 2$. We have  
$$
(r(S),a(S),\delta(S))=(1,1,1),\ (2,2,0),\ (2,2,1),\ (2,0,0). 
$$
We order these cases according to the natural ordering of the corresponding 
quotients $Y=X/\{1,\tau \}$, see below. 

We consider all possible $\theta$ for these cases. 

The case: $(r(S),a(S),\delta(S))=(1,1,1)$. Then $S\cong \langle 2 \rangle$. 
Since $S_-$ should be hyperbolic, $S_-=S$ and $\theta=-1$ on $S$. 
Then $S_+=\{0\}$.  
We consider this case in Sect. \ref{P2} below.  
For this case $Y=X/\{1,\tau\}=\bp^2$. 

The case $(r(S),a(S),\delta(S))=(2,2,0)$. Then $S\cong U(2)$. 
We consider this case 
in Sects. \ref{hyperboloid} and \ref{ellipsoid}. By Proposition  \ref{factorU(2)}, we have 
$Y=\bp^1\times \bp^1$ in the non-degenerate case, which we are 
considering in this paper (in the degenerate case $Y=\bff_2$).  
Let us consider possible $\theta$. Dividing 
form of $S$ by $2$, we get the unimodular lattice $U$. It follows 
that if $\rk S_-=1$, then $S_-\cong \langle 4 \rangle$ and 
$S_+\cong \langle -4 \rangle$. If $\rk S_-=2$, then $\theta = -1$ on 
$S$. Thus, we get two cases 

The case $S\cong U(2)$, the involution $\theta$ is $-1$ on $S$.
We consider this case in Sect. \ref{hyperboloid}.  
Then  $Y=X/\{1,\tau\}=\bp^1\times \bp^1$ in 
the non-degenerate case which we consider. Moreover,  $Y=X/\{1,\tau\}=\bp^1\times \bp^1$ over $\br$ is a hyperboloid, if 
$Y(\br )\not=\emptyset$.

The case $S\cong U(2)$, $S_-\cong \langle 4 \rangle$, 
$S_+\cong \langle -4 \rangle$. We consider this case in Sect. 
\ref{ellipsoid}. For this case, $Y=X/\{1,\tau\}=\bp^1\times \bp^1$ and 
$Y$ over $\br$ is an ellipsoid. Really, if $Y=\bff_2$, then any 
anti-holomorphic involution of $Y$ acts as $-1$ in $H^2(Y,\br)=\br^2$, and 
then $\theta=-1$ in $S\otimes \br$.

The case $(r(S),a(S),\delta(S))=(2,2,1)$. Then 
$S\cong \langle 2 \rangle\oplus \langle -2 \rangle$. 
Assume that $\rk S_-=1$.   
Since the lattice 
$S(2^{-1})\cong \langle 1 \rangle \oplus \langle -1 \rangle$ is 
unimodular and odd, it follows that $S_-=\langle 2 \rangle$ and 
$S_+\cong \langle -2 \rangle$. Again the lattice $S_+$ has elements with 
square $-2$ which is impossible for $\theta$. Thus, $\theta=-1$ on $S$ 
and $S_+=\{0\}$. We consider this case in Sect. \ref{F1}.
For this case $Y=X/\{1,\tau\}=\bff_1$. 

The case $(r(S),a(S),\delta(S))=(2,0,0)$. Then $S\cong U$ where 
$U=\left(\begin{array}{cc}
0&1\\
1&0\end{array}
\right)
$.
Since $S_-$ is hyperbolic, $\rk S_-=1$ or $2$. Let $\rk S_-=1$. 
Since $S$ is unimodular and even, $S_-\cong \langle 2 \rangle$ and 
$S_+\cong \langle -2 \rangle$. Then $S_+$ has elements with square $-2$ which 
is impossible for $\theta$. 
Thus, $\rk S_-=2$ and $\theta=-1$ on $S$. Then $S_-=S$ and 
$S_+=\{0\}$.  
We consider this case in Sect. \ref{F4} below.
For this case, $Y=X/\{1,\tau\}=\bff_4$. 

\bigskip

\section{Connected components of moduli of real non-singular curves of 
degree 6 in $\bp^2$}
\label{P2} 

Here we consider the case $(r(S),a(S),\delta(S))=(1,1,1)$. Then 
$S\cong \langle 2 \rangle$. The involution $\theta=-1$ on $S$. Since 
$S$ has no elements with square $-4$, the group $G$ is 
trivial. 

For this case, $Y=X/\{1,\tau\}$ is 
a rational surface with Picard number $\rk S=1$. It follows that 
$Y\cong \bp^2$. Then $|-2K_Y|$ consists of curves of degree $6$. 
Thus, $X$ is a double covering of $\bp^2$ ramified in a non-singular 
curve of degree $6$ and $\tau$ is the involution of the double covering. 
The lattice $S=\bz h$ where $h$ is the preimage of a line in $\bp^2$, 
evidently $h^2=2$.

The $\bp^2$ has only one real structure which is defined by a homogeneous 
real coordinates $(x_0:x_1:x_2)$ in $\bp^2$. 
Then $\bp^2(\br)=\br\bp^2$ is the real projective plane. 
A curve $A$ is the zero set $P=0$ of 
a real homogeneous polynomial $P(x_0,x_1,x_2)$ of degree $6$. 
A positive curve $A^+$ is equivalent 
to a half $P\ge 0$ of the $\br\bp^2$. It defines the real double covering  
$(X,\tau,\varphi)$ of $\bp^2$ ramified along $A$ by  the condition 
that $A^+$ is the image of the real part $X_\varphi(\br)$. A curve $A$ 
is defined by the polynomial $P$ up to $\lambda P$ where 
$\lambda \in \br^\ast$. A positive curve $A^+$ is defined by the polynomial 
$P$ up to $\lambda P$ where $\lambda \in \br_{++}$. Here $\br_{++}$ is 
the set of positive real numbers.  

Classification of non-singular $A$ and $A^+$ of degree 6 is well-known.  
The isotopy classification was obtained by Gudkov \cite{Gudkov69}.   
Classification of connected components of moduli was obtained in 
\cite{Nikulin79}. See also Arnol'd \cite{Arnold71} and Rokhlin 
\cite{Rokhlin78} for the geometric interpretation of invariants. 
The last classification (for positive curves) 
is equivalent to description of connected components of 
$$
\left( \left(\br^{28}-Discr\right)/\br_{++} \right)/ PGL(3,\br)
$$
where $\br^{28}$ is the space of 
real homogeneous polynomials $P(x_0,x_1,x_2)$ of degree $6$, 
and $Discr$ is the subspace of singular polynomials, 
defining singular curves $P=0$. Since both groups $\br_{++}$ and 
$PGL(3,\br)$ are connected, this classification is equivalent to the 
classification of connected components of $\br^{28}-Discr$ which is called 
{\it the rigid isotopy classification} due to V.A. Rokhlin \cite{Rokhlin78}.   

For completeness, 
we outline general arguments of Sect. \ref{modrealK3}, 
in this case. 
    
We have $S_+=0$, $S_-=S$, $H_+=0$, $\rho=0$. The invariants $s=1$, $p=0$, 
$l(A_{S_+})=0$, $l(A_{S_-})=1$.  The group $H=H_-\subset \bz h/2\bz h$. 
We get two possibilities for the group $H$.  

(1) $H=0$ (equivalently, $\delta_h=1$). Then $a_H=0$, $\delta_H=0$, $k_\rho=0$,
$\mu_\rho=0$, $\sigma_\rho\equiv 0\pmod 8$. 
If $\delta_{\varphi S}=0$, then 
$v=0$ (equivalently, $\delta_\varphi=0$), $c_v\equiv 0\pmod 4$, 
$\delta_{\varphi S_+}=0$, $\delta_{\varphi S_-}=1$, 
$\varepsilon_{v_+}\equiv 0\pmod 2$. 

(2) $H=[h]=\{0,h\}\mod 2S$ (equivalently, $\delta_h=0$).
Then $a_H=1$, $\delta_H=1$, $k_\rho=0$,
$\mu_\rho=0$, $\sigma_\rho\equiv -1\pmod 8$. 
If $\delta_{\varphi S}=0$, then $v=h\mod 2S$ 
(equivalently, $\delta_{\varphi h}=0$, see 
\cite{Nikulin79}), $c_v\equiv -1\pmod 4$, $\delta_{\varphi S_+}=1$, 
$\delta_{\varphi S_-}=0$, 
$\varepsilon_{v_-}\equiv 0\pmod 2$.  

In particular, from these relations, $a_{H_+}=0$, 
$a_H=a_{H_-}=1-\delta_h$.

Since $S_{\pm}$ are 2-elementary, by Proposition \ref{genus-isomor2}, 
Theorem \ref{genus-isomor} and results of Sect. 1, 
we get that {\it connected components of moduli of 
$(X,\tau,\varphi)$ (or $A^+$) are in one-to-one 
correspondence with the genus invariants}
\begin{equation}
(r,a,\delta_h,\delta_{\varphi S}, v) = 
(r(\varphi),a(\varphi),\delta(\varphi)_h, \delta_{\varphi S}, v(\varphi)).
\label{<2>invariants}
\end{equation}
{\it The complete list of the invariants 
(where $r,a (\geq 0) \in \bz$, $\delta_h,\delta_{\varphi S}\in \{0,1\}$, 
the $v$ is defined only if $\delta_{\varphi S}=0$) is given by 
the conditions:} 

\medskip

\noindent
{\bf Inequalities and general relations:}

\noindent
$1\le r\le 20$; $1-\delta_h+\delta_{\varphi S}\le a$;  
$a\le r$; $r+a\le 22-2\delta_h$.

\noindent
$r+a\equiv 0\pmod 2$; 

\noindent
if $\delta_h=1$ and $\delta_{\varphi S}=0$, then $v=0\pmod {2S}$ and 
$2-r\equiv 0\pmod 4$;  

\noindent 
if $\delta_h=0$ and $\delta_{\varphi S}=0$, then $v=h \pmod {2S}$ and  
$2-r\equiv -1\pmod 4$.

\noindent 
{\bf Boundary conditions:}

\noindent
If $a=1-\delta_h$ and $v=0\pmod {2S}$, then $2-r\equiv 0\pmod 8$;

\noindent 
If $a=1-\delta_h$ and $v=h\pmod {2S}$, then $2-r\equiv -1\pmod 8$;

\noindent
If $a=2-\delta_h$ and $\delta_{\varphi S}=1$, then 
$2-r\equiv -1+\delta_h\pm 1 \pmod 8$;

\noindent
If $a=3-\delta_h$ and $\delta_{\varphi S}=1$, then 
$2-r\not\equiv -1+\delta_h+4\pmod 8$.

\noindent
If $a=r$, $\delta_{\varphi S}=0$ and $v=0$, then $2-r\equiv 0\mod 8$.

\noindent
If $r+a=22-2\delta_h$, $\delta_{\varphi S}=0$ and $v=h\pmod {2S}$, 
then $2-r\equiv -1\pmod 8$. 

\medskip

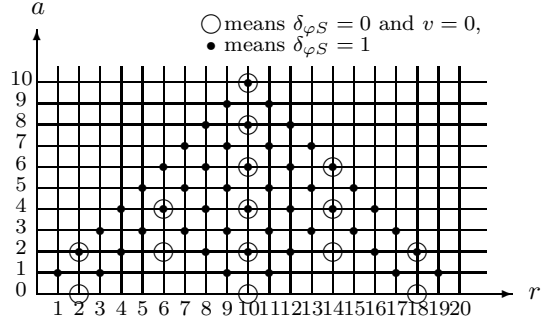
\begin{figure}

\begin{picture}(200,140)
%%non-singular degree 6  \delta_h = 1
%%\put(24,120){$\delta_h = 1$}
\put(66,102){\circle{7}}
\put(71,100){{\tiny means $\delta_{\varphi S}=0$ and $v=0$,}}
%{\circle{5}}
%{{\tiny means $\delta_{\varphi S}=0$ and $v \neq 0$}}
\put(66,94){\circle*{3}}
\put(71,92){{\tiny means $\delta_{\varphi S}=1$}}

\multiput(8,0)(8,0){20}{\line(0,1){86}}
\multiput(0,8)(0,8){10}{\line(1,0){170}}
\put(0,0){\vector(0,1){100}}
\put(0,0){\vector(1,0){180}}
\put(  6,-8){{\tiny $1$}}
\put( 14,-8){{\tiny $2$}}
\put( 22,-8){{\tiny $3$}}
\put( 30,-8){{\tiny $4$}}
\put( 38,-8){{\tiny $5$}}
\put( 46,-8){{\tiny $6$}}
\put( 54,-8){{\tiny $7$}}
\put( 62,-8){{\tiny $8$}}
\put( 70,-8){{\tiny $9$}}
\put( 76,-8){{\tiny $10$}}
\put( 84,-8){{\tiny $11$}}
\put( 92,-8){{\tiny $12$}}
\put(100,-8){{\tiny $13$}}
\put(108,-8){{\tiny $14$}}
\put(116,-8){{\tiny $15$}}
\put(124,-8){{\tiny $16$}}
\put(132,-8){{\tiny $17$}}
\put(140,-8){{\tiny $18$}}
\put(148,-8){{\tiny $19$}}
\put(156,-8){{\tiny $20$}}

\put(-8, -1){{\tiny $0$}}
\put(-8,  7){{\tiny $1$}}
\put(-8, 15){{\tiny $2$}}
\put(-8, 23){{\tiny $3$}}
\put(-8, 31){{\tiny $4$}}
\put(-8, 39){{\tiny $5$}}
\put(-8, 47){{\tiny $6$}}
\put(-8, 55){{\tiny $7$}}
\put(-8, 63){{\tiny $8$}}
\put(-8, 71){{\tiny $9$}}
\put(-10, 79){{\tiny $10$}}

\put( -2,106){{\footnotesize $a$}} %vertical
\put(186, -2){{\footnotesize $r$}} %horizontal

%% Type 0  \circle{7}
%%  ( 8r,8a)
\put( 16, 0){\circle{7}}
\put( 80, 0){\circle{7}}
\put(144, 0){\circle{7}}
\put( 16,16){\circle{7}}
\put( 48,16){\circle{7}}
\put( 80,16){\circle{7}}
\put(112,16){\circle{7}}
\put(144,16){\circle{7}}
\put( 48,32){\circle{7}}
\put( 80,32){\circle{7}}
\put(112,32){\circle{7}}
\put( 80,48){\circle{7}}
\put(112,48){\circle{7}}
\put( 80,64){\circle{7}}
\put( 80,80){\circle{7}}

%% Type Ib  \circle*{3}
%%  ( 8r,8a)
\put(  8, 8){\circle*{3}}
\put( 24, 8){\circle*{3}}
\put( 72, 8){\circle*{3}}
\put( 88, 8){\circle*{3}}
\put(136, 8){\circle*{3}}
\put(152, 8){\circle*{3}}
\put( 16,16){\circle*{3}}
\put( 32,16){\circle*{3}}
\put( 64,16){\circle*{3}}
\put( 80,16){\circle*{3}}
\put( 96,16){\circle*{3}}
\put(128,16){\circle*{3}}
\put(144,16){\circle*{3}}
\put( 24,24){\circle*{3}}
\put( 40,24){\circle*{3}}
\put( 56,24){\circle*{3}}
\put( 72,24){\circle*{3}}
\put( 88,24){\circle*{3}}
\put(104,24){\circle*{3}}
\put(120,24){\circle*{3}}
\put(136,24){\circle*{3}}
\put( 32,32){\circle*{3}}
\put( 48,32){\circle*{3}}
\put( 64,32){\circle*{3}}
\put( 80,32){\circle*{3}}
\put( 96,32){\circle*{3}}
\put(112,32){\circle*{3}}
\put(128,32){\circle*{3}}
\put( 40,40){\circle*{3}}
\put( 56,40){\circle*{3}}
\put( 72,40){\circle*{3}}
\put( 88,40){\circle*{3}}
\put(104,40){\circle*{3}}
\put(120,40){\circle*{3}}
\put( 48,48){\circle*{3}}
\put( 64,48){\circle*{3}}
\put( 80,48){\circle*{3}}
\put( 96,48){\circle*{3}}
\put(112,48){\circle*{3}}
\put( 56,56){\circle*{3}}
\put( 72,56){\circle*{3}}
\put( 88,56){\circle*{3}}
\put(104,56){\circle*{3}}
\put( 64,64){\circle*{3}}
\put( 80,64){\circle*{3}}
\put( 96,64){\circle*{3}}
\put( 72,72){\circle*{3}}
\put( 88,72){\circle*{3}}
\put( 80,80){\circle*{3}}
\end{picture}

\caption{$\bp^2$: All possible $(r,a,\delta_{\varphi S},v)$ 
with $H=0$ (i. e. $\delta_h=1$)}
\label{6_delta_h=1}
\end{figure}

\begin{figure}

\begin{picture}(200,140)
%%non-singular degree 6  \delta_h = 0
%%\put(24,122){$\delta_h = 0$}
%{\circle{7}}
%{{\tiny means $\delta_{\varphi S}=0$ and $v=0$,}}
\put(66,112){\circle{7}}
\put(71,110){{\tiny means $\delta_{\varphi S}=0$ and $v=h$, }}
\put(66,104){\circle*{3}}
\put(71,102){{\tiny means $\delta_{\varphi S}=1$}}
\multiput(8,0)(8,0){20}{\line(0,1){94}}
\multiput(0,8)(0,8){11}{\line(1,0){170}}
\put(0,0){\vector(0,1){108}}
\put(0,0){\vector(1,0){180}}
\put(  6,-8){{\tiny $1$}}
\put( 14,-8){{\tiny $2$}}
\put( 22,-8){{\tiny $3$}}
\put( 30,-8){{\tiny $4$}}
\put( 38,-8){{\tiny $5$}}
\put( 46,-8){{\tiny $6$}}
\put( 54,-8){{\tiny $7$}}
\put( 62,-8){{\tiny $8$}}
\put( 70,-8){{\tiny $9$}}
\put( 76,-8){{\tiny $10$}}
\put( 84,-8){{\tiny $11$}}
\put( 92,-8){{\tiny $12$}}
\put(100,-8){{\tiny $13$}}
\put(108,-8){{\tiny $14$}}
\put(116,-8){{\tiny $15$}}
\put(124,-8){{\tiny $16$}}
\put(132,-8){{\tiny $17$}}
\put(140,-8){{\tiny $18$}}
\put(148,-8){{\tiny $19$}}
\put(156,-8){{\tiny $20$}}

\put(-8, -1){{\tiny $0$}}
\put(-8,  7){{\tiny $1$}}
\put(-8, 15){{\tiny $2$}}
\put(-8, 23){{\tiny $3$}}
\put(-8, 31){{\tiny $4$}}
\put(-8, 39){{\tiny $5$}}
\put(-8, 47){{\tiny $6$}}
\put(-8, 55){{\tiny $7$}}
\put(-8, 63){{\tiny $8$}}
\put(-8, 71){{\tiny $9$}}
\put(-10, 79){{\tiny $10$}}
\put(-10, 87){{\tiny $11$}}

\put( -2,114){{\footnotesize $a$}} %vertical
\put(186, -2){{\footnotesize $r$}} %horizontal

%% Type Ia  \circle{5}
%%  ( 8r,8a)
\put(  24, 8){\circle{7}}
\put( 88, 8){\circle{7}}
\put(152, 8){\circle{7}}
\put( 24,24){\circle{7}}
\put( 56,24){\circle{7}}
\put( 88,24){\circle{7}}
\put( 120,24){\circle{7}}
\put(152,24){\circle{7}}
\put( 56,40){\circle{7}}
\put( 88,40){\circle{7}}
\put(120,40){\circle{7}}
\put( 56,56){\circle{7}}
\put(88,56){\circle{7}}
\put(88,72){\circle{7}}
\put(88,88){\circle{7}}

%% Type Ib  \circle*{3}
%%  ( 8r,8a)
\put( 16,16){\circle*{3}}
\put( 32,16){\circle*{3}}
\put( 80,16){\circle*{3}}
\put( 96,16){\circle*{3}}
\put(144,16){\circle*{3}}
\put(160,16){\circle*{3}}
\put( 24,24){\circle*{3}}
\put( 40,24){\circle*{3}}
\put( 72,24){\circle*{3}}
\put( 88,24){\circle*{3}}
\put(104,24){\circle*{3}}
\put(136,24){\circle*{3}}
\put(152,24){\circle*{3}}
\put( 32,32){\circle*{3}}
\put( 48,32){\circle*{3}}
\put( 64,32){\circle*{3}}
\put( 80,32){\circle*{3}}
\put( 96,32){\circle*{3}}
\put(112,32){\circle*{3}}
\put(128,32){\circle*{3}}
\put(144,32){\circle*{3}}
\put( 40,40){\circle*{3}}
\put( 56,40){\circle*{3}}
\put( 72,40){\circle*{3}}
\put( 88,40){\circle*{3}}
\put(104,40){\circle*{3}}
\put(120,40){\circle*{3}}
\put(136,40){\circle*{3}}
\put( 48,48){\circle*{3}}
\put( 64,48){\circle*{3}}
\put( 80,48){\circle*{3}}
\put( 96,48){\circle*{3}}
\put(112,48){\circle*{3}}
\put(128,48){\circle*{3}}
\put( 56,56){\circle*{3}}
\put( 72,56){\circle*{3}}
\put( 88,56){\circle*{3}}
\put(104,56){\circle*{3}}
\put(120,56){\circle*{3}}
\put( 64,64){\circle*{3}}
\put( 80,64){\circle*{3}}
\put( 96,64){\circle*{3}}
\put(112,64){\circle*{3}}
\put( 72,72){\circle*{3}}
\put( 88,72){\circle*{3}}
\put(104,72){\circle*{3}}
\put( 80,80){\circle*{3}}
\put( 96,80){\circle*{3}}
\put( 88,88){\circle*{3}}
\end{picture}

\caption{$\bp^2$: All possible $(r,a,\delta_{\varphi S}, v)$ with $H=[h]$ 
(i. e. $\delta_h=0$)}
\label{6_delta_h=0}
\end{figure}
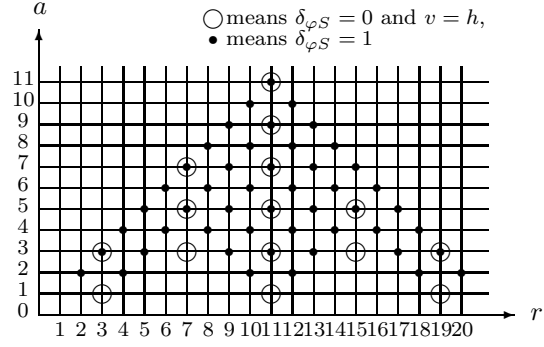

All possible invariants \eqref{<2>invariants} satisfying 
these conditions are given in Figures \ref{6_delta_h=1} and 
\ref{6_delta_h=0}. 

\medskip  

The relations between related involutions are:  
\begin{equation}
r(\varphi)+r(\tau\varphi)=21;\  
a(\tau\varphi)-a(\varphi)=2\delta(\varphi)_h-1;\ 
\delta(\varphi)_h+\delta(\tau\varphi)_h=1;  
\label{<2>related1}
\end{equation}
\begin{equation}
\delta_{\varphi S}=\delta_{\tau\varphi S};\   
v(\varphi)+v(\tau\varphi)\equiv h\pmod {2S}, \text{\ if\ }  
\delta_{\varphi S}=\delta_{\tau\varphi S}=0. 
\label{<2>related2}
\end{equation}

Because of the relation $\delta(\varphi)_h+\delta(\tau\varphi)_h=1$, 
it is sufficient to describe only involutions with $\delta(\varphi)_h=0$ 
or with $\delta(\varphi)_h=1$. 

It follows that  
there are $49$ connected components of moduli of positive curves 
with $\delta_{\varphi S}=1$, and $15$ with 
$\delta_{\varphi S}=0$ (for the fixed $\delta_h=1$ or $0$). 
Thus, there are $128$ connected components of moduli of 
positive curves (or involutions $(X,\tau, \varphi)$) and 
64 connected components of moduli of 
real non-singular curves of degree 6 on the real projective plane. 

\bigskip

According to Arnol'd \cite{Arnold71}, the invariant $\delta_h=0$, if 
$A^+\not=\emptyset$ and the positive curve $A^+$ is orientable. 
Otherwise, $\delta_h=1$.  
Using this, Bezout Theorem and \eqref{realcomponents}, \eqref{realmod2}, 
\eqref{<2>related1}, \eqref{<2>related2},  
one can draw the corresponding pictures.  
It is enough to draw them for $\delta_h=0$ when the positive curve is 
orientable. We have 

\medskip 

\noindent
{\it If $(\delta_{h};r,a,\delta_{\varphi S})=(0;11,11,0)$, 
the positive curve is the whole $\br\bp^2$.} 

\medskip

\begin{figure} 
\centerline{\includegraphics[width=4cm]{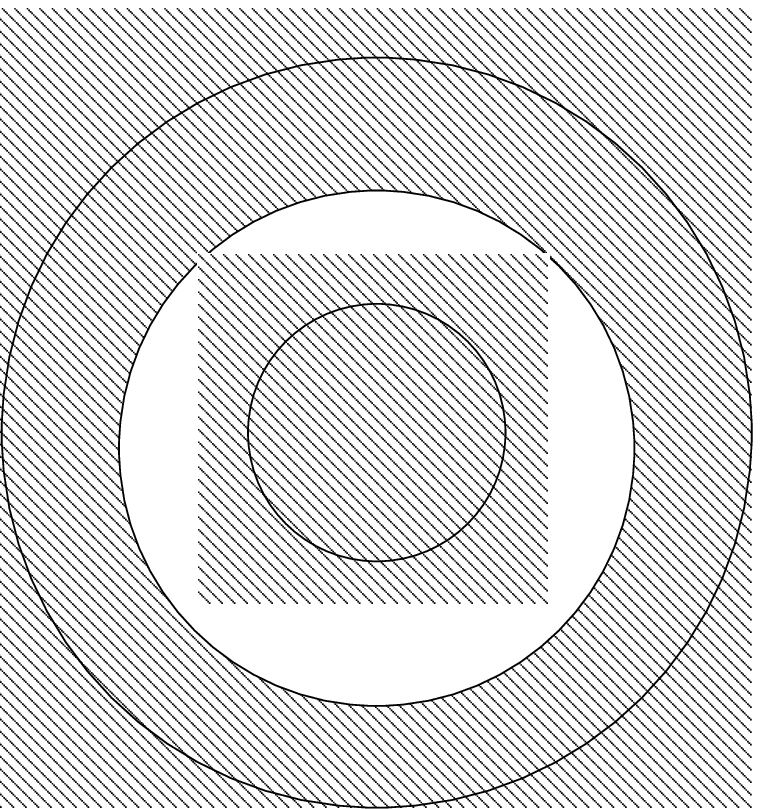}} 
\caption{$\bp^2$: $A^+$ for $H=[h]$ and 
$(r,a,\delta_{\varphi S},v)=(11,9,0,h)$}
\label{<2>graph1}
\end{figure}

\noindent 
{\it If $(\delta_h;r,a,\delta_{\varphi S})=(0;11,9,0)$ the 
positive curve is drawn in Figure \ref{<2>graph1}.}

\medskip
 
\begin{figure} 
\centerline{\includegraphics[width=6cm]{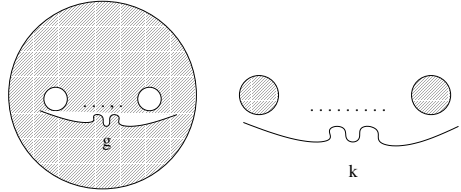}} 
\caption{$\bp^2$: $A^+$ for $H=[h]$ and 
$(r,a,\delta_{\varphi S},v)\not=(11,9,0,h)$, $(11,11,0,h)$.}
\label{<2>graph2}
\end{figure}

\noindent
{\it If $\delta_h=0$ and 
$(r,a,\delta_{\varphi S})$ is different from $(11,11,0)$ 
and $(11,9,0)$ the positive curve is drawn in  
Figure \ref{<2>graph2} where $g=(22-r-a)/2$ and $k=(r-a)/2$.}

\medskip

\noindent   
{\it The invariant $\delta_{\varphi S}=0$, if and only if the real curve 
$A(\br)$ divides $A(\bc )$.} 
See \cite{Arnold71} and \cite{Rokhlin78}. 

More exactly, by \cite{Arnold71}  
the invariant $\delta_{\varphi S}=0$, if the curve $A(\br)$ divides $A(\bc )$. 
By the results above, connected component of moduli is defined by 
the isotopy type and by the invariant $\delta_{\varphi S}$. Thus, to prove 
the opposite statement, it is sufficient to construct a dividing curve 
on the plane for each isotopy type with $\delta_{\varphi S}=0$. This 
construction is explained in \cite{Rokhlin78}.

%\begin{figure}
%\begin{picture}(200,100)
%\put( 60,40){\oval(95,65)}
%\put( 35,40){\circle{15}}
%\put( 85,40){\circle{15}}
%\multiput(50,40)(5,0){5}{\circle*{1.2}}
%\put(34,30){$\underbrace{\ \ \ \ \ \ \ 
%\ \ \ \ }_{(r-a)/2\ \mbox{{\tiny ovals}}%}$}
%
%\put(126,40){\circle{15}}
%\put(176,40){\circle{15}}
%\multiput(142,40)(5,0){5}{\circle*{1.2}}
%\put(117,30){$\underbrace{\ \ \ \ \ \ \ \ \ \ \ 
%}_{(20-r-a)/2\ \mbox{{\tiny ovals}}}$}
%\end{picture}
%\caption{A non-singular real curve of degree $6$ on $\br\bp^2$ 
%with the invariants $(r,a,\delta_h=1, \delta_{\varphi S})$ 
%corresponding to the non-oriented positive curve}
%\label{non-singular-degree6}
%\end{figure}

\medskip

\noindent
{\it The connected component of moduli of 
a positive non-singular real curve $A$ of the degree 6 on the real 
projective plane is defined by its isotopy type and by 
dividiness of $A(\bc)$ by $A(\br)$. All possibilities are described above.}

\bigskip

Applying results of Sect. \ref{modrealK3}, below we get similar results 
for all lattices $S$ of $\rk S=2$.

\bigskip

%%%%%%%%%%%%%%%%%%%%%%%%%%%%%%%%%%%%%%%%%%%%%%%%%%%%%%
\section{Connected components of moduli of real non-singular curves of 
bidegree $(4,4)$ on hyperboloid}
\label{hyperboloid}

Here we consider the case $S \cong U(2)$ and $\theta=-1$ on $S$.  
Then $Y = X/\{1,\tau \} = \bp^1\times \bp^1$  in the non-degenerate 
case which we only consider, and $A = X^\tau$ is a non-singular curve 
of bidegree $(4,4)$ in $Y$. 

Let $S=\bz e_1+\bz e_2$ where $e_1^2=e_2^2=0$ and $(e_1,e_2) = 2$. 
The generators $e_1$ and $e_2$ of $S = H_2(X;\bz)^\tau$ are 
classes of preimages of $\text{pt}\times \bp^1$ and $\bp^1 \times \text{pt}$ 
respectively. 

Since $\theta = -1$ on $S$, then the real structure on $Y$ (defined 
by the anti-holomorphic involution $\varphi$ of $X$) is one of the 
four structures: 
(usual, usual) = {\it hyperboloid}, (usual, spin), (spin, usual) and 
(spin, spin), 
where {\it usual} means 
the usual real structure on $\bp^1$, 
which is given by the anti-holomoprhic involution 
$(z_0:z_1)\mapsto (\overline{z_0}:\overline{z_1})$ 
in some homogeneous coordinates $(z_0:z_1)$ on $\bp^1$, while 
{\it spin} is given by 
$(z_0:z_1)\mapsto (\overline{z_1}:-\overline{z_0})$. 
In hyperpoloid case $Y(\br) = \br\bp^1\times \br\bp^1$ is a $2$-torus. 
In all three spin cases 
$Y(\br)$ and $X(\br)$ are empty, and $A(\br)$ is also empty. 

If $Y$ has the usual real structure (i.e. it is a hyperboloid), 
then $A$ is the zero set of a real bi-homogeneous polynomial 
$P(x_0:x_1;y_0:y_1)$ 
of bidegree $(4,4)$. 
For a positive curve $(A,\varphi)$, equivalently for  
$A^+ = \pi(X_\varphi(\br))$, 
like for real non-singular curves of degree $6$ in $\bp^2$, 
we can choose $P$ by the condition that 
$A = \{ P=0 \}$ in $\bp^1\times \bp^1$ and 
$A^+ = \{ P\ge 0\}$ on $\br\bp^1\times \br\bp^1$. 
The polynomial $P$ is defined up to 
$\br_{++} \times ((PGL(2,\br)\times PGL(2,\br)) \rtimes \bz/2\bz)$. 
Here $\br_{++}$ denote all positive real numbers.  
Thus, classification of connected components of moduli 
of positive curves $A^+$ on a hyperboloid 
is equivalent to the description of connected components of 
$$
\left((\br^{25}-Discr)/\br_{++} \right)
/\left((PGL(2,\br)\times PGL(2,\br)) \rtimes \bz/2\bz \right)
$$
where the discriminant $Discr$ is the set of all polynomials giving 
singular (over $\bc$) curves. 
The group 
$(PGL(2,\br)\times PGL(2,\br)) \rtimes \bz/2\bz$ 
has $8$ connected components 
and has dimension $6$ over $\br$. 

%%%%%%%%%%%% spin case %%%%%%%%%%%%%%%%%%%%%
If $Y$ has one of spin structures, then $Discr$ has codimension two. 
Hence the moduli of positive curves and the moduli of curves are connected. 
(We will also give another proof of that below.) 

%%%%%%%%%%%%%%%%%%%%% enumeration %%%%%%%%%%%%%%%%%%%%%%
\medskip

Now we outline general arguments of Subsect. 
\ref{genclassification} in the $(S,\theta)=(U(2),-1)$ case (for finding 
connected components of moduli).  

Since 
$\Delta(S_+, S)^{(-4)} = \emptyset$ and 
$\Delta(S_-, S)^{(-4)} = \{ \pm(e_1 - e_2) \}$, the group  
$W^{(-4)}(S_+,S)$ is trivial, and the group  
$W^{(-4)}(S_-,S)$ consists of identity and $g$ 
where $g(e_1)=e_2$ and $g(e_2)=e_1$.
By Lemma \ref{GroupG}, we have $G = \{ \mbox{identity},\ g \} \cong \bz/2\bz$. 

We have $s=2$, $p=0$. 
$A_{S_-} = S_-^*/S_- = S^*/S$ is generated by 
$e_1^* = \frac12 e_2$ and $e_2^* = \frac12 e_1$, 
and hence it is isomorphic to $(\bz/2\bz)^2$ and $l(A_{S_-})=2$. 
We have 
$A_{S_+}=0,\ H_+ =0,\ \rho=0$. 
Since $\Gamma_- = 0$, 
$H = H_- \ (\subset S/2S)$ is one of the following 5 subgroups; 
$$\textstyle 
0,\ [ e_1 ],\ [ e_2 ],\ 
[ h ]=[e_1+e_2],\ S/2S=[e_1,e_2],$$
where we set $h = e_1 + e_2$ and 
consider $e_1$, $e_2$ and $h$ mod $2S$. 
Since $H_+ = 0$, we have $q_\rho = (- q_{S_-})|H_-$. 
%%%this means (- q_{S_-})|\frac12H_-
Since $\delta(S)=0$, 
we always have $\delta_H=0$ and 
$\delta_{\varphi S_\pm}=\delta_\varphi$. 
If $\delta_\varphi = 0$, then $\varepsilon_{v_\pm}\equiv 0\pmod 2$.

Thus, we have the following cases: 

(1) $H=0$:  
Then $a_H=0$, $k_\rho=0$, $\mu_\rho=0$, $\sigma_\rho \equiv 0\pmod 8$. 
If $\delta_{\varphi S}=0$, then 
$v=0$, equivalently, $\delta_\varphi=0$. We have $c_v\equiv 0\pmod 4$, 
$\delta_{\varphi S_\pm}=\delta_\varphi=0$. 

(2) $H=[e_1 ]$:  
Then $a_H=1$, $q_\rho \cong z$, 
$k_\rho=1$, $\mu_\rho=0$, $\sigma_\rho \equiv 0\pmod 8$. 

(3) $H=[ e_2 ]$:
Then $a_H=1$, $q_\rho \cong z$, 
$k_\rho=1$, $\mu_\rho=0$, $\sigma_\rho \equiv 0\pmod 8$. 

(4) $H=[ h ]=[e_1+e_2]$:
Then $a_H=1$, $q_\rho \cong w$, 
$k_\rho=1$, $\mu_\rho=1$. 

(5) $H=S/2S=[e_1,e_2]$: 
Then $a_H=2$, $q_\rho \cong u_+(2)$, 
$k_\rho=0$, $\mu_\rho=0$, $\sigma_\rho \equiv 0\pmod 8$. 
If $\delta_{\varphi S}=0$, then $v=0$, and $\delta_{\varphi S}=0$ 
is equivalent to $\delta_\varphi=0$. 

From above, $\mu_\rho=1$ if and only if $H=[ h ]$. 
If $\mu_\rho=0$, then $\sigma_\rho \equiv 0\pmod 8$. 

We see that genus 
(hence, by Theorem \ref{genus-isomor} and Proposition \ref{genus-isomor2}, 
the isomorphism class) of an integral involutions $(L,\varphi,S)$ 
of the type $(U(2),-1)$ satisfying (RSK3) 
is determined by the data 
$(r, a, H, \delta_{\varphi S}, v)$, 
where the $v$ is defined only if $\delta_{\varphi S}=0$. 

By Theorem  \ref{geninvtheorem}, the complete list of the above invariants 
is given by the conditions:\\

\noindent
{\bf Type 0} ($\delta_\varphi=0$)\\
$a_H + k_\rho \leq a \leq r$\ and\ $a \equiv 0 \pmod{2}.$\\
$2 \leq r \leq 18$ and $2-r \equiv 0 \pmod{4}.$\\
$r+a \leq 2a_H +18.$\\
If $a = a_H + k_\rho,\ \mu_\rho=0,$ then $2-r \equiv 0 \pmod{8}.$\\
If $a=r$, then $2-r \equiv 0 \pmod{8}.$\\
If $r+a=2a_H +18$, then $2-r \equiv 0 \pmod{8}.$\\

\noindent
{\bf Type Ia} ($\delta_{\varphi S}=0,\ \delta_\varphi=1$)\\
$a_H=1$,\ $k_\rho =1.$\\
$c_v \equiv 0 \pmod{4}$ if $\mu_\rho=0,$ \ and 
$c_v \equiv 2 \pmod{4}$ if $\mu_\rho=1.$\\
$2 \leq a \leq r$\ and\ $a \equiv 0 \pmod{2}.$\\
$2 \leq r \leq 18$ and $2-r \equiv c_v \pmod{4}.$\\
$r+a \leq 20.$\\
If $a = 2,\ \mu_\rho=0,$ then $2-r \equiv 0 \pmod{8}.$\\

\noindent
{\bf Type Ib} ($\delta_{\varphi S}=1$)\\
$a_H + k_\rho +1 \leq a \leq r,$\\
$1 \leq r \leq 19,$\\
$r+a \equiv 0 \pmod{2},$\ \ $r+a \leq 2a_H +18.$\\
If $a = a_H + k_\rho +1,\ \mu_\rho=0,$ then $2-r \equiv \pm 1 \pmod{8}.$\\
If $a = a_H + k_\rho +2,\ \mu_\rho=0,$ then $2-r \not \equiv 4 \pmod{8}.$\\

\medskip

We see $a_H=1$ in Type Ia case. 
Then $H$ is generated by the unique nonzero element, 
hence it is nothing but $v$. 
Hence we can use the invariant $\delta_\varphi$ instead of $v$ 
in the collection of invariants above. 
Thus the isomorphism class is determined by the data 
$$
(r, a, H, \delta_{\varphi S}, \delta_\varphi).
$$

Since $G = \{ \mbox{identity},\ g \}$, 
the triplets $(r,a,[ e_1 ])$
         and $(r,a,[ e_2 ])$ 
represent the same isomorphism class for each Type (0,\ Ia and Ib). 
The other different triplets represent different isomorphism classes. 

\medskip

All possible data 
\begin{equation}
(r, a, H, \delta_{\varphi S}, v)
\label{invhyperboloid}
\end{equation}
(which are 
equivalent to  $(r, a, H, \delta_{\varphi S}, \delta_\varphi)$)   
are given in  Figures \ref{hyperboloid-0} --- \ref{hyperboloid-h}. 
We have $\delta_\varphi =0$, if and only if 
$\delta_{\varphi S}=v=0$.  

The relations between related involutions are as follows. 
Since $\delta(S)=0$ and $\theta=-1$ on $S$, then $s_\theta=0\mod 2L$. 
Hence, by Theorem \ref{relinvtheorem}, we have:
\begin{equation}
r(\varphi)+r(\tau\varphi)=20;\ \ 
a_{H(\varphi)} + a_{H(\tau\varphi)}=2;\ \ 
a(\varphi)-a_{H(\varphi)} = a(\tau\varphi)-a_{H(\tau\varphi)};
\label{hyperbo-related1}
\end{equation}
\begin{equation}
H(\tau\varphi)=H(\varphi)^\perp\ 
\text{w. r. t.}\ b_{S_-}\ \text{on\ }A_{S_-};\ \ 
\delta_{\varphi S}=\delta_{\tau\varphi S};\ \ 
s_\varphi \equiv s_{\tau\varphi} \mod 2L.
\label{hyperbo-related2}
\end{equation}

Hence, 
involutions of Type 0 (resp. Type Ib) with $H=0$ ($11$ (resp. $39$) classes) are related to 
involutions of Type 0 (resp. Type Ib) with $H=S/2S$. 
Involutions of Type 0 (resp. Type Ia, Type Ib) with $H=[ e_i ]$ 
($9$ (resp. $11$, $30$) classes) are related to 
involutions of Type 0 (resp. Type Ia, Type Ib) with $H=[e_i ]$. 
(More precisely, 
the class $(r, a)$ is related to the class $(20-r, a)$.) 
Involutions of Type 0 (resp. Type Ia, Type Ib) with $H= [ h ]$ 
($11$ (resp. $12$, $36$) classes) are related to 
involutions of Type 0 (resp. Type Ia, Type Ib) with $H=[ h ]$. 
(More precisely, 
the class $(r, a)$ is related to the class with $(20-r, a)$, too.) 

Thus there are 
$42 (= 11\times2 +9+11)$    classes 
(i.e., connected components of moduli of positives curves 
of bidegree $(4,4)$) 
of Type 0, 
$23 (= 11+12)$              classes of Type Ia, and 
$144 (= 39\times2 + 30+36)$ classes of Type Ib.

Moreover, if we identify related involutions, there are 
$26 (=11+7+8)$   classes 
(i.e., connected components of moduli of real non-singular curves 
of bidegree $(4,4)$) 
of Type 0, 
$14 (=8+6)$      classes of Type Ia, and 
$76 (=39+17+20)$ classes of Type Ib.

\begin{figure}
\begin{picture}(200,140)
%%hyperboloid
%%\put(24,120){$H=0$}
%\put(53,102){{\tiny (Here,}}
\put(66,104){\circle{7}}
\put(71,102){{\tiny means $\delta_{\varphi S}=0$ and $v=0$,}}
%{\circle{5}}
%{{\tiny means $\delta_{\varphi S}=0$ and $v \neq 0$}}
\put(66, 96){\circle*{3}}
\put(71, 94){{\tiny means $\delta_{\varphi S}=1$}}
\multiput(8,0)(8,0){20}{\line(0,1){86}}
\multiput(0,8)(0,8){10}{\line(1,0){170}}
\put(0,0){\vector(0,1){100}}
\put(0,0){\vector(1,0){180}}
\put(  6,-8){{\tiny $1$}}
\put( 14,-8){{\tiny $2$}}
\put( 22,-8){{\tiny $3$}}
\put( 30,-8){{\tiny $4$}}
\put( 38,-8){{\tiny $5$}}
\put( 46,-8){{\tiny $6$}}
\put( 54,-8){{\tiny $7$}}
\put( 62,-8){{\tiny $8$}}
\put( 70,-8){{\tiny $9$}}
\put( 76,-8){{\tiny $10$}}
\put( 84,-8){{\tiny $11$}}
\put( 92,-8){{\tiny $12$}}
\put(100,-8){{\tiny $13$}}
\put(108,-8){{\tiny $14$}}
\put(116,-8){{\tiny $15$}}
\put(124,-8){{\tiny $16$}}
\put(132,-8){{\tiny $17$}}
\put(140,-8){{\tiny $18$}}
\put(148,-8){{\tiny $19$}}
\put(156,-8){{\tiny $20$}}

\put(-8, -1){{\tiny $0$}}
\put(-8,  7){{\tiny $1$}}
\put(-8, 15){{\tiny $2$}}
\put(-8, 23){{\tiny $3$}}
\put(-8, 31){{\tiny $4$}}
\put(-8, 39){{\tiny $5$}}
\put(-8, 47){{\tiny $6$}}
\put(-8, 55){{\tiny $7$}}
\put(-8, 63){{\tiny $8$}}
\put(-8, 71){{\tiny $9$}}
\put(-10, 79){{\tiny $10$}}

\put( -2,106){{\footnotesize $a$}} %vertical
\put(186, -2){{\footnotesize $r$}} %horizontal

%% Type 0  \circle{7}
%%  ( 8r,8a)
\put( 16, 0){\circle{7}}
\put( 80, 0){\circle{7}}
\put(144, 0){\circle{7}}
\put( 16,16){\circle{7}}
\put( 48,16){\circle{7}}
\put( 80,16){\circle{7}}
\put(112,16){\circle{7}}
\put( 48,32){\circle{7}}
\put( 80,32){\circle{7}}
\put( 80,48){\circle{7}}
\put( 80,64){\circle{7}}

%% Type Ib  \circle*{3}
%%  ( 8r,8a)
\put(  8,8){\circle*{3}}
\put( 24,8){\circle*{3}}
\put( 72,8){\circle*{3}}
\put( 88,8){\circle*{3}}
\put(136,8){\circle*{3}}
\put( 16,16){\circle*{3}}
\put( 32,16){\circle*{3}}
\put( 64,16){\circle*{3}}
\put( 80,16){\circle*{3}}
\put( 96,16){\circle*{3}}
\put(128,16){\circle*{3}}
\put( 24,24){\circle*{3}}
\put( 40,24){\circle*{3}}
\put( 56,24){\circle*{3}}
\put( 72,24){\circle*{3}}
\put( 88,24){\circle*{3}}
\put(104,24){\circle*{3}}
\put(120,24){\circle*{3}}
\put( 32,32){\circle*{3}}
\put( 48,32){\circle*{3}}
\put( 64,32){\circle*{3}}
\put( 80,32){\circle*{3}}
\put( 96,32){\circle*{3}}
\put(112,32){\circle*{3}}
\put( 40,40){\circle*{3}}
\put( 56,40){\circle*{3}}
\put( 72,40){\circle*{3}}
\put( 88,40){\circle*{3}}
\put(104,40){\circle*{3}}
\put( 48,48){\circle*{3}}
\put( 64,48){\circle*{3}}
\put( 80,48){\circle*{3}}
\put( 96,48){\circle*{3}}
\put( 56,56){\circle*{3}}
\put( 72,56){\circle*{3}}
\put( 88,56){\circle*{3}}
\put( 64,64){\circle*{3}}
\put( 80,64){\circle*{3}}
\put( 72,72){\circle*{3}}
\end{picture}

\caption{$\bhh$: All possible $(r,a,\delta_{\varphi S}, v)$ with $H=0$}
\label{hyperboloid-0}
\end{figure}

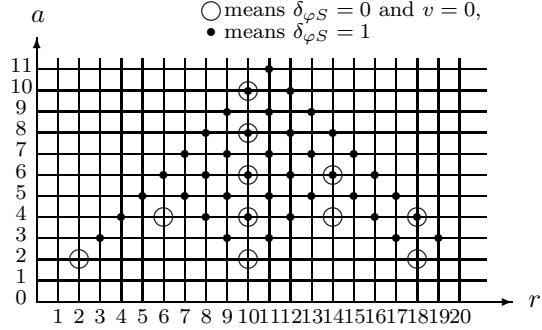
\begin{figure}
\begin{picture}(200,140)
%%hyperboloid
\put(66,110){\circle{7}}
\put(71,108){{\tiny means $\delta_{\varphi S}=0$ and $v=0$,}}
%{\circle{5}}
%{{\tiny means $\delta_{\varphi S}=0$ and $v \neq 0$}}
\put(66, 102){\circle*{3}}
\put(71, 100){{\tiny means $\delta_{\varphi S}=1$}}
\multiput(8,0)(8,0){20}{\line(0,1){92}}
\multiput(0,8)(0,8){11}{\line(1,0){170}}
\put(0,0){\vector(0,1){100}}
\put(0,0){\vector(1,0){180}}
\put(  6,-8){{\tiny $1$}}
\put( 14,-8){{\tiny $2$}}
\put( 22,-8){{\tiny $3$}}
\put( 30,-8){{\tiny $4$}}
\put( 38,-8){{\tiny $5$}}
\put( 46,-8){{\tiny $6$}}
\put( 54,-8){{\tiny $7$}}
\put( 62,-8){{\tiny $8$}}
\put( 70,-8){{\tiny $9$}}
\put( 76,-8){{\tiny $10$}}
\put( 84,-8){{\tiny $11$}}
\put( 92,-8){{\tiny $12$}}
\put(100,-8){{\tiny $13$}}
\put(108,-8){{\tiny $14$}}
\put(116,-8){{\tiny $15$}}
\put(124,-8){{\tiny $16$}}
\put(132,-8){{\tiny $17$}}
\put(140,-8){{\tiny $18$}}
\put(148,-8){{\tiny $19$}}
\put(156,-8){{\tiny $20$}}

\put(-8, -1){{\tiny $0$}}
\put(-8,  7){{\tiny $1$}}
\put(-8, 15){{\tiny $2$}}
\put(-8, 23){{\tiny $3$}}
\put(-8, 31){{\tiny $4$}}
\put(-8, 39){{\tiny $5$}}
\put(-8, 47){{\tiny $6$}}
\put(-8, 55){{\tiny $7$}}
\put(-8, 63){{\tiny $8$}}
\put(-8, 71){{\tiny $9$}}
\put(-10, 79){{\tiny $10$}}
\put(-10, 87){{\tiny $11$}}

\put( -2,106){{\footnotesize $a$}} %vertical
\put(186, -2){{\footnotesize $r$}} %horizontal

%% Type 0  \circle{7}
%%  ( 8r,8a)
\put( 16, 16){\circle{7}}
\put( 80, 16){\circle{7}}
\put(144, 16){\circle{7}}
\put( 48,32){\circle{7}}
\put( 80,32){\circle{7}}
\put(112,32){\circle{7}}
\put(144,32){\circle{7}}
\put( 80,48){\circle{7}}
\put( 112,48){\circle{7}}

\put( 80,64){\circle{7}}
\put( 80,80){\circle{7}}

%% Type Ib  \circle*{3}
%%  ( 8r,8a)
\put( 24,24){\circle*{3}}
\put( 72,24){\circle*{3}}
\put( 88,24){\circle*{3}}
\put(136,24){\circle*{3}}
\put(152,24){\circle*{3}}
\put( 32,32){\circle*{3}}
\put( 64,32){\circle*{3}}
\put( 80,32){\circle*{3}}
\put( 96,32){\circle*{3}}
\put(128,32){\circle*{3}}
\put(144,32){\circle*{3}}
\put( 40,40){\circle*{3}}
\put( 56,40){\circle*{3}}
\put( 72,40){\circle*{3}}
\put( 88,40){\circle*{3}}
\put(104,40){\circle*{3}}
\put(120,40){\circle*{3}}
\put(136,40){\circle*{3}}
\put( 48,48){\circle*{3}}
\put( 64,48){\circle*{3}}
\put( 80,48){\circle*{3}}
\put( 96,48){\circle*{3}}
\put(112,48){\circle*{3}}
\put(128,48){\circle*{3}}

\put( 56,56){\circle*{3}}
\put( 72,56){\circle*{3}}
\put( 88,56){\circle*{3}}
\put(104,56){\circle*{3}}
\put(120,56){\circle*{3}}
\put( 64, 64){\circle*{3}}
\put( 80,64){\circle*{3}}
\put( 96,64){\circle*{3}}
\put( 112,64){\circle*{3}}

\put( 72,72){\circle*{3}}
\put( 88,72){\circle*{3}}
\put( 104,72){\circle*{3}}
\put( 80,80){\circle*{3}}
\put( 96,80){\circle*{3}}
\put( 88,88){\circle*{3}}
\end{picture}

\caption{$\bhh$: All possible $(r,a,\delta_{\varphi S}, v)$ 
with $H=[e_1,e_2]=S/2S$}
\label{hyperboloid-e_1e_2}
\end{figure}

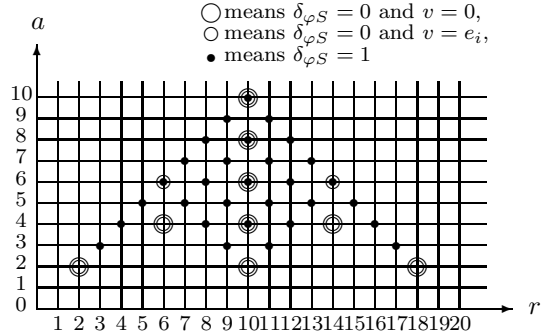
\begin{figure}
\begin{picture}(200,140)
%%hyperboloid
%%\put(24,122){$H=[ e_i ]\ \ (i=1\ \mbox{or}\ 2)$}
\put(66,112){\circle{7}}
\put(71,110){{\tiny means $\delta_{\varphi S}=0$ and $v=0$,}}
\put(66,104){\circle{5}}
\put(71,102){{\tiny means $\delta_{\varphi S}=0$ and $v=e_i$,}}
\put(66, 95){\circle*{3}}
\put(71, 94){{\tiny means $\delta_{\varphi S}=1$}}
\multiput(8,0)(8,0){20}{\line(0,1){86}}
\multiput(0,8)(0,8){10}{\line(1,0){170}}
\put(0,0){\vector(0,1){100}}
\put(0,0){\vector(1,0){180}}
\put(  6,-8){{\tiny $1$}}
\put( 14,-8){{\tiny $2$}}
\put( 22,-8){{\tiny $3$}}
\put( 30,-8){{\tiny $4$}}
\put( 38,-8){{\tiny $5$}}
\put( 46,-8){{\tiny $6$}}
\put( 54,-8){{\tiny $7$}}
\put( 62,-8){{\tiny $8$}}
\put( 70,-8){{\tiny $9$}}
\put( 76,-8){{\tiny $10$}}
\put( 84,-8){{\tiny $11$}}
\put( 92,-8){{\tiny $12$}}
\put(100,-8){{\tiny $13$}}
\put(108,-8){{\tiny $14$}}
\put(116,-8){{\tiny $15$}}
\put(124,-8){{\tiny $16$}}
\put(132,-8){{\tiny $17$}}
\put(140,-8){{\tiny $18$}}
\put(148,-8){{\tiny $19$}}
\put(156,-8){{\tiny $20$}}

\put(-8, -1){{\tiny $0$}}
\put(-8,  7){{\tiny $1$}}
\put(-8, 15){{\tiny $2$}}
\put(-8, 23){{\tiny $3$}}
\put(-8, 31){{\tiny $4$}}
\put(-8, 39){{\tiny $5$}}
\put(-8, 47){{\tiny $6$}}
\put(-8, 55){{\tiny $7$}}
\put(-8, 63){{\tiny $8$}}
\put(-8, 71){{\tiny $9$}}
\put(-10, 79){{\tiny $10$}}

\put( -2,106){{\footnotesize $a$}} %vertical
\put(186, -2){{\footnotesize $r$}} %horizontal

%% Type 0  \circle{7}
%%  ( 8r,8a)
\put( 16,16){\circle{7}}
\put( 80,16){\circle{7}}
\put(144,16){\circle{7}}
\put( 48,32){\circle{7}}
\put( 80,32){\circle{7}}
\put(112,32){\circle{7}}
\put( 80,48){\circle{7}}
\put( 80,64){\circle{7}}
\put( 80,80){\circle{7}}

%% Type Ia  \circle{5}
%%  ( 8r,8a)
\put( 16,16){\circle{5}}
\put( 80,16){\circle{5}}
\put(144,16){\circle{5}}
\put( 48,32){\circle{5}}
\put( 80,32){\circle{5}}
\put(112,32){\circle{5}}
\put( 48,48){\circle{5}}
\put( 80,48){\circle{5}}
\put(112,48){\circle{5}}
\put( 80,64){\circle{5}}
\put( 80,80){\circle{5}}

%% Type Ib  \circle*{3}
%%  ( 8r,8a)
\put( 24,24){\circle*{3}}
\put( 72,24){\circle*{3}}
\put( 88,24){\circle*{3}}
\put(136,24){\circle*{3}}
\put( 32,32){\circle*{3}}
\put( 64,32){\circle*{3}}
\put( 80,32){\circle*{3}}
\put( 96,32){\circle*{3}}
\put(128,32){\circle*{3}}
\put( 40,40){\circle*{3}}
\put( 56,40){\circle*{3}}
\put( 72,40){\circle*{3}}
\put( 88,40){\circle*{3}}
\put(104,40){\circle*{3}}
\put(120,40){\circle*{3}}
\put( 48,48){\circle*{3}}
\put( 64,48){\circle*{3}}
\put( 80,48){\circle*{3}}
\put( 96,48){\circle*{3}}
\put(112,48){\circle*{3}}
\put( 56,56){\circle*{3}}
\put( 72,56){\circle*{3}}
\put( 88,56){\circle*{3}}
\put(104,56){\circle*{3}}
\put( 64,64){\circle*{3}}
\put( 80,64){\circle*{3}}
\put( 96,64){\circle*{3}}
\put( 72,72){\circle*{3}}
\put( 88,72){\circle*{3}}
\put( 80,80){\circle*{3}}
\end{picture}

\caption{$\bhh$: All possible $(r,a,\delta_{\varphi S}, v)$ 
with $H=[e_i]\ \ (i=1\ \mbox{or}\ 2)$}
\label{hyperboloid-ei}
\end{figure}

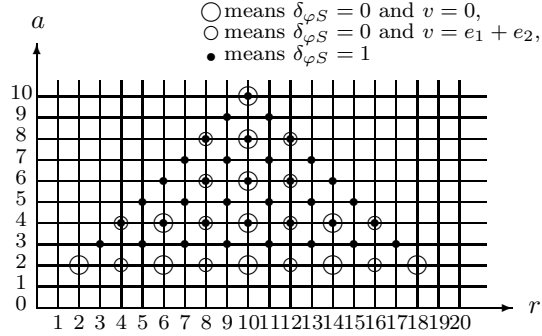
\begin{figure}
\begin{picture}(200,140)
%%hyperboloid
%%\put(24,122){$H=[ e_1+e_2 ]$}
\put(66,112){\circle{7}}
\put(71,110){{\tiny means $\delta_{\varphi S}=0$ and $v=0$,}}
\put(66,104){\circle{5}}
\put(71,102){{\tiny means $\delta_{\varphi S}=0$ and $v=e_1+e_2$,}}
\put(66, 95){\circle*{3}}
\put(71, 94){{\tiny means $\delta_{\varphi S}=1$}}
\multiput(8,0)(8,0){20}{\line(0,1){86}}
\multiput(0,8)(0,8){10}{\line(1,0){170}}
\put(0,0){\vector(0,1){100}}
\put(0,0){\vector(1,0){180}}
\put(  6,-8){{\tiny $1$}}
\put( 14,-8){{\tiny $2$}}
\put( 22,-8){{\tiny $3$}}
\put( 30,-8){{\tiny $4$}}
\put( 38,-8){{\tiny $5$}}
\put( 46,-8){{\tiny $6$}}
\put( 54,-8){{\tiny $7$}}
\put( 62,-8){{\tiny $8$}}
\put( 70,-8){{\tiny $9$}}
\put( 76,-8){{\tiny $10$}}
\put( 84,-8){{\tiny $11$}}
\put( 92,-8){{\tiny $12$}}
\put(100,-8){{\tiny $13$}}
\put(108,-8){{\tiny $14$}}
\put(116,-8){{\tiny $15$}}
\put(124,-8){{\tiny $16$}}
\put(132,-8){{\tiny $17$}}
\put(140,-8){{\tiny $18$}}
\put(148,-8){{\tiny $19$}}
\put(156,-8){{\tiny $20$}}

\put(-8, -1){{\tiny $0$}}
\put(-8,  7){{\tiny $1$}}
\put(-8, 15){{\tiny $2$}}
\put(-8, 23){{\tiny $3$}}
\put(-8, 31){{\tiny $4$}}
\put(-8, 39){{\tiny $5$}}
\put(-8, 47){{\tiny $6$}}
\put(-8, 55){{\tiny $7$}}
\put(-8, 63){{\tiny $8$}}
\put(-8, 71){{\tiny $9$}}
\put(-10, 79){{\tiny $10$}}

\put( -2,106){{\footnotesize $a$}} %vertical
\put(186, -2){{\footnotesize $r$}} %horizontal

%% Type 0  \circle{7}
%%  ( 8r,8a)
\put( 16,16){\circle{7}}
\put( 48,16){\circle{7}}
\put( 80,16){\circle{7}}
\put(112,16){\circle{7}}
\put(144,16){\circle{7}}
\put( 48,32){\circle{7}}
\put( 80,32){\circle{7}}
\put(112,32){\circle{7}}
\put( 80,48){\circle{7}}
\put( 80,64){\circle{7}}
\put( 80,80){\circle{7}}

%% Type Ia  \circle{5}
%%  ( 8r,8a)
\put( 32,16){\circle{5}}
\put( 64,16){\circle{5}}
\put( 96,16){\circle{5}}
\put(128,16){\circle{5}}
\put( 32,32){\circle{5}}
\put( 64,32){\circle{5}}
\put( 96,32){\circle{5}}
\put(128,32){\circle{5}}
\put( 64,48){\circle{5}}
\put( 96,48){\circle{5}}
\put( 64,64){\circle{5}}
\put( 96,64){\circle{5}}

%% Type Ib  \circle*{3}
%%  ( 8r,8a)
\put( 24,24){\circle*{3}}
\put( 40,24){\circle*{3}}
\put( 56,24){\circle*{3}}
\put( 72,24){\circle*{3}}
\put( 88,24){\circle*{3}}
\put(104,24){\circle*{3}}
\put(120,24){\circle*{3}}
\put(136,24){\circle*{3}}
\put( 32,32){\circle*{3}}
\put( 48,32){\circle*{3}}
\put( 64,32){\circle*{3}}
\put( 80,32){\circle*{3}}
\put( 96,32){\circle*{3}}
\put(112,32){\circle*{3}}
\put(128,32){\circle*{3}}
\put( 40,40){\circle*{3}}
\put( 56,40){\circle*{3}}
\put( 72,40){\circle*{3}}
\put( 88,40){\circle*{3}}
\put(104,40){\circle*{3}}
\put(120,40){\circle*{3}}
\put( 48,48){\circle*{3}}
\put( 64,48){\circle*{3}}
\put( 80,48){\circle*{3}}
\put( 96,48){\circle*{3}}
\put(112,48){\circle*{3}}
\put( 56,56){\circle*{3}}
\put( 72,56){\circle*{3}}
\put( 88,56){\circle*{3}}
\put(104,56){\circle*{3}}
\put( 64,64){\circle*{3}}
\put( 80,64){\circle*{3}}
\put( 96,64){\circle*{3}}
\put( 72,72){\circle*{3}}
\put( 88,72){\circle*{3}}
\put( 80,80){\circle*{3}}
\end{picture}

\caption{$\bhh$: All possible $(r,a,\delta_{\varphi S}, v)$ 
with $H=[e_1+e_2]=[h]$}
\label{hyperboloid-h}
\end{figure}

\bigskip

If two positive curves are in one connected component of moduli or 
are related,  the real structure on $Y$ stays the same. By 
\eqref{realcomponents} and \eqref{realmod2}, the $A^+$ is empty, if and only if 
$(r,a,\delta_{\varphi})=(10,10,0)$. It follows that the component  
$(r, a, H, \delta_{\varphi S}, \delta_\varphi) =
(10,10,[ h ],0,0)$
corresponds to the real structure (spin, spin), the component 
$(10,10,[ e_i ],0,0)$ corresponds to the real structure (usual, spin) 
(or (spin, usual)). And all the remaining components 
correspond to the real structure (usual, usual) (namely, hyperboloid). 
The component 
$(10,10,S/2S,0,0)$ consists of 
empty $A^+$ on hyperboloid.  

%%%%%%%%%%%%%%%%%%%%%%%%%%%%%%%%%%%

The isotopy classification of real non-singular curves of 
bidegree $(4,4)$ on a hyperboloid 
was obtained by Gudkov \cite{Gudkov79}. 
Zvonilov \cite{Zvonilov92} clarified all the 
{\it complex schemes} of curves of bidegree $(4,4)$ on hyperboloid and 
ellipsoid, where {\it complex schemes} mean real schemes 
(i.e. real isotopy types) with divideness 
and their complex orientations (if dividing). 
Curves on hyperboloid and ellipsoid were studied in 
\cite{Zvonilov82}, \cite{Zvonilov83}, 
\cite{Matsuoka(Saito)91}, \cite{Matsuoka(Saito)90a}, 
\cite{Matsuoka(Saito)92}.  
In these articles the notions: \ \ 
{\it torsion} $(s,t)\ (\in \bz \times \bz)$ of 
a connected component of $A(\br)$,\ \ 
{\it oval},\ \ {\it non-oval},\ \ {\it odd branch}\ \ and {\it even branch},\ 
are well-known.  
We quote the isotopy classification of non-singular curves of 
bidegree $(4,4)$ on hyperboloid  
from \cite{Zvonilov92} in \textsc{Table} \ref{44isotopy} where we use 
notations due to Viro \cite{Viro80} and Zvonilov \cite{Zvonilov92}, 
and $l_1$ and $l_2$ denote non-ovals with torsions $(1,0)$ and $(0,1)$ 
respectively (e. g. see Figure \ref{hyperbo-4e1}).

\begin{table}
\begin{center}
\begin{tabular}{|l|}  \hline
$\langle 1\langle m \rangle \sqcup n \rangle$,   \\
\ \ where $(m,n) = 
(9,0),\ (5,4),\ (1,8),\ (8,0),\ (5,3),\ (4,4),\ (1,7);$         \\
\ \ \ or $m \geq 1,\ n \geq 0$ and $m + n \leq 7$.     \\ \hline
$\langle m \rangle$,           \\
\ \ where $0 \leq m \leq 9$.                        \\ \hline
$\langle 1\langle 1 \rangle \sqcup 1\langle 1 \rangle \rangle$    \\  \hline
$\langle l_1 + l_2,\ m,\ l_1 + l_2,\ n \rangle$,      \\
\ \ where $0 \leq m \leq n$ and $m + n \leq 8$.   \\ \hline
$\langle 4(l_1 + l_2) \rangle$               \\ \hline
$\langle l_1,\ m,\ l_1,\ n \rangle$,    \\
\ \ where $(m,n) = (0,8),\ (4,4),\ (0,7),\ (3,4)$;\\
\ \ \ or $0 \leq m \leq n$ and $m + n \leq 6$.  \\ \hline
$\langle 4(l_1) \rangle$            \\ \hline
$\langle 2(l_1 + 2l_2) \rangle$     \\ \hline
\end{tabular}
\end{center}
\caption{The list of isotopy types of real non-singular 
curves of bidegree $(4,4)$ 
on a hyperboloid.}
\label{44isotopy}
\end{table}

\begin{figure} 
\centerline{\includegraphics[width=5cm]{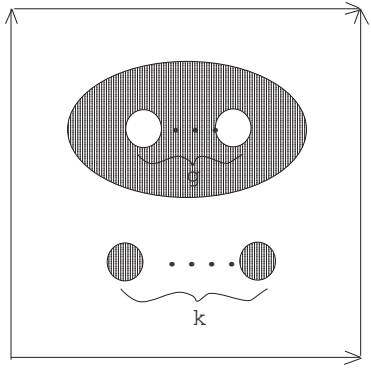}} 
\caption{$\bhh$: $A^+$ for $H=S/2S=[e_1,e_2]$ 
with $(r,a,\delta_{\varphi S},v) \neq (10,10,0,0)\ \text{nor}\ (10,8,0,0)$.}
\label{hyperbo}
\end{figure}

\begin{figure} 
\centerline{\includegraphics[width=5cm]{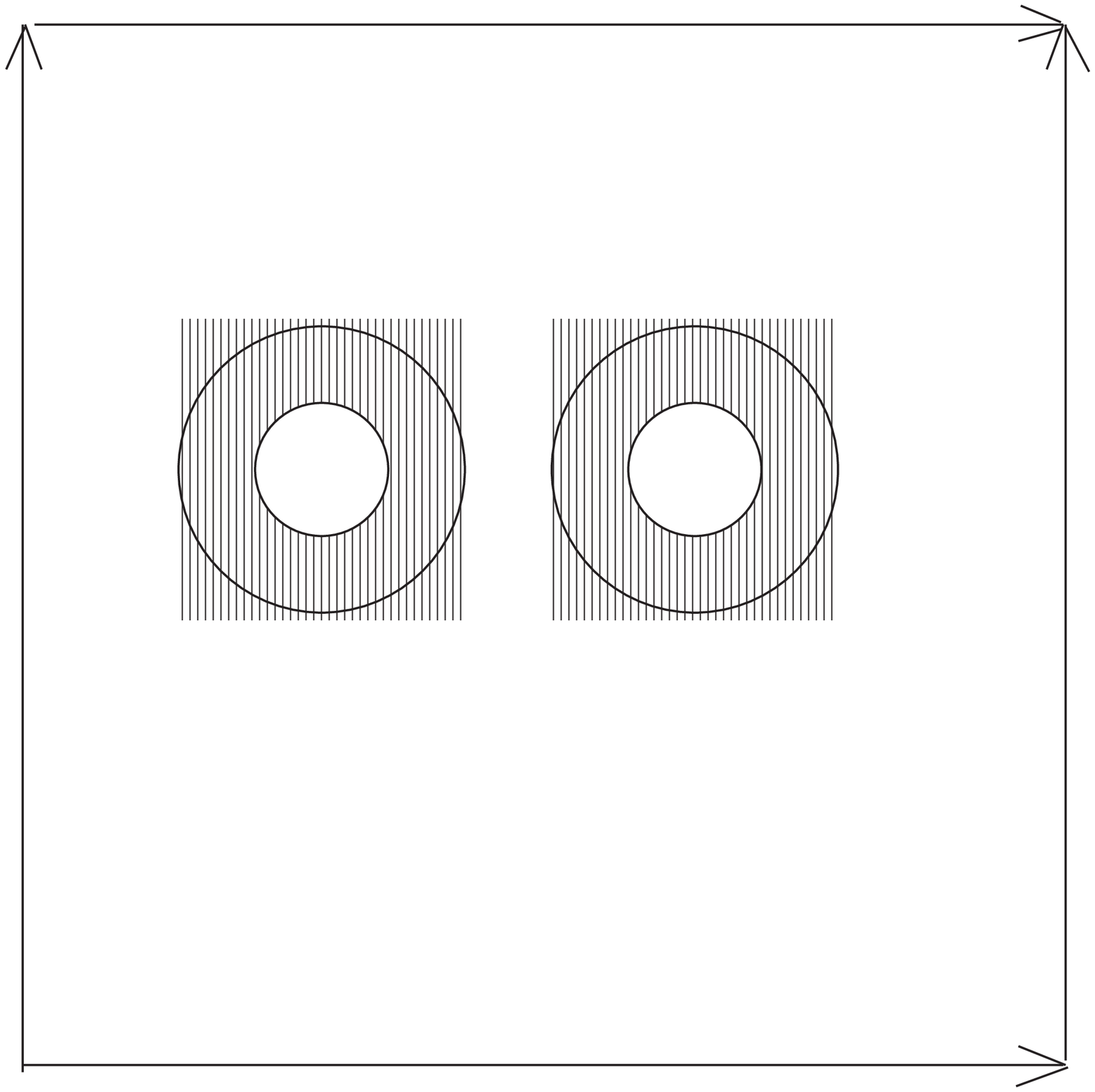}} 
\caption{$\bhh$: $A^+$ for $H=S/2S=[e_1,e_2]$ 
with $(r,a,\delta_{\varphi S},v)=(10,8,0,0)$.}
\label{hyperbo1111}
\end{figure}

\begin{figure} 
\centerline{\includegraphics[width=7cm]{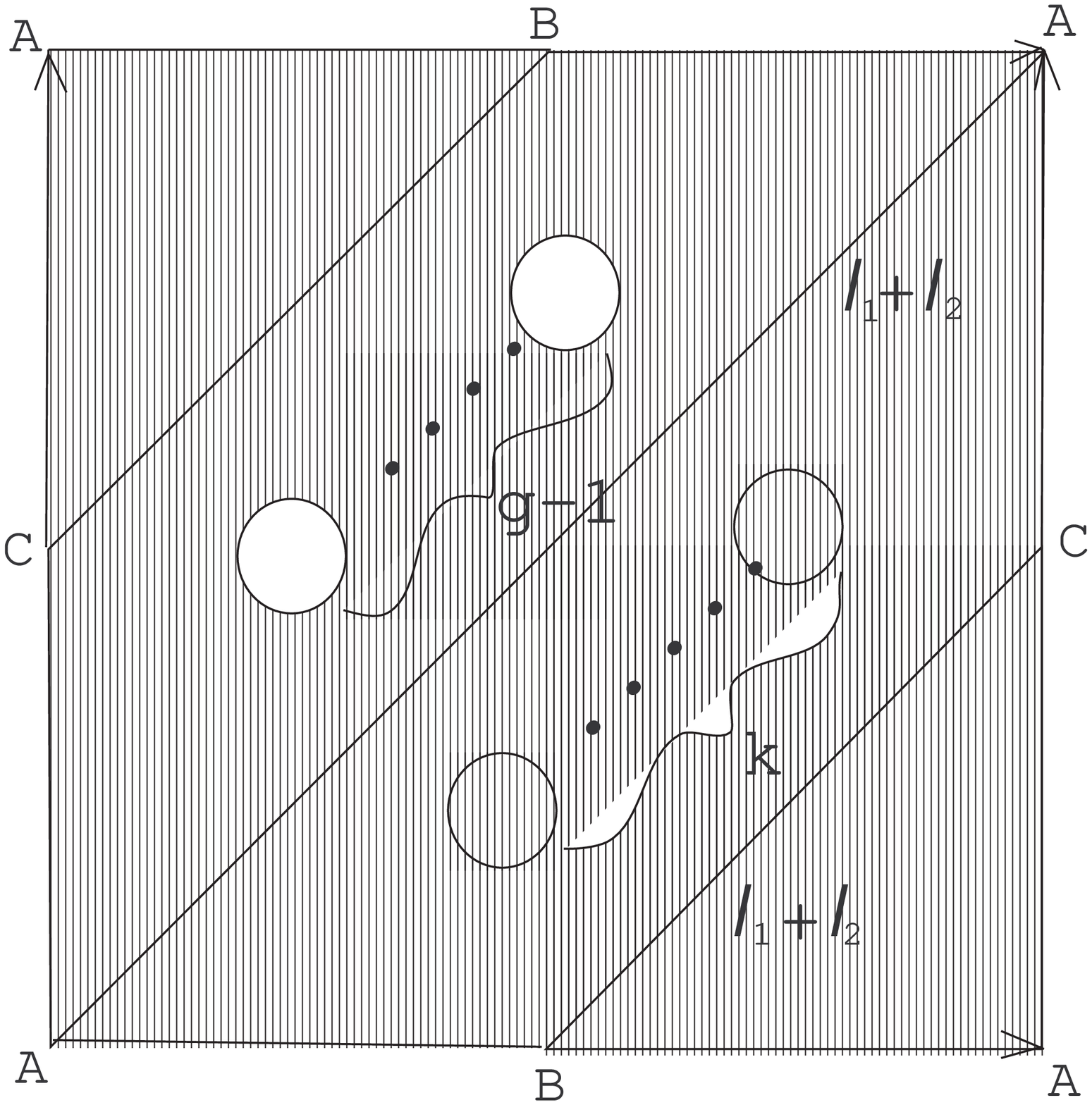}} 
\caption{$\bhh$: $A^+$ for $H=[h]=[e_1+e_2]$ 
with $(r,a,\delta_{\varphi S},v) \neq (10,10,0,0)\ \text{nor}\ (10,8,0,0)$.}
\label{hyperbo-h}
\end{figure}

\begin{figure} 
\centerline{\includegraphics[width=5cm]{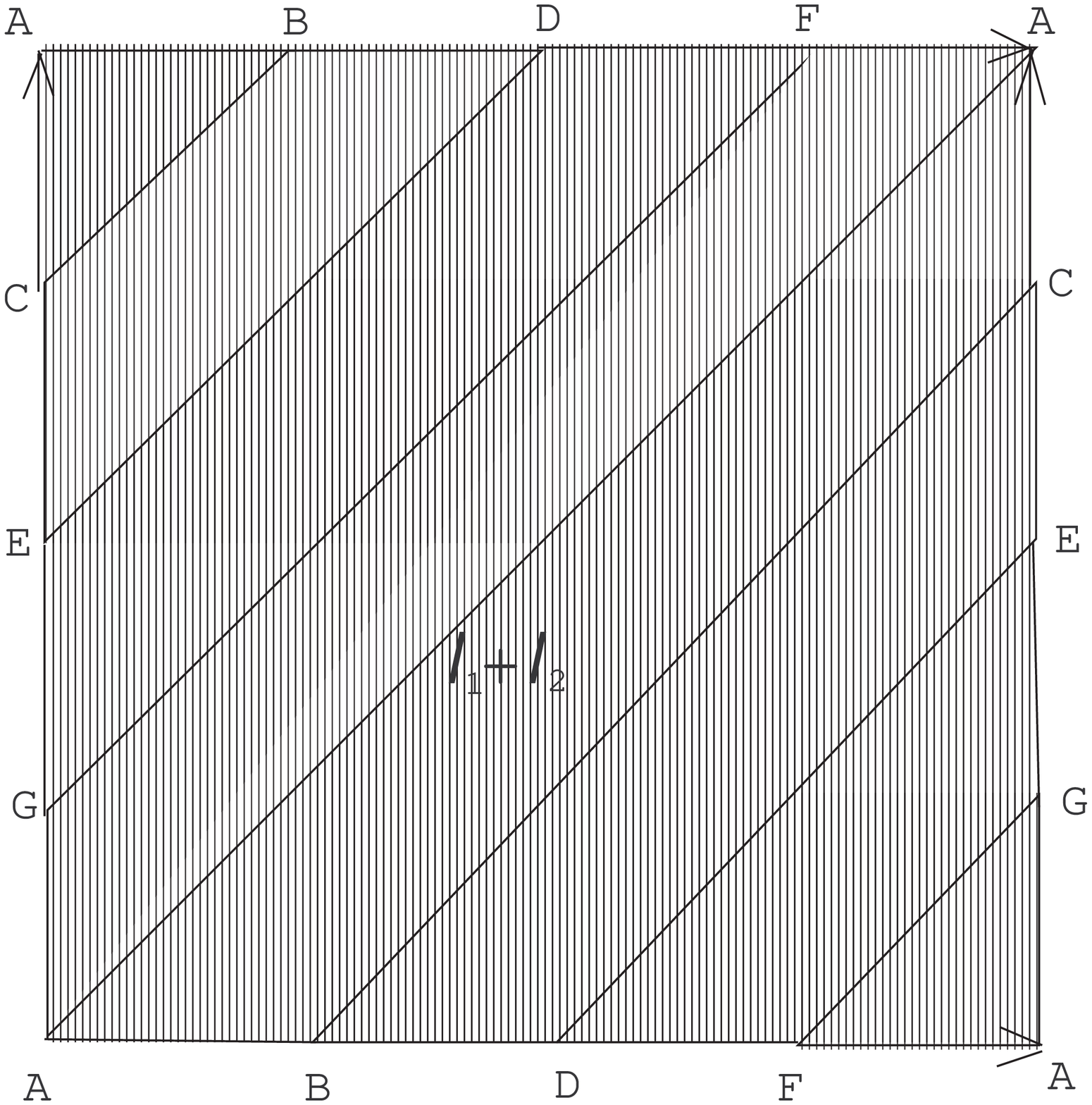}} 
\caption{$\bhh$: $A^+$ for $H=[h]=[e_1+e_2]$ 
with $(r,a,\delta_{\varphi S},v)=(10,8,0,0)$.}
\label{hyperbo-4h}
\end{figure}

\begin{figure} 
\centerline{\includegraphics[width=5cm]{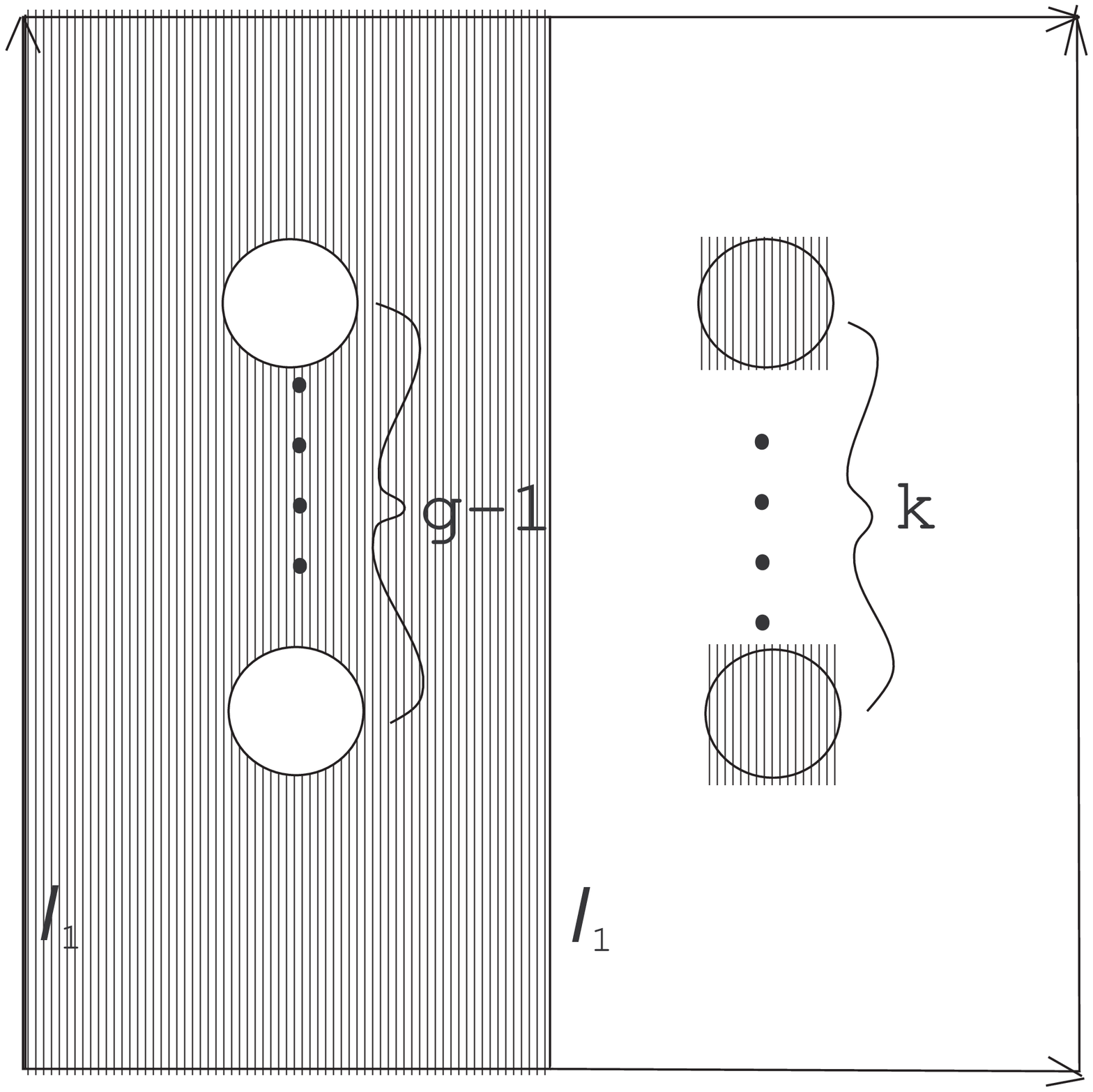}} 
\caption{$\bhh$: $A^+$ for $H=[e_1]$ 
with $(r,a,\delta_{\varphi S},v) \neq (10,10,0,0)\ \text{nor}\ (10,8,0,0)$.}
\label{hyperbo-e1}
\end{figure}

\begin{figure} 
\centerline{\includegraphics[width=5cm]{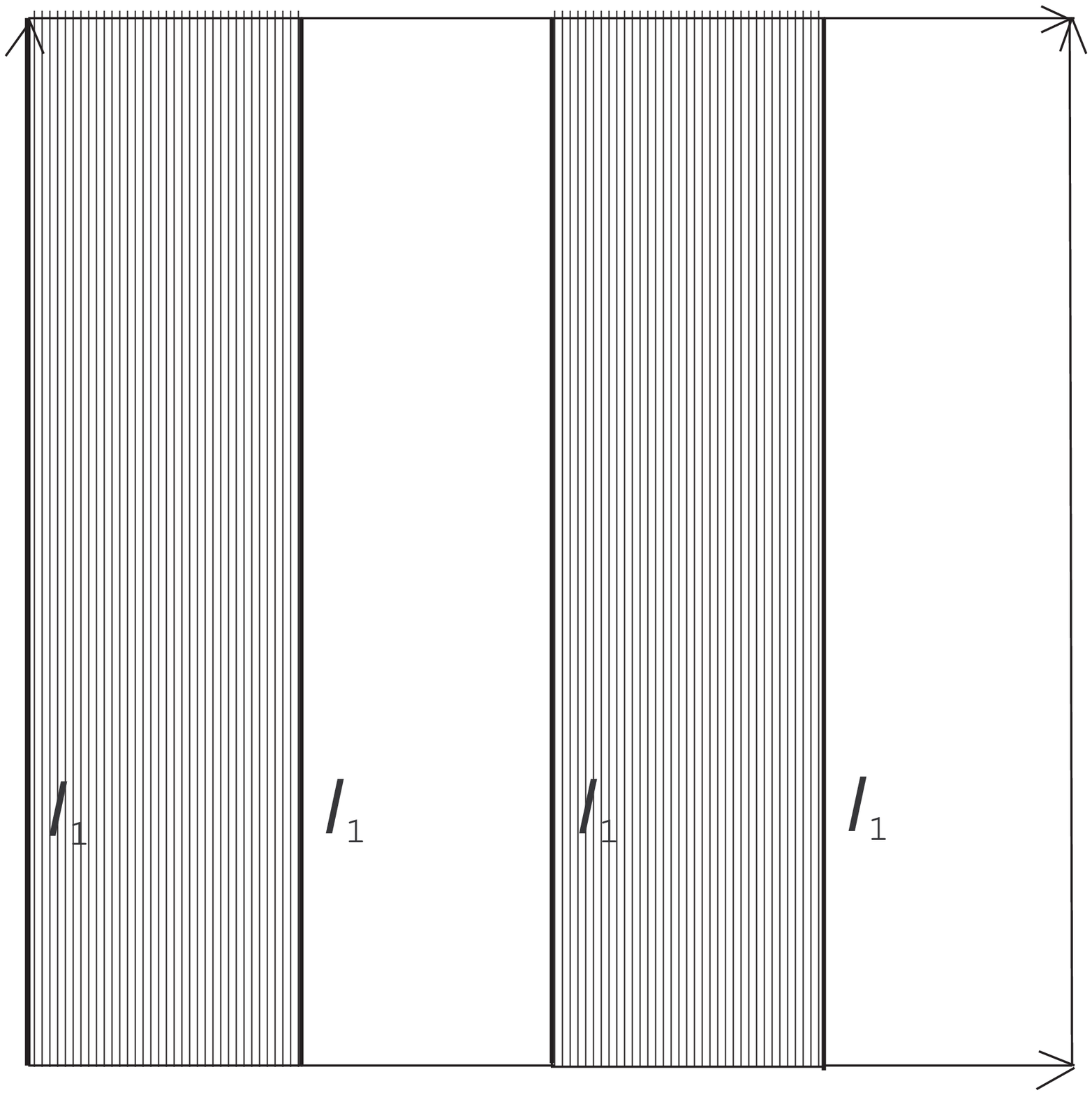}} 
\caption{$\bhh$: $A^+$ for $H=[e_1]$ 
with $(r,a,\delta_{\varphi S},v)=(10,8,0,0)$.}
\label{hyperbo-4e1}
\end{figure}

If a curve $A(\br)$ on hyperboloid has only ovals or $A(\br) = \emptyset$, 
then one of $A^+$ and $A^-$ contains the outermost component. 
Thus we divide isotopy types of positive curves $A^+$ into the 
following 4 cases:\\
\begin{description}
\item[(i)]$A(\br)$ has only ovals or $A(\br) = \emptyset$, and $A^+$ contains 
the outermost component. 
          (In this case we say $A^+$ is {\it outer}.)
\item[(ii)]$A(\br)$ has only ovals or $A(\br) = \emptyset$, and $A^+$ does 
not contain the outermost 
component. (In this case we say $A^+$ is {\it inner}.) E. g. see Figures  \ref{hyperbo}, \ref{hyperbo1111}. 
\item[(iii)]$A(\br)$ has even branches. E. g. see Figures \ref{hyperbo-h}, 
\ref{hyperbo-4h}  
\item[(iv)]$A(\br)$ has odd branches. E. g. see Figures \ref{hyperbo-e1}, 
\ref{hyperbo-4e1}. 
\end{description}

The subgroup $H=H_-$ is determined by the invariants 
$\delta_{e_1},\ \delta_{e_2}$ and $\delta_h$. 
We have the following interpretation of the invariant $\delta_{x_-}$. 

\begin{theorem}\label{realcurve}
Let $X$ be a compact K\"{a}hler surface with an 
anti-holomorphic involution $\varphi$, 
and $X(\br)$ be the fixed point set of $\varphi$. 
We assume that $H_1(X;\bz) = 0$ and $X(\br) \neq \emptyset$. 
We set $L = H_2(X;\bz)$ 
and $L_\varphi = \{ x \in L \ |\ \varphi_* (x) = - x \}$. 
Let $C$ be a $1$-dimensional complex submanifold of $X$ with 
$\varphi(C)=C$, and 
$C(\br)$ be the fixed point set of $\varphi$ on $C$. 
Let $[C]$ and $[C(\br)]$ denote the homology classes represented by 
$C$ and $C(\br)$ respectively. 
We set $x_- := [C] \ \ (\in L_\varphi)$. 
Then we have 
$$\delta_{x_-} = 0\ \ \ \Longleftrightarrow \ \ \ 
[C(\br)] = 0\ \ \mbox{in}\ H_1(X(\br);\bz/2\bz),$$
where we define the symbol $\delta_{x_-}$ in the same way as 
$\delta_{x_-}$ in Subsect. \ref{genclassification}. 
\end{theorem}
\Proof 
Suppose that $[C(\br)] \neq 0\ \mbox{in}\ H_1(X(\br);\bz/2\bz)$. 
Then there exists an embedded circle $D$ on $X(\br)$ 
such that $[C(\br)] \cdot [D] \neq 0$, 
where $[C(\br)] \cdot [D]$ means the $\bz/2\bz$-intersection number in 
$X(\br)$. We set 
$E_- = \{ x \in H^2(X;\bz) \ |\ \varphi^* (x) = - x \}$. 
$E_-$ is the Poincar\'{e} dual to $L_\varphi$. 
By Th\'{e}or\`{e}me 2.4 in \cite{Mangolte97}, 
there exists a surjective canonical morphism 
$$\varphi_X : E_- \rightarrow H^1(X(\br);\bz/2\bz)$$
such that 
$$(\varphi_X(\gamma),\varphi_X(\gamma')) \equiv Q(\gamma,\gamma') \pmod{2} 
\ \ \forall \gamma,\gamma' \in E_-,$$
where $(\ ,\ )$ and $Q$ are the forms induced by the cup products 
on $H^1(X(\br);\bz/2\bz)$ and $H^2(X;\bz)$. 
Moreover, by Th\'{e}or\`{e}me 2.5 in \cite{Mangolte97}, 
the following diagram commutes: 
$$
\begin{array}{ccc}
\Pic(X)^G            &  \stackrel{c_1}{\longrightarrow}  &      E_-        \\
\alpha \downarrow    &                              & \varphi_X \downarrow  \\
H^1_{\alg}(X(\br);\bz/2\bz)  &\stackrel{i}{\longrightarrow} & 
H^1(X(\br);\bz/2\bz), 
\end{array}
$$
where $\alpha$ is defined as in \cite{Mangolte97},p.562 (see also below), 
$c_1$ is the first Chern class map, and 
$i$ is the inclusion map. 
In the sequel, $A^P$ denotes the Poincar\'{e} dual element to a (co)homology 
class $A$. 
Since $\varphi_{X}$ is surjective, there exists 
$\gamma \ (\in E_-)$ such that $\varphi_X(\gamma) = [D]^P \ 
(\in H^1(X(\br);\bz/2\bz))$. 
We consider the {\it divisor class} $[C]$ in $\Pic(X)^G$. 
Then, as is well known, its first Chern class $c_1([C])$ is 
the Poincar\'{e} dual to $x_-$. 
By the definition (see \cite{Mangolte97}) of $\alpha$, we see 
$$\alpha([C])=\eta(C)=[C(\br)]^P\ (\in H^1_{\alg}(X(\br);\bz/2\bz)).$$
On the other hand, by the above commutative diagram, we see 
$$\alpha([C])=\varphi_X(c_1([C]))=\varphi_X(x_-^P).$$
Hence, we have 
$[C(\br)] \cdot [D]
=([C(\br)]^P,[D]^P)
=(\varphi_X(x_-^P),\varphi_X(\gamma))
\equiv Q(x_-^P,\gamma) \pmod{2}
= (x_-, \gamma^P)$. 
Thus we have $\delta_{x_-} =1$. 
This completes the proof of the implication $\Longrightarrow$. 
The converse assertion can be proved by the same argument as 
the proof of Lemma 2 in \cite{Matsuoka(Saito)92}. 
Suppose that $[C(\br)] = 0\ \mbox{in}\ H_1(X(\br);\bz/2\bz)$. 
Then $X(\br)$ and $C(\br)$ satisfy the conditions a) and b) 
of Remark 2.2 in \cite{Kharlamov75a}. By that remark, 
Lemma 2.3 is applicable to the involution $\varphi : X \to X$ and $C$. 
Hence, we see $(x_-)_{\mathrm{mod}\ 2}\ (\in H_2(X;\bz/2\bz))$ is orthogonal 
to $\mathrm{Im} \alpha_2$. 
Since $H_1(X;\bz) = 0$, as in the proof of Lemma 3.7 in \cite{Kharlamov76}, 
we have 
$\mathrm{Im} \alpha_2 = 
\{ x \in H_2(X;\bz/2\bz) \ |\ \varphi_* (x) = x \}$. 
Thus we have $\delta_{x_-} = 0$.  \QED
%%%%%%%%%%%%%%%%%%%

By Theorem \ref{realcurve}, we have 
$H=0$ if and only if $A^+$ is outer, 
$H=S/2S$ if and only if $A^+$ is inner, 
$H=[h]$ if and only if $A(\br)$ has even branches, 
$H=[e_1]$ if and only if $A(\br)$ has odd branches with odd $s$, and 
$H=[e_2]$ if and only if $A(\br)$ has odd branches with odd $t$. 

When $A$ is a dividing curve on hyperboloid with non-ovals,  
we define the number $\hat{l}\ (\in \bz/2\bz)$ 
as follows (see \cite{Matsuoka(Saito)91}): \ 
Fix a complex orientation of $A(\br)$. 
When the number of non-ovals of $A(\br)$ is $2$, we define 
$\hat{l}=0$ 
if the complex orientations of the $2$ non-ovals are different in 
$\br\bp^1\times \br\bp^1$, 
$\hat{l}=1$ if otherwise. 
When the number of non-ovals of $A(\br)$ is $4$, we fix a non-oval $E$ and 
define 
$\hat{l}=0$ 
if the number of non-ovals 
whose complex orientations are the same as $E$ in $\br\bp^1\times \br\bp^1$ 
is even, 
$\hat{l}=1$ if odd.  

We quote the following lemma. 
\begin{lemma}[\cite{Matsuoka(Saito)91}]\label{lemma:X(R)}
For a positive curve $(A,\varphi)$ of bidegree $(4,4)$ on hyperboloid, 
we have the following:\\
{\rm (1)}\ $[X_\varphi(\br)] = [X_{\tau\varphi}(\br)]$ in $H_2(X;\bz/2\bz)$\\
{\rm (2)}\ If $A$ is dividing, then 
$$[X_\varphi(\br)] \ (= [X_{\tau\varphi}(\br)]) = \left\{
\begin{array}{cl}
0           & (if\ \ A(\br)\ has\ only\ ovals)\\
\hat{l} h   & (if\ \ A(\br)\ has\ even\ branches)\\
\hat{l} e_1 & (if\ \ A(\br)\ has\ odd\ branches\ with\ odd\ s)\\
\hat{l} e_2 & (if\ \ A(\br)\ has\ odd\ branches\ with\ odd\ t)
\end{array}
\right.
$$
in $H_2(X;\bz/2\bz).$  \QED
\end{lemma}
Then we have the following:
\begin{proposition}\label{prop:div =>}
Let $A$ be a non-singular real curve of bidegree $(4,4)$ on hyperboloid. 
If $A$ is dividing, 
then 
the positive curves $(A,\varphi)$ satisfies $\delta_{\varphi S} = 0$. 
\end{proposition}
\Proof 
If $A$ is dividing, then by (2) of Lemma \ref{lemma:X(R)}, 
we have $[X_\varphi(\br)] \equiv s_{\varphi}$ (mod $2H_2(X;\bz)$) 
for some $s_{\varphi}\ \in S$. 
Hence we have $\delta_{\varphi S} = 0$.  \QED

\medskip

Moreover, we have the following interpretation of $v$ 
when $A$ is a dividing curve with non-ovals on  hyperboloid:\\
$v = 0$              if and only if $\hat{l} = 0$,\\
$v = h$   (mod $2S$) if and only if $A(\br)$ has even branches and 
$\hat{l} = 1$,\\
$v = e_1$ (mod $2S$) if and only if $A(\br)$ has odd branches with odd $s$ 
and $\hat{l} = 1$,\\
$v = e_2$ (mod $2S$) if and only if $A(\br)$ has odd branches with odd $t$ 
and $\hat{l} = 1$. 

\medskip

By the above geometric interpretations and Zvonilov's classification 
\cite{Zvonilov92}, from the invariants \eqref{invhyperboloid} of 
the connected component of moduli of a positive curve $A^+$ on hyperboloid, 
we can uniquely determine the complex  
scheme of $A$.  
In particular, one can draw the picture of the positive curve $A^+$.  
By \eqref{hyperbo-related1} and \eqref{hyperbo-related2}, 
it is enough to draw them for $H=S/2S$, $[h]$ and $[e_1]$.  
See Figures \ref{hyperbo} --- \ref{hyperbo-4e1}. 
Remind that 
$A(\br)$ is an empty curve on hyperboloid if 
$H=S/2S$ and $(r,a,\delta_{\varphi S},v)=(10,10,0,0)$, 
$Y$ has the (spin,spin) structure if 
$H=[h]$ and $(r,a,\delta_{\varphi S},v)=(10,10,0,0)$, 
and $Y$ has the (usual,spin) (or (spin,usual)) structure if 
$H=[e_i]$ and $(r,a,\delta_{\varphi S},v)=(10,10,0,0)$.

Moreover, there is known (due to Zvonilov \cite{Zvonilov92}) 
a dividing curve on hyperboloid for each connected component with 
$\delta_{\varphi S}=0$. Hence, the opposite statement of Proposition 
\ref{prop:div =>} is true:

\begin{theorem}\label{prop:div}
Let $A$ be a non-singular real curve of bidegree $(4,4)$ on hyperboloid.  
Then, $A$ is dividing or $A(\br)=\emptyset$, 
if and only if the positive curves $(A,\,\varphi)$ (or $A^+$) 
has $\delta_{\varphi S} = 0$. \QED 
\end{theorem}

\medskip

%%If $\delta_{\varphi S}=0$ and $H=0$, then we have $v=0$. 
%%If $\delta_{\varphi S}=0$ and $H=S/2S$, then we also have $v=0$ by 
%%\eqref{hyperbo-related1} and \eqref{hyperbo-related2}. 
It turns out that 
if $\delta_{\varphi S}=0$ and $H=[h]$, then $(r,a)$ determines $v$. 
But when $\delta_{\varphi S}=0$ and $H=[e_i]$, $(r,a)$ does not 
always determine $v$. 
As stated above, when $\delta_{\varphi S}=0$, $H=[e_i]$ and 
$(r,a,v)\neq(10,10,0)$, 
$v = 0$ if and only if $\hat{l} = 0$. 

\medskip

Thus, we finally get

\begin{theorem}
A connected component of moduli of 
a positive  real non-singular curve $A^+$ of bidegree $(4,4)$ on hyperboloid 
is defined by the isotopy type of $A^+ \subset \br\bp^1\times \br\bp^1$ 
(up to the action of $(PGL(2,\br)\times PGL(2,\br)) \rtimes \bz/2\bz$), 
divideness of $A(\bc)$ by $A(\br)$, 
and by the invariant $\hat{l}\mod 2$ defined by the complex orientation 
(if $A(\br)$ has odd branches and is dividing). 
\end{theorem}

\bigskip

\section{Connected components of moduli of real non-singular curves of 
bidegree $(4,4)$ on ellipsoid}
\label{ellipsoid}

Here we consider the case $S\cong U(2)$ and $\rk S_-=1$ (ellipsoid cases). 
Like for hyperboloid, $Y=\bp^1\times\bp^1$, and we use the same notations 
for $S$. Obviously, 
$\theta(e_1)= - e_2$ and $\theta(e_2)= - e_1$, 
we have $S_+ = \bz (e_1-e_2)$ and $S_- = \bz (e_1+e_2)$. 
Thus $S_+\cong \langle -4 \rangle$ and $S_-\cong \langle 4 \rangle$. 

Since 
$\Delta(S_+, S)^{(-4)} = \{ \pm(e_1 - e_2) \}$ and 
$\Delta(S_-, S)^{(-4)} = \emptyset$, 
$W^{(-4)}(S_+,S)$ consists of identity and $g$ 
where $g(e_1)=e_2$ and $g(e_2)=e_1$, and 
$W^{(-4)}(S_-,S)$ is trivial. 
By Lemma \ref{GroupG}, we have $G = \{ \mbox{identity},\ g \} \cong \bz/2\bz$. 

We call these data $(S,\theta, G)$ as {\it the ellipsoid case}.

For the ellipsoid case we have $s=2$, $p=1$. 
$A_{S_+} \cong \bz/4\bz$ is generated by $\frac14(e_1-e_2)$, and 
$A_{S_-} \cong \bz/4\bz$ is generated by $\frac14(e_1+e_2)$. 
We see 
$\Gamma_+ = [ e_1-e_2 ] = 2({S_+}^\ast \cap (\frac12 S_+))/2S_+$ and 
$\Gamma_- = [ e_1+e_2 ]  = 2({S_-}^\ast \cap (\frac12 S_-))/2S_-$. 
Hence we have 
$H_\pm = \Gamma_\pm$. 
Since 
$(e_1-e_2,0)=(0,e_1+e_2)$ in $H=H_+\oplus_\gamma H_-$, we see 
$$H = \{ 0,\ (e_1-e_2,0)=(0,e_1+e_2) \} \cong \bz/2\bz,$$ 
$a_H = 1$ and $q_\rho \cong w$. 
Hence $\delta_H = 0$, $k_\rho = \mu_\rho = 1$. 

If $\delta_\varphi = 0$ (Type 0), then $c_v \equiv 0 \pmod{4}$. 
If $\delta_{\varphi S}=0$ and $\delta_\varphi = 1$ (Type Ia), 
then $v = (e_1-e_2,0)=(0,e_1+e_2)$ and 
$c_v \equiv 2 \pmod{4}$. 

We have $\delta_{\varphi S} = \delta_{\varphi S_\pm}.$ 
Assume that $\delta_{\varphi S} =0$. Then  
$$
\gamma_\pm=
\left\{
\begin{array}{ll}
0         & \text{if $\delta_\varphi=0$},\\
q_1(2)^2  &\text{if $\delta_{\varphi S}=0$ and $\delta_\varphi = 1$.}
\end{array}\right. 
$$
If $\delta_\varphi=0$, then $(q_{S_\pm})_v = q_{S_\pm}$. 
Hence we have $\varepsilon_{v_\pm}=0$. 
If $\delta_{\varphi S}=0$ and $\delta_\varphi = 1$, 
then $(q_{S_\pm})_v$ is defined on $\bz/4\bz$, and 
the square of a generator is $\pm \frac54 \pmod{2\bz}.$ 
Hence $(q_{S_\pm})_v$ is the discriminant quadratic form of 
the $2$-adic lattice 
$K^{(2)}_{\pm 5}(2^2) = \langle \pm 5\cdot 2^2 \rangle$ (\cite{Nikulin79}). 
Thus we have $\varepsilon_{v_\pm}=1$. 

\medskip

We see that the genus of  
an integral involutions $(L,\varphi,S)$ is determined by the data 
$$
(r,a,\delta_{\varphi S},\delta_\varphi).
$$ 
By Theorem \ref{geninvtheorem}, the complete list of the above invariants 
for ellipsoid is given by the conditions:\\

\noindent
{\bf Type 0} ($\delta_\varphi=0$)\\
$2 \leq a \leq r$ and $a$ is even.\\
$r = 2,6,10,14,18.$\\
$r+a \leq 22.$\\
If $a = r$, then $(r,a)=(2,2),\ (10,10).$\\
If $r+a = 22$, then $(r,a)=(18,4).$\\

\noindent
{\bf Type Ia} ($\delta_{\varphi S}=0,\ \delta_\varphi=1$)\\
$2 \leq a \leq r$ and $a$ is even.\\
$r = 4,8,12,16,20.$\\
$r+a \leq 22.$\\
If $a = r$, then $(r,a)=(4,4).$\\
If $r+a = 22$, then $(r,a)=(20,2),\ (12,10).$\\

\noindent
{\bf Type Ib} ($\delta_{\varphi S}=1$)\\
$3 \leq a \leq r.$\\
$r+a$ is even.\\
$2 \leq r \leq 20.$\\
$r+a \leq 22.$\\

\medskip

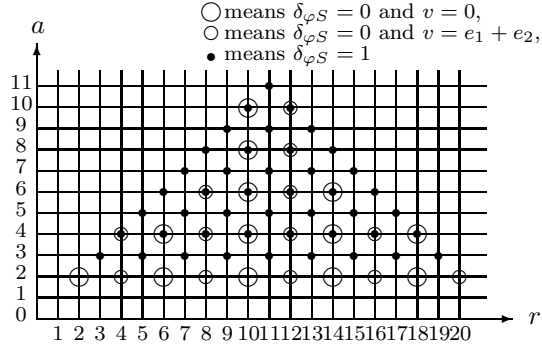
\begin{figure}
\begin{picture}(200,140)
%%ellipsoid
\put(66,116){\circle{7}}
\put(71,114){{\tiny means $\delta_{\varphi S}=0$ and $v=0$,}}
\put(66,108){\circle{5}}
\put(71,106){{\tiny means $\delta_{\varphi S}=0$ and $v=e_1+e_2$,}}
\put(66, 99){\circle*{3}}
\put(71, 98){{\tiny means $\delta_{\varphi S}=1$}}
\multiput(8,0)(8,0){20}{\line(0,1){94}}
\multiput(0,8)(0,8){11}{\line(1,0){170}}
\put(0,0){\vector(0,1){104}}
\put(0,0){\vector(1,0){180}}
\put(  6,-8){{\tiny $1$}}
\put( 14,-8){{\tiny $2$}}
\put( 22,-8){{\tiny $3$}}
\put( 30,-8){{\tiny $4$}}
\put( 38,-8){{\tiny $5$}}
\put( 46,-8){{\tiny $6$}}
\put( 54,-8){{\tiny $7$}}
\put( 62,-8){{\tiny $8$}}
\put( 70,-8){{\tiny $9$}}
\put( 76,-8){{\tiny $10$}}
\put( 84,-8){{\tiny $11$}}
\put( 92,-8){{\tiny $12$}}
\put(100,-8){{\tiny $13$}}
\put(108,-8){{\tiny $14$}}
\put(116,-8){{\tiny $15$}}
\put(124,-8){{\tiny $16$}}
\put(132,-8){{\tiny $17$}}
\put(140,-8){{\tiny $18$}}
\put(148,-8){{\tiny $19$}}
\put(156,-8){{\tiny $20$}}

\put(-8, -1){{\tiny $0$}}
\put(-8,  7){{\tiny $1$}}
\put(-8, 15){{\tiny $2$}}
\put(-8, 23){{\tiny $3$}}
\put(-8, 31){{\tiny $4$}}
\put(-8, 39){{\tiny $5$}}
\put(-8, 47){{\tiny $6$}}
\put(-8, 55){{\tiny $7$}}
\put(-8, 63){{\tiny $8$}}
\put(-8, 71){{\tiny $9$}}
\put(-10, 79){{\tiny $10$}}
\put(-10, 87){{\tiny $11$}}

\put( -2,108){{\footnotesize $a$}} %vertical
\put(186, -2){{\footnotesize $r$}} %horizontal

%% Type 0  \circle{7}
%%  ( 8r,8a)
\put( 16,16){\circle{7}}
\put( 48,16){\circle{7}}
\put( 80,16){\circle{7}}
\put(112,16){\circle{7}}
\put(144,16){\circle{7}}
\put( 48,32){\circle{7}}
\put( 80,32){\circle{7}}
\put(112,32){\circle{7}}
\put(144,32){\circle{7}}
\put( 80,48){\circle{7}}
\put(112,48){\circle{7}}
\put( 80,64){\circle{7}}
\put( 80,80){\circle{7}}

%% Type Ia  \circle{5}
%%  ( 8r,8a)
\put( 32,16){\circle{5}}
\put( 64,16){\circle{5}}
\put( 96,16){\circle{5}}
\put(128,16){\circle{5}}
\put(160,16){\circle{5}}
\put( 32,32){\circle{5}}
\put( 64,32){\circle{5}}
\put( 96,32){\circle{5}}
\put(128,32){\circle{5}}
\put( 64,48){\circle{5}}
\put( 96,48){\circle{5}}
\put( 96,64){\circle{5}}
\put( 96,80){\circle{5}}

%% Type Ib  \circle*{3}
%%  ( 8r,8a)
\put( 24,24){\circle*{3}}
\put( 40,24){\circle*{3}}
\put( 56,24){\circle*{3}}
\put( 72,24){\circle*{3}}
\put( 88,24){\circle*{3}}
\put(104,24){\circle*{3}}
\put(120,24){\circle*{3}}
\put(136,24){\circle*{3}}
\put(152,24){\circle*{3}}
\put( 32,32){\circle*{3}}
\put( 48,32){\circle*{3}}
\put( 64,32){\circle*{3}}
\put( 80,32){\circle*{3}}
\put( 96,32){\circle*{3}}
\put(112,32){\circle*{3}}
\put(128,32){\circle*{3}}
\put(144,32){\circle*{3}}
\put( 40,40){\circle*{3}}
\put( 56,40){\circle*{3}}
\put( 72,40){\circle*{3}}
\put( 88,40){\circle*{3}}
\put(104,40){\circle*{3}}
\put(120,40){\circle*{3}}
\put(136,40){\circle*{3}}
\put( 48,48){\circle*{3}}
\put( 64,48){\circle*{3}}
\put( 80,48){\circle*{3}}
\put( 96,48){\circle*{3}}
\put(112,48){\circle*{3}}
\put(128,48){\circle*{3}}
\put( 56,56){\circle*{3}}
\put( 72,56){\circle*{3}}
\put( 88,56){\circle*{3}}
\put(104,56){\circle*{3}}
\put(120,56){\circle*{3}}
\put( 64,64){\circle*{3}}
\put( 80,64){\circle*{3}}
\put( 96,64){\circle*{3}}
\put(112,64){\circle*{3}}
\put( 72,72){\circle*{3}}
\put( 88,72){\circle*{3}}
\put(104,72){\circle*{3}}
\put( 80,80){\circle*{3}}
\put( 96,80){\circle*{3}}
\put( 88,88){\circle*{3}}
\end{picture}

\caption{$\bee$: All possible $(r,a,\delta_{\varphi S},v)$ 
(here $H=[e_1+e_2]$)} \label{ellipsoid-graph}
\end{figure}

All possible data 
\begin{equation}
(r,a,\delta_{\varphi S},v)
\label{geninvellipsoid}
\end{equation}
 which are 
equivalent to $(r,a,\delta_{\varphi S},\delta_\varphi)$ 
are given in Figure \ref{ellipsoid-graph}. We have $\delta_{\varphi}=0$, 
if and only if $\delta_{\varphi S}=v=0$.  
Thus there are 
13 cases of Type 0, and 13 cases of Type Ia, and 
45 cases of Type Ib.

\bigskip

\bigskip

The relations between related involutions are:  
\begin{equation}
r(\varphi)+r(\tau\varphi)=22,\ \ a(\varphi) = a(\tau\varphi),
\label{ell-related1}
\end{equation}
\begin{equation}
\delta_{\varphi S}=\delta_{\tau\varphi S},\ \  
s_\varphi + s_{\tau\varphi} \equiv e_1 +e_2  \mod 2L.
\label{ell-related2}
\end{equation}
Hence, if a positive curve $A^+$ satisfies $\delta_\varphi = 0$, 
then the related positive curve $A^-$ 
satisfies $(\delta_{\tau\varphi S},\delta_{\tau\varphi}) = (0,1)$.  
If a positive curve $A^+$ satisfies $\delta_{\varphi S}= 1$ 
and $r(\varphi)\le 11$,  
then the related positive curve $A^-$  
satisfies $\delta_{\tau\varphi S}=1$ and $r(\tau \varphi)\ge 11$.

For the ellipsoid case, conditions $(b_+)$ and $(b_-)$ of 
Theorem \ref{genus-isomor} are not valid in all cases. We have to 
check conditions $(a_+)$ and $(a_-)$ in the rest cases directly. 
As the result, we have 

\begin{theorem} 
\label{gencl-ellipsoid}
For the ellipsoid case, genus always defines the class of 
an integral involution. In particular (by Theorems \ref{genmod}, 
\ref{gnmoddpn+} and \ref{geninvtheorem}), the genus invariants 
$$
(r,a,\delta_{\varphi S},\delta_\varphi)
$$ 
define the connected component of muduli of positive non-singular 
curves of bidegree $(4,4)$ on ellipsoid up to the action of 
the automorphism group of ellipsoid. There are exactly 
$13$ connected components of moduli with $\delta_\varphi=0$, exactly 
$13$ connected components of moduli with 
$(\delta_{\varphi S},\delta_\varphi)=(0,1)$, and exactly $45$ connected 
components of moduli with $\delta_{\varphi S}=1$.

Identifying related positive curves,  
we get that there are exactly $38$ connected components of moduli 
of real non-singular curves of bidegree $(4,4)$ on ellipsoid 
up to the action of the automorphism group of ellipsoid. 
There are $13$ of them with $\delta_\varphi=0$ 
(equivalently, $(\delta_{\varphi S},\delta_\varphi)=(0,1)$), and 
$25$ connected components with $\delta_{\varphi S}=1$: they are distinguished 
from $45$ by the inequality $r(\varphi)\le 11$ (equivalently, by 
$r(\varphi)\ge 11$). 
\end{theorem}

\Proof For the ellipsoid case, $a_{H_+}=1$, $l(A_{S_+})=1$, $p=1$. 
Moreover, $\kappa (A_{S_+})=0$. Thus, the 
condition $(b_+)$ of Theorem \ref{genus-isomor} means that 
either $r-a>0$, or $r-a=0$ and $a\ge 5$, or $r-a=0$, $a=4$ and 
$\delta_{\varphi S_+}=0$ (for our case,  
$\delta_{\varphi S_+}=\delta_{\varphi S}$). 
From the list of possible invariants above, we see 
that this condition is not satisfied only in the following cases: 

(a) $\delta_\varphi=0$ and $r=a=2$;

(b) $\delta_{\varphi S}=1$ and $r=a=3$; 

(c) $\delta_{\varphi S}=1$ and $r=a=4$.  

Let us consider all these cases. 

Case (a): In this case, the lattice $L^\varphi\cong U(2)$, and 
the lattice $K_+$ is the orthogonal complement to the primitive 
sublattice $S_+\subset U(2)$ where $S_+\cong \langle -4 \rangle$. 
It follows that $K_+\cong \langle 4 \rangle$. Obviously, the 
condition $(a_+)$ of Theorem \ref{genus-isomor} is valid in this case. 

In cases (b) and (c), the invariant $\delta_{\varphi}=1$, and the 
lattice $L^\varphi\cong \langle 2 \rangle \oplus (r-1)\langle -2 \rangle$.  
The lattice $K_+$ is the orthogonal 
complement to the primitive sublattice 
$S_+\cong \langle -4 \rangle \subset L^\varphi$. Equivalently, 
$K_+$ is the orthogonal complement $\delta^\perp$ to a primitive  
element $\delta\in L^\varphi$ with $\delta^2=-4$. Since 
$\delta_{\varphi S}=\delta_{\varphi S_+}=1$, 
the discriminant form $q_{K_+}$ must be isomorphic to 
$q_{\pm 1}(2)\oplus q^\prime$ in this case. 
Here $q_{\pm 1}(2)$ is a non-degenerate quadratic form 
on the group of order $2$. Equivalently, there exists  
$x\in A_{K_+}$ of order 2 such that 
$q_{K_+}(x)\not\equiv 0\mod 1$. See (\cite{Nikulin83}, page 116).

Considering the odd unimodular lattice 
$L=L^\varphi(1/2)\cong \langle 1 \rangle \oplus (r-1)\langle -1 \rangle$ 
instead of $L^\varphi$, we should prove the following 
for $r=3,\,4$: For all elements  
$\delta \in L$ with $\delta^2=-2$ and with odd $K=\delta^\perp$, 
the lattices $K=\delta^\perp$ are isomorphic, and the 
canonical homomorphism  
\begin{equation}
\xi: O(K)\to O(q_{K(2)})
\label{epimorphic}
\end{equation}
is epimorphic. Now let us consider cases (b) and (c) 
separately. 

Case (b). In this case, by simple calculations, the discriminant form 
$q_{K(2)}\cong q_{1}(2)\oplus q_{-1}(4)$. Its automorphism 
group $O(q_{K(2)})=\pm 1$ is the image of the automorphism  
subgroup $\pm 1$ of $O(K)$.  Thus, \eqref{epimorphic} is valid. 
To prove that all lattices $K=\delta^\perp$ are isomorphic, one 
should argue similarly to the case (c) below which is more 
complicated.

Case (c). In this case, 
$q_{K(2)}\cong q_{1}(2)\oplus q_{-1}(2) \oplus q_{-1}(4)$  
(this is clear, but we shall give the proof in considerations below). 
By simple calculations, $O(q_{K(2)})\cong (\bz/2)^2$. 

Let us consider the group $W(L)$ generated by reflections in all 
elements of $L$ with squares $-1$ and $-2$. The Dynkin diagram of 
the fundamental chamber of $W(L)$ is well-known (e. g. see 
\cite{Vinberg72}). It is (with usual notations) 

\begin{equation}
\centerline{$\bullet$ --- $\bullet$ --- $\circ$ --- $\star$}
\label{diagram1}
\end{equation}

Here the two $\bullet$ correspond to elements $\delta\in L$ 
with $\delta^2=-2$ and with odd orthogonal complement 
$K=\delta^\perp$. 
They (i. e. their $O(L)$-orbits) give all elements of $L$ we are 
looking for. In particular, there exist no more than two possibilities 
for $K$ up to isomorphisms. 
The $\circ$ corresponds to elements $e\in L$ with $e^2=-1$. 
The $\star$ corresponds to elements $v\in L$ with $v^2=-2$ and 
with even $v^\perp$ in $L$. 

Let $e_0, e_1, e_2, e_3$ be the standard basis of $L$ with the 
diagonal matrix $diag(1,-1,-1,-1)$. The diagram \eqref{diagram1} is 
given (from the left) by $\delta_2=-e_2+e_3$, $\delta_1=-e_1+e_2$, 
$e=e_1$ and $v=e_0-e_1-e_2-e_3$. The permutation  
$(e_0,e_1,e_2,e_3)\to (e_0, e_2,e_3,e_1)$ gives the automorphism of 
$L$ such that $\delta_1\mapsto \delta_2$. Thus the lattices $K$ 
corresponding to $K=\delta_1^\perp$ and $K=\delta_2^\perp$ are 
isomorphic. Any element $\delta\in K$ with $\delta^2=-2$ and with 
odd $\delta^\perp$ belongs to $O(L)$-orbit of $\delta_1$ or $\delta_2$.  

The lattice $K=\delta_2^\perp$ has the basis 
$e_0$, $e_1$, $\delta=e_2+e_3$ with the diagonal matrix 
$diag(1,-1,-4)$. It follows  that $q_{K(2)}\cong q_1(2)\oplus 
q_{-1}(2)\oplus q_{-1}(4)$ with the corresponding bases 
$e_0/2$, $e_1/2$, $\delta/4$ $\mod K$. As we have mentioned, 
it is easy to see that $O(q_{K(2)})\cong (\bz /2)^2$. The automorphisms 
$\pm 1$ of the lattice $K$ give the automorphisms $\pm 1$ of $q_{K(2)}$. 
To prove \eqref{epimorphic}, it is sufficient to construct an  
automorphism of $K$ which does not give the automorphism $\pm 1$ in 
$O(q_{K(2)})$. 

This is the automorphism $f$: $e_0\mapsto 2e_0+e_1+\delta$, 
$e_1\mapsto e_0+\delta$, $\delta\mapsto 2e_0+2e_1+\delta$ of the lattice $K$.
One can easily check that this is an automorphism of $K$. We have that 
$f:\delta/4\mapsto e_0/2+e_1/2+\delta/4$. It follows that $f$ does not 
give $\pm 1$ in $O(q_{K(2)})$. Thus, \eqref{epimorphic} is 
epimorphic. 

In the case (b), the Dynkin diagram of 
$W(L)$ is 
\begin{equation}
\centerline{$\bullet$ --- $\circ$ --- $\circ$\ .}
\label{diagram2}
\end{equation}
  
Replacing $\phi$ by $\tau\phi$, we get the statement $(a_-)$ or 
Theorem 13.

It follows the theorem. 

%%%%%%%%%%%%%%%%%%%%%%%%%%%%%%%%%%%%

\bigskip

\begin{figure} 
\centerline{\includegraphics[width=5cm]{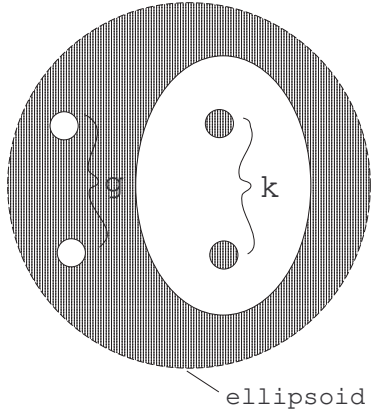}} 
\caption{$\bee$: $A^+$ 
for $(r,a,\delta_{\varphi S})\not=(10,10,0)$, $(12,10,0)$, 
$(10,8,0)$, $(12,8,0)$.}
\label{elli}
\end{figure}

\begin{figure} 
\centerline{\includegraphics[width=5cm]{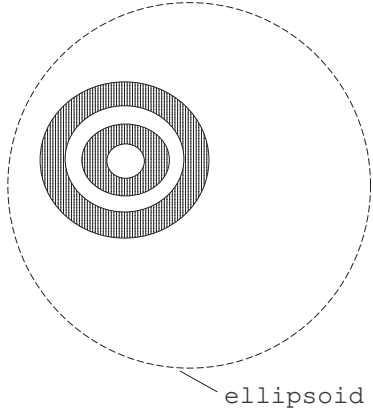}} 
\caption{$\bee$: $A^+$ for $(r,a,\delta_{\varphi S})=(10,8,0)$ \ 
($A^+=\emptyset$ for $(r,a,\delta_{\varphi S})=(10,10,0)$).}
\label{elli-nest}
\end{figure}

Now we consider the geometric interpretation of the above calculations. 
Using \eqref{realcomponents}, \eqref{realmod2} for both positive curves 
$A^+$, $A^-$ and 
\eqref{ell-related1}, \eqref{ell-related2}, we can easily draw the 
picture of $A^+$ on ellipsoid for all invariants \eqref{geninvellipsoid}.  
See Figures \ref{elli} and \ref{elli-nest}. As usual, we set 
$g=(22-r-a)/2$ and $k=(r-a)/2$.   
Thus we get the isotopy classification of real non-singular curves of 
bidegree $(4,4)$ on ellipsoid. We can see that the  
isotopy type is defined by the topology type of $A^+$. 
The isotopy types of $A$ are: 

\medskip

a nest of depth $4$;\ \ $\langle 0 \rangle$ (or $\emptyset$);\ \   
$\langle 1\langle m \rangle \sqcup n \rangle$ 
where $m+n \leq 9,\ m \leq n$. 

\medskip 

This isotopy classification was first obtained by Gudkov and Shustin \cite{GudkovShustin80}. As we have mentioned in Sect. \ref{hyperboloid}, 
Zvonilov \cite{Zvonilov92} classified complex schemes of these isotopy 
types. Below we partly repeat his results.

We have similar statement to Arnold \cite{Arnold71} 

\begin{proposition}\label{prop:div => (elli)}
Let $(A,\varphi)$ (or $A^+$) be a positive curve of bidegree $(4,4)$ on 
an ellipsoid. 
If $A$ is dividing, then 
$$[X_\varphi(\br)]=c (e_1 + e_2)$$
in $H_2(X;\bz/2\bz)$ 
for some $c \ (\in \bz/2\bz).$ 
In particular, we have $\delta_{\varphi S}=0$. 
\end{proposition}

\Proof
We use arguments like Wilson \cite{Wilson}, Lemma 6.7. 
Let $C^+$ and $C^-$ be the connected components of $A \setminus A(\br).$ 
$C^\pm$ have the orientations induced from $A.$ 
$A^+$ and $A^-$ are the halfs of $Y(\br) \setminus A(\br)$. 
$A^+ = \pi(X_\varphi(\br))$ and $A^- = \pi(X_{\tau\varphi}(\br))$. 
$Y(\br)$ (ellipsoid) is homeomorphic to the $2$-sphere $S^2$, 
and orientable. 
We give $A^\pm$ the orientations induced from an orientation of $Y(\br)$. 
We regard $C^\pm$ and $A^\pm$ as $\bz$-chain. 
Then we have $\partial(C^+ + A^+) = 2\partial D$ 
where $D$ is the sum of some $2$-disks on $Y(\br)$. 
Hence $C^+ + A^+ - 2D$ is $\bz$-cycle. 
We consider the anti-holomorphic involution $\theta=\varphi\mod\{1,\tau\}$ on $Y$. 
Then 
$\theta(\text{pt}\times \bp^1) = - (\bp^1 \times \text{pt}),$
$\theta(\bp^1 \times \text{pt})= - (\text{pt}\times \bp^1),$
$\theta(C^+)=-C^-$ and $\theta$ fix $\bz$-chains on $Y(\br)$. 

For any integer $m,n$, we have
$$(1+\theta)(m(\text{pt}\times \bp^1) + n(\bp^1 \times \text{pt}))
=(m-n)(\text{pt}\times \bp^1 - \bp^1 \times \text{pt}).$$
Hence 
$(1+\theta)(C^+ + A^+ - 2D) \sim c(\text{pt}\times \bp^1 - \bp^1 \times \text{pt})$ 
for some integer $c$, 
where $\sim$ means $\bz$-homologous. 
Since $A$ is of bidegree $(4,4)$, 
$C^+ + C^- \sim 4(\text{pt}\times \bp^1 + \bp^1 \times \text{pt})$. 
Thus we have 
$$
2[C^+ + A^+ - 2D - 2(\text{pt}\times \bp^1 + \bp^1 \times \text{pt})]
=
c[\text{pt}\times \bp^1 - \bp^1 \times \text{pt}]
$$
in $H_2(\bp^1\times \bp^1;\bz)$. 
Since $H_2(\bp^1\times \bp^1;\bz)$ is free, 
we have 
$
[C^+ + A^+ - 2D - 2(\text{pt}\times \bp^1 + \bp^1 \times \text{pt})]
=
c[\text{pt}\times \bp^1 - \bp^1 \times \text{pt}]
$
for some integer $c$. 
Hence we have 
$$C^+ + A^+ \sim c(\text{pt}\times \bp^1 + \bp^1 \times \text{pt})$$
as $\bz/2\bz$-cycles for some $c \ (\in \bz/2\bz)$. 
Taking the transfer (like \cite{Wilson}, Lemma 6.7) gives 
$[X_\varphi(\br)]=c (e_1 + e_2)$ 
in $H_2(X;\bz/2\bz)$. \QED

Suppose that $A$ is dividing. 
Then $\delta_{\varphi S}=0$ by Proposition \ref{prop:div => (elli)}. 
Since 
$s_\varphi \equiv [X_\varphi(\br)] \mod 2H_2(X;\bz)$ 
and 
$s_\varphi + s_{\tau\varphi} \equiv (e_1 - e_2) \mod 2H_2(X;\bz)$ (see above), 
we have 
$$
[X_\varphi(\br)] + [X_{\tau\varphi}(\br)] = (e_1 + e_2)
$$
in $H_2(X;\bz/2\bz)$. 
Hence, by Proposition \ref{prop:div => (elli)}, we see 
$[X_{\tau\varphi}(\br)]=(e_1 + e_2)$ if $[X_\varphi(\br)]=0$, and 
$[X_{\tau\varphi}(\br)]=0$ if $[X_\varphi(\br)]=(e_1 + e_2)$ 
in $H_2(X;\bz/2\bz)$.

By Proposition \ref{prop:div => (elli)}, 
if $\delta_{\varphi S}=1$, then the curves are not dividing. 
Hence, the following $6$ isotopy types 
$\langle 7 \sqcup 1\langle 0 \rangle \rangle$,\ 
$\langle 6 \sqcup 1\langle 1 \rangle \rangle$,\ 
$\langle 5 \sqcup 1\langle 2 \rangle \rangle$,\ 
$\langle 4 \sqcup 1\langle 3 \rangle \rangle$,\ 
$\langle 4 \sqcup 1\langle 1 \rangle \rangle$ and 
$\langle 3 \sqcup 1\langle 2 \rangle \rangle$ 
with $\delta_{\varphi S}=1$ 
are not dividing. 
By Zvonilov \cite{Zvonilov92}, these $6$ isotopy types are of {\it indeterminate type}, i.e., there exist dividing curves of these isotopy types. Then these dividing curves must realize some connected components with 
$\delta_{\varphi S}=0$. As a result, it turns out that 
if $\delta_{\varphi S}=0$ and $A(\br) \neq \emptyset$, 
then $A$ is dividing. 
Thus we have the following characterization of the invariant $\delta_{\varphi S}$.
 
\begin{theorem}
Let $A^+$ be a positive non-singular real curve of bidegree $(4,4)$ on 
ellipsoid. Then, $A$ is dividing or $A(\br)=\emptyset$, if and only if  
$\delta_{\varphi S} = 0$. \QED
\end{theorem}

Finally we get

\begin{theorem}
A connected component of moduli of 
a positive  real non-singular curve $A^+$ of bidegree $(4,4)$ on ellipsoid 
is defined by the topology type (or isotopy type) of $A^+$ and 
by divideness of $A(\bc)$ by $A(\br)$ (or the invariant $\delta_{\varphi S}$).
All these possibilities are presented in figures $\bee$ above. 
\end{theorem}

\bigskip 

%%%%%%%%%%%%%%%%%%%%%%%%%%%%%%%%%%%%%%%%%%%%%%%%%%%%%%%%%%
\section{Connected components of moduli of real non-singular curves in 
$|-2K_{\bff_1}|$
\label{F1}}
%<2>+<-2>

Here we consider the case 
$S \cong \langle 2 \rangle \oplus \langle -2 \rangle$ and 
$\theta=-1$ on $S$. 

Let us show that $Y=X/\{1,\tau\}=\bff_1$ in this case. We apply general 
results about right DPN-pairs which we have cited in Sect. \ref{modDPN}.  
By \eqref{fixedtau}, \eqref{fixedtau1}, $A=X^\tau$ is a non-singular 
irreducible curve of genus $9$. Since $A\not=\emptyset$,  
$Y$ is a rational surface with the Picard number $2$. 
Then $Y\cong \bff_n$ with $n\ge 0$. Since $A\subset Y$ is 
a right DPN-pair, exceptional curves of $Y$ are rational with square 
$-1$, $-2$ or $-4$. Thus, $Y\cong \bff_n$ with $n=0$, $1$, $2$ or $4$. 
If $n=0$, then $Y=\bp^1\times \bp^1$, and we have seen that $S=U(2)$ in 
this case. If $n=2$, then preimage of the exceptional section gives 
an element of $S$ with square $-4$ (it is also a root of $S$), but 
$S$ does not have elements with square $-4$. If $n=4$, then the 
exceptional rational curve with square $-4$ must be a component of $A$, but $A$ is irreducible of genus $9$. Thus, $Y=\bff_1$. 

This case is ($\da$)-non-degenerate because $S$ has no elements with square $-4$.  

For the exceptional section of $\bff_1$ with square $-1$, we denote by $e$ its 
preimage on $X$. Then $e^2=-2$. We consider the contraction $\bff_1\to \bp^2$ of 
the exceptional section and denote by $h$ the preimage on $X$ of  
a line $l\subset \bp^2$. Then $h^2=2$, $h\cdot e=0$. We have $S=\bz h+\bz e$.  

Since $\theta=-1$, we have $S_+=\{ 0 \}$. Obviously, $s=2$, $p=0$. 

Since 
$\Delta(S_+, S)^{(-4)} = \Delta(S_-, S)^{(-4)} = \emptyset$, 
by Lemma \ref{GroupG}, 
$G$ is trivial. 

The group $A_{S_-}$ is generated by 
$h^* = \frac12 h$ and $e^* = -\frac12 e$, 
and hence $A_{S_-} \cong (\bz/2\bz)^2$ and $l(A_{S_-})=2$. 
The characteristic element of $q_{S_-}$ is $h+e$ (mod $2S$). 

We have $A_{S_+}=0,\ H_+ =0,\ \rho=0$. 
Since $\Gamma_- = 0$, 
$H = H_- \ (\subset S/2S)$ is one of the following 5 subgroups: 
$$
\textstyle 
0,\ [ h ],\ [ e ],\ [ h+e ],\ S/2S=[h,e],
$$
where we consider $h$, $e$ mod $2S$. 
Since $H_+ = 0$, we have $q_\rho = (- q_{S_-})|H_-$. For the groups $H$ 
we have:

(1) $H=0$. 
Then $a_H=0$, $q_\rho=0$, 
$\delta_H=0$, $k_\rho=0$, $\mu_\rho=0$, $\sigma_\rho \equiv 0\pmod 8$. 
If $\delta_{\varphi S}=0$, 
then $v=0$, equivalently, $\delta_\varphi=0$, and $c_v\equiv 0\pmod 4$. 

(2) $H=[ h ]$. 
Then $a_H=1$, $q_\rho \cong q_{-1}(2)$, 
$\delta_H=1$, $k_\rho=0$, $\mu_\rho=0$, $\sigma_\rho \equiv -1\pmod 8$. 
If $\delta_{\varphi S}=0$, 
then $v=h$, and hence $c_v\equiv -1\pmod 4$. 

(3) $H=[ e ]$. 
Then $a_H=1$, $q_\rho \cong q_1(2)$, 
$\delta_H=1$, $k_\rho=0$, $\mu_\rho=0$, $\sigma_\rho \equiv 1\pmod 8$. 
If $\delta_{\varphi S}=0$, 
then $v=e$, and hence $c_v\equiv 1\pmod 4$. 

(4) $H=[ h+e ]$.
Then $a_H=1$, $q_\rho \cong z$, 
$\delta_H=0$, $k_\rho=1$, $\mu_\rho=0$, $\sigma_\rho \equiv 0\pmod 8$. 
If $\delta_{\varphi S}=0$, 
then $v=0\ \text{or}\ h+e$, and hence $c_v\equiv 0\pmod 4$. 

(5) $H=S/2S$. 
Then $a_H=2$, $q_\rho \cong q_1(2)\oplus q_{-1}(2)$, 
$\delta_H=1$, $k_\rho=0$, $\mu_\rho=0$, $\sigma_\rho \equiv 0\pmod 8$. 
If $\delta_\varphi=0$, then $c_v\equiv 0\pmod 4$. 
If $\delta_{\varphi S}=0$ and $\delta_\varphi=1$ (Type Ia), 
then we also have $c_v\equiv 0\pmod 4$ as we mention below. \\

\medskip

From above we always have $\mu_\rho=0$. 
We see $\delta_{\varphi S_+}=\delta_\varphi$. 
We also see $\delta_{\varphi S_-}=0$ if and only if 
$\delta_{\varphi S}=0$,\ $v = h+e$ and 
$H=[ h+e ] \ \text{or}\ S/2S$. 
If $\delta_{\varphi S_+}=\delta_\varphi =0$, 
then we have $\varepsilon_{v_+} \equiv 0 \pmod 2$. 
If $\delta_{\varphi S_-}=0$, then we have 
$c_v\equiv 0\pmod 4$ (see (4) and (5) above), 
and hence $\gamma_- = q_1(2)\oplus q_{-1}(2)$. 
We see $(q_{S_-})_v = u_+(2)$, and hence 
$\varepsilon_{v_-} \equiv 0 \pmod 2$.

We see that the genus 
(hence, by Theorem \ref{genus-isomor} and Proposition \ref{genus-isomor2},  isomorphism class) of an integral involutions $(L,\varphi,S)$ of 
the type $(\langle 2 \rangle \oplus \langle -2 \rangle,-1)$ 
satisfying (RSK3) is determined by the data 
\begin{equation}
(r, a, H, \delta_{\varphi S}, v),
\label{invF_1}
\end{equation}
where the $v\ (\in H)$ is defined only if $\delta_{\varphi S}=0$. We remind 
that $\delta_\varphi=0$ if and only if $\delta_{\varphi S}=v=0$. \\

\begin{figure}
\begin{picture}(200,140)
%%<2>+<-2>   $H=0$
%%\put(24,120){$H=0$}
\put(66,102){\circle{7}}
\put(71,100){{\tiny means $\delta_{\varphi S}=0$ and $v=0$,}}
%{\circle{5}}
%{{\tiny means $\delta_{\varphi S}=0$ and $\delta_\varphi=1$}}
\put(66,94){\circle*{3}}
\put(71,92){{\tiny means $\delta_{\varphi S}=1$}}

\multiput(8,0)(8,0){20}{\line(0,1){86}}
\multiput(0,8)(0,8){10}{\line(1,0){170}}
\put(0,0){\vector(0,1){100}}
\put(0,0){\vector(1,0){180}}
\put(  6,-10){{\tiny $1$}}
\put( 14,-10){{\tiny $2$}}
\put( 22,-10){{\tiny $3$}}
\put( 30,-10){{\tiny $4$}}
\put( 38,-10){{\tiny $5$}}
\put( 46,-10){{\tiny $6$}}
\put( 54,-10){{\tiny $7$}}
\put( 62,-10){{\tiny $8$}}
\put( 70,-10){{\tiny $9$}}
\put( 76,-10){{\tiny $10$}}
\put( 84,-10){{\tiny $11$}}
\put( 92,-10){{\tiny $12$}}
\put(100,-10){{\tiny $13$}}
\put(108,-10){{\tiny $14$}}
\put(116,-10){{\tiny $15$}}
\put(124,-10){{\tiny $16$}}
\put(132,-10){{\tiny $17$}}
\put(140,-10){{\tiny $18$}}
\put(148,-10){{\tiny $19$}}
\put(156,-10){{\tiny $20$}}

\put(-8, -1){{\tiny $0$}}
\put(-8,  7){{\tiny $1$}}
\put(-8, 15){{\tiny $2$}}
\put(-8, 23){{\tiny $3$}}
\put(-8, 31){{\tiny $4$}}
\put(-8, 39){{\tiny $5$}}
\put(-8, 47){{\tiny $6$}}
\put(-8, 55){{\tiny $7$}}
\put(-8, 63){{\tiny $8$}}
\put(-8, 71){{\tiny $9$}}
\put(-10, 79){{\tiny $10$}}

\put( -2,106){{\footnotesize $a$}} %vertical
\put(186, -2){{\footnotesize $r$}} %horizontal

%% Type 0  \circle{7}
%%  ( 8r,8a)
\put( 16, 0){\circle{7}}
\put( 80, 0){\circle{7}}
\put(144, 0){\circle{7}}
\put( 16,16){\circle{7}}
\put( 48,16){\circle{7}}
\put( 80,16){\circle{7}}
\put(112,16){\circle{7}}
%%%\put(144,16){\circle{7}}
\put( 48,32){\circle{7}}
\put( 80,32){\circle{7}}
\put(112,32){\circle{7}}
\put( 80,48){\circle{7}}
%%%\put(112,48){\circle{7}}
\put( 80,64){\circle{7}}
%%%\put( 80,80){\circle{7}}    %%12data

%% Type Ib  \circle*{3}
%%  ( 8r,8a)
\put(  8, 8){\circle*{3}}
\put( 24, 8){\circle*{3}}
\put( 72, 8){\circle*{3}}
\put( 88, 8){\circle*{3}}
\put(136, 8){\circle*{3}}
%%%\put(152, 8){\circle*{3}}
\put( 16,16){\circle*{3}}
\put( 32,16){\circle*{3}}
\put( 64,16){\circle*{3}}
\put( 80,16){\circle*{3}}
\put( 96,16){\circle*{3}}
\put(128,16){\circle*{3}}
%%%\put(144,16){\circle*{3}}
\put( 24,24){\circle*{3}}
\put( 40,24){\circle*{3}}
\put( 56,24){\circle*{3}}
\put( 72,24){\circle*{3}}
\put( 88,24){\circle*{3}}
\put(104,24){\circle*{3}}
\put(120,24){\circle*{3}}
%%%\put(136,24){\circle*{3}}
\put( 32,32){\circle*{3}}
\put( 48,32){\circle*{3}}
\put( 64,32){\circle*{3}}
\put( 80,32){\circle*{3}}
\put( 96,32){\circle*{3}}
\put(112,32){\circle*{3}}
%%%\put(128,32){\circle*{3}}
\put( 40,40){\circle*{3}}
\put( 56,40){\circle*{3}}
\put( 72,40){\circle*{3}}
\put( 88,40){\circle*{3}}
\put(104,40){\circle*{3}}
%%%\put(120,40){\circle*{3}}
\put( 48,48){\circle*{3}}
\put( 64,48){\circle*{3}}
\put( 80,48){\circle*{3}}
\put( 96,48){\circle*{3}}
%%%\put(112,48){\circle*{3}}
\put( 56,56){\circle*{3}}
\put( 72,56){\circle*{3}}
\put( 88,56){\circle*{3}}
%%%\put(104,56){\circle*{3}}
\put( 64,64){\circle*{3}}
\put( 80,64){\circle*{3}}
%%%\put( 96,64){\circle*{3}}
\put( 72,72){\circle*{3}}
%%%\put( 88,72){\circle*{3}}
%%%\put( 80,80){\circle*{3}}   %%39data
\end{picture}

\caption{$\bff_1$: All possible $(r,a,\delta_{\varphi S}, v)$ with $H=0$}
\label{<2>+<-2>-0}
\end{figure}

\begin{figure}
\begin{picture}(200,140)
%%<2>+<-2>   $H=[h,e]$
%%\put(24,120){$H=0$}
\put(66,110){\circle{7}}
\put(71,108){{\tiny means $\delta_{\varphi S}=0$ and $v=h+e$,}}
%{\circle{5}}
%{{\tiny means $\delta_{\varphi S}=0$ and $\delta_\varphi=1$}}
\put(66,102){\circle*{3}}
\put(71,100){{\tiny means $\delta_{\varphi S}=1$}}

\multiput(8,0)(8,0){20}{\line(0,1){94}}
\multiput(0,8)(0,8){11}{\line(1,0){170}}
\put(0,0){\vector(0,1){100}}
\put(0,0){\vector(1,0){180}}
\put(  6,-10){{\tiny $1$}}
\put( 14,-10){{\tiny $2$}}
\put( 22,-10){{\tiny $3$}}
\put( 30,-10){{\tiny $4$}}
\put( 38,-10){{\tiny $5$}}
\put( 46,-10){{\tiny $6$}}
\put( 54,-10){{\tiny $7$}}
\put( 62,-10){{\tiny $8$}}
\put( 70,-10){{\tiny $9$}}
\put( 76,-10){{\tiny $10$}}
\put( 84,-10){{\tiny $11$}}
\put( 92,-10){{\tiny $12$}}
\put(100,-10){{\tiny $13$}}
\put(108,-10){{\tiny $14$}}
\put(116,-10){{\tiny $15$}}
\put(124,-10){{\tiny $16$}}
\put(132,-10){{\tiny $17$}}
\put(140,-10){{\tiny $18$}}
\put(148,-10){{\tiny $19$}}
\put(156,-10){{\tiny $20$}}

\put(-8, -1){{\tiny $0$}}
\put(-8,  7){{\tiny $1$}}
\put(-8, 15){{\tiny $2$}}
\put(-8, 23){{\tiny $3$}}
\put(-8, 31){{\tiny $4$}}
\put(-8, 39){{\tiny $5$}}
\put(-8, 47){{\tiny $6$}}
\put(-8, 55){{\tiny $7$}}
\put(-8, 63){{\tiny $8$}}
\put(-8, 71){{\tiny $9$}}
\put(-10, 79){{\tiny $10$}}
\put(-10, 88){{\tiny $11$}}

\put( -2,106){{\footnotesize $a$}} %vertical
\put(186, -2){{\footnotesize $r$}} %horizontal

%% Type 0  \circle{7}
%%  ( 8r,8a)

\put( 16, 16){\circle{7}}
\put( 80, 16){\circle{7}}
\put(144, 16){\circle{7}}
\put( 48,32){\circle{7}}
\put( 80,32){\circle{7}}
\put(112,32){\circle{7}}
\put(144,32){\circle{7}}
\put( 48,48){\circle{7}}
\put( 80,48){\circle{7}}
\put(112,48){\circle{7}}
\put( 80,64){\circle{7}}
%%%\put(112,48){\circle{7}}

\put( 80,80){\circle{7}}

%%%\put( 80,80){\circle{7}}    %%12data

%% Type Ib  \circle*{3}
%%  ( 8r,8a)
\put( 24, 24){\circle*{3}}
\put( 72, 24){\circle*{3}}
\put( 88, 24){\circle*{3}}
\put(136, 24){\circle*{3}}
\put(152, 24){\circle*{3}}
\put( 32,32){\circle*{3}}
\put( 64,32){\circle*{3}}
\put( 80,32){\circle*{3}}
\put( 96,32){\circle*{3}}
\put(128,32){\circle*{3}}
\put(144,32){\circle*{3}}
\put( 40,40){\circle*{3}}
\put( 56,40){\circle*{3}}
\put( 72,40){\circle*{3}}
\put( 88,40){\circle*{3}}
\put(104,40){\circle*{3}}
\put(120,40){\circle*{3}}
\put(136,40){\circle*{3}}
\put( 48,48){\circle*{3}}
\put( 64,48){\circle*{3}}
\put( 80,48){\circle*{3}}
\put( 96,48){\circle*{3}}
\put(112,48){\circle*{3}}
\put(128,48){\circle*{3}}
\put( 56,56){\circle*{3}}
\put( 72,56){\circle*{3}}
\put( 88,56){\circle*{3}}
\put(104,56){\circle*{3}}
\put(120,56){\circle*{3}}
\put( 64,64){\circle*{3}}
\put( 80,64){\circle*{3}}
\put( 96,64){\circle*{3}}
\put(112,64){\circle*{3}}
\put( 72,72){\circle*{3}}
\put( 88,72){\circle*{3}}
\put(104,72){\circle*{3}}
\put( 80,80){\circle*{3}}
\put( 96,80){\circle*{3}}
\put( 88,88){\circle*{3}}
%%%\put( 80,96){\circle*{3}}   %%39data
\end{picture}

\caption{$\bff_1$: All possible $(r,a,\delta_{\varphi S}, v)$ with 
$H=S/2S=[h,e]$}
\label{<2>+<-2>-[h,e]}
\end{figure}
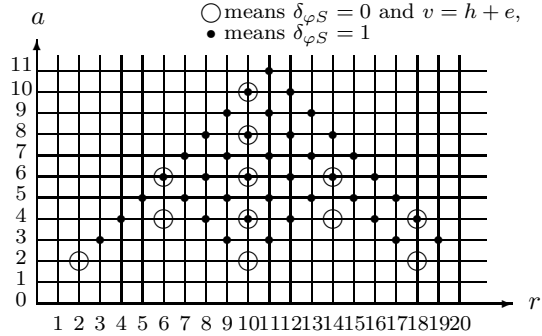

\begin{figure}

\begin{picture}(200,140)
%%<2>+<-2>   $H=[ h ]$
%%\put(24,122){$H=[ h ]$}
%{\circle{7}}
%{{\tiny means $\delta_\varphi=0$,}}
\put(66,102){\circle{6}}
\put(71,100){{\tiny means $\delta_{\varphi S}=0$ and $v=h$,}}
\put(66, 94){\circle*{3}}
\put(71, 92){{\tiny means $\delta_{\varphi S}=1$}}
\multiput(8,0)(8,0){20}{\line(0,1){86}}
\multiput(0,8)(0,8){10}{\line(1,0){170}}
\put(0,0){\vector(0,1){100}}
\put(0,0){\vector(1,0){180}}
\put(  6,-8){{\tiny $1$}}
\put( 14,-8){{\tiny $2$}}
\put( 22,-8){{\tiny $3$}}
\put( 30,-8){{\tiny $4$}}
\put( 38,-8){{\tiny $5$}}
\put( 46,-8){{\tiny $6$}}
\put( 54,-8){{\tiny $7$}}
\put( 62,-8){{\tiny $8$}}
\put( 70,-8){{\tiny $9$}}
\put( 76,-8){{\tiny $10$}}
\put( 84,-8){{\tiny $11$}}
\put( 92,-8){{\tiny $12$}}
\put(100,-8){{\tiny $13$}}
\put(108,-8){{\tiny $14$}}
\put(116,-8){{\tiny $15$}}
\put(124,-8){{\tiny $16$}}
\put(132,-8){{\tiny $17$}}
\put(140,-8){{\tiny $18$}}
\put(148,-8){{\tiny $19$}}
\put(156,-8){{\tiny $20$}}

\put(-8, -1){{\tiny $0$}}
\put(-8,  7){{\tiny $1$}}
\put(-8, 15){{\tiny $2$}}
\put(-8, 23){{\tiny $3$}}
\put(-8, 31){{\tiny $4$}}
\put(-8, 39){{\tiny $5$}}
\put(-8, 47){{\tiny $6$}}
\put(-8, 55){{\tiny $7$}}
\put(-8, 63){{\tiny $8$}}
\put(-8, 71){{\tiny $9$}}
\put(-10, 79){{\tiny $10$}}

\put( -2,106){{\footnotesize $a$}} %vertical
\put(186, -2){{\footnotesize $r$}} %horizontal

%% Type Ia  \circle{6}
%%  ( 8r,8a)
\put( 24, 8){\circle{6}}
\put( 88, 8){\circle{6}}
\put(152, 8){\circle{6}}
\put( 24,24){\circle{6}}
\put( 56,24){\circle{6}}
\put( 88,24){\circle{6}}
\put(120,24){\circle{6}}
\put( 56,40){\circle{6}}
\put( 88,40){\circle{6}}
\put(120,40){\circle{6}}
\put( 56,56){\circle{6}}
\put( 88,56){\circle{6}}
\put( 88,72){\circle{6}}   %%13 data

%% Type Ib  \circle*{3}
%%  ( 8r,8a)
\put( 16,16){\circle*{3}}
\put( 32,16){\circle*{3}}
\put( 80,16){\circle*{3}}
\put( 96,16){\circle*{3}}
\put(144,16){\circle*{3}}
\put( 24,24){\circle*{3}}
\put( 40,24){\circle*{3}}
\put( 72,24){\circle*{3}}
\put( 88,24){\circle*{3}}
\put(104,24){\circle*{3}}
\put(136,24){\circle*{3}}
\put( 32,32){\circle*{3}}
\put( 48,32){\circle*{3}}
\put( 64,32){\circle*{3}}
\put( 80,32){\circle*{3}}
\put( 96,32){\circle*{3}}
\put(112,32){\circle*{3}}
\put(128,32){\circle*{3}}
\put( 40,40){\circle*{3}}
\put( 56,40){\circle*{3}}
\put( 72,40){\circle*{3}}
\put( 88,40){\circle*{3}}
\put(104,40){\circle*{3}}
\put(120,40){\circle*{3}}
\put( 48,48){\circle*{3}}
\put( 64,48){\circle*{3}}
\put( 80,48){\circle*{3}}
\put( 96,48){\circle*{3}}
\put(112,48){\circle*{3}}
\put( 56,56){\circle*{3}}
\put( 72,56){\circle*{3}}
\put( 88,56){\circle*{3}}
\put(104,56){\circle*{3}}
\put( 64,64){\circle*{3}}
\put( 80,64){\circle*{3}}
\put( 96,64){\circle*{3}}
\put( 72,72){\circle*{3}}
\put( 88,72){\circle*{3}}
\put( 80,80){\circle*{3}}  %%39 data
\end{picture}

\caption{$\bff_1$: All possible $(r,a,\delta_{\varphi S}, v)$ with $H=[ h ]$}
\label{<2>+<-2>-h}
\end{figure}

\begin{figure}

\begin{picture}(200,140)
%%<2>+<-2>   $H=[ h ]$
%%\put(24,122){$H=[ h ]$}
%{\circle{7}}
%{{\tiny means $\delta_\varphi=0$,}}
\put(66,102){\circle{6}}
\put(71,100){{\tiny means $\delta_{\varphi S}=0$ and $v=e$,}}
\put(66, 94){\circle*{3}}
\put(71, 92){{\tiny means $\delta_{\varphi S}=1$}}
\multiput(8,0)(8,0){20}{\line(0,1){86}}
\multiput(0,8)(0,8){10}{\line(1,0){170}}
\put(0,0){\vector(0,1){100}}
\put(0,0){\vector(1,0){180}}
\put(  6,-8){{\tiny $1$}}
\put( 14,-8){{\tiny $2$}}
\put( 22,-8){{\tiny $3$}}
\put( 30,-8){{\tiny $4$}}
\put( 38,-8){{\tiny $5$}}
\put( 46,-8){{\tiny $6$}}
\put( 54,-8){{\tiny $7$}}
\put( 62,-8){{\tiny $8$}}
\put( 70,-8){{\tiny $9$}}
\put( 76,-8){{\tiny $10$}}
\put( 84,-8){{\tiny $11$}}
\put( 92,-8){{\tiny $12$}}
\put(100,-8){{\tiny $13$}}
\put(108,-8){{\tiny $14$}}
\put(116,-8){{\tiny $15$}}
\put(124,-8){{\tiny $16$}}
\put(132,-8){{\tiny $17$}}
\put(140,-8){{\tiny $18$}}
\put(148,-8){{\tiny $19$}}
\put(156,-8){{\tiny $20$}}

\put(-8, -1){{\tiny $0$}}
\put(-8,  7){{\tiny $1$}}
\put(-8, 15){{\tiny $2$}}
\put(-8, 23){{\tiny $3$}}
\put(-8, 31){{\tiny $4$}}
\put(-8, 39){{\tiny $5$}}
\put(-8, 47){{\tiny $6$}}
\put(-8, 55){{\tiny $7$}}
\put(-8, 63){{\tiny $8$}}
\put(-8, 71){{\tiny $9$}}
\put(-10, 79){{\tiny $10$}}

\put( -2,106){{\footnotesize $a$}} %vertical
\put(186, -2){{\footnotesize $r$}} %horizontal

%% Type Ia  \circle{6}
%%  ( 8r,8a)
\put( 136, 8){\circle{6}}
\put( 72, 8){\circle{6}}
\put(8, 8){\circle{6}}
\put( 136,24){\circle{6}}
\put( 104,24){\circle{6}}
\put( 72,24){\circle{6}}
\put(40,24){\circle{6}}
\put(104,40){\circle{6}}
\put( 72,40){\circle{6}}
\put(40,40){\circle{6}}
\put( 72,56){\circle{6}}
\put( 104,56){\circle{6}}
\put( 72,72){\circle{6}}   
%%13 data

%% Type Ib  \circle*{3}
%%  ( 8r,8a)

\put( 16,16){\circle*{3}}
\put (64,16){\circle*{3}}
\put( 80,16){\circle*{3}}
\put(128,16){\circle*{3}}
\put(144,16){\circle*{3}}
\put( 24,24){\circle*{3}}
\put( 56,24){\circle*{3}}
\put( 72,24){\circle*{3}}
\put( 88,24){\circle*{3}}
\put (120,24){\circle*{3}}
\put( 136,24){\circle*{3}}
\put( 32,32){\circle*{3}}
\put( 48,32){\circle*{3}}
\put( 64,32){\circle*{3}}
\put( 80,32){\circle*{3}}
\put( 96,32){\circle*{3}}
\put(112,32){\circle*{3}}
\put(128,32){\circle*{3}}
\put( 40,40){\circle*{3}}
\put( 56,40){\circle*{3}}
\put( 72,40){\circle*{3}}
\put( 88,40){\circle*{3}}
\put(104,40){\circle*{3}}
\put(120,40){\circle*{3}}
\put( 48,48){\circle*{3}}
\put( 64,48){\circle*{3}}
\put( 80,48){\circle*{3}}
\put( 96,48){\circle*{3}}
\put(112,48){\circle*{3}}
\put( 56,56){\circle*{3}}
\put( 72,56){\circle*{3}}
\put( 88,56){\circle*{3}}
\put(104,56){\circle*{3}}
\put( 64,64){\circle*{3}}
\put( 80,64){\circle*{3}}
\put( 96,64){\circle*{3}}
\put( 72,72){\circle*{3}}
\put( 88,72){\circle*{3}}
\put( 80,80){\circle*{3}}  %%39 data
\end{picture}

\caption{$\bff_1$: All possible $(r,a,\delta_{\varphi S}, v)$ with $H=[ e ]$}
\label{<2>+<-2>-e}
\end{figure}

\begin{figure}
\begin{picture}(200,140)
%%<2>+<-2>   $H= [ h+s ]$
\put(66,112){\circle{7}}
\put(71,110){{\tiny means $\delta_{\varphi S}=0$ and $v=0$,}}
\put(66,102){\circle{5}}
\put(71,100){{\tiny means $\delta_{\varphi S}=0$ and $v=h+e$,}}
\put(66,94){\circle*{3}}
\put(71,92){{\tiny means $\delta_{\varphi S}=1$}}

\multiput(8,0)(8,0){20}{\line(0,1){86}}
\multiput(0,8)(0,8){10}{\line(1,0){170}}
\put(0,0){\vector(0,1){100}}
\put(0,0){\vector(1,0){180}}
\put(  6,-8){{\tiny $1$}}
\put( 14,-8){{\tiny $2$}}
\put( 22,-8){{\tiny $3$}}
\put( 30,-8){{\tiny $4$}}
\put( 38,-8){{\tiny $5$}}
\put( 46,-8){{\tiny $6$}}
\put( 54,-8){{\tiny $7$}}
\put( 62,-8){{\tiny $8$}}
\put( 70,-8){{\tiny $9$}}
\put( 76,-8){{\tiny $10$}}
\put( 84,-8){{\tiny $11$}}
\put( 92,-8){{\tiny $12$}}
\put(100,-8){{\tiny $13$}}
\put(108,-8){{\tiny $14$}}
\put(116,-8){{\tiny $15$}}
\put(124,-8){{\tiny $16$}}
\put(132,-8){{\tiny $17$}}
\put(140,-8){{\tiny $18$}}
\put(148,-8){{\tiny $19$}}
\put(156,-8){{\tiny $20$}}

\put(-8, -1){{\tiny $0$}}
\put(-8,  7){{\tiny $1$}}
\put(-8, 15){{\tiny $2$}}
\put(-8, 23){{\tiny $3$}}
\put(-8, 31){{\tiny $4$}}
\put(-8, 39){{\tiny $5$}}
\put(-8, 47){{\tiny $6$}}
\put(-8, 55){{\tiny $7$}}
\put(-8, 63){{\tiny $8$}}
\put(-8, 71){{\tiny $9$}}
\put(-10, 79){{\tiny $10$}}

\put( -2,106){{\footnotesize $a$}} %vertical
\put(186, -2){{\footnotesize $r$}} %horizontal

%% Type 0  \circle{7}
%%  ( 8r,8a)
\put( 16,16){\circle{7}}
\put( 80,16){\circle{7}}
\put(144,16){\circle{7}}
\put( 48,32){\circle{7}}
\put( 80,32){\circle{7}}
\put(112,32){\circle{7}}
\put( 80,48){\circle{7}}
\put(112,48){\circle{7}}
\put( 80,64){\circle{7}}
\put( 80,80){\circle{7}}    %%10data

%% Type Ia  \circle{5}
%%  ( 8r,8a)
\put( 16,16){\circle{5}}
\put( 80,16){\circle{5}}
\put(144,16){\circle{5}}
\put( 48,32){\circle{5}}
\put( 80,32){\circle{5}}
\put(112,32){\circle{5}}
\put( 48,48){\circle{5}}
\put( 80,48){\circle{5}}
\put( 80,64){\circle{5}}
\put( 80,80){\circle{5}}    %%10data

%% Type Ib  \circle*{3}
%%  ( 8r,8a)
\put( 24,24){\circle*{3}}
\put( 72,24){\circle*{3}}
\put( 88,24){\circle*{3}}
\put(136,24){\circle*{3}}
\put( 32,32){\circle*{3}}
\put( 64,32){\circle*{3}}
\put( 80,32){\circle*{3}}
\put( 96,32){\circle*{3}}
\put(128,32){\circle*{3}}
\put( 40,40){\circle*{3}}
\put( 56,40){\circle*{3}}
\put( 72,40){\circle*{3}}
\put( 88,40){\circle*{3}}
\put(104,40){\circle*{3}}
\put(120,40){\circle*{3}}
\put( 48,48){\circle*{3}}
\put( 64,48){\circle*{3}}
\put( 80,48){\circle*{3}}
\put( 96,48){\circle*{3}}
\put(112,48){\circle*{3}}
\put( 56,56){\circle*{3}}
\put( 72,56){\circle*{3}}
\put( 88,56){\circle*{3}}
\put(104,56){\circle*{3}}
\put( 64,64){\circle*{3}}
\put( 80,64){\circle*{3}}
\put( 96,64){\circle*{3}}
\put( 72,72){\circle*{3}}
\put( 88,72){\circle*{3}}
\put( 80,80){\circle*{3}}   %%30data
\end{picture}

\caption{$\bff_1$: All possible $(r,a,\delta_{\varphi S}, v)$ 
with $H= [ h+e ]$}
\label{<2>+<-2>-he}
\end{figure}
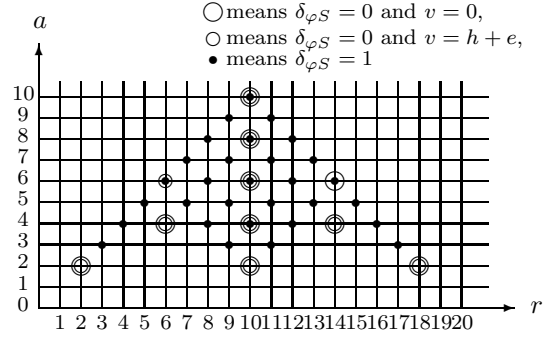

By Theorem \ref{geninvtheorem}, 
the complete list of the above invariants is given by the conditions:

\noindent
{\bf Type 0} ($\delta_\varphi=0$)\\
$H=0\ \text{or}\ [ h+e ]$.\\
$r=2,\,6,\,10,\,14,\,18$.\\
$a \equiv  0 \pmod{2}.$\\
$a \leq r$.\\
$r+a \leq 2a_H +18$.\\
If $H=0$, then $a \geq 0$.\\
If $H=[ h+e ]$, then $a \geq 2$.\\
If $H=0$ and $a=0$, then $2-r \equiv 0 \pmod{8}.$\\
If $H=[ h+e ]$ and $a=2$, then $2-r \equiv 0 \pmod{8}.$\\
If $a=r$, then $2-r \equiv 0 \pmod 8$.\\

\noindent
{\bf Type Ia} ($\delta_{\varphi S}=0,\ \delta_\varphi=1$)\\
$a_H \geq 1$.\\
$c_v \equiv \sigma_\rho \pmod{4}$. 
Hence, $c_v \equiv 0 \pmod{4}$ if $H=S/2S$.\\
If $H=[ h ]$, 
then $a \geq 1$, $r=3,\,7,\,11,\,15$ and $a \equiv  1 \pmod{2}$. \\
If $H=[ e ]$, 
then $a \geq 1$, $r=1,\,5,\,9,\,13,\,17$ and $a \equiv  1 \pmod{2}$. \\
If $H=[ h+e ]$, 
then $a \geq 2$, $r=2,\,6,\,10,\,14,\,18$ and $a \equiv  0 \pmod{2}$. \\
If $H=S/2S$, 
then $a \geq 2$, $r=2,\,6,\,10,\,14,\,18$ and $a \equiv  0 \pmod{2}$. \\
If $H=[ h ]$ or $H=[ e ]$ and $a=1$, 
then $2-r \equiv \sigma_\rho \pmod{8}.$\\
If $H=[ h+e ]$ or $H=S/2S$ and $a=2$, 
then $2-r \equiv 0 \pmod{8}.$

\noindent
$a \leq r$, \ $r+a \leq 2a_H +18$.\\
If $\delta_{\varphi S}=0$,\ $v = h+e$,\ 
$H=[ h+e ] \ \text{or}\ S/2S$ and $r+a = 2a_H +18$, \\
then $2-r \equiv 0 \pmod{8}.$\\

\noindent
{\bf Type Ib} ($\delta_{\varphi S}=1$)\\
$a \geq a_H + k_\rho +1$.\\
$r+a \equiv 0 \pmod{2}.$\\
If $a = a_H + k_\rho +1$, then $2-r \equiv \sigma_\rho \pm 1 \pmod{8}.$\\
If $a = a_H + k_\rho +2$, then $2-r \not \equiv \sigma_\rho +4 \pmod{8}.$

\noindent
$1 \leq r \leq 19$, \ $a \leq r$, \ $r+a \leq 2a_H +18$.\\
If $\delta_{\varphi S}=0$,\ $v = h+e$,\ 
$H=[ h+e ] \ \text{or}\ S/2S$ and $r+a = 2a_H +18$,\\
then $2-r \equiv 0 \pmod{8}.$\\

\medskip

In Type Ia case, if $H = S/2S$, then $c_v \equiv 0 \pmod{4}$. 
Thus $v = h+e$. 
Hence, $H$ determines $v$ in Type Ia case. 
Thus, by the above, the isomorphism class of an integral 
involutions $(L,\varphi,S)$ of 
the type $(\langle 2 \rangle \oplus \langle -2 \rangle,-1)$ satisfying (RSK3) 
is determined by the data 
$$
(r, a, H, \delta_{\varphi S}, \delta_\varphi).
$$
We give all these possible data in Figures \ref{<2>+<-2>-0} --- \ref{<2>+<-2>-he}. 

By Theorem \ref{relinvtheorem} about related involutions, we have 
\begin{equation}\begin{split} 
&r(\varphi)+r(\tau\varphi)=20,\ \ 
a(\varphi) - a_{H(\varphi)}=a(\tau\varphi)-a_{H(\tau\varphi)},\\
&a_{H(\varphi)} + a_{H(\tau\varphi)} = 2,\ \ 
H(\tau\varphi)=H(\varphi)^\perp\ \text{w.r.t.}\  b_{S_-},\\
&\delta_{\varphi S}=\delta_{(\tau\varphi) S},\ \ 
s_\varphi+s_{\tau\varphi}=h+e\mod 2L. 
\label{relatedF_1} 
\end{split}\end{equation}
Thus, integral involutions of Type 0 with $H=0$ ($12$ classes) are related to 
involutions of Type Ia with $H=S/2S$. 
Involutions of Type 0  with $H= [ h+e ]$ 
($10$ classes) are related to 
involutions of Type Ia with $H= [ h+e ]$. 
Involutions of Type Ia with $H= [ h ]$ ($13$ classes) 
are related to 
involutions of Type Ia with $H= [ e ]$. 
Involutions of Type Ib with $H=0$ ($39$ classes) are related to 
involutions of Type Ib with $H=S/2S$. 
Involutions of Type Ib with $H= [ h ]$ ($39$ classes) 
are related to 
involutions of Type Ib with $H= [ e ]$. 
Finally, the class 
$(r, a, H=[ h+e ], \text{Type\,Ib})$ is related to 
$(20-r, a, H=[ h+e ], \text{Type\,Ib})$. 
(There are $30$ classes of Type Ib with $H=[ h+e ]$.)

Moreover, if we identify related involutions, there are 
$35 (= 12+10+13)$ classes with $\delta_{\varphi S}=0$ and 
$95 (= 39+39+17)$ classes with $\delta_{\varphi S}=1$. 

\begin{figure} 
\centerline{\includegraphics[width=5cm]{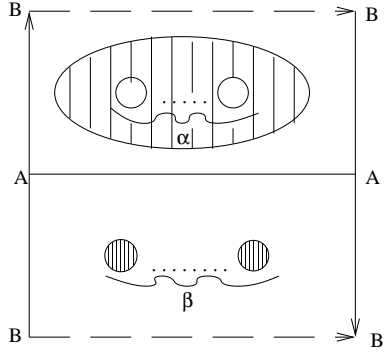}} 
\caption{$\bff_1$: $A^+$ for $H=[h,e]$ with $(r,a,\delta_{\varphi S},v)\not=(10,10,0,h+e)$; here $\alpha=g$,  $\beta=k$.}
\label{F1he1}
\end{figure}

\begin{figure} 
\centerline{\includegraphics[width=5cm]{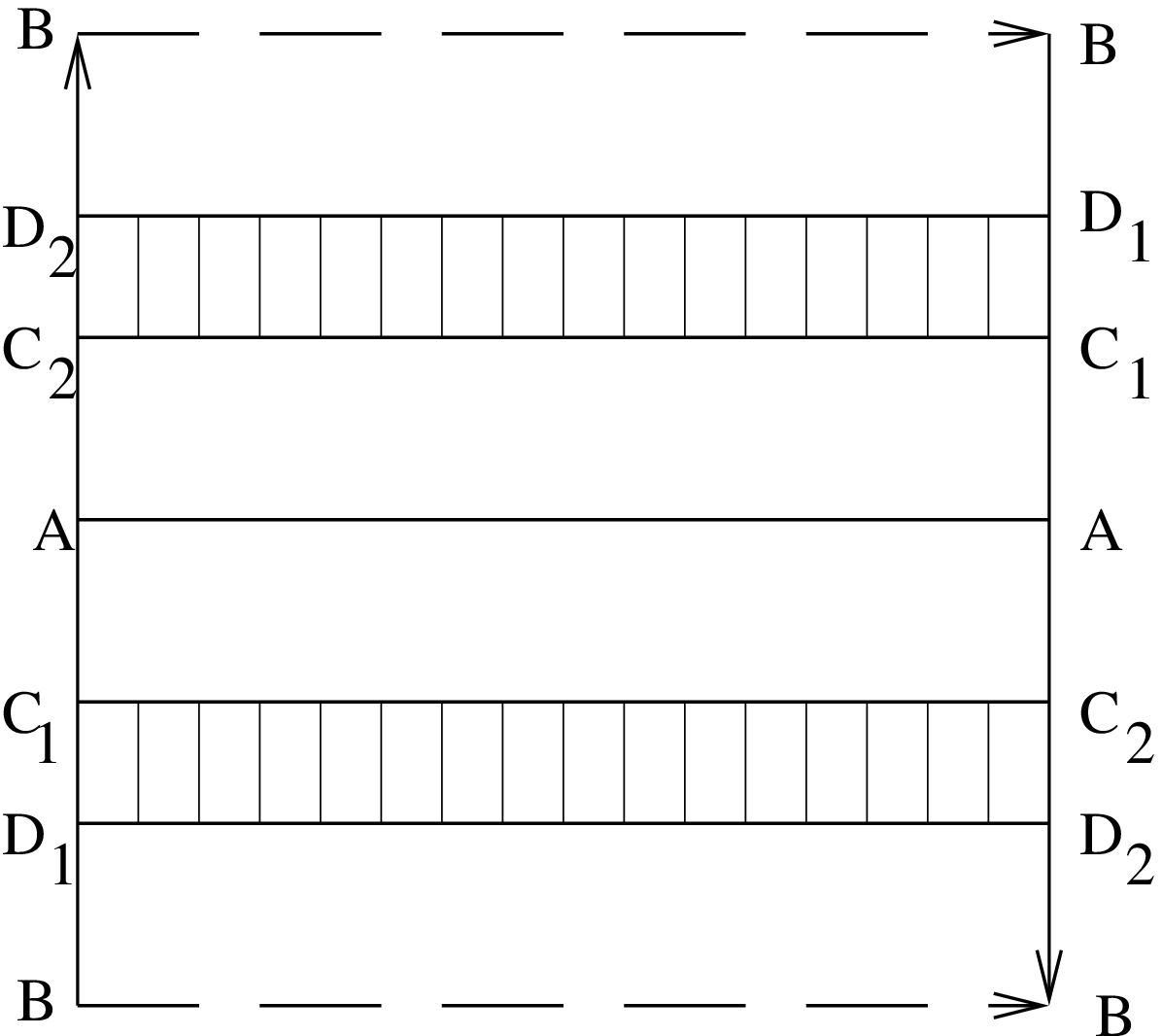}} 
\caption{$\bff_1$: $A^+$ for $H=[h,e]$ with 
$(r,a,\delta_{\varphi S},v)=(10,10,0,h+e)$.}
\label{F1he2}
\end{figure}

\begin{figure} 
\centerline{\includegraphics[width=5cm]{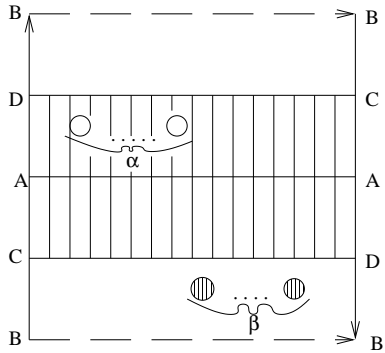}} 
\caption{$\bff_1$: $A^+$ for $H=[h]$; here $\alpha=g-1$, $\beta=k$.}
\label{F1he3}
\end{figure}

\begin{figure} 
\centerline{\includegraphics[width=5cm]{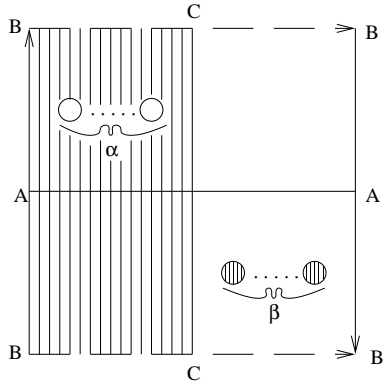}} 
\caption{$\bff_1$: $A^+$ for $H=[h+e]$ with  
$(r,a,\delta_{\varphi S})\not=(10,8,0),\ (10,10,0)$; here $\alpha=g-1$,  $\beta=k$.}
\label{F1h+e1}
\end{figure}

\begin{figure} 
\centerline{\includegraphics[width=5cm]{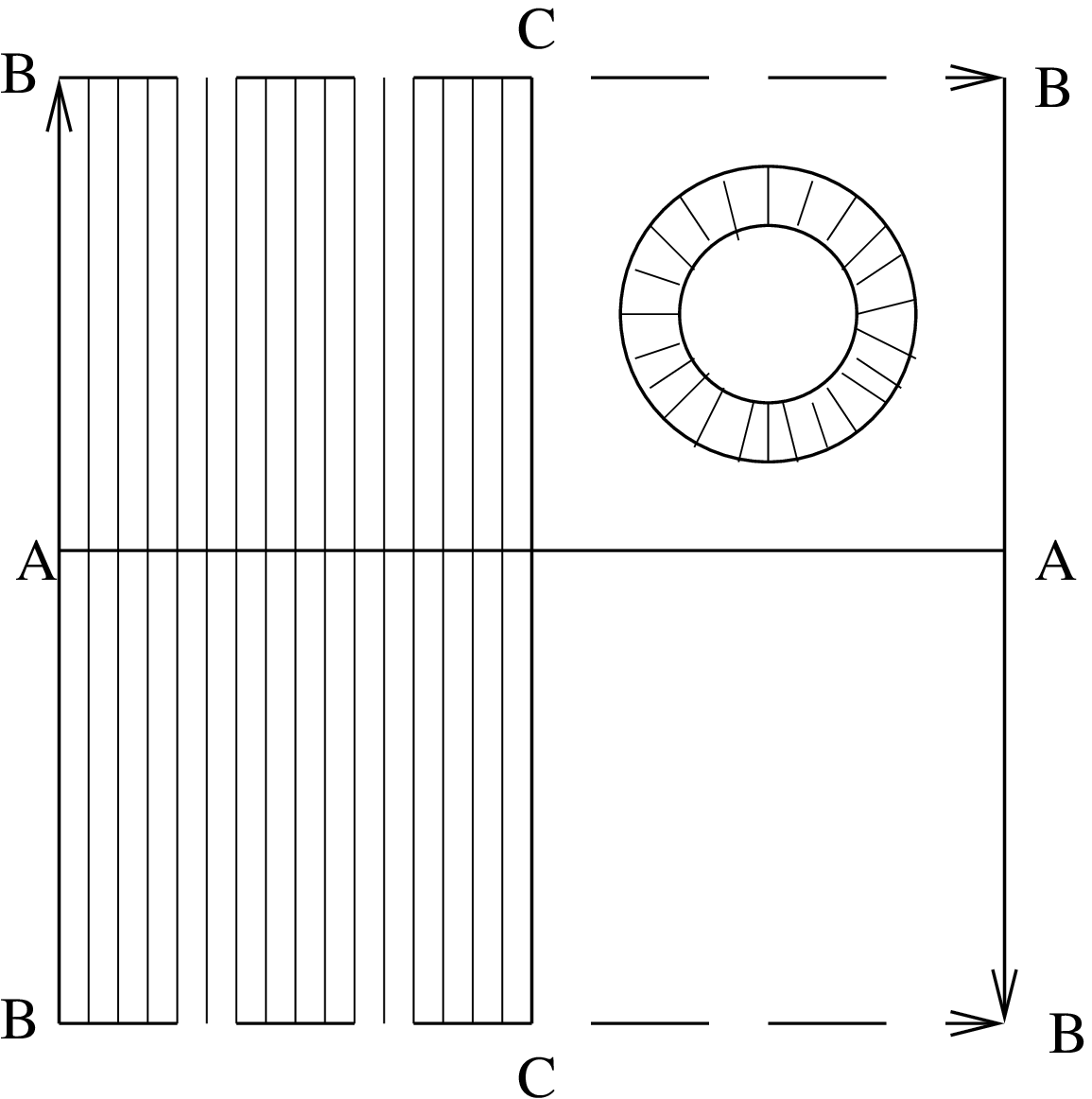}} 
\caption{$\bff_1$: $A^+$ for $H=[h+e]$ with  
$(r,a,\delta_{\varphi S},v)=(10,8,0,0)$.}
\label{F1h+e2}
\end{figure}

\begin{figure} 
\centerline{\includegraphics[width=5cm]{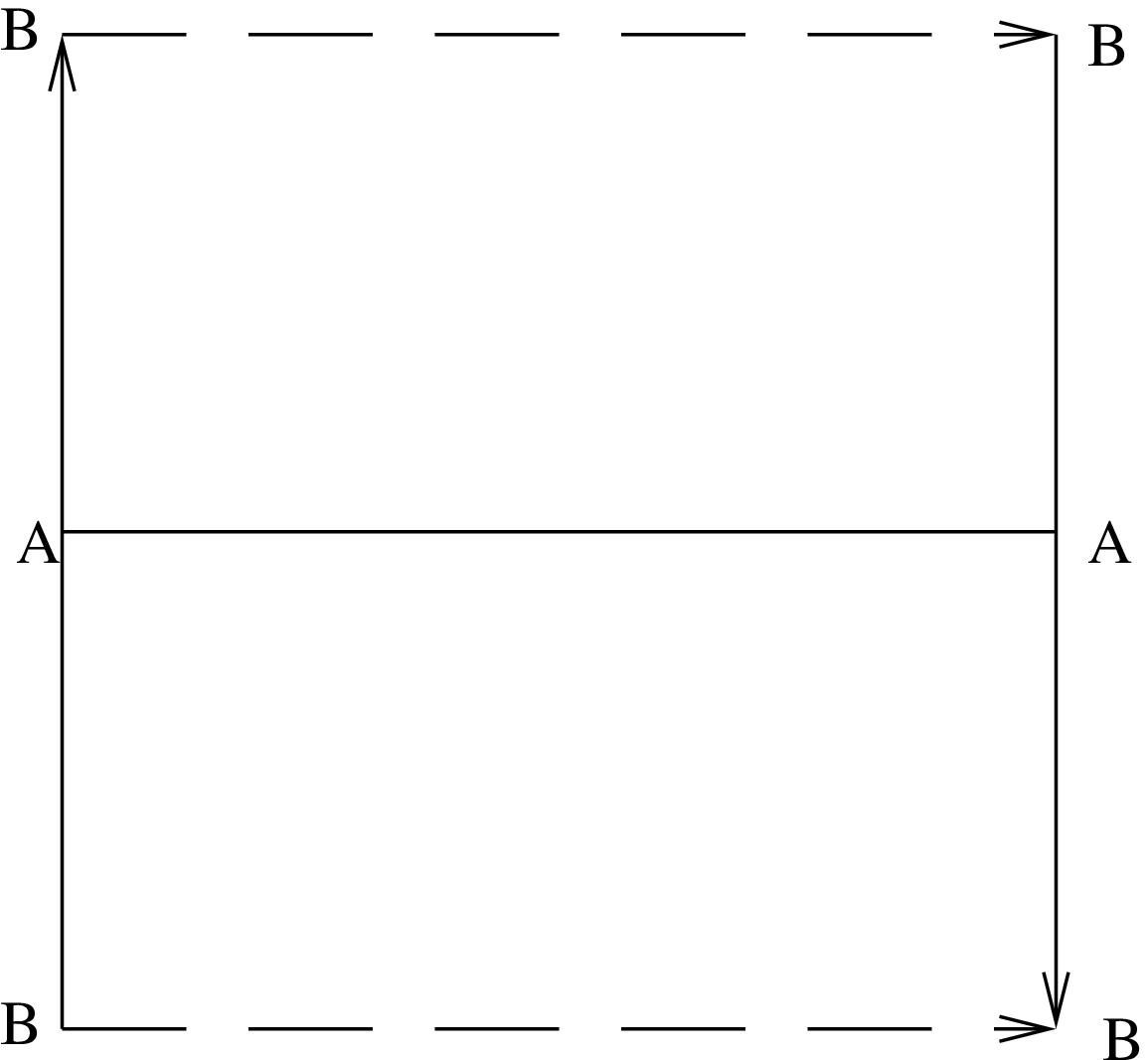}} 
\caption{$\bff_1$: $A^+=\emptyset$ for $H=[h+e]$ with  
$(r,a,\delta_{\varphi S},v)=(10,10,0,0)$.}
\label{F1h+e3}
\end{figure}

\medskip
%5620

Let us consider the geometric interpretation of the above calculations. 
We denote by $s$ the exceptional section of $\bff_1$ with $s^2=-1$ and 
by $c$ the fiber of the natural fibration $\pi:\bff_1\to \bp^1$. The 
contraction of $s$ (as an exceptional curve of the first kind) gives 
the natural morphism $p:\bff_1\to \bp^2$, and we denote $P=p(s)$. 
A non-singular curve $A\in |-2K_{\bff_1}|$ gives then a curve $A_1=p(A)$ 
of degree $6$ in $\bp^2$ with only one singular point $P$ which is a 
quadratic singular point resolving by one blow-up. A small 
deformation of $A$ or $A_1$ (in the same connected component of moduli) 
makes $P$ non-degenerate. Using Bezoute theorem, one can easily 
draw all in principle possible pictures of $A_1$.   
For example, one can find these pictures in Figure 1 of 
\cite{Itenberg94}. Lifting these pictures to $\bff_1$ 
and using \eqref{realcomponents}, \eqref{realmod2} and Theorem 
\ref{realcurve} applied to both positive curves 
$A^+$ and $A^-$, we get from  Figures \ref{<2>+<-2>-0} --- \ref{<2>+<-2>-he} 
all pictures of $A^+$ on $\bff_1$ up to isotopy.   
They are all presented in Figures \ref{F1he1}---\ref{F1h+e3} where 
we denote $g=(22-r-a)/2$ and $k=(r-a)/2$, as usual. The interval $AA$ 
denotes $s(\br)$, the real part of the exceptional section $s$, 
and the interval $BB$ denotes the $p^{-1}(l(\br))$ 
of a real projective line $l\subset \bp^2$  which does not contain $P$.  

We have the following geometric interpretation of the invariants  
$\delta_h$ and $\delta_e$ (see \eqref{deltax1}) of $A^+$ if 
$A^+\not=\emptyset$. 
We have $\delta_h=0$ (equivalently, $h\mod2\in H$), if and only if 
homotopically $l(\br)\subset A^-$ (i. e. some deformation of 
$l(\br)$ is contained in $A^-$). Similarly, $\delta_e=0$, if and only if 
homotopically $s(\br)\subset A^-$. The invariants $\delta_h$ and $\delta_e$ 
for both positive curves $A^+$ and $A^-$ are sufficient to find the group 
$H$. 

Thus, we get 
the isotopy classification of real non-singular curves $A\in |-2K_{\bff_1}|$. 
In this classification we don't care about position of $A(\br)$ with 
respect to the real part $s(\br)$ of the exceptional section (one 
can  see that more delicate classification in \cite{Itenberg94}).  

We have the following interpretation of the invariant $\delta_{\varphi S}$: 
{\it one has $\delta_{\varphi S}=0$, if and only if the curve $A$ is dividing: 
$A(\br)$ divides $A(\bc)$ in two connected parts or $A(\br)=\emptyset$.  
Equivalently, $A(\br)=0$ in $H_1(A(\bc),\bz/2)$.} 

In one direction it follows from considerations similar to Kharlamov 
in \cite{Kharlamov75a} and \cite{Kharlamov76}. We have  

\begin{proposition}\label{div => delta_varphiS=0}
Let $(X,\tau,\varphi)$ be a real K3 surface with a non-symplectic involution 
of the type $(S,\theta)$ with non empty $A=X^\tau$. 
%%%%%We assume $Y(\br) \neq \emptyset$, 
%%%%%where we use the same symbol $\theta$ for the anti-holomorphic involution 
%%%%%$\varphi_{{\rm mod}\ \tau}$ on $Y:=X/\{1,\tau\}$, and 
%%%%%$Y(\br)$ denotes the fixed point set of $\theta$ in $Y$. 
Then, if $(A, \theta)$ is dividing, then $\delta_{\varphi S}=0$. 
\end{proposition}

\Proof
Suppose that $[A(\br)]=0\ \text{in}\ H_1(A;\bz/2\bz)$. 
Then $A(\br)$ and $A$ satisfy the condition a) in Sect. 2.2 of  
\cite{Kharlamov75a}. Since $A(\br)$ is a disjoint union of circles, 
$w_1(A(\br))=0$, and the condition b) in Sect. 2.2 of 
\cite{Kharlamov75a} is also satisfied. 
Using the notation of \cite{Kharlamov75a}, we obtain that 
$l(P_{A(\br)})$ realizes the nulls of the group $H_2(P_A;\bz/2\bz)$. 
Since $\tau\varphi=\varphi\tau$, we have $\tau (X(\br)) = X(\br)$, and  
$A(\br)=A \cap X(\br)$. 

The tangent bundle $T(A)$, which is real $2$-dimensional, 
is isomorphic to the normal bundle $N(A)$ of $A$ in $X$ (because $X$ is 
a K3 surface), and the tangent bundle $T(A(\br))$ is isomorphic to the 
normal bundle $N(A(\br))$ of $A(\br)$ in $X(\br)$. 
Since $l(P_{A(\br)})$ realizes the nulls of the group $H_2(P_A,\bz/2\bz)$, 
by Lemma 1 in Sec. 2.3 of \cite{Kharlamov75a}, the class $[X(\br)]$ in 
$H_2(X,\bz/2\bz)$ is orthogonal to $\mathrm{Im}\,\alpha_2$ in 
$H_2(X,\bz/2\bz)$ with respect to the intersection pairing  
where $\alpha_2$ is the homomorphism in the Smith exact sequence for 
$(X,\tau)$ (see \cite{Kharlamov75a}). 
Since $H_1(X,\bz) = 0$, as in the proof of Lemma 3.7 in \cite{Kharlamov76}, 
we have 
$\mathrm{Im}\,\alpha_2 = \{ x \in H_2(X,\bz/2\bz) \ |\ \tau_* (x) = x \}$. 
It follows, $\delta_{\varphi S}=0$. 
\QED

Unfortunately, we don't know a simple proof of the opposite statement: 
if $\delta_{\varphi S}=0$, then the curve $A$ is dividing or 
$A(\br)=\emptyset$, 

Like in previous cases 
(because we know that invariants \eqref{invF_1} define the connected 
component of moduli),    
it would be enough to construct a dividing curve $A$ on $\bff_1$ 
in each case when $\delta_{\varphi S}=0$, and it should follow from known 
results about real curves of degree 6 (e. g. see \cite{Itenberg92} and 
\cite{Itenberg94}). Another general argument is to  
apply Donaldson's trick \cite{D}, like in \cite{DIK2000}.  We 
can consider $\varphi$ as holomorphic involution and $\tau$ as 
anti-holomorphic. Then the statement follows from Theorem 
\ref{realcurve}. The same argument works in all previous cases, 
but we did not need it. Thus, we finally get 

\begin{theorem}
A connected component of moduli of a positive real non-singular curve 
$A\in |-2K_{\bff_1}|$ is defined by the isotopy type of 
$A^+\subset \bff_1(\br)$ and by the divideness by $A(\br)$ of $A(\bc)$ (equivalently, by the invariant $\delta_{\varphi S}$). 
All these possibilities are presented in 
in the Figures $\bff_1$ above and in \eqref{relatedF_1}.
\label{componentsF_1}
\end{theorem}

\bigskip 

\section{Connected components of moduli of real non-singular curves in 
$|-2K_{\bff_4}|$}
\label{F4}

%U case

Here we consider the case $S \cong U$ and $\theta = -1$. 

Like in previous cases, we can prove that $Y=X/\{1,\tau\}=\bff_4$ 
in this case. 
This case is ($\da$)-non-degenerate because $S$ has no roots with square 
$-4$.

We have $S_+=\{ 0 \}$ and $S_- = S$, $s=2$, $p=0$. 

Since 
$\Delta(S_+, S)^{(-4)} = \Delta(S_-, S)^{(-4)} = \emptyset$, 
by Lemma \ref{GroupG} the group  
$G$ is trivial. 

%Let $e$ and $e'$ be generators of $U$ with $e^2 = e'^2 = 0$ and $(e,e')=1$. 

Since $U$ is unimodular, $A_{S_-} = \{ 0 \}$.  Obviously, $A_{S_+} = \{ 0 \}$. 
Hence, $H = H_+ = H_- = 0$, $\rho=0$ and $\delta_{\varphi S}=\delta_\varphi$. 

By Theorem \ref{genus-isomor} and Proposition \ref{genus-isomor2}, 
the genus defines the isomorphism class of an integral involutions $(L,\varphi,S)$ of the type $(U,-1)$ satisfying (RSK3), and it 
defines the connected component of moduli. It  is determined by the data 
\begin{equation}
(r, a, \delta_{\varphi S}=\delta_\varphi,v)
\label{invF_4}
\end{equation}
where $v=0$ if $\delta_{\varphi S}=0$ (otherwise, $v$ is not defined). 
The complete list of the invariants \eqref{invF_4} is given by the conditions:

\medskip

\noindent
$\delta_\varphi=0$ case

$r=2,6,10,14,18$,\ \ $a \equiv  0 \pmod{2}$,\ \ $a \leq r$,\ \ $r+a \leq 20$.

If $a=0$, $a=r$ or $r+a=20$, then $2-r \equiv 0 \pmod{8}.$\\

\noindent
$\delta_\varphi=1$ case

$a \geq 1$,\ \ $r+a \equiv 0 \pmod{2}$,\ \ $1 \leq r \leq 19$,\ \ $a \leq r$,\ \ $r+a \leq 20$.

If $a=1$, then $2-r \equiv \pm 1 \pmod{8}.$

If $a=2$, then $2-r \not \equiv 4 \pmod{8}.$\\

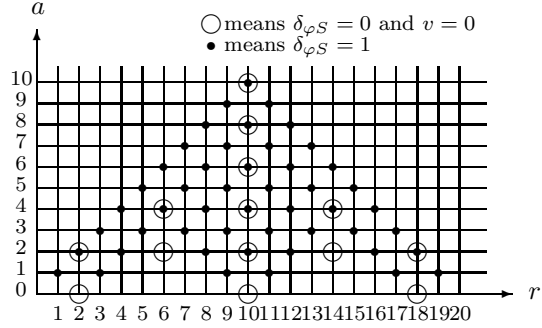
\begin{figure}
\begin{picture}(200,140)
%% (U,-1)
\put(66,102){\circle{7}}
\put(71,100){{\tiny means $\delta_{\varphi S}=0$ and $v=0$}}
\put(66,94){\circle*{3}}
\put(71,92){{\tiny means $\delta_{\varphi S}=1$}}

\multiput(8,0)(8,0){20}{\line(0,1){86}}
\multiput(0,8)(0,8){10}{\line(1,0){170}}
\put(0,0){\vector(0,1){100}}
\put(0,0){\vector(1,0){180}}
\put(  6,-10){{\tiny $1$}}
\put( 14,-10){{\tiny $2$}}
\put( 22,-10){{\tiny $3$}}
\put( 30,-10){{\tiny $4$}}
\put( 38,-10){{\tiny $5$}}
\put( 46,-10){{\tiny $6$}}
\put( 54,-10){{\tiny $7$}}
\put( 62,-10){{\tiny $8$}}
\put( 70,-10){{\tiny $9$}}
\put( 76,-10){{\tiny $10$}}
\put( 84,-10){{\tiny $11$}}
\put( 92,-10){{\tiny $12$}}
\put(100,-10){{\tiny $13$}}
\put(108,-10){{\tiny $14$}}
\put(116,-10){{\tiny $15$}}
\put(124,-10){{\tiny $16$}}
\put(132,-10){{\tiny $17$}}
\put(140,-10){{\tiny $18$}}
\put(148,-10){{\tiny $19$}}
\put(156,-10){{\tiny $20$}}

\put(-8, -1){{\tiny $0$}}
\put(-8,  7){{\tiny $1$}}
\put(-8, 15){{\tiny $2$}}
\put(-8, 23){{\tiny $3$}}
\put(-8, 31){{\tiny $4$}}
\put(-8, 39){{\tiny $5$}}
\put(-8, 47){{\tiny $6$}}
\put(-8, 55){{\tiny $7$}}
\put(-8, 63){{\tiny $8$}}
\put(-8, 71){{\tiny $9$}}
\put(-10, 79){{\tiny $10$}}

\put( -2,106){{\footnotesize $a$}} %vertical
\put(186, -2){{\footnotesize $r$}} %horizontal

%% \delta_\varphi = 0
%%  ( 8r,8a)
\put( 16, 0){\circle{7}}
\put( 80, 0){\circle{7}}
\put(144, 0){\circle{7}}
\put( 16,16){\circle{7}}
\put( 48,16){\circle{7}}
\put( 80,16){\circle{7}}
\put(112,16){\circle{7}}
\put(144,16){\circle{7}}
\put( 48,32){\circle{7}}
\put( 80,32){\circle{7}}
\put(112,32){\circle{7}}
\put( 80,48){\circle{7}}
\put( 80,64){\circle{7}}
\put( 80,80){\circle{7}}    %%14data

%% \delta_\varphi = 1
%%  ( 8r,8a)
\put(  8, 8){\circle*{3}}
\put( 24, 8){\circle*{3}}
\put( 72, 8){\circle*{3}}
\put( 88, 8){\circle*{3}}
\put(136, 8){\circle*{3}}
\put(152, 8){\circle*{3}}
\put( 16,16){\circle*{3}}
\put( 32,16){\circle*{3}}
\put( 64,16){\circle*{3}}
\put( 80,16){\circle*{3}}
\put( 96,16){\circle*{3}}
\put(128,16){\circle*{3}}
\put(144,16){\circle*{3}}
\put( 24,24){\circle*{3}}
\put( 40,24){\circle*{3}}
\put( 56,24){\circle*{3}}
\put( 72,24){\circle*{3}}
\put( 88,24){\circle*{3}}
\put(104,24){\circle*{3}}
\put(120,24){\circle*{3}}
\put(136,24){\circle*{3}}
\put( 32,32){\circle*{3}}
\put( 48,32){\circle*{3}}
\put( 64,32){\circle*{3}}
\put( 80,32){\circle*{3}}
\put( 96,32){\circle*{3}}
\put(112,32){\circle*{3}}
\put(128,32){\circle*{3}}
\put( 40,40){\circle*{3}}
\put( 56,40){\circle*{3}}
\put( 72,40){\circle*{3}}
\put( 88,40){\circle*{3}}
\put(104,40){\circle*{3}}
\put(120,40){\circle*{3}}
\put( 48,48){\circle*{3}}
\put( 64,48){\circle*{3}}
\put( 80,48){\circle*{3}}
\put( 96,48){\circle*{3}}
\put(112,48){\circle*{3}}
\put( 56,56){\circle*{3}}
\put( 72,56){\circle*{3}}
\put( 88,56){\circle*{3}}
\put(104,56){\circle*{3}}
\put( 64,64){\circle*{3}}
\put( 80,64){\circle*{3}}
\put( 96,64){\circle*{3}}
\put( 72,72){\circle*{3}}
\put( 88,72){\circle*{3}}
\put( 80,80){\circle*{3}}   %%49data
\end{picture}

\caption{$\bff_4$: All possible $(r, a, \delta_{\varphi S}, v)$ (here $H=0$)}
\label{Ugraph}
\end{figure}

All possible data \eqref{invF_4} are given on Figure \ref{Ugraph}. 
There are 
$14$ isomorphism classes with $\delta_\varphi = 0$ and 
$49$ isomorphism classes with $\delta_\varphi = 1$. \\

By Theorem \ref{relinvtheorem}, we have for related involutions  
\begin{equation}
r(\varphi)+r(\tau\varphi)=20,\ \ a(\varphi) = a(\tau\varphi)\ \ \text{and}\ \ 
\delta_{\varphi S}= \delta_{\tau\varphi S},\ s_\varphi=s_{\tau\varphi}. 
\label{relatedF4}
\end{equation}
It corresponds to the symmetry with respect to the line $r=10$ of 
Figure \ref{Ugraph}. 
Thus, if we identify related involutions, there are 
$10$ isomorphism classes with $\delta_\varphi = 0$ and 
$27$ isomorphism classes with $\delta_\varphi = 1$. \\

\medskip

Now let us consider the geometric interpretation of the above results. 
Denote by $s$ the exceptional rational section with $s^2=-4$ of 
$\bff=\bff_4$, and by $c$ the fiber of the natural fibration $f:\bff_4\to s$. 
One has $c^2=0$ and $p_g(c)=0$. We have $-2K_\bff=12c+4s$. Thus, for 
$A\in |-2K_\bff|$, one has $A\cdot s=-4$. It follows  
$A=s+A_1$ where $A_1\in |12c+3c|$. We have $c\cdot A_1=3$ and $s\cdot A_1=0$. 
It follows that a non-singular $A\in |-2K_\bff|$ has two irreducible 
components: $A=s+A_1$. Thus, we are describing connected components of 
moduli of non-singular curves $A_1$ in the linear system $|12c+3s|$. 
Any such $A_1$ gives a non-singular $A=s+A_1$ from $|-2K_\bff|$. 

We mention that the lattice $S\cong U$ is $\bz C+\bz E$ where 
$C=\pi^\ast (c)$ and $E=\pi^\ast(s)/2$. We have $C^2=0$, $E^2=-2$ and 
$C\cdot E=1$. As usual, we denote $\pi:X\to Y=\bff_4$ the quotient morphism. 

If $\bff(\br)$ is not empty, then $\bff(\br)$ is a torus and 
$s(\br)$ is a circle giving a generator of the torus. The curves 
$c(\br)$ where $f(c)\in s(\br)$ give another generator of the torus. 
It follows that a real curve $A_1(\br)$ belongs to the open cylinder 
$\bff(\br)-s(\br)$ with its infinity identified with two copies of 
$s(\br)$.

\begin{figure} 
\centerline{\includegraphics[width=3cm]{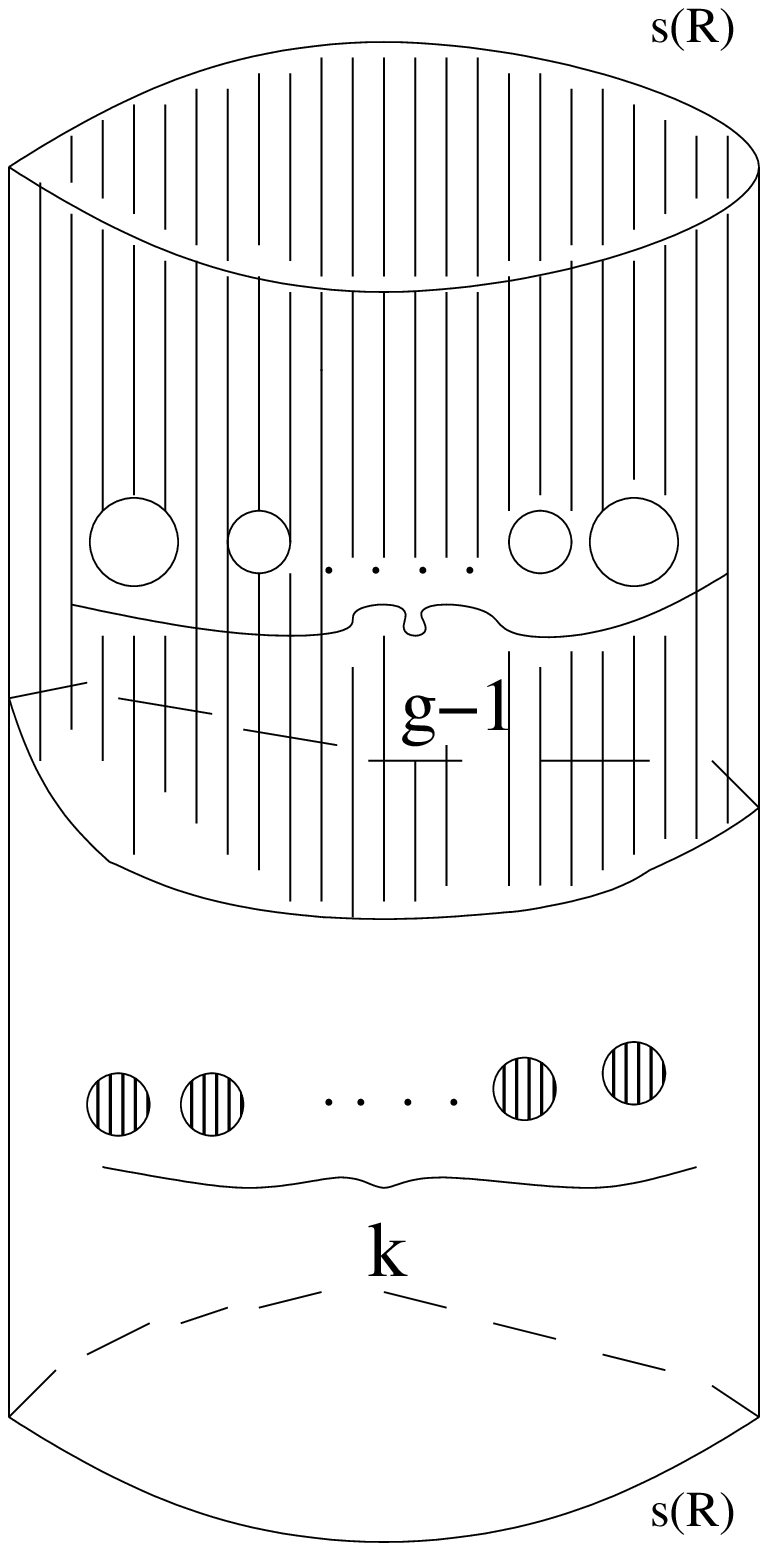}} 
\caption{$\bff_4$: $A^+$ with  
$(r,a,\delta_{\varphi S},v)\not=(10,8,0,0),\ (10,10,0,0)$.}
\label{Ugraph1}
\end{figure}

\begin{figure} 
\centerline{\includegraphics[width=3cm]{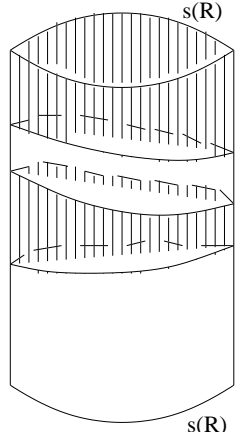}} 
\caption{$\bff_4$: $A^+$ with  
$(r,a,\delta_{\varphi S},v)=(10,8,0,0)$ ($\bff_4(\br)=\emptyset$ if 
$(r,a,\delta_{\varphi S},v)=(10,10,0,0)$).}
\label{Ugraph2}
\end{figure}

Using $c\cdot A_1=3$ and \eqref{realcomponents}, \eqref{realmod2} 
for both involutions $\varphi$ and $\tau\varphi$, we get at once that 
a positive curve $A^+$ with the invariants \eqref{invF_4} has the isotopy 
type given in Figures \ref{Ugraph1} and \ref{Ugraph2}. In Figures 
\ref{Ugraph1} and \ref{Ugraph2} we assume that 
$(r,a,\delta_\varphi)\not=(10,10,0)$. Then 
 $\bff(\br)$ is a torus. If $(r,a,\delta_\varphi)=(10,10,0)$, then   
$\bff(\br)$ is empty. As usual, $g=(22-r-a)/2$ and $k=(r-a)/2$.

Like in Sect. \ref{F1}, we obtain the following interpretation of 
the invariant $\delta_{\varphi}$: 
{\it one has $\delta_{\varphi}=0$, if and only if the curve $A$ is dividing: 
$A(\br)$ divides $A(\bc)$ in two connected parts or $A(\br)=\emptyset$. 
Equivalently, $A(\br)=0$ in $H_1(A(\bc),\bz/2)$.} 

Thus, finally, we get 
 
\begin{theorem}
\label{theoremF_4}
The connected component of moduli of positive real non-singular curves $A^+$,  
$A\in |-2K_{\bff_4}|$ up to the action of the automorphism group of 
$\bff_4$ is defined by the isotopy type and by divideness of 
$A(\bc)$ by $A(\br)$. All these possibilities are presented in 
Figures $\bff_4$ above. 
\end{theorem} 

We remark that the corresponding real K3 surfaces $(X,\tau,\varphi)$ 
give an important class of K3 surfaces. They are elliptic K3 surfaces  
(preimages of $c$) with section (the preimage of $s$). 
The holomorphic involution $\tau$ gives the inverse map of the elliptic pencil.

\bigskip

\section{Application: Real polarized K3 surfaces as deformations of 
general K3 double  rational scrolls}
\label{polK3}

Consider a non-degenerate real K3 surfaces $(X,\tau,\varphi)$ 
with a non-symplectic involution of the type $(S,\theta)$ where 
$(S,\theta)$ is of one of 5 types with $\rk S\le 2$ described in 
Sect. \ref{typesStheta}. Assume that $P\in S$ is a primitive 
$nef$ element with $P^2=n$, $n\in \bn$ is even, and such that 
$\theta(P)=-P$. Moreover, assume that the complete linear system 
$|P|$ or $|2P|$ gives the quotient map $\pi:X\to Y=X/\{1,\tau\}$. 
By standard results about linear systems on K3 surfaces, e. g. 
see \cite{Saint-Donat74}, it is valid in the following and only the 
following cases:  

{\it The Case of $\bp^2$, i. e. $S=\langle 2 \rangle$ and $\theta=-1$:} 
Then $P=h$ and $P^2=n=2$, both $|P|$ and $|2P|$ give the quotient 
map $\pi:X\to Y=\bp^2$. Thus, $n=2$.  

{\it The Case of ellipsoid $\bee$, i. e. $S=U(2)$ and $\theta(e_1)=-e_2$:} 
Then $P=e_1+e_2$, $P^2=n=4$ and $|P|$ gives the 
quotient map $\pi$ to the ellipsoid. Thus, $n=4$.  

{\it The Case of hyperboloid $\bhh$, i. e. $S=U(2)$ and $\theta=-1$:}
Then either $P=n_1e_1+e_2$ (we denote this case as $\bhh_{(1)}$) or 
$P=e_1+n_1e_2$ (we denote this case as $\bhh_{(2)}$), 
and $P^2=n=4n_1$, $n_1\ge 1$. Then $|P|$ gives the quotient map $\pi$ 
to the hyperboloid. Thus, $n=4n_1\equiv 0\mod 4$ and $n\ge 4$.  

{\it The Case of $\bff_1$, i. e. $S=\langle 2\rangle \oplus \langle -2 
\rangle$ and $\theta=-1$.} We have  $P=n_1{c}+e=n_1h+e(1-n_1)$ where 
$c=h-e$ and $P^2=n=4n_1-2$, $n_1\ge 2$. 
Then $|P|$ gives the quotient map to $\bff_1$. We have $n\equiv 2\mod 4$ 
and $n\ge 6$.  

{\it The case of $\bff_4$, i. e. $S=U$ and $\theta=-1$.} 
We denote $C=\pi^\ast(c)$ and $E=\pi^\ast(s)/2$ where $C^2=0$, 
$E^2=-2$ and $C\cdot E=1$.  We have two possibilities: 

{\it The case $(\bff_4)^{(1)}$:} $P=n_1C+E$ where $n_1\ge 3$ 
and $P^2=n=2n_1-2$, then $|2P|$ gives the quotient map and $|P|$ gives the composition of $\pi$ with the natural fibration $\bff_4\to \bp^1$. 
We have $n\equiv 0\mod 2$ and $n\ge 4$ for this case

{\it The case $(\bff_4)^{(2)}$:} $P=n_1C+2E$ where $n_1\ge 5$ is odd and 
$P^2=n=4n_1-8$,  then $|P|$ gives the quotient map $\pi$. We have 
$n \equiv 4\mod 8$ and $n\ge 12$ for this case. 

For each of these types: 
\begin{equation}
\Ta=\{\bp^2, \,\bee, \,\bhh_{(1)},\,\bhh_{(2)},\, 
\bff_1,\,(\bff_4)^{(1)},\,(\bff_4)^{(2)}\}
\label{typeT}
\end{equation}
the element $P$ is 
uniquely defined by $n=P^2$, if it exists. If for $t\in\Ta$ such $P$ 
does exist we write $n\in t$. 
For $t\in \Ta$ we denote by $\pi_0(\M_t)$ the set of connected 
components of moduli of non-degenerate $(X,\tau, \phi)$ of the type 
$t\in\Ta$ which we had described in Sects  \ref{P2} --- \ref{F4}. 

Since the branch curve is non-singular in all these cases, the $P$ 
is ample. In all these cases $|3P|$ gives an embedding of $X$ to a projective 
space. If $n\ge 4$, a small deformation of the pair $(X,P)$ gives a K3 surface 
with the very ample $P$, i. e. $|P|:X\subset \bp^{n/2+1}$. 

We denote by $\M_{n,k}$ where $n,\,k\in \bn$ and $n$ is even, 
the moduli space of real polarized K3 surfaces $(X, kP)$ with the polarization $kP$ where $P$ is primitive and 
$P^2=n$, and we require that the linear 
system $|kP|$ gives an embedding $|kP|:X\subset \bp^{k^2n/2+1}$, i. e. 
$kP$ is very ample.  Thus, in all cases above, for $k\ge 3$ 
the correspondence $(X,\tau,\phi)\mapsto (X,kP)$ gives an embedding of the moduli of $(X,\tau,\varphi)$ to the moduli $\M_{n,k}$, $k\ge 3$, of real polarized K3 surfaces $(X,kP)$ with the polarization $kP$ of the degree 
$nk^2$ where $P$ is primitive. Thus, for $k\ge 3$ we get the natural  
{\it deformation map}  
\begin{equation} 
D_n: \bigcup_{t\in \Ta:\ n\in t}{\pi_0(\M_t})\to \pi_0(\M_{n,k}). 
\label{deformdoubleK3}
\end{equation} 
A polarized K3 surface $(X,kP)$ which belongs to a connected component 
of moduli from the image of $D$ can be obtained by a deformation  
inside $\M_{n,k}$ from a real K3 surface which is a  
general double covering (i. e. with a non-singular branch curve) 
of a {\it rational scroll.}  All these scrolls and 
the types of their double coverings are given by the types $t\in \Ta$ 
such that $n\in t$. We call such K3 surfaces as 
{\it general K3 double rational scrolls}.  

One of consequences of  Theorem 3.10.1 in \cite{Nikulin79} (we shall 
discuss it in details later) is that the 
natural embedding   
\begin{equation}
\M_{n}=\M_{n,1}\to \M_{n,k},\ \ (X,P)\mapsto (X,kP), 
\label{multpolarization}
\end{equation}
gives the isomorphism of the set of connected components of $\M_{n}=\M_{n,1}$ 
and the set of connected components of $\M_{n,k}$ for $k\ge 2$. A difference 
between $\M_{n}$ and $\M_{n,k}$ for $k\ge 2$ is only on the 
boundary of $\M_{n}$ in $\M_{n,k}$, and this boundary (in spite of it 
has the real codimension one in $\M_{n,k}$) does not divide connected 
components of $\M_{n,k}$, $k\ge 2$.  

Thus, the deformation map $D_n$ in \eqref{deformdoubleK3} is 
defined for all $k\ge 1$. {\it The image of $D$ gives connected 
components of moduli of real polarized K3 surfaces $(X,kP)$ for any 
$k\ge 1$, which can be obtained as a deformation of a general real  
K3 double rational scroll.} 
Using results of previous sections and \cite{Nikulin79}, 
the deformation map $D_n$ can be computed explicitly and gives all 
{\it real polarized K3 surfaces $(X,kP)$, $k\ge 1$,    
which are deformations of general real K3 double rational scrolls. 
Moreover, we can even enumerate (classify) all types of these deformations.}   

A polarized real K3 surface $(X,\text{cong},kP)$, where $P$ is primitive 
and $k\ge 1$, defines an integral involution with a condition 
$(L=H_2(X,\bz),\varphi=\text{cong}_\ast, S_n=[P])$ where $S_n=[P]=\bz P\cong 
\langle n\rangle$ and $\theta=\phi|S_n=-1$. Theorem 3.10.1 in \cite{Nikulin79} 
gives that the connected component of moduli of 
$(X,\text{cong},kP)$ is defined by the isomorphism class of its 
integral involution $(L,\varphi,S_n)$ with the condition $(S_n,-1)$. 
All invariants of the genus of $(L,\varphi,S_n)$ are given by the 
invariants  
\begin{equation}
(n;r,a,\delta_P,\delta_\varphi, \delta_{\varphi P}) 
\label{genK3}
\end{equation}
which are the same as in \eqref{genusinvariants}. All possibilities for 
the genus invariants \eqref{genK3} are given in Theorem 3.4.3 
from \cite{Nikulin79} where one should put $l_{(+)}=3$, $l_{(-)}=19$, 
$t_{(+)}=1$ and $t_{(-)}=r-1$.  The genus \eqref{genK3} 
defines the class (Theorem 3.3.1 in \cite{Nikulin79}), 
if $r\le 18-\delta_{P \varphi}$ (actually, similarly to Theorem \ref{genus-isomor} one can give much better sufficient conditions).

Assume that the real K3 surface $(X,\tau,\varphi)$ has the 
integral involution $(L,\varphi, S)$ of the type $(S,\theta)$, and 
$(S,\theta)$ and $P\in S$ is of one of types $t\in \Ta$ from 
\eqref{typeT} such that $n\in t$. 
Assume that $(L,\varphi, S)$ has the invariants 
\eqref{geninv2} which are 
\begin{equation}
(r,a, H_+, H_-, \delta_{\varphi S},v)\ \ \ \text{where}\ \ 
\delta_\varphi=\left\{
\begin{array}{cl}
0 &\mbox{if}\ (\delta_{\varphi S},v)=(0,0)\\
1 &\mbox{otherwise}
\end{array}
\right. 
 \label{geninvgendouble}
\end{equation}
(we don't need $q_\rho$, if  $\rk S\le 2$). Then the integral involution 
of the polarized K3 $(X,P)$ is equal to $(L,\varphi, S_n=\bz P)$ where 
$S_n=\bz P\subset S$ is the sublattice generated by $P$. It follows that 
the invariants $(r,a,\delta_\varphi)$ of $(X,P)$ in \eqref{genK3} are 
same as in \eqref{geninvgendouble}. We have 
\begin{equation}
\delta_P=\left\{
\begin{array}{cl}
0 &\mbox{if}\ P\mod 2S_- \in H_-\\
1 &\mbox{otherwise}
\end{array}
\right. ,
\label{deltaP}
\end{equation}
and 
\begin{equation}
\delta_{\varphi P}=\left\{
\begin{array}{cl}
0 &\mbox{if}\ P\equiv v\mod 2S_-\\
1 &\mbox{otherwise}
\end{array}
\right. . 
\label{deltaphiP}
\end{equation}
Symbolically we write in this case 
\begin{equation}
D_n:(t;r,a, H_+, H_-, \delta_{\varphi S},v)\Longrightarrow
(n;r,a,\delta_P,\delta_\varphi, \delta_{\varphi P}). 
\label{DdoubleK3}
\end{equation}

Applying results of previous Sections (especially Theorem \ref{genmod}), 
we get  

\begin{theorem} Let us fix a type $t\in\Ta$ of a rational scroll. 
There exists a real polarized K3 surface $(X,kP)$ with the invariants 
$$
(n;r,a,\delta_P,\delta_\varphi, \delta_{\varphi P})
$$ 
which  is a deformation of a general real K3 double rational scroll 
$t$, if and only if $n\in t$ and   
$$
D_n:(t;r,a, H_+, H_-, \delta_{\varphi S},v)\Longrightarrow
(n;r,a,\delta_P,\delta_\varphi, \delta_{\varphi P}). 
$$
All these polarized K3 surfaces $(X,kP)$ belong to one connected component 
of moduli of real polarized K3 surfaces.  
\label{doubleK3}
\end{theorem}

Thus, to find all real polarized K3 surfaces which are deformations of 
general real K3 double rational scrolls, we need to find out when 
one has \eqref{DdoubleK3}. These calculations are routine, and we only 
formulate the result. Using notations of Sects. \ref{P2} --- \ref{F4}, we get:

\medskip 

\begin{equation}
\text{\bf Types of deformations of general real K3 double rational scrolls:} 
\label{typesdef}
\end{equation} 

\medskip 

{\it The Case of $\bp^2$:\ \ } $n\in \bp^2$, if and only if $n=2$; $P=h$;   
$$
D_2:(\bp^2;r,a,H=0,\delta_{\varphi S},v)
\Longrightarrow 
(2;r,a,\delta_P=1,\delta_{\varphi}=\delta_{\varphi S},\delta_{\varphi P}=1);
$$
$$
D_2:(\bp^2;r,a,H=[h],\delta_{\varphi S},v)
\Longrightarrow 
(2;r,a,\delta_P=0,\delta_{\varphi}=1,\delta_{\varphi P}=\delta_{\varphi S}).
$$
See Figures \ref{6_delta_h=1} and \ref{6_delta_h=0} for the list of all possible 
invariants $(\bp^2;r,a,H,\delta_{\varphi S},v)$.  

\bigskip

{\it The Case of ellipsoid $\bee$:\ \ } $n\in \bee$, if and inly if $n=4$; 
$P=e_1+e_2$;  
$$
D_4:(\bee;r,a,H=[e_1+e_2],\delta_{\varphi S},v)
\Longrightarrow
(4;r,a,\delta_P=0,\delta_\varphi, \delta_{\varphi P})\ \ \text{where} 
$$
$$
\delta_{\varphi P}=\left\{
\begin{array}{cl}
0 &\mbox{if}\ (\delta_{\varphi S},v)=(0,e_1+e_2)\\ 
1 &\mbox{otherwise} 
\end{array}\right .  \ . 
$$
See Figure \ref{ellipsoid-graph} for the list of all possible invariants 
$(\bee;r,a,H=[e_1+e_2],\delta_{\varphi S},v)$. 

\bigskip

{\it The Case of hyperboloid $\bhh_{(1)}$.} $n\in \bhh_{(1)}$, 
if and only if $n\equiv 0\mod 4$ and $n\ge 4$; $P=(n/4)e_1+e_2$;   
$$
D_n:(\bhh_{(1)};r,a,H=0,\delta_{\varphi S},v)
\Longrightarrow
(n;r,a,\delta_P=1,\delta_\varphi=\delta_{\varphi S}, \delta_{\varphi P}=1);
$$ 
$$
D_n:(\bhh_{(1)};r,a,H=[e_1,e_2],\delta_{\varphi S},v)
\Longrightarrow
(n;r,a,\delta_P=0,\delta_\varphi=\delta_{\varphi S}, 
\delta_{\varphi P}=1);
$$ 
$$
D_n:(\bhh_{(1)};r,a,H=[e_1],\delta_{\varphi S},v)
\Longrightarrow
(n;r,a,\delta_P=1,\delta_\varphi, \delta_{\varphi P}=1);
$$ 
$$
D_n:(\bhh_{(1)};r,a,H=[e_2],\delta_{\varphi S},v)
\Longrightarrow
(n;r,a,\delta_P,\delta_\varphi, \delta_{\varphi P})\ \ \text{where}
$$ 
$$
\delta_P=
\left\{\begin{array}{cl}
0 &\mbox{if}\ n\equiv 0\mod 8\\
1 &\mbox{otherwise}
\end{array}\right . ,\ \ 
\delta_{\varphi P}=\left\{
\begin{array}{cl}
0 &\mbox{if}\ (n~\mod 8, \delta_{\varphi S},v)=
(0~\mod 8,0,e_2)\\
1 &\mbox{otherwise} 
\end{array}\right . ; 
$$
$$
D_n:(\bhh_{(1)};r,a,H=[e_1+e_2],\delta_{\varphi S},v)
\Longrightarrow
(n;r,a,\delta_P,\delta_\varphi, \delta_{\varphi P}) \ \ \text{where} 
$$
$$
\delta_P=
\left\{\begin{array}{cl}
0 &\mbox{if}\ n\equiv 4\mod 8\\
1 &\mbox{otherwise}
\end{array}\right . ,\   
\delta_{\varphi P}=\left\{
\begin{array}{cl}
0 &\mbox{if}\ (n\mod 8,\delta_{\varphi S},v)=
(4\mod 8,0,e_1+e_2)\\
1 &\mbox{otherwise} 
\end{array}\right . . 
$$
See figures \ref{hyperboloid-0} --- \ref{hyperboloid-h} for the 
list of all possible invariants 
$(\bhh;r,a,H,\delta_{\varphi S},v)$.

\bigskip

{\it The Case of hyperboloid $\bhh_{(2)}$:} Change $e_1$ and $e_2$ 
places in the previous case.

\bigskip

{\it The Case of $\bff_1$:} $n\in \bff_1$, if and only if 
$n\equiv 2\mod 4$ and $n\ge 6$; $P=\left((n+2)/4\right)h+\left((2-n)/4\right)e$;  
$$
D_n:(\bff_1;r,a,H=0,\delta_{\varphi S},v)
\Longrightarrow
(n;r,a,\delta_P=1,\delta_\varphi=\delta_{\varphi S}, \delta_{\varphi P}=1);
$$ 
$$
D_n:(\bff_1;r,a,H=[h,e],\delta_{\varphi S},v)
\Longrightarrow
(n;r,a,\delta_P=0,\delta_\varphi=1, \delta_{\varphi P}=1);  
$$
$$
D_n:(\bff_1;r,a,H=[h],\delta_{\varphi S},v)
\Longrightarrow
(n;r,a,\delta_P,\delta_\varphi=1, \delta_{\varphi P})\ \text{where}
$$ 
$$
\delta_P=
\left\{\begin{array}{cl}
0 &\mbox{if}\ n\equiv 2\mod 8\\
1 &\mbox{otherwise}
\end{array}\right .\ ,\ \ 
\delta_{\varphi P}=\left\{
\begin{array}{cl}
0 &\mbox{if}\ (n~\mod 8, \delta_{\varphi S})=
(2~\mod 8,0)\\
1 &\mbox{otherwise} 
\end{array}\right . ; 
$$
$$
D_n:(\bff_1;r,a,H=[e],\delta_{\varphi S},v)
\Longrightarrow
(n;r,a,\delta_P,\delta_\varphi=1, \delta_{\varphi P})\ \ \text{where}
$$ 
$$
\delta_P=
\left\{\begin{array}{cl}
0 &\mbox{if}\ n\equiv -2\mod 8\\
1 &\mbox{otherwise}
\end{array}\right .\ ,\ \ 
\delta_{\varphi P}=\left\{
\begin{array}{cl}
0 &\mbox{if}\ (n~\mod 8, \delta_{\varphi S})=
(-2~\mod 8,0)\\
1 &\mbox{otherwise} 
\end{array}\right . ; 
$$
$$
D_n:(\bff_1;r,a,H=[h+e],\delta_{\varphi S},v)
\Longrightarrow
(n;r,a,\delta_P=1,\delta_\varphi, \delta_{\varphi P}=1).
$$
See the list of all possible invariants 
$(\bff_1;r,a,H,\delta_{\varphi S},v)$ in 
figures \ref{<2>+<-2>-0} --- \ref{<2>+<-2>-he}.

\bigskip

{\it The Case of $(\bff_4)^{(1)}$:} $n\in (\bff_4)^{(1)}$, if and only if 
$n\equiv 0\mod 2$ and $n\ge 4$; $P=(n/2+1)C+E$;   
$$
D_n:((\bff_4)^{(1)};r,a,H=0,\delta_{\varphi S},v)
\Longrightarrow
(n;r,a,\delta_P=1,\delta_\varphi=\delta_{\varphi S}, \delta_{\varphi P}=1). 
$$

\bigskip

{\it The Case of $(\bff_4)^{(2)}$:} $n\in (\bff_4)^{(2)}$, if and only if 
$n\equiv 4\mod 8$ and $n\ge 12$; \linebreak $P=(n/4+2)C+2E$;   
$$
D_n:((\bff_4)^{(2)};r,a,H=0,\delta_{\varphi S},v)
\Longrightarrow
(n;r,a,\delta_P=1,\delta_\varphi=\delta_{\varphi S}, \delta_{\varphi P}=1). 
$$

See the figure \ref{Ugraph} for the list of all possible invariants 
$(\bff_4;r,a,H,\delta_{\varphi S},v)$.

\bigskip

We can apply this calculations to answer the following interesting question: 
What are real K3 surfaces which are deformations of general real K3 
double rational scrolls?

We have 

\begin{proposition} A real polarized K3 surface $(X,kP)$ where $P$ is 
primitive and $n=P^2\ge 6$ cannot be a deformation of a general real K3 double rational scroll, if either $(r,a)=(20,2)$ (equivalently $X(\br)=(T_0)^{10}$) 
or $r+a=22$ and  $\delta_{\varphi P}=0$ (equivalently, 
$X(\br)=(T_0)^m$ and $X(\br)\sim P\mod 2$ in $H_2(X,\bz)$).  
\label{nondoubleK3}
\end{proposition}  

\Proof If $r=20$, then $L_\varphi=L_-$ is positive definite or the rank 2. 
It cannot contain $S_-$ of the rank two because $S_-$ is hyperbolic. 
Thus, $r=20$ may happen only for $\bp^2$ when $n=2$, and for ellipsoid 
when $n=4$. From the list of possible 
invariants $(r,a,\delta_\varphi)$ of real K3 surfaces 
(see Figure \ref{Sgraph}) we get that $a=2$ and $\delta_\varphi=1$, if 
$r=20$. By \eqref{realcomponents}, this is equivalent to $X(\br)=T_0^{10}$. 

Assume that $r+a=22$ and $n\ge 6$. This  is possible only for 
$(\bhh;H=[e_1,e_2])$ and $(\bff_1;H=[h,e])$, but all these cases 
give $\delta_{\varphi P}=1$ by deformations \eqref{typesdef}. 
By \eqref{realcomponents}, $r+a=22$ is equivalent to 
$X(\br)=(T_0)^m$, and $\delta_{\varphi P}=0$ is 
equivalent to $X(\br)\sim P\mod 2$ in $H_2(X,\bz)$, by \eqref{realmod2}.

This proves the statement  
\medskip

We have 

\begin{theorem} A genus invariant   
$(n;r,a,\delta_P,\delta_\varphi,\delta_{\varphi P})$ 
of a real polarized K3 surface $(X,kP)$  
can be obtained by the deformation $D_n$ of   
some general real K3 double rational scroll $t\in \Ta$ from  
\eqref{typeT}, if and only if the following condition is valid:   
\begin{equation} 
n\le  4,\  \text{if either}\ (r,a)=(20,2)\ \text{or}\ r+a=22\ \text{and}\  
\delta_{\varphi P}=0 
\label{exceptions} 
\end{equation} 
(we just exclude cases of Proposition \ref{nondoubleK3}). All possible types of 
these deformations are given in \eqref{typesdef}.   
\end{theorem}

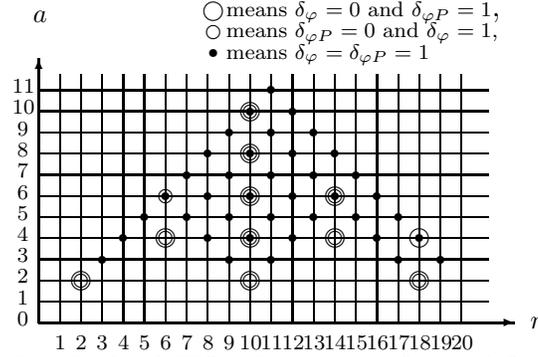
\begin{figure}
\begin{picture}(200,140)
%% K3- n equiv 0mod8, delta_P=0

\put(66,118){\circle{7}}
\put(71,116){{\tiny means $\delta_\varphi=0$ and $\delta_{\varphi P}=1$},}
\put(66,110){\circle{5}}
\put(71,108){{\tiny means $\delta_{\varphi P}=0$ and $\delta_\varphi=1$,}}
\put(66,102){\circle*{3}}
\put(71,100){{\tiny means $\delta_\varphi=\delta_{\varphi P}=1$}}
\multiput(8,0)(8,0){20}{\line(0,1){94}}
\multiput(0,8)(0,8){11}{\line(1,0){170}}
\put(0,0){\vector(0,1){100}}
\put(0,0){\vector(1,0){180}}
\put(  6,-10){{\tiny $1$}}
\put( 14,-10){{\tiny $2$}}
\put( 22,-10){{\tiny $3$}}
\put( 30,-10){{\tiny $4$}}
\put( 38,-10){{\tiny $5$}}
\put( 46,-10){{\tiny $6$}}
\put( 54,-10){{\tiny $7$}}
\put( 62,-10){{\tiny $8$}}
\put( 70,-10){{\tiny $9$}}
\put( 76,-10){{\tiny $10$}}
\put( 84,-10){{\tiny $11$}}
\put( 92,-10){{\tiny $12$}}
\put(100,-10){{\tiny $13$}}
\put(108,-10){{\tiny $14$}}
\put(116,-10){{\tiny $15$}}
\put(124,-10){{\tiny $16$}}
\put(132,-10){{\tiny $17$}}
\put(140,-10){{\tiny $18$}}
\put(148,-10){{\tiny $19$}}
\put(156,-10){{\tiny $20$}}

\put(-8, -1){{\tiny $0$}}
\put(-8,  7){{\tiny $1$}}
\put(-8, 15){{\tiny $2$}}
\put(-8, 23){{\tiny $3$}}
\put(-8, 31){{\tiny $4$}}
\put(-8, 39){{\tiny $5$}}
\put(-8, 47){{\tiny $6$}}
\put(-8, 55){{\tiny $7$}}
\put(-8, 63){{\tiny $8$}}
\put(-8, 71){{\tiny $9$}}
\put(-10, 79){{\tiny $10$}}
\put(-10, 87){{\tiny $11$}}
\put( -2,114){{\footnotesize $a$}} %vertical
\put(186, -2){{\footnotesize $r$}} %horizontal

%% \delta_\varphi = 0
%%  ( 8r,8a)
%\put( 16, 0){\circle{7}}
%\put( 80, 0){\circle{7}}
%\put(144, 0){\circle{7}}
\put( 16,16){\circle{7}}
%\put( 48,16){\circle{7}}
\put( 80,16){\circle{7}}
%\put(112,16){\circle{7}}
\put(144,16){\circle{7}}
\put( 48,32){\circle{7}}
\put( 80,32){\circle{7}}
\put(112,32){\circle{7}}
\put( 80,48){\circle{7}}
\put( 112,48){\circle{7}}
\put( 80,64){\circle{7}}
\put( 80,80){\circle{7}}    %%14data

%% \delta_\varphi = 1
%%  ( 8r,8a)
%\put(  8, 8){\circle*{3}}
%\put( 24, 8){\circle*{3}}
%\put( 72, 8){\circle*{3}}
%\put( 88, 8){\circle*{3}}
%\put(136, 8){\circle*{3}}
%\put(152, 8){\circle*{3}}
%\put( 16,16){\circle*{3}}
\put( 16,16){\circle{5}}
%\put( 32,16){\circle*{3}}
%\put( 64,16){\circle*{3}}
%\put( 80,16){\circle*{3}}
\put( 80,16){\circle{5}}
\put( 80,32){\circle{5}}
\put( 80,48){\circle{5}}
\put( 80,64){\circle{5}}
\put( 80,80){\circle{5}}
%\put( 96,16){\circle*{3}}
%\put(128,16){\circle*{3}}
%\put(144,16){\circle*{3}}
\put(144,16){\circle{5}}
\put( 24,24){\circle*{3}}
%\put( 40,24){\circle*{3}}
%\put( 56,24){\circle*{3}}
\put( 72,24){\circle*{3}}
\put( 88,24){\circle*{3}}
%\put(104,24){\circle*{3}}
%\put(120,24){\circle*{3}}
\put(136,24){\circle*{3}}
\put( 32,32){\circle*{3}}
%\put( 48,32){\circle*{3}}
\put( 48,32){\circle{5}}
\put( 48,48){\circle{5}}
\put( 64,32){\circle*{3}}
\put( 80,32){\circle*{3}}
\put( 96,32){\circle*{3}}
%\put(112,32){\circle*{3}}
\put(112,32){\circle{5}}
\put(112,48){\circle{5}}
\put(128,32){\circle*{3}}
\put( 40,40){\circle*{3}}
\put( 56,40){\circle*{3}}
\put( 72,40){\circle*{3}}
\put( 88,40){\circle*{3}}
\put(104,40){\circle*{3}}
\put(120,40){\circle*{3}}
\put( 48,48){\circle*{3}}
\put( 64,48){\circle*{3}}
\put( 80,48){\circle*{3}}
\put( 96,48){\circle*{3}}
\put(112,48){\circle*{3}}
\put( 56,56){\circle*{3}}
\put( 72,56){\circle*{3}}
\put( 88,56){\circle*{3}}
\put(104,56){\circle*{3}}
\put( 64,64){\circle*{3}}
\put( 80,64){\circle*{3}}
\put( 96,64){\circle*{3}}
\put( 72,72){\circle*{3}}
\put( 88,72){\circle*{3}}
\put( 80,80){\circle*{3}}   %%49data
\put( 88,88){\circle*{3}}
\put( 96,80){\circle*{3}}
\put( 104,72){\circle*{3}}
\put( 112,64){\circle*{3}}
\put( 120,56){\circle*{3}}
\put( 128,48){\circle*{3}}
\put( 136,40){\circle*{3}}
\put( 144,32){\circle*{3}}
\put( 144,32){\circle{7}}
\put( 152,24){\circle*{3}}
%%\put( 160,16){\circle*{3}}
\end{picture}

\caption{Real polarized K3 when $n\equiv 0\mod 8$ and $\delta_P=0$ 
without $(r,a)=(20,2)$ and $(r+a=22, \delta_{\varphi P}=0)$}
\label{K3-0mod8-0}
\end{figure}

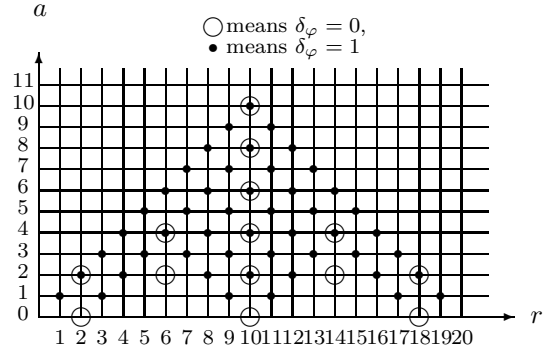
\begin{figure}
\begin{picture}(200,140)
%% K3-n equiv 0mod8 delta_P=1
\put(66,110){\circle{7}}
\put(71,108){{\tiny means $\delta_\varphi=0$,}}
\put(66,102){\circle*{3}}
\put(71,100){{\tiny means $\delta_\varphi=1$}}

\multiput(8,0)(8,0){20}{\line(0,1){94}}
\multiput(0,8)(0,8){11}{\line(1,0){170}}
\put(0,0){\vector(0,1){100}}
\put(0,0){\vector(1,0){180}}
\put(  6,-10){{\tiny $1$}}
\put( 14,-10){{\tiny $2$}}
\put( 22,-10){{\tiny $3$}}
\put( 30,-10){{\tiny $4$}}
\put( 38,-10){{\tiny $5$}}
\put( 46,-10){{\tiny $6$}}
\put( 54,-10){{\tiny $7$}}
\put( 62,-10){{\tiny $8$}}
\put( 70,-10){{\tiny $9$}}
\put( 76,-10){{\tiny $10$}}
\put( 84,-10){{\tiny $11$}}
\put( 92,-10){{\tiny $12$}}
\put(100,-10){{\tiny $13$}}
\put(108,-10){{\tiny $14$}}
\put(116,-10){{\tiny $15$}}
\put(124,-10){{\tiny $16$}}
\put(132,-10){{\tiny $17$}}
\put(140,-10){{\tiny $18$}}
\put(148,-10){{\tiny $19$}}
\put(156,-10){{\tiny $20$}}

\put(-8, -1){{\tiny $0$}}
\put(-8,  7){{\tiny $1$}}
\put(-8, 15){{\tiny $2$}}
\put(-8, 23){{\tiny $3$}}
\put(-8, 31){{\tiny $4$}}
\put(-8, 39){{\tiny $5$}}
\put(-8, 47){{\tiny $6$}}
\put(-8, 55){{\tiny $7$}}
\put(-8, 63){{\tiny $8$}}
\put(-8, 71){{\tiny $9$}}
\put(-10, 79){{\tiny $10$}}
\put(-10, 87){{\tiny $11$}}
\put( -2,114){{\footnotesize $a$}} %vertical
\put(186, -2){{\footnotesize $r$}} %horizontal

%% \delta_\varphi = 0
%%  ( 8r,8a)
\put( 16, 0){\circle{7}}
\put( 80, 0){\circle{7}}
\put(144, 0){\circle{7}}
\put( 16,16){\circle{7}}
\put( 48,16){\circle{7}}
\put( 80,16){\circle{7}}
\put(112,16){\circle{7}}
\put(144,16){\circle{7}}
\put( 48,32){\circle{7}}
\put( 80,32){\circle{7}}
\put(112,32){\circle{7}}
\put( 80,48){\circle{7}}
%\put( 112,48){\circle{7}}
\put( 80,64){\circle{7}}
\put( 80,80){\circle{7}}    %%14data

%% \delta_\varphi = 1
%%  ( 8r,8a)
\put(  8, 8){\circle*{3}}
\put( 24, 8){\circle*{3}}
\put( 72, 8){\circle*{3}}
\put( 88, 8){\circle*{3}}
\put(136, 8){\circle*{3}}
\put(152, 8){\circle*{3}}
\put( 16,16){\circle*{3}}
\put( 32,16){\circle*{3}}
\put( 64,16){\circle*{3}}
\put( 80,16){\circle*{3}}
\put( 96,16){\circle*{3}}
\put(128,16){\circle*{3}}
\put(144,16){\circle*{3}}
\put( 24,24){\circle*{3}}
\put( 40,24){\circle*{3}}
\put( 56,24){\circle*{3}}
\put( 72,24){\circle*{3}}
\put( 88,24){\circle*{3}}
\put(104,24){\circle*{3}}
\put(120,24){\circle*{3}}
\put(136,24){\circle*{3}}
\put( 32,32){\circle*{3}}
\put( 48,32){\circle*{3}}
\put( 64,32){\circle*{3}}
\put( 80,32){\circle*{3}}
\put( 96,32){\circle*{3}}
\put(112,32){\circle*{3}}
\put(128,32){\circle*{3}}
\put( 40,40){\circle*{3}}
\put( 56,40){\circle*{3}}
\put( 72,40){\circle*{3}}
\put( 88,40){\circle*{3}}
\put(104,40){\circle*{3}}
\put(120,40){\circle*{3}}
\put( 48,48){\circle*{3}}
\put( 64,48){\circle*{3}}
\put( 80,48){\circle*{3}}
\put( 96,48){\circle*{3}}
\put(112,48){\circle*{3}}
\put( 56,56){\circle*{3}}
\put( 72,56){\circle*{3}}
\put( 88,56){\circle*{3}}
\put(104,56){\circle*{3}}
\put( 64,64){\circle*{3}}
\put( 80,64){\circle*{3}}
\put( 96,64){\circle*{3}}
\put( 72,72){\circle*{3}}
\put( 88,72){\circle*{3}}
\put( 80,80){\circle*{3}}   %%49data
\end{picture}

\caption{Real polarized K3 when $n\equiv 0\mod 8$ and $\delta_P=1$ 
(then $\delta_{\varphi P}=1$)}
\label{K3-0mod8-1}
\end{figure}

\begin{figure}
\begin{picture}(200,140)
%% K3- n equiv 4mod8, delta_P=0
\put(66,118){\circle{7}}
\put(71,116){{\tiny means $\delta_\varphi=0$ and $\delta_{\varphi P}=1$},}
\put(66,110){\circle{5}}
\put(71,108){{\tiny means $\delta_{\varphi P}=0$ and $\delta_\varphi=1$,}}
\put(66,102){\circle*{3}}
\put(71,100){{\tiny means $\delta_\varphi=\delta_{\varphi P}=1$}}

\multiput(8,0)(8,0){20}{\line(0,1){94}}
\multiput(0,8)(0,8){11}{\line(1,0){170}}
\put(0,0){\vector(0,1){100}}
\put(0,0){\vector(1,0){180}}
\put(  6,-10){{\tiny $1$}}
\put( 14,-10){{\tiny $2$}}
\put( 22,-10){{\tiny $3$}}
\put( 30,-10){{\tiny $4$}}
\put( 38,-10){{\tiny $5$}}
\put( 46,-10){{\tiny $6$}}
\put( 54,-10){{\tiny $7$}}
\put( 62,-10){{\tiny $8$}}
\put( 70,-10){{\tiny $9$}}
\put( 76,-10){{\tiny $10$}}
\put( 84,-10){{\tiny $11$}}
\put( 92,-10){{\tiny $12$}}
\put(100,-10){{\tiny $13$}}
\put(108,-10){{\tiny $14$}}
\put(116,-10){{\tiny $15$}}
\put(124,-10){{\tiny $16$}}
\put(132,-10){{\tiny $17$}}
\put(140,-10){{\tiny $18$}}
\put(148,-10){{\tiny $19$}}
\put(156,-10){{\tiny $20$}}

\put(-8, -1){{\tiny $0$}}
\put(-8,  7){{\tiny $1$}}
\put(-8, 15){{\tiny $2$}}
\put(-8, 23){{\tiny $3$}}
\put(-8, 31){{\tiny $4$}}
\put(-8, 39){{\tiny $5$}}
\put(-8, 47){{\tiny $6$}}
\put(-8, 55){{\tiny $7$}}
\put(-8, 63){{\tiny $8$}}
\put(-8, 71){{\tiny $9$}}
\put(-10, 79){{\tiny $10$}}
\put(-10, 87){{\tiny $11$}}
\put( -2,114){{\footnotesize $a$}} %vertical
\put(186, -2){{\footnotesize $r$}} %horizontal

%% \delta_\varphi = 0
%%  ( 8r,8a)
%\put( 16, 0){\circle{7}}
%\put( 80, 0){\circle{7}}
%\put(144, 0){\circle{7}}
\put( 16,16){\circle{7}}
\put( 48,16){\circle{7}}
\put( 80,16){\circle{7}}
\put(112,16){\circle{7}}
\put(144,16){\circle{7}}
\put( 48,32){\circle{7}}
\put( 80,32){\circle{7}}
\put(112,32){\circle{7}}
\put( 80,48){\circle{7}}
\put( 112,48){\circle{7}}
\put( 80,64){\circle{7}}
\put( 80,80){\circle{7}}    %%14data
%% \delta_\varphi = 1
%%  ( 8r,8a)
%\put(  8, 8){\circle*{3}}
%\put( 24, 8){\circle*{3}}
%\put( 72, 8){\circle*{3}}
%\put( 88, 8){\circle*{3}}
%\put(136, 8){\circle*{3}}
%\put(152, 8){\circle*{3}}
%\put( 16,16){\circle*{3}}
%\put( 32,16){\circle*{3}}
\put( 32,16){\circle{5}}
\put( 32,32){\circle{5}}
%\put( 64,16){\circle*{3}}
\put( 64,16){\circle{5}}
\put( 64,32){\circle{5}}
\put( 64,48){\circle{5}}
\put( 64,64){\circle{5}}
%\put( 80,16){\circle*{3}}
%\put( 96,16){\circle*{3}}
\put( 96,16){\circle{5}}
\put( 96,32){\circle{5}}
\put( 96,48){\circle{5}}
\put( 96,64){\circle{5}}
%\put(128,16){\circle*{3}}
\put(128,16){\circle{5}}
\put(128,32){\circle{5}}
%\put(144,16){\circle*{3}}
\put( 24,24){\circle*{3}}
\put( 40,24){\circle*{3}}
\put( 56,24){\circle*{3}}
\put( 72,24){\circle*{3}}
\put( 88,24){\circle*{3}}
\put(104,24){\circle*{3}}
\put(120,24){\circle*{3}}
\put(136,24){\circle*{3}}
\put( 32,32){\circle*{3}}
\put( 48,32){\circle*{3}}
\put( 64,32){\circle*{3}}
\put( 80,32){\circle*{3}}
\put( 96,32){\circle*{3}}
\put(112,32){\circle*{3}}
\put(128,32){\circle*{3}}
\put( 40,40){\circle*{3}}
\put( 56,40){\circle*{3}}
\put( 72,40){\circle*{3}}
\put( 88,40){\circle*{3}}
\put(104,40){\circle*{3}}
\put(120,40){\circle*{3}}
\put( 48,48){\circle*{3}}
\put( 64,48){\circle*{3}}
\put( 80,48){\circle*{3}}
\put( 96,48){\circle*{3}}
\put(112,48){\circle*{3}}
\put( 56,56){\circle*{3}}
\put( 72,56){\circle*{3}}
\put( 88,56){\circle*{3}}
\put(104,56){\circle*{3}}
\put( 64,64){\circle*{3}}
\put( 80,64){\circle*{3}}
\put( 96,64){\circle*{3}}
\put( 72,72){\circle*{3}}
\put( 88,72){\circle*{3}}
\put( 80,80){\circle*{3}}   %%49data
\put( 88,88){\circle*{3}}
\put( 96,80){\circle*{3}}
\put( 104,72){\circle*{3}}
\put( 112,64){\circle*{3}}
\put( 120,56){\circle*{3}}
\put( 128,48){\circle*{3}}
\put( 136,40){\circle*{3}}
\put( 144,32){\circle*{3}}
\put( 144,32){\circle{7}}
\put( 152,24){\circle*{3}}
%%\put( 160,16){\circle*{3}}
\end{picture}

\caption{Real polarized K3 when $n\equiv 4\mod 8$ and $\delta_P=0$ without 
$(r,a)=(20,2)$ and $(r+a=22, \delta_{\varphi P}=0)$}
\label{K3-4mod8-0}
\end{figure}
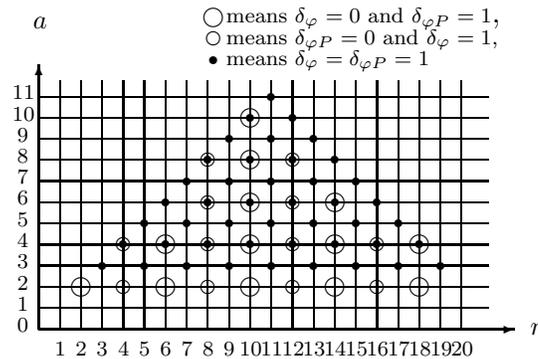

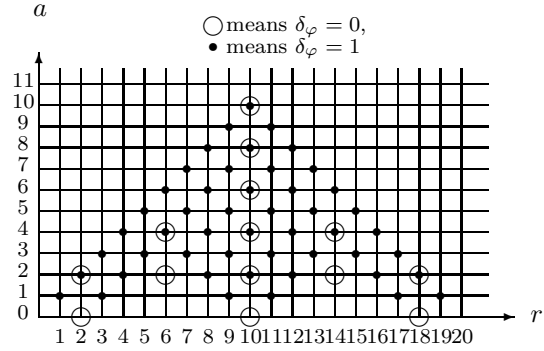
\begin{figure}
\begin{picture}(200,140)
%% K3-n equiv 4mod8 delta_P=1
\put(66,110){\circle{7}}
\put(71,108){{\tiny means $\delta_\varphi=0$,}}
\put(66,102){\circle*{3}}
\put(71,100){{\tiny means $\delta_\varphi=1$}}

\multiput(8,0)(8,0){20}{\line(0,1){94}}
\multiput(0,8)(0,8){11}{\line(1,0){170}}
\put(0,0){\vector(0,1){100}}
\put(0,0){\vector(1,0){180}}
\put(  6,-10){{\tiny $1$}}
\put( 14,-10){{\tiny $2$}}
\put( 22,-10){{\tiny $3$}}
\put( 30,-10){{\tiny $4$}}
\put( 38,-10){{\tiny $5$}}
\put( 46,-10){{\tiny $6$}}
\put( 54,-10){{\tiny $7$}}
\put( 62,-10){{\tiny $8$}}
\put( 70,-10){{\tiny $9$}}
\put( 76,-10){{\tiny $10$}}
\put( 84,-10){{\tiny $11$}}
\put( 92,-10){{\tiny $12$}}
\put(100,-10){{\tiny $13$}}
\put(108,-10){{\tiny $14$}}
\put(116,-10){{\tiny $15$}}
\put(124,-10){{\tiny $16$}}
\put(132,-10){{\tiny $17$}}
\put(140,-10){{\tiny $18$}}
\put(148,-10){{\tiny $19$}}
\put(156,-10){{\tiny $20$}}

\put(-8, -1){{\tiny $0$}}
\put(-8,  7){{\tiny $1$}}
\put(-8, 15){{\tiny $2$}}
\put(-8, 23){{\tiny $3$}}
\put(-8, 31){{\tiny $4$}}
\put(-8, 39){{\tiny $5$}}
\put(-8, 47){{\tiny $6$}}
\put(-8, 55){{\tiny $7$}}
\put(-8, 63){{\tiny $8$}}
\put(-8, 71){{\tiny $9$}}
\put(-10, 79){{\tiny $10$}}
\put(-10, 87){{\tiny $11$}}
\put( -2,114){{\footnotesize $a$}} %vertical
\put(186, -2){{\footnotesize $r$}} %horizontal

%% \delta_\varphi = 0
%%  ( 8r,8a)
\put( 16, 0){\circle{7}}
\put( 80, 0){\circle{7}}
\put(144, 0){\circle{7}}
\put( 16,16){\circle{7}}
\put( 48,16){\circle{7}}
\put( 80,16){\circle{7}}
\put(112,16){\circle{7}}
\put(144,16){\circle{7}}
\put( 48,32){\circle{7}}
\put( 80,32){\circle{7}}
\put(112,32){\circle{7}}
\put( 80,48){\circle{7}}
%\put( 112,48){\circle{7}}
\put( 80,64){\circle{7}}
\put( 80,80){\circle{7}}    %%14data

%% \delta_\varphi = 1
%%  ( 8r,8a)
\put(  8, 8){\circle*{3}}
\put( 24, 8){\circle*{3}}
\put( 72, 8){\circle*{3}}
\put( 88, 8){\circle*{3}}
\put(136, 8){\circle*{3}}
\put(152, 8){\circle*{3}}
\put( 16,16){\circle*{3}}
\put( 32,16){\circle*{3}}
\put( 64,16){\circle*{3}}
\put( 80,16){\circle*{3}}
\put( 96,16){\circle*{3}}
\put(128,16){\circle*{3}}
\put(144,16){\circle*{3}}
\put( 24,24){\circle*{3}}
\put( 40,24){\circle*{3}}
\put( 56,24){\circle*{3}}
\put( 72,24){\circle*{3}}
\put( 88,24){\circle*{3}}
\put(104,24){\circle*{3}}
\put(120,24){\circle*{3}}
\put(136,24){\circle*{3}}
\put( 32,32){\circle*{3}}
\put( 48,32){\circle*{3}}
\put( 64,32){\circle*{3}}
\put( 80,32){\circle*{3}}
\put( 96,32){\circle*{3}}
\put(112,32){\circle*{3}}
\put(128,32){\circle*{3}}
\put( 40,40){\circle*{3}}
\put( 56,40){\circle*{3}}
\put( 72,40){\circle*{3}}
\put( 88,40){\circle*{3}}
\put(104,40){\circle*{3}}
\put(120,40){\circle*{3}}
\put( 48,48){\circle*{3}}
\put( 64,48){\circle*{3}}
\put( 80,48){\circle*{3}}
\put( 96,48){\circle*{3}}
\put(112,48){\circle*{3}}
\put( 56,56){\circle*{3}}
\put( 72,56){\circle*{3}}
\put( 88,56){\circle*{3}}
\put(104,56){\circle*{3}}
\put( 64,64){\circle*{3}}
\put( 80,64){\circle*{3}}
\put( 96,64){\circle*{3}}
\put( 72,72){\circle*{3}}
\put( 88,72){\circle*{3}}
\put( 80,80){\circle*{3}}   %%49data
\end{picture}

\caption{Real polarized K3 when $n\equiv 4\mod 8$ and $\delta_P=1$ 
(then $\delta_{\varphi P}=1$)}
\label{K3-4mod8-1}
\end{figure}

\begin{figure}
\begin{picture}(200,140)
%% K3- n equiv 2mod8, delta_P=0
\put(66,110){\circle{5}}
\put(71,108){{\tiny means $\delta_{\varphi P}=0$ and $\delta_\varphi=1$,}}
\put(66,102){\circle*{3}}
\put(71,100){{\tiny means $\delta_\varphi=\delta_{\varphi P}=1$}}

\multiput(8,0)(8,0){20}{\line(0,1){94}}
\multiput(0,8)(0,8){11}{\line(1,0){170}}
\put(0,0){\vector(0,1){100}}
\put(0,0){\vector(1,0){180}}
\put(  6,-10){{\tiny $1$}}
\put( 14,-10){{\tiny $2$}}
\put( 22,-10){{\tiny $3$}}
\put( 30,-10){{\tiny $4$}}
\put( 38,-10){{\tiny $5$}}
\put( 46,-10){{\tiny $6$}}
\put( 54,-10){{\tiny $7$}}
\put( 62,-10){{\tiny $8$}}
\put( 70,-10){{\tiny $9$}}
\put( 76,-10){{\tiny $10$}}
\put( 84,-10){{\tiny $11$}}
\put( 92,-10){{\tiny $12$}}
\put(100,-10){{\tiny $13$}}
\put(108,-10){{\tiny $14$}}
\put(116,-10){{\tiny $15$}}
\put(124,-10){{\tiny $16$}}
\put(132,-10){{\tiny $17$}}
\put(140,-10){{\tiny $18$}}
\put(148,-10){{\tiny $19$}}
\put(156,-10){{\tiny $20$}}

\put(-8, -1){{\tiny $0$}}
\put(-8,  7){{\tiny $1$}}
\put(-8, 15){{\tiny $2$}}
\put(-8, 23){{\tiny $3$}}
\put(-8, 31){{\tiny $4$}}
\put(-8, 39){{\tiny $5$}}
\put(-8, 47){{\tiny $6$}}
\put(-8, 55){{\tiny $7$}}
\put(-8, 63){{\tiny $8$}}
\put(-8, 71){{\tiny $9$}}
\put(-10, 79){{\tiny $10$}}
\put(-10, 87){{\tiny $11$}}
\put( -2,114){{\footnotesize $a$}} %vertical
\put(186, -2){{\footnotesize $r$}} %horizontal

%% \delta_\varphi = 0
%%  ( 8r,8a)
%\put( 16, 0){\circle{7}}
%\put( 80, 0){\circle{7}}
%\put(144, 0){\circle{7}}
%\put( 16,16){\circle{7}}
%\put( 48,16){\circle{7}}
%\put( 80,16){\circle{7}}
%\put(112,16){\circle{7}}
%\put(144,16){\circle{7}}
%\put( 48,32){\circle{7}}
%\put( 80,32){\circle{7}}
%\put(112,32){\circle{7}}
%\put( 80,48){\circle{7}}
%\put( 112,48){\circle{7}}
%\put( 80,64){\circle{7}}
%\put( 80,80){\circle{7}}    %%14data
%% \delta_\varphi = 1
%%  ( 8r,8a)
%\put(  8, 8){\circle*{3}}
%\put( 24, 8){\circle*{3}}
\put( 24, 8){\circle{5}}
\put( 24, 24){\circle{5}}
%\put( 72, 8){\circle*{3}}
%\put( 88, 8){\circle*{3}}
\put (88,8){\circle{5}}
\put( 88, 24){\circle{5}}
\put( 88, 40){\circle{5}}
\put( 88, 56){\circle{5}}
\put( 88, 72){\circle{5}}
%\put(136, 8){\circle*{3}}
%\put(152, 8){\circle*{3}}
\put(152, 8){\circle{5}}
\put( 16,16){\circle*{3}}
\put( 32,16){\circle*{3}}
%\put( 64,16){\circle*{3}}
\put( 80,16){\circle*{3}}
\put( 96,16){\circle*{3}}
%\put(128,16){\circle*{3}}
\put(144,16){\circle*{3}}
\put( 24,24){\circle*{3}}
\put( 40,24){\circle*{3}}
%\put( 56,24){\circle*{3}}
\put( 56,24){\circle{5}}
\put( 56,40){\circle{5}}
\put( 56,56){\circle{5}}
\put( 72,24){\circle*{3}}
\put( 88,24){\circle*{3}}
\put(104,24){\circle*{3}}
%\put(120,24){\circle*{3}}
\put(120,24){\circle{5}}
\put(120,40){\circle{5}}
\put(136,24){\circle*{3}}
\put( 32,32){\circle*{3}}
\put( 48,32){\circle*{3}}
\put( 64,32){\circle*{3}}
\put( 80,32){\circle*{3}}
\put( 96,32){\circle*{3}}
\put(112,32){\circle*{3}}
\put(128,32){\circle*{3}}
\put( 40,40){\circle*{3}}
\put( 56,40){\circle*{3}}
\put( 72,40){\circle*{3}}
\put( 88,40){\circle*{3}}
\put(104,40){\circle*{3}}
\put(120,40){\circle*{3}}
\put( 48,48){\circle*{3}}
\put( 64,48){\circle*{3}}
\put( 80,48){\circle*{3}}
\put( 96,48){\circle*{3}}
\put(112,48){\circle*{3}}
\put( 56,56){\circle*{3}}
\put( 72,56){\circle*{3}}
\put( 88,56){\circle*{3}}
\put(104,56){\circle*{3}}
\put( 64,64){\circle*{3}}
\put( 80,64){\circle*{3}}
\put( 96,64){\circle*{3}}
\put( 72,72){\circle*{3}}
\put( 88,72){\circle*{3}}
\put( 80,80){\circle*{3}}   %%49data
\put( 88,88){\circle*{3}}
\put( 96,80){\circle*{3}}
\put( 104,72){\circle*{3}}
\put( 112,64){\circle*{3}}
\put( 120,56){\circle*{3}}
\put( 128,48){\circle*{3}}
\put( 136,40){\circle*{3}}
\put( 144,32){\circle*{3}}
%\put( 144,32){\circle{7}}
\put( 152,24){\circle*{3}}
%%\put( 160,16){\circle*{3}}
\end{picture}

\caption{Real polarized K3 when $n\equiv 2\mod 8$ and $\delta_P=0$ (then  
$\delta_{\varphi}=1$) without $(r,a)=(20,2)$ and 
$(r+a=22, \delta_{\varphi P}=0$)}
\label{K3-2mod8-0}
\end{figure}
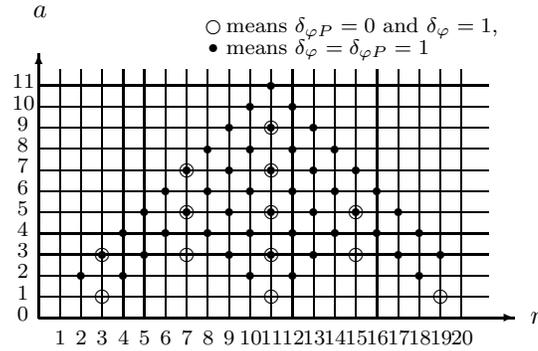

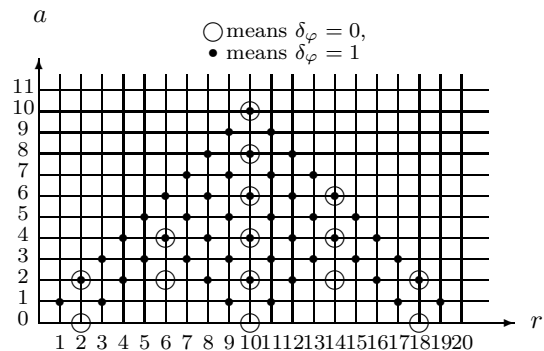
\begin{figure}
\begin{picture}(200,140)
%% K3-n equiv 2mod8 delta_P=1
\put(66,110){\circle{7}}
\put(71,108){{\tiny means $\delta_\varphi=0$,}}
\put(66,102){\circle*{3}}
\put(71,100){{\tiny means $\delta_\varphi=1$}}

\multiput(8,0)(8,0){20}{\line(0,1){94}}
\multiput(0,8)(0,8){11}{\line(1,0){170}}
\put(0,0){\vector(0,1){100}}
\put(0,0){\vector(1,0){180}}
\put(  6,-10){{\tiny $1$}}
\put( 14,-10){{\tiny $2$}}
\put( 22,-10){{\tiny $3$}}
\put( 30,-10){{\tiny $4$}}
\put( 38,-10){{\tiny $5$}}
\put( 46,-10){{\tiny $6$}}
\put( 54,-10){{\tiny $7$}}
\put( 62,-10){{\tiny $8$}}
\put( 70,-10){{\tiny $9$}}
\put( 76,-10){{\tiny $10$}}
\put( 84,-10){{\tiny $11$}}
\put( 92,-10){{\tiny $12$}}
\put(100,-10){{\tiny $13$}}
\put(108,-10){{\tiny $14$}}
\put(116,-10){{\tiny $15$}}
\put(124,-10){{\tiny $16$}}
\put(132,-10){{\tiny $17$}}
\put(140,-10){{\tiny $18$}}
\put(148,-10){{\tiny $19$}}
\put(156,-10){{\tiny $20$}}

\put(-8, -1){{\tiny $0$}}
\put(-8,  7){{\tiny $1$}}
\put(-8, 15){{\tiny $2$}}
\put(-8, 23){{\tiny $3$}}
\put(-8, 31){{\tiny $4$}}
\put(-8, 39){{\tiny $5$}}
\put(-8, 47){{\tiny $6$}}
\put(-8, 55){{\tiny $7$}}
\put(-8, 63){{\tiny $8$}}
\put(-8, 71){{\tiny $9$}}
\put(-10, 79){{\tiny $10$}}
\put(-10, 87){{\tiny $11$}}
\put( -2,114){{\footnotesize $a$}} %vertical
\put(186, -2){{\footnotesize $r$}} %horizontal

%% \delta_\varphi = 0
%%  ( 8r,8a)

\put( 16, 0){\circle{7}}
\put( 80, 0){\circle{7}}
\put(144, 0){\circle{7}}
\put( 16,16){\circle{7}}
\put( 48,16){\circle{7}}
\put( 80,16){\circle{7}}
\put(112,16){\circle{7}}
\put(144,16){\circle{7}}
\put( 48,32){\circle{7}}
\put( 80,32){\circle{7}}
\put(112,32){\circle{7}}
\put( 80,48){\circle{7}}
\put( 112,48){\circle{7}}
\put( 80,64){\circle{7}}
\put( 80,80){\circle{7}}    %%14data

%% \delta_\varphi = 1
%%  ( 8r,8a)
\put(  8, 8){\circle*{3}}
\put( 24, 8){\circle*{3}}
\put( 72, 8){\circle*{3}}
\put( 88, 8){\circle*{3}}
\put(136, 8){\circle*{3}}
\put(152, 8){\circle*{3}}
\put( 16,16){\circle*{3}}
\put( 32,16){\circle*{3}}
\put( 64,16){\circle*{3}}
\put( 80,16){\circle*{3}}
\put( 96,16){\circle*{3}}
\put(128,16){\circle*{3}}
\put(144,16){\circle*{3}}
\put( 24,24){\circle*{3}}
\put( 40,24){\circle*{3}}
\put( 56,24){\circle*{3}}
\put( 72,24){\circle*{3}}
\put( 88,24){\circle*{3}}
\put(104,24){\circle*{3}}
\put(120,24){\circle*{3}}
\put(136,24){\circle*{3}}
\put( 32,32){\circle*{3}}
\put( 48,32){\circle*{3}}
\put( 64,32){\circle*{3}}
\put( 80,32){\circle*{3}}
\put( 96,32){\circle*{3}}
\put(112,32){\circle*{3}}
\put(128,32){\circle*{3}}
\put( 40,40){\circle*{3}}
\put( 56,40){\circle*{3}}
\put( 72,40){\circle*{3}}
\put( 88,40){\circle*{3}}
\put(104,40){\circle*{3}}
\put(120,40){\circle*{3}}
\put( 48,48){\circle*{3}}
\put( 64,48){\circle*{3}}
\put( 80,48){\circle*{3}}
\put( 96,48){\circle*{3}}
\put(112,48){\circle*{3}}
\put( 56,56){\circle*{3}}
\put( 72,56){\circle*{3}}
\put( 88,56){\circle*{3}}
\put(104,56){\circle*{3}}
\put( 64,64){\circle*{3}}
\put( 80,64){\circle*{3}}
\put( 96,64){\circle*{3}}
\put( 72,72){\circle*{3}}
\put( 88,72){\circle*{3}}
\put( 80,80){\circle*{3}}   %%49data
\end{picture}

\caption{Real polarized K3 when $n\equiv 2\mod 8$ and $\delta_P=1$ 
(then $\delta_{\varphi P}=1$)}
\label{K3-2mod8-1}
\end{figure}

\begin{figure}
\begin{picture}(200,140)
%% K3- n equiv -2mod8, delta_P=0
\put(66,110){\circle{5}}
\put(71,108){{\tiny means $\delta_{\varphi P}=0$ and $\delta_\varphi=1$,}}
\put(66,102){\circle*{3}}
\put(71,100){{\tiny means $\delta_\varphi=\delta_{\varphi P}=1$}}

\multiput(8,0)(8,0){20}{\line(0,1){94}}
\multiput(0,8)(0,8){11}{\line(1,0){170}}
\put(0,0){\vector(0,1){100}}
\put(0,0){\vector(1,0){180}}
\put(  6,-10){{\tiny $1$}}
\put( 14,-10){{\tiny $2$}}
\put( 22,-10){{\tiny $3$}}
\put( 30,-10){{\tiny $4$}}
\put( 38,-10){{\tiny $5$}}
\put( 46,-10){{\tiny $6$}}
\put( 54,-10){{\tiny $7$}}
\put( 62,-10){{\tiny $8$}}
\put( 70,-10){{\tiny $9$}}
\put( 76,-10){{\tiny $10$}}
\put( 84,-10){{\tiny $11$}}
\put( 92,-10){{\tiny $12$}}
\put(100,-10){{\tiny $13$}}
\put(108,-10){{\tiny $14$}}
\put(116,-10){{\tiny $15$}}
\put(124,-10){{\tiny $16$}}
\put(132,-10){{\tiny $17$}}
\put(140,-10){{\tiny $18$}}
\put(148,-10){{\tiny $19$}}
\put(156,-10){{\tiny $20$}}

\put(-8, -1){{\tiny $0$}}
\put(-8,  7){{\tiny $1$}}
\put(-8, 15){{\tiny $2$}}
\put(-8, 23){{\tiny $3$}}
\put(-8, 31){{\tiny $4$}}
\put(-8, 39){{\tiny $5$}}
\put(-8, 47){{\tiny $6$}}
\put(-8, 55){{\tiny $7$}}
\put(-8, 63){{\tiny $8$}}
\put(-8, 71){{\tiny $9$}}
\put(-10, 79){{\tiny $10$}}
\put(-10, 87){{\tiny $11$}}
\put( -2,114){{\footnotesize $a$}} %vertical
\put(186, -2){{\footnotesize $r$}} %horizontal

%% \delta_\varphi = 0
%%  ( 8r,8a)
%\put( 16, 0){\circle{7}}
%\put( 80, 0){\circle{7}}
%\put(144, 0){\circle{7}}
%\put( 16,16){\circle{7}}
%\put( 48,16){\circle{7}}
%\put( 80,16){\circle{7}}
%\put(112,16){\circle{7}}
%\put(144,16){\circle{7}}
%\put( 48,32){\circle{7}}
%\put( 80,32){\circle{7}}
%\put(112,32){\circle{7}}
%\put( 80,48){\circle{7}}
%\put( 112,48){\circle{7}}
%\put( 80,64){\circle{7}}
%\put( 80,80){\circle{7}}    %%14data

%% \delta_\varphi = 1
%%  ( 8r,8a)
%\put(  8, 8){\circle*{3}}
\put(  8, 8){\circle{5}}
%\put( 24, 8){\circle*{3}}
%\put( 72, 8){\circle*{3}}
\put( 72, 8){\circle{5}}
\put( 72, 24){\circle{5}}
\put( 72, 40){\circle{5}}
\put( 72, 56){\circle{5}}
\put( 72, 72){\circle{5}}
%\put( 88, 8){\circle*{3}}
%\put(136, 8){\circle*{3}}
\put(136, 8){\circle{5}}
\put(136, 24){\circle{5}}
%\put(152, 8){\circle*{3}}
\put( 16,16){\circle*{3}}
%\put( 32,16){\circle*{3}}
\put( 64,16){\circle*{3}}
\put( 80,16){\circle*{3}}
%\put( 96,16){\circle*{3}}
\put(128,16){\circle*{3}}
\put(144,16){\circle*{3}}
\put( 24,24){\circle*{3}}
%\put( 40,24){\circle*{3}}
\put( 40,24){\circle{5}}
\put( 40,40){\circle{5}}
\put( 56,24){\circle*{3}}
\put( 72,24){\circle*{3}}
\put( 88,24){\circle*{3}}
%\put(104,24){\circle*{3}}
\put(104,24){\circle{5}}
\put(104,40){\circle{5}}
\put(104,56){\circle{5}}
\put(120,24){\circle*{3}}
\put(136,24){\circle*{3}}
\put( 32,32){\circle*{3}}
\put( 48,32){\circle*{3}}
\put( 64,32){\circle*{3}}
\put( 80,32){\circle*{3}}
\put( 96,32){\circle*{3}}
\put(112,32){\circle*{3}}
\put(128,32){\circle*{3}}
\put( 40,40){\circle*{3}}
\put( 56,40){\circle*{3}}
\put( 72,40){\circle*{3}}
\put( 88,40){\circle*{3}}
\put(104,40){\circle*{3}}
\put(120,40){\circle*{3}}
\put( 48,48){\circle*{3}}
\put( 64,48){\circle*{3}}
\put( 80,48){\circle*{3}}
\put( 96,48){\circle*{3}}
\put(112,48){\circle*{3}}
\put( 56,56){\circle*{3}}
\put( 72,56){\circle*{3}}
\put( 88,56){\circle*{3}}
\put(104,56){\circle*{3}}
\put( 64,64){\circle*{3}}
\put( 80,64){\circle*{3}}
\put( 96,64){\circle*{3}}
\put( 72,72){\circle*{3}}
\put( 88,72){\circle*{3}}
\put( 80,80){\circle*{3}}   %%49data
\put( 88,88){\circle*{3}}
\put( 96,80){\circle*{3}}
\put( 104,72){\circle*{3}}
\put( 112,64){\circle*{3}}
\put( 120,56){\circle*{3}}
\put( 128,48){\circle*{3}}
\put( 136,40){\circle*{3}}
\put( 144,32){\circle*{3}}
%\put( 144,32){\circle{7}}
\put( 152,24){\circle*{3}}
%%\put( 160,16){\circle*{3}}
\end{picture}

\caption{Real polarized K3 when $n\equiv -2\mod 8$ and $\delta_P=0$ (then  
$\delta_{\varphi}=1$) without $(r,a)=(20,2)$ and $(r+a=22, \delta_{\varphi P}=0$)}
\label{K3-6mod8-0}
\end{figure}
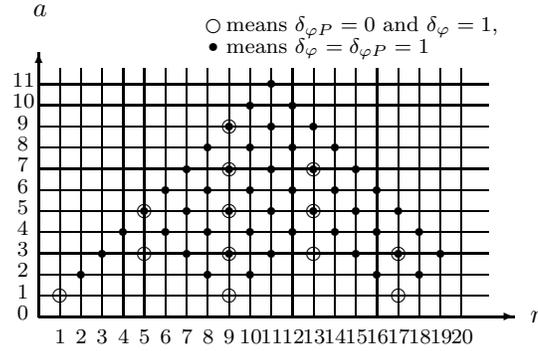

\Proof Using Theorem 3.4.3 from \cite{Nikulin79}, we can get the 
list of all possible invariants $(n; r, a, \delta_\varphi,$ 
$\delta_{\varphi P})$ 
of real polarized K3 when $(r,a)\not=(20,2)$ and $r+a<22$, 
if $\delta_{\varphi P}=0$. They depend on $n\mod 8$ and are presented in 
figures \ref{K3-0mod8-0} --- \ref{K3-6mod8-0} (when $r+a=22$ and 
$\delta_{\varphi P}=0$, the result of \cite{Nikulin79} is more 
complicated: it depends on $n\mod 16$ when $(r,a)\not=(20,2)$; 
if $(r,a)=(20,2)$, 
it depends even on prime decomposition of $n$).  
By inspection (which requires some time and patience), one can see that 
all these invariants can be obtained applying all possible deformations 
\eqref{typesdef} from  
$\bp^2$, $\bhh$, $\bff_1$ and $\bff_4$.    
Again applying Theorem 3.4.3 from \cite{Nikulin79} when $n\le 4$ and 
either $(r,a)=(20,2)$ or $r+a=22$ and $\delta_{\varphi P}=0$,  
we get all of them as deformations from $\bp^2$ and $\bee$. The list of  
all possible invariants for $n=4$ is given in figures \ref{K3-4mod8-1}  
and \ref{K3-4-0}. 
 
It proves the statement. 

\bigskip 

As a particular result, we get that any non-singular real quartic in 
$\bp^3$ ($n=4$) can be obtained as a deformation of a general real 
double K3 cover of ellipsoid, hyperboloid or $\bff_4$. These deformations 
can be given explicitly \cite{Shah81} and give an effective method of 
construction of all (i. e. representatives of all connected components 
of moduli) examples of real quartics in $\bp^3$. Actually, only one 
connected component for $n=4$: with 
$(r=19,\, a=1,\, \delta_P=\delta_\varphi=\delta_{\varphi P}=1)$,   
requires $\bff_4$. All other can be obtained from ellipsoid and 
hyperboloid.

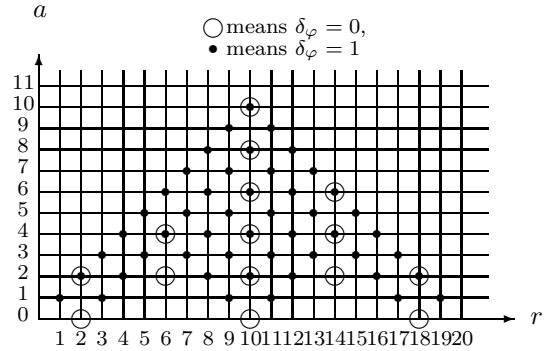
\begin{figure}
\begin{picture}(200,140)
%% K3-n equiv -2mod8 delta_P=1
\put(66,110){\circle{7}}
\put(71,108){{\tiny means $\delta_\varphi=0$,}}
\put(66,102){\circle*{3}}
\put(71,100){{\tiny means $\delta_\varphi=1$}}

\multiput(8,0)(8,0){20}{\line(0,1){94}}
\multiput(0,8)(0,8){11}{\line(1,0){170}}
\put(0,0){\vector(0,1){100}}
\put(0,0){\vector(1,0){180}}
\put(  6,-10){{\tiny $1$}}
\put( 14,-10){{\tiny $2$}}
\put( 22,-10){{\tiny $3$}}
\put( 30,-10){{\tiny $4$}}
\put( 38,-10){{\tiny $5$}}
\put( 46,-10){{\tiny $6$}}
\put( 54,-10){{\tiny $7$}}
\put( 62,-10){{\tiny $8$}}
\put( 70,-10){{\tiny $9$}}
\put( 76,-10){{\tiny $10$}}
\put( 84,-10){{\tiny $11$}}
\put( 92,-10){{\tiny $12$}}
\put(100,-10){{\tiny $13$}}
\put(108,-10){{\tiny $14$}}
\put(116,-10){{\tiny $15$}}
\put(124,-10){{\tiny $16$}}
\put(132,-10){{\tiny $17$}}
\put(140,-10){{\tiny $18$}}
\put(148,-10){{\tiny $19$}}
\put(156,-10){{\tiny $20$}}

\put(-8, -1){{\tiny $0$}}
\put(-8,  7){{\tiny $1$}}
\put(-8, 15){{\tiny $2$}}
\put(-8, 23){{\tiny $3$}}
\put(-8, 31){{\tiny $4$}}
\put(-8, 39){{\tiny $5$}}
\put(-8, 47){{\tiny $6$}}
\put(-8, 55){{\tiny $7$}}
\put(-8, 63){{\tiny $8$}}
\put(-8, 71){{\tiny $9$}}
\put(-10, 79){{\tiny $10$}}
\put(-10, 87){{\tiny $11$}}
\put( -2,114){{\footnotesize $a$}} %vertical
\put(186, -2){{\footnotesize $r$}} %horizontal

%% \delta_\varphi = 0
%%  ( 8r,8a)
\put( 16, 0){\circle{7}}
\put( 80, 0){\circle{7}}
\put(144, 0){\circle{7}}
\put( 16,16){\circle{7}}
\put( 48,16){\circle{7}}
\put( 80,16){\circle{7}}
\put(112,16){\circle{7}}
\put(144,16){\circle{7}}
\put( 48,32){\circle{7}}
\put( 80,32){\circle{7}}
\put(112,32){\circle{7}}
\put( 80,48){\circle{7}}
\put( 112,48){\circle{7}}
\put( 80,64){\circle{7}}
\put( 80,80){\circle{7}}    %%14data

%% \delta_\varphi = 1
%%  ( 8r,8a)
\put(  8, 8){\circle*{3}}
\put( 24, 8){\circle*{3}}
\put( 72, 8){\circle*{3}}
\put( 88, 8){\circle*{3}}
\put(136, 8){\circle*{3}}
\put(152, 8){\circle*{3}}
\put( 16,16){\circle*{3}}
\put( 32,16){\circle*{3}}
\put( 64,16){\circle*{3}}
\put( 80,16){\circle*{3}}
\put( 96,16){\circle*{3}}
\put(128,16){\circle*{3}}
\put(144,16){\circle*{3}}
\put( 24,24){\circle*{3}}
\put( 40,24){\circle*{3}}
\put( 56,24){\circle*{3}}
\put( 72,24){\circle*{3}}
\put( 88,24){\circle*{3}}
\put(104,24){\circle*{3}}
\put(120,24){\circle*{3}}
\put(136,24){\circle*{3}}
\put( 32,32){\circle*{3}}
\put( 48,32){\circle*{3}}
\put( 64,32){\circle*{3}}
\put( 80,32){\circle*{3}}
\put( 96,32){\circle*{3}}
\put(112,32){\circle*{3}}
\put(128,32){\circle*{3}}
\put( 40,40){\circle*{3}}
\put( 56,40){\circle*{3}}
\put( 72,40){\circle*{3}}
\put( 88,40){\circle*{3}}
\put(104,40){\circle*{3}}
\put(120,40){\circle*{3}}
\put( 48,48){\circle*{3}}
\put( 64,48){\circle*{3}}
\put( 80,48){\circle*{3}}
\put( 96,48){\circle*{3}}
\put(112,48){\circle*{3}}
\put( 56,56){\circle*{3}}
\put( 72,56){\circle*{3}}
\put( 88,56){\circle*{3}}
\put(104,56){\circle*{3}}
\put( 64,64){\circle*{3}}
\put( 80,64){\circle*{3}}
\put( 96,64){\circle*{3}}
\put( 72,72){\circle*{3}}
\put( 88,72){\circle*{3}}
\put( 80,80){\circle*{3}}   %%49data
\end{picture}

\caption{Real polarized K3 when $n\equiv -2\mod 8$ and $\delta_P=1$ 
(then $\delta_{\varphi P}=1$)}
\label{K3-6mod8-1}
\end{figure}

\begin{figure}
\begin{picture}(200,140)
%% K3- n equiv 4mod8, delta_P=0
\put(66,118){\circle{7}}
\put(71,116){{\tiny means $\delta_\varphi=0$ and $\delta_{\varphi P}=1$},}
\put(66,110){\circle{5}}
\put(71,108){{\tiny means $\delta_{\varphi P}=0$ and $\delta_\varphi=1$,}}
\put(66,102){\circle*{3}}
\put(71,100){{\tiny means $\delta_\varphi=\delta_{\varphi P}=1$}}

\multiput(8,0)(8,0){20}{\line(0,1){94}}
\multiput(0,8)(0,8){11}{\line(1,0){170}}
\put(0,0){\vector(0,1){100}}
\put(0,0){\vector(1,0){180}}
\put(  6,-10){{\tiny $1$}}
\put( 14,-10){{\tiny $2$}}
\put( 22,-10){{\tiny $3$}}
\put( 30,-10){{\tiny $4$}}
\put( 38,-10){{\tiny $5$}}
\put( 46,-10){{\tiny $6$}}
\put( 54,-10){{\tiny $7$}}
\put( 62,-10){{\tiny $8$}}
\put( 70,-10){{\tiny $9$}}
\put( 76,-10){{\tiny $10$}}
\put( 84,-10){{\tiny $11$}}
\put( 92,-10){{\tiny $12$}}
\put(100,-10){{\tiny $13$}}
\put(108,-10){{\tiny $14$}}
\put(116,-10){{\tiny $15$}}
\put(124,-10){{\tiny $16$}}
\put(132,-10){{\tiny $17$}}
\put(140,-10){{\tiny $18$}}
\put(148,-10){{\tiny $19$}}
\put(156,-10){{\tiny $20$}}

\put(-8, -1){{\tiny $0$}}
\put(-8,  7){{\tiny $1$}}
\put(-8, 15){{\tiny $2$}}
\put(-8, 23){{\tiny $3$}}
\put(-8, 31){{\tiny $4$}}
\put(-8, 39){{\tiny $5$}}
\put(-8, 47){{\tiny $6$}}
\put(-8, 55){{\tiny $7$}}
\put(-8, 63){{\tiny $8$}}
\put(-8, 71){{\tiny $9$}}
\put(-10, 79){{\tiny $10$}}
\put(-10, 87){{\tiny $11$}}
\put( -2,114){{\footnotesize $a$}} %vertical
\put(186, -2){{\footnotesize $r$}} %horizontal

%% \delta_\varphi = 0
%%  ( 8r,8a)
%\put( 16, 0){\circle{7}}
%\put( 80, 0){\circle{7}}
%\put(144, 0){\circle{7}}
\put( 16,16){\circle{7}}
\put( 48,16){\circle{7}}
\put( 80,16){\circle{7}}
\put(112,16){\circle{7}}
\put(144,16){\circle{7}}
\put( 48,32){\circle{7}}
\put( 80,32){\circle{7}}
\put(112,32){\circle{7}}
\put( 80,48){\circle{7}}
\put( 112,48){\circle{7}}
\put( 80,64){\circle{7}}
\put( 80,80){\circle{7}}    %%14data
%% \delta_\varphi = 1
%%  ( 8r,8a)
%\put(  8, 8){\circle*{3}}
%\put( 24, 8){\circle*{3}}
%\put( 72, 8){\circle*{3}}
%\put( 88, 8){\circle*{3}}
%\put(136, 8){\circle*{3}}
%\put(152, 8){\circle*{3}}
%\put( 16,16){\circle*{3}}
%\put( 32,16){\circle*{3}}
\put( 32,16){\circle{5}}
\put( 32,32){\circle{5}}
%\put( 64,16){\circle*{3}}
\put( 64,16){\circle{5}}
\put( 64,32){\circle{5}}
\put( 64,48){\circle{5}}
\put( 64,64){\circle{5}}
%\put( 80,16){\circle*{3}}
%\put( 96,16){\circle*{3}}
\put( 96,16){\circle{5}}
\put( 96,32){\circle{5}}
\put( 96,48){\circle{5}}
\put( 96,64){\circle{5}}
\put( 96,80){\circle{5}}

%\put(128,16){\circle*{3}}
\put(128,16){\circle{5}}
\put(128,32){\circle{5}}
%\put(144,16){\circle*{3}}
\put( 24,24){\circle*{3}}
\put( 40,24){\circle*{3}}
\put( 56,24){\circle*{3}}
\put( 72,24){\circle*{3}}
\put( 88,24){\circle*{3}}
\put(104,24){\circle*{3}}
\put(120,24){\circle*{3}}
\put(136,24){\circle*{3}}
\put( 32,32){\circle*{3}}
\put( 48,32){\circle*{3}}
\put( 64,32){\circle*{3}}
\put( 80,32){\circle*{3}}
\put( 96,32){\circle*{3}}
\put(112,32){\circle*{3}}
\put(128,32){\circle*{3}}
\put( 40,40){\circle*{3}}
\put( 56,40){\circle*{3}}
\put( 72,40){\circle*{3}}
\put( 88,40){\circle*{3}}
\put(104,40){\circle*{3}}
\put(120,40){\circle*{3}}
\put( 48,48){\circle*{3}}
\put( 64,48){\circle*{3}}
\put( 80,48){\circle*{3}}
\put( 96,48){\circle*{3}}
\put(112,48){\circle*{3}}
\put( 56,56){\circle*{3}}
\put( 72,56){\circle*{3}}
\put( 88,56){\circle*{3}}
\put(104,56){\circle*{3}}
\put( 64,64){\circle*{3}}
\put( 80,64){\circle*{3}}
\put( 96,64){\circle*{3}}
\put( 72,72){\circle*{3}}
\put( 88,72){\circle*{3}}
\put( 80,80){\circle*{3}}   %%49data
\put( 88,88){\circle*{3}}
\put( 96,80){\circle*{3}}
\put( 104,72){\circle*{3}}
\put( 112,64){\circle*{3}}
\put( 120,56){\circle*{3}}
\put( 128,48){\circle*{3}}
\put( 136,40){\circle*{3}}
\put( 144,32){\circle*{3}}
\put( 144,32){\circle{7}}
\put( 152,24){\circle*{3}}
%%\put( 160,16){\circle*{3}}
\put (160,16){\circle{5}}
\end{picture}

\caption{Real polarized K3 when $n=4$ and $\delta_P=0$}
\label{K3-4-0}
\end{figure}
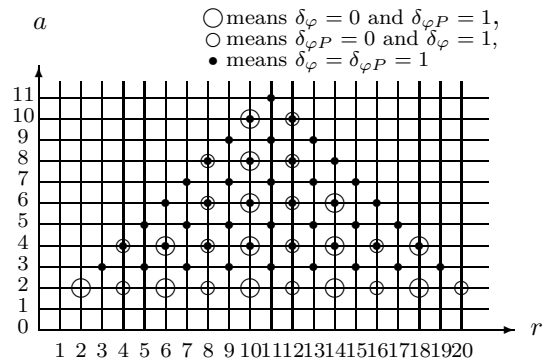

%%%%%%%%%%%%%%%%%%%%%%%%%%%%%%%%%%%%%%%%%%%%%%%%%%%%%%

\end{document}